\documentclass[11pt,a4paper]{article}
\usepackage[T1]{fontenc}
\usepackage[utf8]{inputenc}
\usepackage{lmodern,amsmath,amssymb,amsthm,mathtools}
\usepackage[a4paper,margin=25mm,headheight=15pt]{geometry}
\usepackage{microtype,booktabs,array,longtable,enumitem}
\usepackage{float}
\usepackage{xcolor,tikz,placeins}
\usetikzlibrary{arrows.meta,positioning,calc}
\usepackage[hyphens]{url}
\usepackage{hyperref}
\definecolor{ink}{HTML}{163B4B}
\definecolor{accent}{HTML}{28677B}
\definecolor{pale}{HTML}{EEF4F6}
\hypersetup{colorlinks=true,linkcolor=blue,citecolor=blue,urlcolor=blue,
 pdftitle={The Brown--Shields Conjecture: Two Approximation Theorems and a Counterexample},
 pdfsubject={Cyclicity and radial boundary zeros in the classical Dirichlet space}}
\setlist{leftmargin=1.5em,itemsep=2pt,topsep=3pt}
\numberwithin{equation}{section}
\newtheorem{theorem}{Theorem}[section]
\newtheorem{lemma}[theorem]{Lemma}
\newtheorem{proposition}[theorem]{Proposition}
\newtheorem{corollary}[theorem]{Corollary}
\newtheorem{published}[theorem]{Theorem}
\theoremstyle{definition}

\theoremstyle{remark}
\newcolumntype{L}[1]{>{\raggedright\arraybackslash}p{#1}}
\newcommand{\D}{\mathcal D}\newcommand{\DD}{\mathbb D}
\newcommand{\T}{\mathbb T}\newcommand{\X}{\mathcal X}
\newcommand{\Xpr}{\mathcal X_{\mathrm{pr}}}
\newcommand{\CR}{C_R}\newcommand{\dlog}{d_{\log}}
\newcommand{\Zrad}{Z_{\mathrm{rad}}}\newcommand{\Zapp}{Z_{\mathrm{app}}}
\DeclareMathOperator{\supp}{supp}
\DeclareMathOperator{\dist}{dist}

\begin{document}
\pagenumbering{roman}
\begin{center}
{\LARGE\bfseries
The Brown--Shields Conjecture: Two Approximation Theorems and a Counterexample \par}

\vspace{1.2em}
{\large
 \newcommand{\authoraffiliation}[2]{%
   \footnote{%
     \hspace{0.4em}%
     \parbox[t]{\dimexpr\linewidth-2.2em\relax}{%
       \raggedright
       #1\par
       \mbox{Email:\enspace\href{mailto:#2}{\nolinkurl{#2}}}\strut
     }%
   }%
 }%
 Bonan Chen\authoraffiliation{School of Mathematical Sciences,
 Soochow University, Suzhou 215006, P. R. China.}{bnchen@suda.edu.cn}
 \qquad
 Xiang Fang\authoraffiliation{Department of Applied Mathematics,
 National Yang Ming Chiao Tung University, Hsinchu, Taiwan.}{xfang@nycu.edu.tw}
 \qquad
 Feng Guo\authoraffiliation{School of Mathematics,
 South China University of Technology,
 Guangzhou 510640, P. R. China.}{70207994@nuaa.edu.cn}
 \par}
\end{center}

\vspace{1.5em}
 
\begin{abstract}
We construct a noncyclic outer function in the classical Dirichlet
space whose radial zero set has logarithmic capacity zero, answering
the question of Brown and Shields in the negative.
Our first approximation theorem approximates any given proper outer
canonical generator in the Dirichlet norm by the canonical generator of an invariant subspace associated with a nonzero finitely supported positive measure of arbitrarily small total mass.
For a fixed nonzero finitely supported positive measure on the circle,
the second theorem approximates the canonical function of the corresponding invariant subspace in Dirichlet norm by outer canonical functions of trace-zero
spaces. Each approximant is continuous on the closed disk, with zero
set equal to the defining compact set. The theorem gives explicit
bounds on the location and capacity of these sets. The projections
onto the trace-zero spaces converge strongly to the projection onto 	that invariant subspace.
A nested-ball construction then gives counterexamples arbitrarily
close to every proper outer canonical generator, both in Dirichlet
norm and in \(L^1\) distance between the logarithms of the boundary
moduli.
\end{abstract}

\vspace{1.5em}

{\small
\par

\par
}

\bigskip

\noindent \emph{This is an expository version with navigational aids and an extended review of the literature, not intended for publication.} 

\tableofcontents
\clearpage
\pagenumbering{arabic}

\section{Introduction}\label{sec:01}

\subsection{Cyclicity and boundary zeros}\label{subsec:01-problem}

Let \(\DD=\{z\in\mathbb C:|z|<1\}\) be the unit disk and
\(\T=\partial\DD\) its boundary. Write \(A_0\) for area on \(\DD\)
divided by \(\pi\). For \(f(z)=\sum_{n\ge0}a_nz^n\) analytic on
\(\DD\), set
\begin{equation}\label{eq:full-norm}
 \|f\|_{\D}^{2}=\|f\|_{H^2}^{2}+D(f)
 =\sum_{n\ge0}(n+1)|a_n|^2,
 \qquad D(f)=\int_{\DD}|f'(z)|^2\,dA_0(z).
\end{equation}
The Hardy space \(H^2\) consists of the analytic functions with finite
norm \(\|f\|_{H^2}^2=\sum_{n\ge0}|a_n|^2\).
The Dirichlet space \(\D\) consists of the analytic functions for which
\eqref{eq:full-norm} is finite. For \(f\in\D\), write
\begin{equation}\label{eq:module}
 [f]=\overline{\{pf:p\in\mathbb C[z]\}}^{\,\D}
\end{equation}
for the invariant subspace generated by $f$.
We call \(f\) \emph{cyclic} if \([f]=\D\). If \(1\in[f]\), then
every monomial lies in \([f]\), since multiplication by \(z\) is
bounded in \(\D\). Taylor polynomials are dense by
\eqref{eq:full-norm}, so \(1\in[f]\) suffices for cyclicity.

It is easy to see that an interior zero of \(f\) prevents cyclicity. Indeed, evaluation at each point
of \(\DD\) is a continuous functional in the Dirichlet norm, so all functions in
\([f]\) vanish at that point. On the other hand, a boundary zero can coexist with
cyclicity. For example, let \(H_N=\sum_{k=0}^N 1/(k+1)\) and put
\[
 h_N(z)=H_N^{-1}\sum_{k=0}^N\frac{z^k}{k+1}.
\]
Then \(h_N(1)=1\), so \(1-h_N\) is a polynomial multiple of \(1-z\).
By \eqref{eq:full-norm}, \(\|h_N\|_{\D}^2=1/H_N\to0\).
The optimal-error formulas of Fricain, Mashreghi, and Seco
\cite[Theorem~3.9 and Corollary~3.10]{FMS14} give a broader setting
for this calculation. Here the weights are \(k+1\); for \(N\ge1\),
the polynomial \((1-h_N)/(1-z)\) multiplying \(1-z\) has degree
\(N-1\).
Thus \(1-z\) is cyclic even though it vanishes at \(1\).

There are two classical necessary conditions for cyclicity.
The first comes from the Hardy space. Since Dirichlet convergence
implies Hardy convergence, a Dirichlet-cyclic function must also
be cyclic in \(H^2\). Beurling's invariant-subspace theorem \cite{Beurling49}
shows that every nonzero closed shift-invariant subspace of the Hardy
space is an inner function times the whole Hardy space. Here an inner
function is bounded and analytic in the disk, with boundary modulus one
almost everywhere. In particular, a nonzero Hardy function is cyclic
exactly when its inner factor is constant. Thus the Dirichlet-cyclic
function is \emph{Hardy outer}; Duren gives a modern account
\cite[Theorem~7.4, p.~114]{Duren70}.
An outer function is recovered, up to a constant
of modulus one, from its boundary modulus. More precisely, let \(m\)
be normalized arclength on \(\T\), with \(dm(e^{it})=dt/(2\pi)\),
and let \(f^*\) denote the radial
boundary value of \(f\), defined \(m\)-almost everywhere for \(f\in H^2\).
A nonzero \(f\in H^2\) is outer if \(\log|f^*|\) is integrable and
\begin{equation}\label{eq:outer-intro}
 f(z)=\omega\exp\!\left(\int_{\T}\frac{\zeta+z}{\zeta-z}
                    \log|f^*(\zeta)|\,dm(\zeta)\right),\qquad |\omega|=1.
\end{equation}

The second condition comes from the Dirichlet boundary trace.
We use the radial zero set
\begin{equation}\label{eq:radial-zero}
 \Zrad(f)=\{\zeta\in\T:\lim_{r\uparrow1}f(r\zeta)=0\}.
\end{equation}
Thus \(\Zrad(f)\) is determined by the analytic function
\(f\); assigning values to an almost-everywhere representative of
\(f^*\) does not define this set. Logarithmic capacity measures the
size of such a set.
For a compact set \(K\subset\T\), positive logarithmic capacity means
that \(K\) supports a probability measure \(\mu\) of finite logarithmic
energy, that is,
\[
 \iint_{K\times K}\log\frac{1}{|\zeta-\eta|}
       \,d\mu(\zeta)\,d\mu(\eta)<\infty.
\]
A set of capacity zero is called \emph{polar}. A singleton is polar.
Section~\ref{subsec:03-pairing} gives the capacity convention for
arbitrary sets. Some Cantor sets of arclength zero have positive
capacity, as the construction in Section~\ref{subsec:03-seed} illustrates.

Dirichlet boundary traces agree with radial limits outside a polar
set. We say that a statement holds \emph{quasi-everywhere} if it fails
only on a polar set. Vanishing on a given boundary set, with this
exception, defines a closed subspace of \(\D\).
Consequently, if \(\Zrad(f)\) has positive capacity, every limit of
polynomial multiples of \(f\) still vanishes there in this sense, and
the constant \(1\) is excluded. Section~\ref{sec:03} explains the trace
and proves the two necessary conditions in the conventions used here.
Brown and Shields asked whether they also suffice
\cite[Question~12, p.~296]{BS84}. Must a Hardy outer
function in \(\D\) be cyclic if its radial zero set has logarithmic
capacity zero?
Their preceding question asks whether the generated space of every
outer \(f\in\D\) has the form \(\D_E\) for some boundary set \(E\),
where membership means radial vanishing quasi-everywhere on \(E\).
The example below answers both questions. The deduction for
Question~11 follows the proof of Theorem~\ref{thm:main} in
Section~\ref{sec:13}.

\subsection{The counterexample and approximation using finite interior zero sets}\label{subsec:01-main}

We prove that the two necessary conditions can hold while polynomial
multiples still fail to approximate \(1\). The precise statement is the
following; Section~\ref{sec:13} gives its proof and refinements.

\begin{theorem}\label{thm:main}
There exists a Hardy outer function $f\in\D$ such that $\Zrad(f)$ has
logarithmic capacity zero and $[f]\ne\D$.
\end{theorem}

The construction makes the approximation error explicit by using
\emph{canonical} functions. Here \(f\) is canonical if \(\|f\|_{\D}=1\), \(f(0)>0\), and \(f\) maximizes \(\lvert g(0)\rvert\) over the unit ball of its invariant subspace \([f]\). This extremal property identifies the
best approximation of \(1\) by polynomial multiples of \(f\).
For the counterexample in canonical form, with \(a=f(0)\), it gives
\[
 \inf_{p\in\mathbb C[z]}\|1-pf\|_{\D}^2=1-a^2>0.
\]
The minimum is attained by the constant polynomial \(p=a\).
For unit extremals, B\'en\'eteau, Fleeman, Khavinson, Seco, and
Sola proved both the constant-approximant formula and this exact
squared distance \cite[Theorem~2.3(ii)--(iii)]{BFKSS18}. Here the
extremal realizing the positive distance also has a polar radial zero set.
In this canonical class, keeping the value at zero bounded away from
one therefore preserves a positive approximation error.

\textbf{Local refinement.}
Theorem~\ref{thm:local} starts with any outer canonical function \(f_0\)
whose invariant subspace \([f_0]\) is proper, meaning that \([f_0]\ne\D\).
It finds counterexamples arbitrarily close to \(f_0\) in the Dirichlet
norm and in the integrated absolute difference of the logarithms of
their boundary moduli. Their values at zero can be made arbitrarily
close to \(f_0(0)\). For these counterexamples the squared errors
stay bounded away from zero as \(f\) approaches \(f_0\).
This density assertion is relative to the class of noncyclic outer canonical functions, not to all of \(\D\).

Corollary~\ref{cor:dense-radial-zero-set} gives a counterexample whose
radial zero set is a dense polar \(G_\delta\) subset of the circle.

\textbf{Numerical consequence.}
Corollary~\ref{cor:numerical} gives one choice with
\(\tfrac12 e^{-1/16}\le a\le\tfrac12 e^{1/16}\). The squared error
for every polynomial is then at least
\[
 1-\tfrac14 e^{1/8}>\tfrac{175}{256}.
\]

\textbf{Approximation by finite zero sets.}
Let \(E\subset\T\) be compact with positive logarithmic capacity.
Its \emph{trace-zero space} consists of the Dirichlet functions whose
radial limits vanish on \(E\) outside a polar set. Suppose that this
space is nonzero, and consider its unit extremal for evaluation at
zero, normalized to have positive value there.
Corollary~\ref{cor:trace-zero-finite-approx}, obtained from the
finite-zero theorem in Section~\ref{sec:06}, approximates this
extremal in the Dirichlet norm by the corresponding normalized
extremals for spaces of functions vanishing at finitely many points
of the disk.

Luo and Richter asked whether such trace-zero extremals admit
locally uniform approximation by extremals associated with interior
zero sets \cite[Question~6.3]{LR15}. Immediately after the question
they observe that an affirmative answer already permits finite sets.
Luo and Richter \cite[Corollary~4.5, p.~431]{LR15} show that if
the normalized extremals converge pointwise to the normalized
extremal of the target subspace, then the associated orthogonal
projections converge in the strong operator topology.
The approximating points supplied by
Theorem~\ref{thm:finite-zero-approximation} lie on one circle at each
stage, and the corresponding sums of \(-\log|z|\) tend to zero.
This answers Question~6.3.

Corollary~\ref{cor:finite-codimensional-approximation} answers
Question~6.2. It shows that the orthogonal projection onto
any nonzero closed shift-invariant subspace of the classical
Dirichlet space is a strong operator limit of projections onto
finite-codimensional invariant subspaces defined by finitely many
interior vanishing conditions, with multiplicities allowed.

Theorem~\ref{thm:locally-uniform-extremals} answers
\cite[Question~6.4]{LR15} by characterizing the locally uniform limits
of unimodular multiples of Dirichlet extremals. Extremals associated
with finitely many interior vanishing conditions suffice to obtain all
these limits.

\subsection{Boundary values, capacity, and cyclicity}\label{subsec:01-boundary-history}

The Brown--Shields problem has a long history and connects several areas of function theory. Since the present version is intended to be expository, we devote this subsection and the next four to an extended review of these connections and the surrounding literature. We discuss boundary values and capacity, canonical generators and polynomial approximation, logarithmic criteria and factorization, stability under approximation, and related function spaces. Our purpose is not to give an exhaustive survey, but to explain how these developments relate to the problem, what they contribute, and where important distinctions arise. Readers primarily interested in the construction may proceed directly to Subsection~\ref{subsec:01-mechanism}.

The definition of a trace-zero space uses boundary values defined
quasi-everywhere. A useful way to
choose a representative is through \emph{quasi-continuity}: for every
\(\varepsilon>0\), the representative is continuous after an open set of capacity less than \(\varepsilon\) has been removed.
The open set may depend on \(\varepsilon\).
This is a continuity property; the corresponding uniqueness statement is that two such representatives of the same trace agree quasi-everywhere, meaning outside a polar set. Richter and Yilmaz recall these properties of quasi-continuous representatives in \cite[Section~2]{RY19};  the book of El-Fallah, Kellay,
Mashreghi, and Ransford \cite{Primer14} provides a broader account of
Dirichlet boundary theory. Section~\ref{subsec:03-trace} states the
precise representative theorem used here. The following example
shows why this choice matters even for the constant analytic function.

\textbf{Boundary representatives and radial zeros.}
Let \(K\subset\T\) be a Cantor-type compactum of arclength zero and positive
logarithmic capacity. The functions \(1\) and \(1-\mathbf{1}_K\)
represent the same element of \(L^2(m)\). The second representative
vanishes on \(K\), although the analytic constant function \(1\)
has radial limit \(1\) at every point of \(\T\). Thus changing an
arclength representative can create a zero set of positive capacity.

The continuous function \(1\) is a quasi-continuous representative
of this trace. By uniqueness quasi-everywhere, every quasi-continuous
representative of the same trace agrees with \(1\) outside a polar set.
In particular, the Hardy boundary class is specified up to sets of
arclength zero, and its quasi-continuous representative is specified
up to polar sets. The  radial limits are determined pointwise
by the analytic function wherever those limits exist.

The appearance of capacity in this example reflects a basic feature
of finite Dirichlet energy. Beurling's exceptional-set theorem
\cite{Beurling40} guarantees finite radial limits outside a set of
logarithmic capacity zero. Thus the energy controls pointwise boundary
behavior on a scale finer than arclength alone: a set of arclength
zero can still be too large to discard from a Dirichlet boundary
statement. El-Fallah, Kellay, Mashreghi, and Ransford give a
self-contained treatment of the strong-type capacitary inequality
\cite{EFKMR12}. The inequality also appears in their book
\cite[Theorem~3.3.4, p.~37]{Primer14}.
Their estimate controls capacities of boundary level sets across
the range of levels by the Dirichlet norm. This explains why
capacity enters both approximation and trace arguments. These
results supply a boundary theory for finite-energy functions;
the additional cyclicity question asks whether polynomial
multiples of one given function recover the entire space.

Capacity also governs which analytic functions act by bounded
multiplication on \(\D\). Such a function \(\varphi\) is called a
\emph{multiplier}: the map \(g\mapsto\varphi g\) is a bounded
operator from \(\D\) to itself. The derivative identity
\((\varphi g)'=\varphi'g+\varphi g'\) indicates why this asks more
than boundedness of \(\varphi\). Besides controlling \(\varphi\),
one must control the contribution of \(\varphi'\) against every
Dirichlet function. A positive measure \(\nu\) on the disk is a
\emph{Dirichlet Carleson measure} when
\[
 \int_{\DD}|g|^2\,d\nu\le C\|g\|_{\D}^2
 \qquad(g\in\D)
\]
for some constant \(C\). Stegenga's multiplier theorem
\cite{Stegenga80} identifies the multipliers as the bounded analytic
functions for which \(|\varphi'|^2dA_0\) is such a measure. His
capacitary criterion tests the measure above finite unions of boundary
arcs, so the interaction between separated boundary pieces remains
part of the test. Arcozzi, Rochberg, and Sawyer give a discrete formulation
of Carleson embedding questions for analytic Besov spaces:
they reduce these questions to Hardy inequalities on trees
and obtain applications to multipliers and interpolation
\cite{ARS02}.
The survey of Arcozzi, Rochberg, Sawyer, and Wick
\cite[Sections~3.1 and 6.1]{ARSW11} explains these connections and
records Stegenga's characterization. The capacities in these accounts
have their own conventions; the present discussion uses their
geometric meaning without identifying their numerical normalizations
with the obstacle capacity defined in Section~\ref{subsec:03-norm}.

For cyclicity itself, further results based on the boundary
zero set have led to substantial positive conclusions. El-Fallah,
Kellay, and Ransford \cite[Theorem~1.3]{EFKR06} study the
\emph{unrestricted approach-zero set}
\[
 \Zapp(f)=\left\{\zeta\in\T:
       \liminf_{\substack{z\to\zeta\\z\in\DD}}|f(z)|=0\right\}.
\]
A  radial zero belongs to this set, which allows approach
from any direction in the disk. Hedenmalm and Shields, together with Richter and Sundberg, prove cyclicity when this unrestricted set is countable and \(f\) is outer \cite[Theorem~1.2]{EFKR06}.
Thus our counterexample has uncountably many unrestricted approach zeros. El-Fallah, Kellay, and Ransford go further. Suppose that
\(f\in\D\) is outer and that
\(\Zapp(f)\) is the union of a countable set and a perfect set \(F\).
For \(t>0\), let \(F_t\) be the set of points of \(\T\)
at angular distance less than \(t\) from \(F\), and let
\(|F_t|\) denote its arclength. Their theorem proves cyclicity under
the condition
\[
 \int_0^\delta
 \frac{|F_t|}{\bigl(t\log(1/t)\log\log(1/t)\bigr)^2}\,dt<\infty,
 \qquad 0<\delta<e^{-e}.
\]
The condition concerns the whole family of shrinking neighborhoods
of the perfect part. It therefore imposes an integral bound on the
arclengths of these neighborhoods, while allowing a countable part
to be treated separately. This gives a concrete class in which
outerness and sufficiently sparse boundary geometry do ensure
approximation of \(1\) by polynomial multiples.

El-Fallah, Kellay, and Ransford develop these
neighborhood conditions further in \cite{EFKR09}.
Write \(E=\Zapp(f)\), define \(E_t\)
using the same angular-distance convention, and again use \(|E_t|\)
for arclength. Their Theorem~1.1 \cite[p.~2197]{EFKR09} shows that an outer function
\(f\in\D\) is cyclic if, for some \(\mu>0\),
\[
 |E_t|=O(t^\mu)\quad(t\downarrow0),
 \qquad
 \int_0^\pi\frac{dt}{|E_t|}=\infty.
\]
The first hypothesis gives a power bound on the size of the small
neighborhoods; the second expresses further sparsity through the divergence of \(\int_0^\pi dt/|E_t|\). Their Corollary~1.2 \cite[p.~2197]{EFKR09} applies the resulting theory
to a geometric class of generalized Cantor sets. These are obtained
by middle removals leaving \(2^n\) equal arcs of length \(\ell_n\)
at stage \(n\), with
\(\sup_n\ell_{n+1}/\ell_n<1/2\). If \(|f|\) extends continuously
to the closed disk and its boundary zero set is such a set, cyclicity
is equivalent to outerness together with zero logarithmic capacity
of that set. The continuity of the modulus identifies the boundary
zeros with the approach-zero set, and the Cantor geometry makes the
neighborhood sizes accessible. This is a positive resolution of the
cyclicity criterion on a specified geometric class, with both the
regularity and the construction of the set contributing to the result.

A different positive question starts with a closed polar set and
allows the function to be chosen: can one construct a cyclic
Dirichlet function that vanishes on that set? Brown and Cohn
\cite[Theorem~B, p.~43]{BC85} produced continuous cyclic functions with this prescribed
boundary vanishing, as recalled by Elmadani and Labghail
\cite{EL19}. The strength of this construction is that the boundary
constraint is imposed first, and cyclicity survives the choice of
a function satisfying it. The statement is existential. It supplies
at least one function for the prescribed set, and its content does
not require every outer function vanishing there to behave in the
same way. In the approximation language of the opening example, a
whole closed polar set can be accommodated while suitable polynomial
multiples still approach \(1\).

Elmadani and Labghail extend the prescribed-zero question to
weighted Dirichlet spaces \cite[Theorem~1.1]{EL19}. In their setting
the derivative energy is weighted by the Poisson extension of a
finite positive boundary measure \(\mu\); the space, denoted
\(\D(\mu)\), has its corresponding capacity \(c_\mu\).
For every closed set \(E\) with \(c_\mu(E)=0\), they construct a
function that is cyclic in \(\D(\mu)\), is continuous on the closed
disk, and has closed-disk zero set exactly \(E\). Thus the
construction controls the zero set as well as cyclicity.
Their separate sufficient criterion \cite[Theorem~1.2]{EL19}
starts with an outer function in \(\D(\mu)\) that is continuous
on the closed disk. If \(E\) is its boundary zero set, a sufficient
condition for cyclicity is
\[
 \int_0^\delta c_\mu(E_t)\log(1/t)\,\frac{dt}{t}<\infty
\]
for a sufficiently small positive \(\delta\). Here the size of a
neighborhood is measured by the capacity associated with the
weighted energy. The two theorems answer different questions:
one constructs a function with a specified zero set, and the other
tests a given function using the capacities of neighborhoods of its
zero set. Their norm uses the value at zero together with the
weighted energy, so its capacity and extremal normalizations belong
to that setting. The norm \eqref{eq:full-norm} remains fixed for the
classical Dirichlet space in this paper.

There is also a longstanding zero-set problem inside the disk:
which sequences of points of \(\DD\) can be zeros of a nonzero
Dirichlet function? An interior zero already prevents that function
from being cyclic, so this question concerns the existence and
geometry of finite-energy functions with zeros, before any demand
of cyclicity is imposed. Carleson's paper on zeros of functions with
bounded Dirichlet integral \cite{Carleson52Zeros} is an early
contribution to this problem. It is distinct from his paper on
boundary uniqueness sets \cite{Carleson52}, which enters the
smooth-boundary discussion below. Shapiro and Shields
\cite{SS62} developed the interior-zero theory through Hilbert-space
constructions. B\'en\'eteau, Fleeman, Khavinson, Seco, and Sola
revisit that construction in \cite[Section~3]{BFKSS18}, explaining
how orthogonality and reproducing kernels produce distinguished
functions associated with prescribed interior zeros. This construction
connects zero sets with extremal functions, while keeping the
interior-zero problem separate from the radial boundary-zero
criterion for outer functions.

Later work gives several complementary ways to study those
interior sequences. Kellay and Mashreghi \cite{KM12} investigate
sequences accumulating on prescribed closed polar subsets of the
circle, with growth restrictions that depend on the chosen boundary
set. The boundary set here describes where interior zeros may
accumulate. Richter, Ross, and Sundberg \cite{RRS04} refine a result
of Nagel, Rudin, and Shapiro
\cite{NRS82} concerning zeros of functions with
finite Dirichlet integral, sharpening the earlier interior-zero
theory. Their work belongs to the same effort to determine how the
distribution of zeros is constrained by finite energy. In a recent
preprint, Borichev, El-Fallah, and Kellay \cite{BEK26} study an
extremal characterization of interior zero sequences through kernels
of measure-dependent harmonic Dirichlet spaces. This brings
extremal and kernel methods to the recognition of zero sequences.
These results relate the location and distribution of interior
zeros to finite-energy conditions and extremal problems. They
concern prescribed interior zeros in approximation questions,
including the finite-zero spaces mentioned above, and do not
alter the distinction between an accumulation set of interior
zeros and the radial zeros of a given outer function.

To build the boundary-defined spaces used in this paper, one needs
a nonzero function that vanishes on the chosen compact set. Smooth
outer functions vanishing on the compactum are particularly useful. A
\emph{Carleson set} is a nonempty compactum \(K\subset\T\)
for which the logarithm of reciprocal distance is integrable:
\[
 \int_{\T}\log\frac{2}{\dist(\zeta,K)}\,dm(\zeta)<\infty.
\]
The distance here is Euclidean; Section~\ref{subsec:03-ideal}
also gives the equivalent condition in terms of the lengths of the complementary arcs. The condition, associated with Carleson's boundary uniqueness work \cite{Carleson52}, forces arclength zero but permits positive logarithmic capacity. Taylor and Williams \cite[Theorem~3.3]{TW70} construct an outer function whose closed-disk zero set is exactly a given Carleson set and whose derivatives of every order extend continuously to the closed disk and vanish on that compactum. The class of analytic functions with these continuous derivatives is denoted \(A^\infty\); vanishing
of all derivatives on the set is often described as flatness there.
Wells records that the existence of such a flat outer function was proved independently by Novinger and by Taylor and Williams
\cite[Theorem~A, p.~346]{Wells70}.
See also Novinger's published account \cite[Theorem~4.3 and its proof, pp.~85--89]{Novinger71}.

Richter and Yilmaz discuss such functions together with the quasi-continuous boundary formulation in their study of Dirichlet
spaces \cite{RY19}. Their account is useful here because it brings together the boundary vanishing condition and a smooth outer function that satisfies it. In particular, they record the Korenblum formulation \cite{Kor72} used in
Theorem~\ref{input:flat-outer}. For a Carleson set \(K\),
write \(\D_K\) for the functions whose quasi-continuous traces
vanish quasi-everywhere on \(K\). A smooth outer function \(F_K\)
belongs to this space and satisfies \(F_K(0)\ne0\). Evaluation at
zero is therefore a nonzero functional on \(\D_K\), so the
projection of \(1\) onto \(\D_K\) is nonzero and can be
normalized to give its evaluation extremal. If \(K\) has positive
capacity, \(1\) does not belong to \(\D_K\), and the space is
proper. Thus a nonpolar Carleson set gives both a nontrivial
boundary condition and a nonzero extremal to study. The function \(F_K\) shows that evaluation at zero on \(\D_K\)
is nonzero; the choice of the canonical generator of \(\D_K\)
and the approximation of such generators are further steps. Sections~\ref{subsec:03-ideal} and~\ref{subsec:03-seed}
give these starting constructions in the precise norm used here.

Compact polar subsets of \(\T\) also occur as closed-disk zero
sets of multipliers. Write \(\operatorname{Mult}(\D)\) for the space of
analytic multipliers of \(\D\). Its norm is the operator norm
of multiplication, and its weak-* topology agrees, on bounded
sets, with pointwise convergence inside the disk
\cite[Section~2.1]{DH23}. Chalmoukis and
Hartz \cite{CH22} identify the compact polar subsets of \(\T\)
as exactly the compact sets that are \emph{totally null} for this
algebra. Here total nullity has a specific measure-theoretic
meaning: \(|\nu|(E)=0\) for every finite complex boundary measure
\(\nu\) whose polynomial integration functional
\(p\mapsto\int_{\T}p\,d\nu\) extends weak-* continuously to
\(\operatorname{Mult}(\D)\). Thus their theorem connects
logarithmic capacity with the boundary measures detected by the
operator algebra. Combining this connection with interpolation
theory, they show that every compact polar set is exactly the
closed-disk zero set of some cyclic multiplier that is a limit of
polynomials in multiplier norm. This strengthens the prescribed-zero
construction by giving the chosen function bounded action on the
whole Dirichlet space and approximation in that operator norm.
The assertion again constructs a function for each compact set;
its additional multiplier structure is part of what is obtained.

The interpolation theory used here is developed by Davidson and
Hartz \cite{DH23}. Their ambient algebra is the closure of
polynomials in multiplier norm. Since multiplier norm
controls the supremum norm, its elements are continuous on the
closed disk; the closure is taken in the multiplication-operator
norm, rather than the supremum norm defining the disk algebra.
Davidson and Hartz prove that every continuous function on a compact totally null set has an extension in this algebra with multiplier norm equal to the supremum norm of the given function and, if
the given function is not identically zero, with strictly smaller modulus off the set. They also
describe a multiplier-norm closed ideal using its weak-* closure
in \(\operatorname{Mult}(\D)\) together with the condition of
vanishing on its boundary zero set. The topology and ambient
algebra are part of these conclusions: such ideals are different objects from arbitrary Dirichlet-norm closed invariant subspaces.
Noell's survey \cite{Noell20} provides a broader guide to zero
sets, peak sets, and boundary interpolation. Prescribing where a
function is zero, requiring it to take its maximal modulus on a
set, and extending arbitrary continuous functions from that set ask
different questions, whose answers depend on the analytic class
and its boundary regularity. For the invariant subspaces considered
here, the next question is how the whole space can be represented
by one distinguished analytic generator.

\subsection{Canonical generators and polynomial approximation}\label{subsec:01-generation-history}

The distinguished function of a closed invariant subspace can be
chosen by a simple Hilbert-space construction. In the norm
\eqref{eq:full-norm}, the constant \(1\) represents evaluation at
zero: \(\langle h,1\rangle_{\D}=h(0)\). If \(M\) is a closed
subspace and \(P_M\) denotes orthogonal projection onto it, then
\(P_M1\) represents the same evaluation on \(M\). Assume that
\(P_M1\ne0\), or equivalently that some function in \(M\) is
nonzero at zero. Normalizing this projection gives
\[
 e=G_M=\frac{P_M1}{\|P_M1\|_{\D}},
 \qquad e(0)=\|P_M1\|_{\D}>0.
\]
The function \(e\) has unit norm and maximizes \(|h(0)|\) on the
unit ball of \(M\); the positive value at zero fixes its phase.
This is the canonical normalization used in the opening statement.
It also gives a concrete reproducing-kernel interpretation:
\(P_M1\) is the kernel vector for evaluation at zero in \(M\).
The formulation in Guo and Zhao \cite{GZ10} places such extremals
within the invariant-subspace theory for this norm. Kernel identities
and extremal functions are also central to Shimorin's work on
Dirichlet-type spaces \cite{Shimorin98}, where they provide a
related way to organize the geometry of distinguished analytic
functions. For the present construction, the next issue is whether
this one extremal vector represents the whole space through
polynomial multiplication.

That conclusion uses two structural theorems about the Dirichlet
shift. Recall that an invariant subspace \(M\) is one which is shift-invariant when \(zh\in M\)
whenever \(h\in M\). Its \emph{wandering space} consists of the
vectors in \(M\) orthogonal to the shifted subspace \(zM\), and
is written \(M\ominus zM\). Richter and Shields \cite[Theorem~2(c)]{RSh88} proved that
\(M\ominus zM\) is one-dimensional for every nonzero closed invariant subspace \(M\) of the Dirichlet shift. Richter \cite{Richter88} proved the further generation statement: every nonzero vector in that line generates all of \(M\) under polynomial multiplication. The
dimension statement identifies one direction left after the shift; the generation theorem says that this direction and its successive shifts recover the entire subspace. Guo and Zhao
\cite[Theorem~2.1]{GZ10} record the two conclusions separately in the norm used here, and El-Fallah, Kellay, and Ransford \cite{EFKR10} give an accessible account of their place in Dirichlet invariant-subspace theory. For the projected vector above, \(\langle zh,P_M1\rangle_{\D}=(zh)(0)=0\) for \(h\in M\).
Thus, when \(P_M1\ne0\), its normalized version lies in this one-dimensional wandering space, and Richter's theorem gives
\([e]=M\). This is the substantive reason that an extremal defined
by evaluation at one point can serve as a generator of an infinite-dimensional space.

Richter and Sundberg show that the Dirichlet norm of a function in this subspace equals a weighted norm of its quotient by the generator \cite{RS92}. To describe their formula, we now use a weighted space with the \emph{\(H^2\)
norm} as its first term. This differs from the value-at-zero term in the Elmadani--Labghail convention discussed above. Let \(w=|e^*|^2\), and let \(P[w]\) be its Poisson extension with
respect to the normalized arclength \(m\):
\[
 P[w](z)=\int_{\T}\frac{1-|z|^2}{|\zeta-z|^2}\,
                  w(\zeta)\,dm(\zeta).
\]
With the normalized area \(A_0\) fixed in the opening, define
\[
 \begin{gathered}
 D_w(F)=\int_{\DD}P[w](z)|F'(z)|^2\,dA_0(z),
 \qquad
 \mathcal H_w=\{F\in H^2:D_w(F)<\infty\},\\
 \|F\|_{\mathcal H_w}^2=\|F\|_{H^2}^2+D_w(F).
 \end{gathered}
\]
For a unit canonical generator \(e\), the product identity begins
with the polynomial formula
\[
 \|ep\|_{\D}^2=\|p\|_{H^2}^2+D_w(p)
             =\|p\|_{\mathcal H_w}^2
 \qquad(p\in\mathbb C[z]).
\]
The generator's boundary modulus therefore determines the weight
used to measure the derivative of the quotient. Richter and Sundberg's
work on the local Dirichlet integral \cite{RS91}, recalled in
\cite[p.~168]{RS92}, explains why this boundary modulus appears
in the weighted derivative energy. Their invariant-subspace
representation \cite[(3.1), p.~175]{RS92}
gives the displayed product norm in the stated convention. In particular, the exact equality
uses a unit generator; its normalization and the Hardy term
are both part of the identity.

The result is not restricted to polynomial multiples of the
generator. Richter and Sundberg
\cite[Theorem~3.2]{RS92} show that multiplication by the
unit generator is an onto isometry
\[
 \mathcal H_w\longrightarrow M=[e],\qquad F\longmapsto eF.
\]
Consequently, dividing any function in \(M\) by \(e\) gives a function
in the weighted space. The weighted norm of the quotient equals
the Dirichlet norm of the original function. Wu's treatment
\cite{Wu98} gives an operator-theoretic account of this weighted
range and explains the role of polynomial density in obtaining
the onto statement. This range theorem permits estimates for the
whole invariant subspace, including vectors that are limits of
polynomial multiples. In the outer canonical case used here,
\(e\) has no zero in the disk, so the quotient \(h/e\) is
analytic for every analytic \(h\). Analyticity alone does not
put that quotient in the domain of the representation: it must
belong to \(H^2\) and have finite \(D_w\)-energy. In the later
argument using a weighted boundary integral, these are precisely the
properties that must be proved before using the range theorem to conclude
that a candidate \(h\) belongs to \([e]\).

There is a wider operator-theoretic reason for the appearance of
weighted Dirichlet spaces. Richter's representation theory for
cyclic analytic two-isometries \cite{Richter91}, described in
Richter and Sundberg \cite{RS92} and in Wu's account
\cite{Wu98}, places multiplication by \(z\) in a class of shift
models whose norms are determined by boundary measures. The defining property of 
a \emph{two-isometry}, in the classical Dirichlet space, follows
directly from the coefficient weights \(n+1\). Cyclicity says
that one vector and its successive images span a dense subspace;
analyticity means that the ranges of all powers of the shift
have trivial intersection. Richter's model gives the
boundary measure a structural role in describing such an
operator. The weighted space \(\mathcal H_w\) above is a concrete instance of that representation for an invariant
subspace. Its weight varies with the chosen generator, while
the ambient classical Dirichlet norm in this paper stays fixed.

A complementary approach measures cyclicity at each polynomial
degree. Fix a function \(f\in\D\) with \(f(0)\ne0\), and, for
each \(n\ge0\), choose the polynomial \(p_n\) of degree at most
\(n\) minimizing
\[
 \min_{\substack{p\in\mathbb C[z]\\\deg p\le n}}
       \|1-pf\|_{\D}
       =\|1-p_nf\|_{\D}.
\]
The terminology \emph{optimal polynomial approximant} refers to
\(p_n\). The error is measured after multiplying by \(f\);
the residual is the different function \(1-p_nf\).
This formulation does not require \(1/f\) to belong to \(\D\).
It is a finite-dimensional projection problem: \(p_nf\) is
the projection of \(1\) onto the span of
\(f,zf,\ldots,z^nf\). These spaces increase with \(n\), and
their union is dense in \([f]\), so
\(\|1-p_nf\|_{\D}\) decreases to the distance from \(1\) to
\([f]\). In particular, the errors tend to zero exactly when
\(f\) is cyclic. This is the quantitative approximation question
developed in the work of Fricain, Mashreghi, and Seco
\cite{FMS14} and B\'en\'eteau, Condori, Liaw, Seco, and Sola
\cite{BCLSS15}. It retains the original cyclicity problem while
making the degree needed to achieve an approximation visible.

Fricain, Mashreghi, and Seco \cite{FMS14} connect this minimization
with the reproducing-kernel structure of coefficient Hilbert
spaces. Reciprocal coefficient weights enter the kernel, and
their sums give explicit optimal errors in the model calculation
for \(1-z\). Thus the harmonic sum in the opening example is
part of a general kernel calculation: its divergence provides
the vanishing error, and its size measures the quality of the
approximation at a fixed degree. B\'en\'eteau, Condori, Liaw,
Seco, and Sola \cite{BCLSS15} develop the extremal-polynomial
approach to cyclicity in Dirichlet-type spaces. They turn the
degree-by-degree minimization into constructive coefficient and
moment calculations, and obtain error estimates for specified
classes of functions. The common starting point is the
orthogonality of the residual to every permitted polynomial
multiple of \(f\); solving these finite systems gives explicit
approximants as well as a way to test cyclicity. Their work
therefore connects the qualitative question \([f]=\D\) with
concrete approximation problems. When numerical formulas are
compared, the norm remains important: parts of their constructive
analysis use an equivalent integral norm, whereas the computations
in this paper use \eqref{eq:full-norm}.

B\'en\'eteau, Khavinson, Sola, Liaw, and Seco \cite[Proposition~3.1 and Section~4]{BKLSS16}
develop the connection with orthogonal polynomials and reproducing
kernels. For a fixed \(f\), equip the polynomial space with the
inner product
\[
 (p,q)\longmapsto\langle pf,qf\rangle_{\D}.
\]
Orthogonalizing the polynomials in this inner product gives
orthonormal bases adapted to the subspaces of polynomial
multiples. The finite-dimensional reproducing kernels built from
these bases give formulas for the corresponding projections, so approximation
errors can be expressed through kernel quantities. This supplies
a second way to understand the same optimal approximants and
opens questions about their convergence and zeros. B\'en\'eteau
and Centner \cite{BC21} bring these connections together in
their survey, also explaining the earlier least-squares inverse
methods from digital filtering and their coincidence with optimal
polynomial approximants in the Hardy-space setting. That broader
history helps distinguish the approximation problem from any
particular technique for solving it. For the canonical generators
used here, the problem becomes especially simple: the earlier
constant optimal approximant suggests a direct link between the
extremal normalization and the entire sequence of finite-degree
approximation problems.

For a canonical generator, the finite-degree approximation problem
has the same solution at every degree. Let \(e\) have unit Dirichlet
norm, positive value \(e(0)=a\), and the shifted orthogonality
\(\langle z^je,e\rangle_{\D}=0\) for \(j\ge1\).
The residual \(1-ae\) is then orthogonal to \(e\) and to each
positive shift \(z^je\). It is consequently orthogonal to every
polynomial multiple of \(e\). Since the product \(ae\) is already
available at degree zero, the optimal approximating polynomial is
the constant \(a\) for every degree \(n\). B\'en\'eteau, Fleeman, Khavinson, Seco, and Sola explain this connection between
normalized projection extremals, shifted moments, and constant
optimal approximants in \cite[Theorems~2.2 and 2.3]{BFKSS18}.
Equations~\eqref{eq:certificate} and~\eqref{eq:polynomial-error}
give the projection identity and its orthogonal decomposition in
the present norm. Richter's generation theorem identifies the
space generated by the extremal; this later approximation
interpretation explains why its distance from \(1\) can be read
from the value at zero.

This constant case is also a particularly simple instance of
\emph{stabilization}. An optimal-approximant sequence stabilizes
when there is a degree \(N\) after which the polynomials no longer
change: \(p_n=p_N\) for all \(n\ge N\). Felder
\cite{Felder23} studies this question through projections and
also extends the approximation problem to other target functions.
For the fixed-function problem with target \(1\), stabilization
means that the projection onto \([f]\) is already realized by
a polynomial multiple of finite degree. The canonical case
realizes it at degree zero because of the shifted moments
\cite{BFKSS18}. Unit norm and a positive value at zero alone do
not supply those moments: the normalized \(1-z\) example in
Section~\ref{subsec:02-certificate} illustrates the distinction.
Each fixed shifted moment is continuous under strong Dirichlet
convergence, since multiplication by a fixed power of \(z\) is
bounded. Thus the moments can be retained in a limiting
construction. Once positivity at zero is retained,
the identity \eqref{eq:certificate} shows that keeping
the values at zero bounded away from one preserves
a positive approximation error.
This explains why the manuscript follows canonical functions
while controlling their limit; a stabilization result for one
fixed function is a different assertion from convergence of
extremals for changing subspaces.

The same orthogonality has a useful consequence for multiplication.
A \emph{Dirichlet-inner function} is a unit vector \(u\in\D\)
orthogonal to all its positive shifts:
\[
 \|u\|_{\D}=1,\qquad
 \langle z^ju,u\rangle_{\D}=0\quad(j\ge1).
\]
The term refers to this Dirichlet-space condition. The canonical
functions considered here have, in addition, positive value at
zero; a general Dirichlet-inner function may vanish there.
The multiplier norm of \(u\) is the supremum of \(\|ug\|_{\D}\)
over the unit ball of \(\D\), whenever multiplication is bounded.
Richter and Sundberg \cite{RS92} proved the forward
contractivity theorem: Dirichlet-inner functions satisfy
\(\|ug\|_{\D}\le\|g\|_{\D}\) for every \(g\in\D\).
Testing at \(g=1\) shows that their multiplier norm is exactly
one. Thus the multiplier norm controls multiplication on all
of \(\D\), whereas \(\|u\|_{\D}\) measures only the function itself.
In his account of this theorem, Seco \cite[p.~1655]{Seco19}
credits Aleman with another proof. Shimorin
\cite{Shimorin98} develops kernel and extremal-function identities
that give a further proof of Aleman's contractive-multiplier
theorem and permit a vector-valued extension. These contributions
provide several proofs of the contractive-multiplier theorem. Contractivity allows the norm to decrease; it is not
an assertion that multiplication preserves the norm of every vector in \( \D \).

Seco proves a converse in the classical Dirichlet space
\cite[Theorem~2.1, p.~1656]{Seco19}. If a function \(f\) has both Dirichlet norm one
and multiplier norm one, then it has the shifted orthogonality
required of a Dirichlet-inner function. His argument exposes
the moments using very simple multiplication tests. For a fixed
positive integer \(k\), apply contractivity to
\(z^k+\lambda\), where \(\lambda\) ranges over the complex
numbers. In the expanded squared norms, the terms involving
\(|\lambda|^2\) agree because \(\|f\|_{\D}=1\). The remaining
cross term depends on \(\langle z^kf,f\rangle_{\D}\).
Varying the size and argument of \(\lambda\) forces this moment
to vanish. Repeating the test for each \(k\) recovers the whole
orthogonality condition. Together with the earlier
Richter--Sundberg direction, Seco's result characterizes
Dirichlet-inner functions by the equality of the two norms to
one. The characterization concerns the fixed norm of the
classical Dirichlet space and a function already acting as a
bounded multiplier. It offers an operator-norm test for the same
objects that the projection and optimal-approximation descriptions
identify.

Le gives another operator formulation in weighted Hardy
spaces \cite{Le20}. In that setting an analytic series
\(\sum_{n\ge0}a_nz^n\) has squared norm
\(\sum_{n\ge0}\omega_n|a_n|^2\), where the weights are positive,
\(\omega_0=1\), and \(\omega_{n+1}/\omega_n\to1\).
The choice \(\omega_n=n+1\) is the exact Dirichlet norm used
here. Suppose that \(f\) is already a multiplier, write \(M_f\)
for multiplication by \(f\), and let \(M_f^*\) be its
Hilbert-space adjoint. Le's characterization is
\[
 f\text{ is inner in the weighted space}
 \quad\Longleftrightarrow\quad M_f^*M_f1=1.
\]
Here innerness again means unit norm and orthogonality to the
positive shifts. The identity tests the operator on the
constant vector: its matrix coefficients recover the norm of
\(f\) and its shifted moments. It does not assert
\(M_f^*M_f=I\) on the whole space, and boundedness of \(M_f\)
is a hypothesis. Le also describes inner functions through
normalized projections and develops finite constructions from
evaluation kernels. These descriptions extend the ways of
recognizing an extremal function and lead naturally to concrete
examples with prescribed interior zeros.

For those examples, choose a finite set \(A\) of distinct
nonzero points of the disk. Let \(k_\lambda\) be the evaluation
kernel at \(\lambda\), characterized by
\(\langle h,k_\lambda\rangle_{\D}=h(\lambda)\), and put
\[
 M_A=\{h\in\D:h(\lambda)=0\text{ for every }\lambda\in A\}.
\]
The orthogonal complement of \(M_A\) is the finite-dimensional
span of the kernels \(k_\lambda\), \(\lambda\in A\).
Projecting \(1\) away from this span gives \(P_{M_A}1\).
This projection is nonzero: the polynomial
\(\prod_{\lambda\in A}(z-\lambda)\) belongs to \(M_A\) and
does not vanish at the origin. Normalizing \(P_{M_A}1\) therefore gives the canonical function
\(e_A=G_{M_A}\). The Gram matrix
of the evaluation kernels makes the projection computable
through a finite linear system. Shapiro and Shields \cite{SS62} studied zeros of functions
with finite Dirichlet integral; B\'en\'eteau, Fleeman, Khavinson, Seco, and Sola give the normalized
finite-zero and Gram-matrix formulations used for this
comparison \cite[Theorems~3.4 and 3.5]{BFKSS18}.
The finite-zero theorem in Section~\ref{sec:06} takes an
additional step: it chooses such finite sets of points to approximate a specified canonical target while controlling the sum of the
negative logarithms of their moduli.

Felder and Le study analogues of finite
Blaschke products \cite{FL22}. They connect polynomial-generated
invariant subspaces with finite combinations of evaluation
kernels and, for repeated zeros, derivative-evaluation kernels. A derivative kernel here
means the vector representing the functional
\(h\mapsto h^{(j)}(\lambda)\); its existence requires that
this functional be bounded in the space under consideration.
Thus, the zeros that can be prescribed by this method, including
their multiplicities, depend on which evaluation and derivative-evaluation functionals are bounded. The resulting functions may also have zeros that were not prescribed. In the classical Dirichlet space, evaluation at interior points is bounded,
whereas point evaluation at a boundary point is not; the
finite kernel constraints here are therefore interior
constraints. Felder \cite[Theorem~5.1, p.~324]{Felder23} expresses the
projection of \(1\) onto the invariant subspace generated by
a polynomial nonzero at the origin in a weighted Hardy space.
The projection is \(1\) plus a finite linear combination of the
evaluation and derivative-evaluation kernels that exist at the
zeros of the polynomial, with derivative orders below the
corresponding multiplicities.
These results make concrete projections accessible and clarify
which zero conditions they satisfy, without asserting that
every inner function is a finite kernel combination.

It is useful at this point to keep three zero questions in
view. For a fixed function \(f\), one can study zeros of its
optimal approximating polynomial \(p_n\), or zeros of the
residual \(1-p_nf\). At a zero of \(p_n\), the residual equals
\(1\); at a residual zero, the product \(p_nf\) equals \(1\).
The two questions therefore track different parts of the
approximation. The points of \(A\) used above have a
third role: they impose vanishing on every member of \(M_A\)
and determine a canonical function \(e_A\) that changes with
the set. The manuscript chooses these points to control both
norm approximation and the sum of the
negative logarithms of their moduli; see Theorem~\ref{thm:finite-zero-approximation}. They are not automatically
zeros from either optimal-approximation sequence for one
fixed function. 

B\'en\'eteau, Khavinson, Liaw, Seco, and Simanek \cite[Section~2 and Theorem~4.2]{BKLSim19} study the
geometry of zeros of the optimal approximating polynomials.
One of their questions asks how close such a zero can lie to
the origin as the function being approximated varies. They
connect this extremal problem to nonlinear optimization and
the norms of Jacobi matrices, the tridiagonal matrices that
encode the relevant orthogonal-polynomial recurrences.
Their account also credits McDougall with an earlier
Jacobi-matrix connection. The paper develops Jentzsch-type results on the accumulation of zeros
of optimal approximants on the unit circle. For \(f\in H_\omega^2\)
with \(f(0)\ne0\), let \(W_f\) be the orthogonal projection of \(1\)
onto the closure in \(H_\omega^2\) of the polynomial multiples of
\(f\). The results require additional hypotheses on the weights and
assume that \(W_f/f\) has a singularity on the unit
circle~\cite[Theorems~6.1 and~6.2]{BKLSim19}.
Those results use quasi-integral spaces: the weighted norm
has an equivalent integral description against a regular
rotation-invariant measure on the closed disk. Within the
specified hypotheses, the boundary-accumulation conclusions
show how singular behavior is reflected in the zeros of
finite-degree approximants. This part of the theory belongs
to a neighboring class of weights. The classical Dirichlet
weights \(n+1\) do not give a quasi-integral norm of that
kind: their monomial norms grow without bound, while the
monomial norms for a finite measure on the closed disk stay
bounded. The extremal problem and the Jacobi-matrix methods also
bear on the study of approximant zeros. The accumulation
results, however, apply only under the hypotheses and
in the spaces specified in the cited works.

Acuaviva and Seco \cite[Definition~1.2 and Theorem~1.3]{AcuSeco22}
study the zeros of the residual polynomials. Their result
has a precise specialization to the classical Dirichlet norm.
Let \(f\) be a polynomial with simple zeros, none in the
open disk and at least one on the unit circle. Such a
polynomial is called critical in this setting. Let \(p_n\)
be its optimal approximants and set \(r_n=1-p_nf\).
The theorem gives a subsequence along which the zeros of
the residuals become equidistributed on the unit circle.
More precisely, put equal mass \(1/\deg r_n\) at each zero of \(r_n\),
counting multiplicity. Along some subsequence \(n_j\), these probability
measures converge weakly to the normalized arclength measure:
\[
 \frac{1}{\deg r_{n_j}}
 \sum_{r_{n_j}(\zeta)=0}\delta_\zeta
 \ \longrightarrow\ m .
\]
Thus averages of bounded continuous test functions over the
residual zeros tend to their averages around the circle.
The classical weights \(n+1\) satisfy the weight conditions
of their theorem; the simple-zero assumption and the
subsequence conclusion remain part of this specialization.
The result reveals a boundary distribution in the geometry
of the approximation error itself. The construction in this paper also tracks the convergence of outer canonical functions, so it requires control of their boundary logarithms as the functions converge. This leads to the
role of logarithms and their boundary behavior.

\subsection{Logarithms, cyclicity, and factorization}\label{subsec:01-logarithm-history}

Brown \cite[Section~3, pp.~421--422]{Brown90} proved that
\(f\) is cyclic in \(\D\) whenever both \(f\) and \(1/f\)
belong to \(\D\).
A related sufficient condition involves an analytic logarithm of \(f\).
If \(f\) has no zeros in the disk, an analytic branch
\(L=\log f\) satisfies \(e^L=f\). Its real part is
\(\log|f|\), while its imaginary part is a harmonic conjugate
of \(\log|f|\). Membership of \(L\) in a function space also
controls its analytic derivative \(L'=f'/f\).

B\'en\'eteau, Condori, Liaw, Seco, and Sola prove the following
sufficient condition in their work on Dirichlet-type spaces
\cite[Theorem~4.4, p.~281]{BCLSS15}:
if \(f\in\D\) and an analytic logarithm of \(f\) belongs
to \(\D\), then \(f\) is cyclic. Here the logarithmic
assumption also supplies outerness. Indeed, \(L\in H^2\),
so its real part is the Poisson extension of its boundary
real part; exponentiating this representation gives the
outer representation of \(f\). The additional Dirichlet
regularity of \(L\) then helps construct approximations
to \(1\). This criterion asks for control of an analytic
function in the Dirichlet norm, beyond the integrability
of the real boundary function \(\log|f^*|\).

The small-power argument makes that extra control concrete.
For the fixed branch \(L\), put \(f^\tau=e^{\tau L}\)
with \(0<\tau\le1\). Since
\(|f|^{2\tau}\le1+|f|^2\), differentiation gives the estimate
\[
 D(f^\tau)\le
 \tau^2\bigl(D(L)+D(f)\bigr),\qquad 0<\tau\le1.
\]
For every analytic \(h\) of finite Dirichlet energy, its Taylor
coefficients give \(\|h-h(0)\|_{H^2}^2\le D(h)\). Hence
\[
 \|f^\tau-1\|_{\D}^2
 \le |f^\tau(0)-1|^2+2D(f^\tau)\longrightarrow0
 \qquad(\tau\downarrow0).
\] B\'en\'eteau, Condori, Liaw, Seco, and Sola
combine this approximation with Richter and Sundberg's power
invariance \cite{BCLSS15,RS92}: for an outer function, a
positive power that belongs to \(\D\) generates the same
invariant subspace. Hence the small powers lie in \([f]\),
and their limit puts \(1\) in that space. This proof explains
why controlling \(f'/f\) is effective for cyclicity.
Aleman and Richter \cite[pp.~785 and~791]{AR26} trace an
earlier criterion to Aleman's 1993 habilitation thesis
\cite{Aleman93}. In their weighted Dirichlet setting, an outer
function is cyclic if its analytic logarithm belongs to the same
space. That antecedent and the small-power
explanation place the logarithmic condition within a
longer constructive approach to the problem.

More recent logarithmic criteria use the complete Pick
property. A reproducing kernel has this property
when positivity of the associated Pick matrices suffices
for contractive multiplier interpolation, including
matrix-valued interpolation of every finite size. The
normalization point here is \(0\), and the kernel takes
the value \(1\) when one variable is \(0\). For the exact
Dirichlet norm of this paper, the kernel is
\[
 k(z,\lambda)=
 \sum_{n=0}^{\infty}\frac{(z\overline\lambda)^n}{n+1}.
\]
This particular normalized kernel has the complete Pick
property; it is among the explicitly described kernels
in Aleman, Hartz, McCarthy, and Richter \cite{AHMR17}.
Shimorin's work \cite[Theorems~1.1 and~1.2]{Shimorin02} established the complete
Pick property for Dirichlet-type spaces with superharmonic
derivative weights, connecting weighted energy spaces with
multiplier interpolation and commutant lifting. Thus this
complete Pick property also holds for the reproducing kernel
of the norm fixed here, with its particular orthogonality and extremal vectors.

Luo and Rydhe \cite[Theorem~1.2]{LuoRydhe26} give a broader
operator criterion for the complete Pick property.
A nonzero reproducing kernel Hilbert space of analytic functions
on \(\DD\), invariant under both the forward shift \(M_zf=zf\)
and the backward shift \(Lf=(f-f(0))/z\), has a complete Pick
kernel whenever \(M_z\) is concave. In the norm
\eqref{eq:full-norm}, both shifts are bounded and
\[
 \|z^2f\|_{\D}^2-2\|zf\|_{\D}^2+\|f\|_{\D}^2=0
 \qquad(f\in\D).
\]
Thus the Dirichlet shift is a two-isometry, and their criterion
also yields the complete Pick property for the exact norm used here.

For a normalized complete Pick space \(\mathcal H\), its
multiplier Smirnov class \(N^+(\mathcal H)\) consists of
quotients \(b/a\), where \(a\) and \(b\) are multipliers
and multiplication by \(a\) has dense range \(a\mathcal H\)
in \(\mathcal H\) \cite{AHMR17}. This is the cyclic
denominator condition. In the general kernel setting,
cyclicity is described through multiplier actions; in
the classical Dirichlet space it agrees with the
polynomial cyclicity used here. For a multiplier
denominator, polynomial density already identifies the
closure of \(a\D\) with \([a]\). More generally, every
multiplier multiple \(\varphi f\) is a norm limit of
polynomial multiples of \(f\): the analytic Fej\'er means
of \(\varphi\), the arithmetic means of its Taylor partial
sums, are polynomial multipliers, uniformly
bounded in multiplier norm by averaging unitary
rotations, and convergence on polynomial test vectors
extends by density to strong convergence on every
vector of \(\D\). The convergence is strong operator
convergence, not convergence in multiplier norm.
A cyclic denominator is necessarily nonvanishing, but
nonvanishing alone is a weaker condition. Quotients
with merely nonvanishing multiplier denominators form
the Pick--Nevanlinna class \(N(\mathcal H)\);
the cyclic-denominator requirement is essential in the
logarithmic criteria below \cite{AR26}.

Aleman, Hartz, McCarthy, and Richter prove that every function in a normalized
complete Pick space belongs to this multiplier Smirnov
class \cite{AHMR17}. Their representation has additional
structure. If \(\|f\|_{\mathcal H}\le1\), they find
contractive multipliers \(\varphi,\psi\), with
\(\psi(0)=0\), such that
\[
 f=\frac{\varphi}{1-\psi},
 \qquad 1-\psi\text{ is a cyclic multiplier}.
\]
The theorem expresses an arbitrary Hilbert-space
function through functions that act boundedly on
the whole space. The density condition on the
denominator makes this more than a formal pointwise
division. It is the inclusion
\(\mathcal H\subset N^+(\mathcal H)\) that later allows
Dirichlet membership to imply a Smirnov-class
hypothesis. In this quotient theorem the denominator is chosen to suit \(f\);
this is different from our setting, where the denominator is the
fixed canonical function \(e\).

Aleman, Hartz, McCarthy, and Richter prove a further
factorization theorem when a reproducing kernel
contains a complete Pick factor \cite{AHMR18}.
Suppose \(k=sq\), where \(s\) is a normalized complete
Pick kernel and \(q\) is a positive definite kernel.
The factorization can then involve multipliers from
the space \(\mathcal H_s\) into the different space
\(\mathcal H_k\). For a unit-ball function \(F\) in
\(\mathcal H_k\), their representation has the form
\(F=\Phi/(1-\psi)\), with \(\psi\) a multiplier of
\(\mathcal H_s\), \(\psi(0)=0\), and \(\Phi\) a
multiplier from \(\mathcal H_s\) to \(\mathcal H_k\).
These factors satisfy the joint estimate
\[
 \|\psi h\|_{\mathcal H_s}^{\,2}
 +\|\Phi h\|_{\mathcal H_k}^{\,2}
 \le \|h\|_{\mathcal H_s}^{\,2}
 \qquad(h\in\mathcal H_s).
\]
Their construction uses Sarason functions, whose real
parts give harmonic or pluriharmonic majorants in the
analytic examples, to prove the multiplier estimates. This turns a factorization
of kernels into a controlled factorization of
functions. The spaces on the two sides of
multiplication are part of the conclusion; the
numerator is not simply being declared a multiplier
of \(\mathcal H_k\).

An earlier quotient representation appears in Aleman's work on Hilbert spaces between the Hardy and
Dirichlet spaces \cite{Aleman92}. He considers radial
derivative weights that are positive, integrable,
twice continuously differentiable, decreasing and
concave, and tend to zero at the boundary. For
these spaces he represents each function as a
quotient of two bounded functions belonging to
the same space. The construction uses an outer
factor with boundary modulus
\(\min(1,1/|f^*|)\): multiplying \(f\) by this factor
truncates its large boundary values, while the
weighted estimates keep the resulting numerator
and denominator in the space. This provides bounded
vectors inside invariant subspaces and supports
the intersection and codimension questions
developed from Richter and Shields' work.
Aleman's norm in that paper uses the value at zero
and weighted area energy.  Its bounded factors are
functions in the weighted space; that conclusion
does not identify them as multipliers or supply
the cyclic multiplier denominator of the later
complete Pick theorem.

Multiplier interpolation on the ball supplies
another antecedent. Alpay, Bolotnikov, and Kaptano\u{g}lu
\cite[Sections~2--3 and Theorem~10.3]{ABK02} develop Schur-algorithm methods using
reproducing kernels on the unit ball. Their work
organizes interpolation through the relation
between kernel positivity and contractive
multipliers on the Drury--Arveson space. Aleman, Hartz, McCarthy, and Richter
explicitly identify the finite-dimensional
Drury--Arveson result of Alpay, Bolotnikov, and
Kaptano\u{g}lu as an antecedent of their general
multiplier-quotient inclusion \cite{AHMR17}.
The historical connection is therefore between
a concrete interpolation setting on the ball
and the later theorem for normalized complete
Pick spaces. It gives the quotient representation
a place in interpolation theory as well as in
factorization theory.

Jury and Martin obtain Smirnov quotients by a different route
\cite[Theorem~1.4 and Corollary~3.4]{JM21}. They lift functions to the free
Fock setting, where power-series terms are
indexed by words in noncommuting variables,
and use factorization there before returning
to complete Pick spaces. Their quotient
\(h=b/a\) retains a cyclic multiplier
denominator and can be chosen with the
additional property \(1/a\in\mathcal H\).
This gives Hilbert-space control of the
reciprocal denominator itself. Jury and Martin also obtain factorizations of nonzero functions
into a contractive multiplier and a cyclic function of the same
norm as the original function \cite[Theorem~1.1]{JM19}. 
Aleman, Hartz, McCarthy, and Richter
\cite{AHMR23} expressly acknowledge this
earlier Jury--Martin work while developing
intrinsic characterizations and uniqueness
for the factors. The lift is thus useful
both as a method of constructing factors
and as a reason to distinguish the resulting
notions from classical Hardy factorization.

The free-outer factorization uses a different extremal
problem from the canonical projection onto a fixed
invariant subspace. In the normalized complete Pick   
space \(\mathcal H\), associate to \(f\) the functional
\[
 \Lambda_f(\varphi)=\langle\varphi f,f\rangle_{\mathcal H},
 \qquad \varphi\in\operatorname{Mult}(\mathcal H).
\]
Compare \(f\) with all \(g\in\mathcal H\) satisfying
\(\Lambda_g=\Lambda_f\). Testing with the
constant multiplier \(1\) shows that these functions
have the same norm. A nonzero function is \emph{free outer}
when it maximizes the modulus of its value at \(0\)
within this comparison class. A \emph{subinner}
multiplier has multiplier norm one and attains that
norm on some nonzero vector. Aleman, Hartz, McCarthy,
and Richter use these intrinsic descriptions in
\cite{AHMR23}. The functional \(\Lambda_f\) is fixed in the free-outer
extremal problem, whereas
the canonical function maximizes evaluation on
the unit ball of one specified invariant subspace.
The two extremal questions therefore need not
choose the same function.

Aleman, Hartz, McCarthy,
and Richter prove that every nonzero
\(f\in\mathcal H\) has a unique normalized
subinner/free-outer factorization \cite[Theorems~1.4 and~6.3]{AHMR23}:
\[
 f=\varphi h,\qquad h(0)>0,\qquad
 \|\varphi h\|_{\mathcal H}=\|h\|_{\mathcal H},
\]
where \(h\) is free outer and \(\varphi\) is
subinner. The free outer factor is cyclic.
Its norm equals the norm of \(f\); it is not
being normalized to one. Positivity at \(0\)
fixes the phase, and the uniqueness identifies
the factors independently of a particular
construction. Building on the existence result of
Jury and Martin \cite[Theorem~1.1]{JM19}, Aleman, Hartz,
McCarthy, and Richter establish intrinsic
characterizations and normalized uniqueness
for this factorization in complete Pick spaces. The subinner factor
is characterized by its action on a vector,
and is not thereby the canonical evaluation
extremal of an invariant subspace.

The three notions can now be compared through their
definitions. Hardy outerness concerns
boundary factorization in \(H^2\); free
outerness concerns an extremal problem with \(\Lambda_f\) fixed
in a normalized complete Pick space; a
Dirichlet canonical function concerns
evaluation on one invariant subspace.
Cheng, Mashreghi, and Ross \cite{CMR19}, and Le
\cite{Le20}, place the shift-orthogonality
notions in a broader reproducing-kernel
setting. Their vocabulary also helps keep
operator innerness distinct from a Hardy
boundary-modulus condition. Aleman, Hartz, McCarthy, and Richter prove that some Dirichlet
functions admit no factorization into a Dirichlet extremal times
a cyclic Dirichlet function \cite[Theorem~10.1]{AHMR23}. Their proof
uses the extremal function for the invariant subspace associated with a one-atom singular inner function.
This obstruction concerns the choice of the extremal factor in a factorization of a given function. Each row below records a different defining problem and its role.

\begin{center}
\small
\begin{tabular}{@{}L{24mm}L{62mm}L{61mm}@{}}
\toprule
Notion & Defining condition & Role\\
\midrule
Hardy outer &
The \(H^2\) outer representation determined by the boundary
logarithm, up to phase. &
Reconstruction from the boundary modulus and absence of a
nonconstant Hardy inner factor.\\
\addlinespace
Free outer &
Maximal \(|h(0)|\) among functions with the same
functional \(\Lambda_h\), in a normalized complete Pick space. &
The cyclic factor in the normalized subinner/free-outer
factorization.\\
\addlinespace
Dirichlet canonical &
A unit vector maximizing evaluation at \(0\) on a fixed
invariant subspace, with positive value there. &
A distinguished generator in the Dirichlet norm,
with its approximation error determined by the value at zero.\\
\bottomrule
\end{tabular}
\end{center}

Aleman and Richter use
Smirnov-class logarithms to obtain further
cyclicity criteria \cite[Theorem~1.1]{AR26}. In the
classical Dirichlet space, their sufficient
condition says that an outer \(g\in\D\)
is cyclic whenever an analytic logarithm
of \(g\) belongs to \(N^+(\D)\).
The outer hypothesis is part of this
statement. The conclusion is density
of the generated subspace, not merely
the existence of a quotient representation.
Since \(\D\subset N^+(\D)\), this
allows the analytic logarithm to be
represented through a cyclic multiplier
denominator instead of requiring its
own finite Dirichlet energy. It extends
the earlier logarithmic sufficient
conditions within the same exact
classical norm. Replacing the cyclic
denominator by a merely nonvanishing
one would change the function class
and would not give this criterion.

Iterating logarithms provides further
tests in Aleman and Richter's account \cite[Corollaries~1.2 and~1.3]{AR26}. Let \(g\in\D\)
be bounded and outer, with
\(\|g\|_\infty\le1\), and choose
an analytic \(L=\log(1/g)\).
Then \(\operatorname{Re}L\ge0\).
Use their indexing
\[
 G_0(w)=w,\qquad
 G_{n+1}(w)=\log(1+G_n(w)),
\]
where each subsequent logarithm is
the principal branch. For
\(\operatorname{Re}w\ge0\), the
quantity \(1+G_n(w)\) stays in the
right half-plane, so these branches
are well defined; their real parts
also remain nonnegative. If
\(G_n(L)\in N^+(\D)\) for some
\(n\ge1\), then \(g\) is cyclic.
In particular, membership of one
such iterate in \(\D\) is a
concrete sufficient condition.
Each new logarithm contributes
a denominator when differentiated,
which explains why iteration can
make derivative-energy estimates
easier even when the original
analytic logarithm has insufficient
regularity.

Aleman and Richter also identify
additional conditions under which
these tests are necessary
\cite[Theorem~1.4]{AR26}. Continue to assume
that \(g\in\D\) is bounded and
outer with \(\|g\|_\infty\le1\).
It has \emph{bounded argument}
when the imaginary part of a
chosen analytic logarithm is
bounded on the disk. Alternatively,
their classical result applies
if, for some integer \(\ell\ge1\),
\[
 \int_{\DD}|g'(z)|^2
 \left|G_\ell\!\left(\frac{1}{1-|z|^2}\right)\right|^2
 \,dA_0(z)<\infty.
\]
Under either of these additional
hypotheses, cyclicity of \(g\),
membership of \(\log g\) in
\(N^+(\D)\), membership of some
\(G_n(\log(1/g))\), \(n\ge1\),
in \(N^+(\D)\), and membership
of every such iterate in that
class are equivalent. These are
equivalences in the cyclic-denominator
Smirnov class. Dirichlet membership
of an iterate remains a useful
stronger sufficient condition,
and is not being substituted for
\(N^+(\D)\) in the converse.

Aleman, Perfekt, Richter, Sundberg, and Sunkes
\cite[Theorem~1.1, p.~39]{APRSS26} study radially
weighted Besov spaces on the ball
with complete Pick kernels. For
zero-free functions with bounded
argument in their setting,
cyclicity is equivalent to
Smirnov-class membership of the
analytic logarithm. For functions
also bounded by one, they obtain
further equivalences involving
iterates of the logarithm. Their
earlier work \cite{APRSS24}
investigates cyclicity through
classes of multipliers, comparisons
based on bounded ratios, stable
polynomials, and boundary-zero
geometry. The dimension and
weighted norm remain part of
each result; the boundary geometry
of these ball spaces is distinct
from the logarithmic capacity on
the circle used here.

Datta and Richter continue this direction
by distinguishing norm cyclicity
from weak-* sequential cyclicity
of multipliers in radially
weighted Besov spaces
\cite{DR26}. The first concerns
density in the norm of the
function space; the second
concerns sequential approximation
by multiplier products in the
weak-* topology of the multiplier
algebra. Their work compares
these notions in complete Pick
settings and beyond, and develops
logarithmic sufficient conditions.
The distinction makes the topology
of approximation an explicit part
of the question. It is especially
useful when placing the present
argument in context: convergence
in the Dirichlet norm, convergence
of multiplication operators, and
convergence in a multiplier-algebra
weak-* topology are distinct notions of convergence.

For the limiting argument of this
paper, the logarithm has another
role. The analytic-log criteria
of B\'en\'eteau, Condori, Liaw,
Seco, and Sola \cite{BCLSS15}
and of Aleman and Richter
\cite{AR26} test membership of
an analytic logarithm in a
function class. Here we also
follow the real functions
\(\log|f^*|\) in \(L^1(m)\),
because these boundary logarithms
reconstruct an outer function
through \eqref{eq:outer-intro}.
For an outer canonical function, the negative boundary logarithm is
nonnegative and its integral is determined by the value at zero.
For a fixed outer canonical function \(f\),
Proposition~\ref{prop:outer-continuity} combines this identity with an
estimate of the integral of \(-\log|f^*|\) where \(|f^*|\) is small.
This gives an \(L^1(m)\) bound on the difference between the boundary
logarithms of \(f\) and of an outer canonical function \(g\) sufficiently
close to \(f\) in the Dirichlet norm. The proof uses Scheff\'e's
positive-part argument \cite{Scheffe47}.

For the fixed singular target \(G_\mu\), let \(H_{\mu,\ell}\) be the
outer canonical approximants described below. The integral of
\(-\log|H_{\mu,\ell}^*|\) tends to the integral of
\(-\log|G_\mu^*|\) plus \(\mu(\T)\).
Proposition~\ref{prop:finite-recovery-logarithmic-defect} also shows
that the \(L^1(m)\) distance between these logarithms tends to
\(\mu(\T)\).

We also need the limit to be outer. Theorem~\ref{thm:B}, the second approximation theorem explained below, shows
that a Dirichlet-norm limit of outer canonical functions can have
a nonconstant inner factor. To keep track of outerness, we use the
following distance between outer canonical functions:
\[
 \dlog(f,g)=\|f-g\|_{\D}
   +\|\log|f^*|-\log|g^*|\|_{L^1(m)} .
\]
The second term measures the integrated absolute difference of the
boundary logarithms. If a sequence is Cauchy for this distance, the
limiting boundary logarithm reconstructs the outer function by
\eqref{eq:outer-intro}. Convergence in the first term keeps the norm
and the orthogonality identities that characterize canonical
functions. This is the completeness result of
Section~\ref{subsec:04-complete}.

\subsection{Stability, boundary approximation, and changing subspaces}\label{subsec:01-stability-history}

When the generating function varies, convergence in the ambient norm is
only one part of the approximation problem. Aguilera and Seco study how cyclicity can be preserved under such variation
\cite[Theorem~3.1, pp.~509--510]{AS24}. Their sufficient criterion in the coefficient-weight spaces
\(\D_\alpha\), with squared norm \(\sum_{j\ge0}(j+1)^\alpha|a_j|^2\),
keeps explicit control of the quotients. Suppose cyclic functions \(f_n\)
converge to \(f\) in that norm, each \(f/f_n\) is a multiplier, their
multiplier norms are uniformly bounded, and
\(\|1-f/f_n\|_{\D_\alpha}\to0\). Then \(f\) is cyclic. The case
\(\alpha=1\) uses exactly the norm of this paper. This follows from the decomposition
\[
 1-pf=\left(1-\frac{f}{f_n}\right)
       +\frac{f}{f_n}(1-pf_n).
\]
We first choose \(n\) so that \(\|1-f/f_n\|_{\D_\alpha}\) is small.
For this fixed \(n\), cyclicity of \(f_n\) lets us choose a polynomial
\(p\) so that \(\|1-pf_n\|_{\D_\alpha}\) is small enough for the
multiplier bound of \(f/f_n\) to control the second term.
Thus \(\|1-pf\|_{\D_\alpha}\) can be made arbitrarily small.

Aguilera and Seco also examine a finite-dimensional stability question
\cite[Theorem~4.1 and Corollary~4.2, p.~513]{AS24}. Fix \(\alpha\) and a degree \(d\), and let nonzero functions
\(f_n\) converge in \(\D_\alpha\) to a nonzero \(f\). The degree-\(d\)
optimal approximants then converge to the degree-\(d\) approximant for
\(f\). Here no condition \(f(0)\ne0\) is needed. At a fixed degree, the
Gram matrix formed from \(f,zf,\ldots,z^df\) and the corresponding
right-hand side vary continuously with \(f\). The matrix is invertible
when \(f\ne0\), because these polynomial multiples are linearly
independent; continuity of matrix inversion therefore gives continuity of
the coefficients and of the approximating polynomial in \(\D_\alpha\). The
local estimates can depend on \(d,\alpha\), and \(f\). Aguilera and Seco
also exhibit normalized functions for which inverse-Gram bounds
deteriorate already at degree one
\cite[Proposition~4.5, p.~516]{AS24}. Their estimates describe the
conditioning of a fixed finite system. Projection onto the entire
generated subspace in Theorem~\ref{thm:B} involves every degree:
Section~\ref{sec:12} proves approximation of every target vector
and containment of all subsequential weak limits.

For \(\alpha=0\), Aguilera and Seco give explicit cyclic approximants
\(g_n\) in \(H^2\) to the singular inner function of a unit point
mass at \(1\) \cite[proof of Theorem~3.1, pp.~510--511]{AS24}.
Put \(r=n/(n+1)\). For \(\zeta\in\T\), their formula gives
\[
 -\log|g_n^*(\zeta)|=\frac{1+r}{2}P_r(\zeta)+\frac{r-1}{2},
 \qquad P_r(\zeta)=\frac{1-r^2}{|\zeta-r|^2}.
\]
As \(n\to\infty\), \(r\uparrow1\), so \(P_r\,dm\) converges weak-* to
\(\delta_1\). The coefficient of \(P_r\) tends to \(1\) and the
constant term tends to \(0\). Hence \((-\log|g_n^*|)\,dm\) also
converges weak-* to \(\delta_1\). Examples using Cantor sets
and capacity in Dirichlet-type spaces can be found in
\cite{EFKR09,EFKR10Cantor}.

The present construction proceeds in two approximation steps. Approximation~A selects a subspace defined by a finite atomic singular inner function whose canonical generator approximates a prescribed proper outer canonical function,
while the total atomic mass can be made arbitrarily small.
Approximation~B recovers the canonical generator of this fixed subspace by outer canonical functions of trace-zero
spaces. The construction using Cantor copies on equal subarcs in Section~\ref{sec:10} provides the quantitative trace and
logarithmic-mass estimates needed for the geometric
construction in Approximation~B. For a single atom,
Corollary~\ref{cor:R} applies the criterion in Theorem~\ref{thm:C} to
obtain strong Dirichlet convergence. The measures with density given
by the negative boundary logarithms converge weak-* to the sum of
the measure determined by the outer factor and the prescribed point
mass. The limiting invariant subspace and the canonical value at
zero identify the mass. Section~\ref{sec:12} combines the single-atom
sequences and proves the corresponding limit for the finite atomic
measure in Theorem~\ref{thm:B}.

A different approximation problem combines norm control with uniform
approximation of an independently prescribed function on a boundary set,
while avoiding zeros. B\'en\'eteau, Ivrii,
Manolaki, and Seco \cite[Theorem~4.2]{BIMS20} prove a simultaneous
approximation theorem in the Dirichlet space. Let \(E\subset\T\) be closed
with logarithmic capacity zero, let \(h\in\D\) be nonvanishing in the
disk, and let \(\eta\in C(E)\). One can approximate \(h\) in the Dirichlet norm and \(\eta\) uniformly on \(E\) at
the same time by polynomials having no zeros on the closed disk. The
boundary target is independent of the interior norm target. Their
construction uses boundary cutoffs with small energy and exponentials: the
small capacity permits localized changes with small error in the Dirichlet norm, while
the exponential construction retains nonvanishing. This gives a positive
way to combine three requirements that need not be supplied by ordinary
polynomial density: control in the Hilbert-space norm, control on a
specified exceptional boundary set, and absence of zeros.

B\'en\'eteau, Ivrii, Manolaki, and Seco use this simultaneous construction
in a Baire-category argument for universal optimal approximants
\cite[Theorem~4.1]{BIMS20}. For the fixed closed polar set \(E\), universality means
that every continuous function on \(E\) is a uniform limit there of some
subsequence of the optimal approximating polynomials. The subsequence can
depend on the chosen boundary target. The genericity statement is  made in the complete metric space consisting of the zero-free Dirichlet
functions together with the zero function, with the Dirichlet metric. The
universal functions form a dense \(G_\delta\) subset of this space and are
nonzero. A residual subset is one containing a countable intersection of
dense open subsets of that space. Here \emph{residual} is a topological
term, distinct from the approximation residual \(1-p_nf\). The theorem
shows how rich the boundary behavior of optimal approximants can be on a
polar compact set, while fixing both the ambient metric space and the meaning
of a universal subsequence.

B\'en\'eteau, Manolaki, and Seco
\cite[Corollary~1.2 and Remark~1.3, p.~160]{BMS21} compute the
coefficients of \(1-p_n f\) for a polynomial \(f\) with simple
nonzero zeros using a matrix of truncated reproducing kernels
evaluated at those zeros. The matrix depends on the zero set and
the approximation degree. When all the zeros are simple and lie on
the unit circle, their theorem \cite[Theorem~6.1, pp.~175--176]{BMS21} gives
uniform convergence of \(p_n\) to \(1/f\) on compact subsets of the
closed disk that avoid the zeros. The theorem assumes that the
weights are monotone and that \(\sum_{n\ge0}\omega_n^{-1}=\infty\).
It also assumes that \(C^{-1}\omega_n\le\omega_{n+t}\le C\omega_n\)
for some \(C>0\), every \(n\ge1\), and every integer \(0\le t\le n+1\). The classical Dirichlet
weights \(\omega_j=j+1\) satisfy these assumptions. For these weights
the residual \(1-p_n f\) is \(O(1/\log(n+2))\) uniformly on each
such compact set. Their Hardy and Bergman results
\cite[Theorem~1.8, p.~162]{BMS21} also allow simple zeros outside the
closed disk. For the classical Dirichlet weights, Acuaviva and Seco
\cite[Definition~1.2 and Theorem~1.5]{AcuSeco22} obtain the same
uniform convergence of \(p_n\) to \(1/f\) on compact subsets of the
closed disk away from the zeros when \(f\) has simple zeros and none
in the open disk. At a zero of \(f\), the residual \(1-p_n f\)
equals \(1\) for every \(n\).

For closed subspaces of one fixed Hilbert space, convergence of
projections has an abstract description that separates two tasks. Strong
convergence means convergence in norm; weak convergence means convergence
of the inner products against every fixed test vector. Let \(M_n\) and
\(M\) be closed subspaces of that Hilbert space. The first task is to
recover every vector of \(M\) in norm from the spaces \(M_n\). The second
is to ensure that no weak limit from those spaces lies outside \(M\). In
precise form the conditions are
\[
 \begin{gathered}
 x\in M\ \Longrightarrow\
   \text{there exist }x_n\in M_n
   \text{ with }x_n\to x\text{ in norm},\\
 n_j\uparrow\infty,\quad x_j\in M_{n_j},\quad
 x_j\rightharpoonup x
 \ \Longrightarrow\ x\in M .
 \end{gathered}
\]
The second condition applies to every increasing subsequence of
indices. The first requires a vector \(x_n\in M_n\) for every \(n\).
 These are the
Mosco convergence conditions, originating in Mosco's work on
convex sets and variational inequalities
\cite[Definition~1.1, pp.~519--520]{Mosco69}. Chandler-Wilde, Hewett, Moiola, and Besson give a
useful Hilbert-space formulation and a more general variational
convergence result \cite[Definition~2.2 and Lemmas~2.5 and 2.7]{CHMB21}.
Its inner-product specialization identifies these conditions with
\(P_{M_n}h\to P_Mh\) in norm for every fixed \(h\). This is strong
convergence of the projections, without an assertion of convergence in
operator norm. Kato's monograph
\cite[Chapter~VI, Sections~1--3, and Chapter~VIII, Section~3]{Kato1995}
provides broader
background on forms and operator perturbation within which such
convergence questions sit.

In the unweighted Bergman space, Shimorin proved that every
closed \(z\)-invariant subspace \(M\) with
\(\dim(M\ominus zM)=1\) can be approximated by subspaces
defined by finitely many interior zeros, with multiplicities,
in the sense of strong convergence of their orthogonal projections
\cite[Theorems~1, 1A, and 1B and Lemma~2]{Shimorin00}.
The corresponding unit extremal functions can be chosen with
suitable phases so that they converge in the Bergman norm to a
unit extremal function for \(M\).
His Dirichlet-space results obtained through the Cauchy pairing
concern the backward shift \(f\mapsto(f-f(0))/z\)
\cite[Theorem~2]{Shimorin00}.
Corollary~\ref{cor:finite-codimensional-approximation} concerns
multiplication by \(z\) on the classical Dirichlet space.

Luo and Richter study approximation questions for
Dirichlet invariant subspaces using Hankel operators
\cite{LR15}. They describe invariant spaces as kernels of operators,
which leads to the question of how a space can be recovered from
simpler zero conditions. Their Question~6.2 asks whether the
orthogonal projection onto every nonzero closed shift-invariant
subspace is a strong operator limit of projections onto
finite-codimensional invariant subspaces. Question~6.3 asks
for locally uniform approximation of the extremal for a nonzero
trace-zero space associated with a compact set of positive capacity
by extremals for spaces defined by interior zero sets.
Corollary~\ref{cor:trace-zero-finite-approx} gives
Dirichlet-norm convergence for the trace-zero extremals in
Question~6.3, and Corollary~\ref{cor:finite-codimensional-approximation}
answers Question~6.2.

Yilmaz studies finite-codimensional approximation in a more regular
member of the coefficient-weight scale \cite{Yilmaz17,Yilmaz18}.
The space \(\D_2\) can be described as the analytic functions whose
derivatives belong to \(H^2\). 
In that space, he proves that the orthogonal projection onto every nonzero closed invariant subspace
is a strong operator limit of orthogonal projections onto
finite-codimensional invariant subspaces
\cite[Theorems~1.1 and 3.5]{Yilmaz18}.
For the classical Dirichlet space, let \(E\) be a Carleson set
of positive logarithmic capacity and let \(f_E\) be Korenblum's
corresponding outer function. Yilmaz constructs a Blaschke
sequence for which the closed linear span in \(\D\) of the
subspaces of functions vanishing on its successive tails
contains \([f_E]\) and is contained in \(\D_E\)
\cite[Theorems~1.3 and~4.2]{Yilmaz18}.
The paper explicitly notes that identifying \([f_E]\) with
\(\D_E\) would yield the corresponding finite-codimensional
approximation. The 2017 thesis develops this partial result and
also discusses tangential regularity of extremals
\cite{Yilmaz17}. 

Yilmaz's later work concerns a fixed harmonically weighted Dirichlet space
whose defining boundary measure is a finite positive sum of point masses
\cite{Yilmaz23}. The measure supplies the harmonic weight in the
derivative energy and remains fixed throughout the approximation. In that
setting, Theorems~1.1 and~3.3 of \cite{Yilmaz23} state
that every nontrivial invariant subspace admits finite-codimensional
approximation, with strong convergence of the associated orthogonal
projections. A separate theorem treats trace-zero spaces using spaces defined
by finite sets of interior zeros. The distinction between these two
conclusions matters: finite codimension is the general conclusion, while
the more concrete construction using finitely many prescribed zeros belongs
to the trace-zero theorem. Strong convergence again concerns each fixed
vector in one fixed Hilbert space. These statements concern a fixed harmonically weighted Dirichlet
space with a finite-atomic defining measure; they do not assert the
corresponding conclusions for arbitrary weights.

Richter and Yilmaz obtain a generator in the disk algebra for a
trace-zero space associated with a Carleson set
\cite[Theorem~1.3]{RY19}. The disk algebra consists of functions
analytic in the disk and continuous on the closed disk. They obtain
their generator by multiplying the canonical extremal by a cyclic
function. They also prove that the limit of the modulus of a
Dirichlet extremal in a polynomially tangential approach region
equals its unrestricted upper limit at the boundary point
\cite[Theorem~1.4]{RY19}. For the compact sets constructed here,
Section~\ref{subsec:11-realization} finds a sequence in such a region
along which the canonical extremal tends to zero at each point of
the compact set. The resulting canonical function is continuous on
the closed disk, and its zero set there is exactly that compact set.

Richter and Yilmaz also make explicit use of Richter and Sundberg's work on extremals and pseudocontinuations \cite{RS94,RY19}.
Pseudocontinuation relates interior and exterior representations through
matching boundary values; it need not give one analytic function across
the boundary. The Richter--Yilmaz argument does, however, give analytic
extension across appropriate complementary arcs of the invariant space's
boundary zero set. This is an individual regularity
statement. When the compact set and the extremal vary, one must also find
a common neighborhood across the arc and bounds valid for the whole
family. The later continuation argument in this paper supplies such common
neighborhoods and uniform estimates. An extension neighborhood for each function
separately would leave that additional task unresolved.

The exceptional sets in these boundary statements are described by a
capacitary trace theory. In their work on weak products, Richter and
Sundberg develop the boundary trace theory for harmonic weighted
Dirichlet spaces \cite[Section~5]{RS14}. An obstacle capacity measures
the energy needed for a function to be at least one near the set.
The associated boundary theory describes traces and limits outside
sets of capacity zero. Their account builds on general potential
theory and credits Chac\'on's weighted boundary-limit work.
In the classical harmonic case, comparison with logarithmic capacity
identifies the same polar sets. This gives a common language for
the exceptional sets without identifying different numerical
normalizations of capacity. 

Analytic continuation to the whole plane asks a substantially different
regularity question. Cobos and Seco \cite{CS20} study
entire operator-inner functions, where \emph{entire} means analytic on all
of \(\mathbb C\), and operator-innerness is the unit-norm,
shift-orthogonality condition. Their weighted-Hardy setting has positive
weights \(\omega_n\), \(\omega_0=1\), \(\omega_{n+1}/\omega_n\to1\), and
the doubling bound
\[
 \sup_{n\ge1}\ \sup_{n\le k\le2n}
       \frac{\omega_k}{\omega_n}<\infty .
\]
This prevents the weights from growing by more than a fixed factor as the
index runs from \(n\) to \(2n\). Under their additional
no-extraneous-zeros assumption, the only entire operator-inner functions
are normalized monomials, up to unimodular constants. The no-extraneous-zeros assumption controls both the zero set and
the local order of vanishing of the prescribed-zero inner constructions.  Felder and Le subsequently proved the same conclusion for all
positive weight sequences satisfying \(\omega_{n+1}/\omega_n\to1\),
without the doubling condition or the no-extraneous-zeros assumption
\cite[Proposition~4.7]{FL22}. The conclusion identifies how restrictive global analytic extension
can be. The local extension across selected arcs used here asks for less,
and still requires its own estimates when a family of functions is
involved.

These discussions involve three different ways of varying an approximation
problem. For one fixed function, increasing the polynomial degree follows
the projection onto its generated space. At one fixed degree, varying a
nonzero function changes a finite Gram system. For varying invariant
subspaces, the target is instead convergence of projections on every fixed
Hilbert-space vector. The table records the quantity and topology in each
question, with representative sources. It also indicates why a change of
ambient space or norm requires a fresh comparison.

\begin{center}
\small
\begin{tabular}{@{}L{50mm}L{74mm}L{30mm}@{}}
\toprule
What is fixed and what varies & Target quantity and topology & Representative sources\\
\midrule
Fixed \(f\); degree \(n\) increases &
\(p_nf\to P_{[f]}1\) in  the Dirichlet norm. &
\cite{BFKSS18}\\
\addlinespace
Fixed degree \(d\); nonzero \(f_n\to f\ne0\) &
Degree-\(d\) approximants converge in coefficients and polynomial norm. &
\cite{AS24}\\
\addlinespace
Fixed Hilbert space; closed subspaces \(M_n\) vary &
\(P_{M_n}h\to P_Mh\) in norm for each fixed \(h\), under the two Mosco conditions. &
\cite{Mosco69,CHMB21}\\
\bottomrule
\end{tabular}
\end{center}

\subsection{Related function spaces and approximation questions}\label{subsec:01-wider-field}

The related literature uses weights and measures in several distinct ways.
A radial coefficient weight changes the contribution of each Taylor
coefficient to the squared norm. A harmonic or superharmonic derivative weight changes the
contribution of different parts of the disk to the energy.  These choices
define different ambient Hilbert spaces. Richter's representation \cite{Richter91}
and Richter--Sundberg's weighted-space representation \cite{RS92} explain why measures naturally enter invariant-subspace theory; Guillot's work \cite[Theorem~4.6]{Guillot12} develops the boundary theory for
another family of measure-defined spaces. In the weighted-space representation discussed above, the \(H^2\) term and the weighted derivative term
were specified together. A measure can also have a different role: it can
prescribe a singular inner factor and thereby select a subspace inside the
fixed classical Dirichlet space. In that use the ambient norm does not
change. Keeping the role of the measure explicit is necessary when
comparing cyclicity and approximation statements.

Guillot \cite[Theorem~4.6 and Theorems~5.3, 5.8--5.9]{Guillot12} studies capacities, fine boundary behavior, and invariant subspaces for harmonically weighted Dirichlet
spaces. In particular, his work obtains a Brown--Shields-type
characterization when the harmonic weight is the Poisson integral of
finitely many point masses. The capacity associated with the space therefore depends on the defining
measure.  El-Fallah, Elmadani, and
Kellay \cite{EEK16} study cyclicity and describe invariant
subspaces for defining measures whose topological support is countable.
In this setting, the support of the measure enters explicitly into
the weighted cyclicity problem. 

El-Fallah, Elmadani, and Labghail
\cite[Theorem~1 and Corollary~1.1]{EEL22} use
extremal functions to study a further class of weighted Dirichlet spaces.
Their class of defining measures includes canonical Cantor measures, and
they obtain positive weighted Brown--Shields and invariant-subspace
results in that setting. Extremals provide a way to examine the subspace
through a distinguished function, complementing descriptions based
directly on boundary sets. Their paper also considers recovery of an outer
extremal from its modulus on Carleson support
\cite[Theorem~2 and Lemma~5.1]{EEL22}. This recovery question asks
how the boundary modulus on that support determines the chosen extremal. It is separate from
asking whether a given function generates the entire weighted space. 

For positive superharmonic weights, Bahajji-El Idrissi,
El-Fallah, Elmadani, and Hanine extend the inner--outer
description of invariant subspaces and the construction of a
cyclic function vanishing on a prescribed closed polar set
\cite[Theorems~1 and~6.1]{BEIEH26}.
The capacity in these statements is attached to the weighted space;
the cyclic-function construction is existential for each
prescribed closed polar set.

Generalized Cantor geometry also gives positive results for radially
weighted energies. El-Fallah, Kellay, and Ransford
\cite{EFKR10Cantor} combine outerness, continuity of the modulus on the
closed disk, and a capacity-zero condition on a generalized Cantor set
containing the boundary zeros. The geometric structure of the containing
set makes quantitative estimates available, while continuity identifies
the boundary vanishing under discussion. Their parameter describes a
radial derivative weight, with the classical Dirichlet space at parameter
zero. This convention differs from the coefficient-power scale used above,
where \(\D\) corresponds to exponent one. The work is a distinct
Cantor-set study from the authors' invariant-subspace survey, and its
regularity and geometric hypotheses are part of the positive cyclicity
result.

The reproducing kernel gives another way to measure the influence of a
weight. Its diagonal value is the squared norm of point evaluation, so
large diagonal values indicate that evaluation can be large for a function
with controlled norm. El-Fallah, Elmadani, and Kellay
\cite{EEK19Kernel} develop sharp diagonal-kernel estimates and associated
capacity estimates for points and arcs in weighted Dirichlet spaces,
together with criteria for polar sets. This connects analytic evaluation
estimates with the energy needed to impose boundary constraints. The
connection is useful precisely because both quantities depend on the
chosen weighted space; it relates diagonal-kernel estimates to the
potential-theoretic size of sets.

Local energy estimates make this relationship constructive for particular
outer functions. El-Fallah, Kellay, and Mahzouli
\cite{EFKM26} study distance-based outer functions, whose boundary moduli
are built from the distance to a chosen set. They estimate the local
Dirichlet energy of these functions and apply the estimates to
multipliers, polarity, and weighted cyclicity. Local energy separates the
contribution associated with a boundary point before it is integrated
against the defining measure. The resulting estimates can therefore relate
the geometry of the set and the growth of the measure to analytic
membership and multiplication questions. Their spaces consist of analytic
\(H^2\) functions with weighted local energy; they are distinct from the
harmonic trace spaces used to formulate obstacle capacities. The
applications retain the geometric and measure-growth assumptions of the
particular estimates.

Changing the coefficient exponent can also change whether a fixed function
is cyclic. Vavitsas, Wu, and Zarvalis \cite{VWZ26} use critical Riesz capacities of Cantor sets to study this
transition. Their first preprint version gives the following endpoint
construction: for each \(0<\alpha_*\le1\), there is a function in
\(\D_{\alpha_*}\) that is cyclic in every lower-exponent space
\(\D_\alpha\), \(\alpha<\alpha_*\), but is not cyclic at \(\alpha_*\).
Here \(\D_\alpha\) uses the coefficient weights \((n+1)^\alpha\), so the
classical Dirichlet space occurs at exponent one. The critical capacity
records how the boundary set interacts with the changing norm. This is an
across-exponent result: the ambient approximation space changes, and the
statement does not assert the two Brown--Shields boundary conditions at
the endpoint in a fixed classical space.

\paragraph*{Bergman generation and orthogonality.}

The Bergman space offers a different Hilbert-space model for the relation
between extremals and invariant subspaces. Its unweighted norm is
\[
 \|f\|_{A^2}^{\,2}
   =\int_{\DD}|f(z)|^2\,dA_0(z),
\]
with normalized area as above. Hedenmalm's factorization work
\cite{Hedenmalm91} supplies an important predecessor to the Bergman
wandering-subspace theory. As Aleman, Richter, and Sundberg explain in
their account \cite{ARS96}, suitable normalized extremals for zero-set
subspaces serve as contractive divisors: division by the extremal controls
the Bergman norm of functions with the associated zeros. Their discussion
credits a development across Hedenmalm's early factorization papers. The
positive idea is to use an extremal function to remove prescribed zeros
while retaining norm control. Both the direction of the division estimate
and the norm are specific to the Bergman setting.

Aleman, Richter, and Sundberg prove the
wandering-subspace generation theorem for the unweighted Bergman space
\cite{ARS96}. For every closed shift-invariant subspace \(M\subset A^2\),
the vectors orthogonal to \(zM\), together with their successive shifts,
generate the whole subspace:
\[
 M=\overline{\operatorname{span}
 \{z^nw:n\ge0,\ w\in M\ominus zM\}}^{\,A^2}.
\]
The wandering space here can have arbitrary finite or infinite dimension.
Thus the theorem recovers the space from a distinguished orthogonal
complement without reducing that complement to one line. This differs from
Richter's one-dimensional Dirichlet generation theorem while preserving
the same structural question. Aleman, Richter, and Sundberg's introduction
also credits contractive-divisor approaches of Duren, Khavinson, Shapiro, and Sundberg alongside Hedenmalm's
extremal work. These antecedents help explain how function-theoretic
division estimates enter the operator-theoretic generation problem.

Shimorin places wandering-subspace questions within a wider study
of Wold-type decompositions for operators close to isometries
\cite[Theorems~3.6 and~4.1]{Shimorin01Wold}. The question is which part of an operator's
action can be reconstructed from successive images of wandering vectors
and what additional structure is needed for the remaining part. The
Bergman and Dirichlet theorems make this question concrete in two analytic
Hilbert spaces with different norms. Shimorin's results on Wold-type
decompositions provide a broader setting for comparing these generation
phenomena.

Beyond these Hilbert models, even the language of orthogonality requires
care. Hedenmalm and Shields \cite{HS90} study
invariant-subspace questions in Banach spaces. Cheng, Mashreghi, and Ross \cite{CMR19} survey and
develop inner-function ideas in reproducing-kernel spaces, including the
role of shift orthogonality and zero constructions. In a Hilbert space,
orthogonal projection gives a vector and an exact decomposition of squared
norms. In a Banach setting, best approximation and norm inequalities can
provide different notions of orthogonality. The resulting inner-function
terminology describes the geometry of the space and need not mean
unimodular Hardy boundary values.

\paragraph*{Cyclicity in several variables.}

B\'en\'eteau, Condori, Liaw, Seco, and Sola extend optimal approximation to the bidisk
\cite{BCLSS15II}. The domain is the product \(\DD^2\), and their
coefficient norm for \(f(z_1,z_2)= \sum_{j,k\ge0}a_{j,k}z_1^jz_2^k\) is
\[
 \|f\|_{\alpha}^{\,2}
 =\sum_{j,k\ge0}(j+1)^\alpha(k+1)^\alpha
                  |a_{j,k}|^2 .
\]
There are now two coordinate shifts, and cyclicity asks whether
polynomials in both variables times \(f\) are dense in the chosen space.
The authors obtain approximation rates for specified classes and relate
cyclicity to geometry on the distinguished boundary \(\T^2\). The product
domain and product weights make this a two-variable approximation problem,
with a richer range of boundary zero configurations than in the disk.

B\'en\'eteau, Knese, Kosi\'nski, Liaw, Seco, and Sola \cite{BKKLSS16} study cyclic polynomials in
these bidisk Dirichlet-type settings. For polynomials with no zeros in the bidisk, their classification
depends on the weight exponent and on the zero sets of the
irreducible factors on the distinguished boundary.
This makes several-variable geometry a decisive part of the problem:
boundary zeros can occur in configurations that have no one-variable
counterpart. The result is a polynomial classification in its stated
dimension and norm. It shows how the same approximation question leads to
different geometric criteria when both coordinate shifts act.
Cyclicity on the ball has also been studied in the Drury--Arveson space and other weighted Besov spaces \cite{APRSS24}.

\paragraph*{Further approximation settings and singular factors.}

Spectral methods also study the positive operator
\(M_f^*M_f\) associated with multiplication. Monsalve-L\'opez
and Seco \cite{MS25} study this operator for linear factors
\(f(z)=a-z\) in weighted Hardy spaces, with detailed spectral
descriptions in Bergman-type cases. The quadratic form of this
operator records the squared norm of a product,
\(\langle M_f^*M_fh,h\rangle =\|fh\|^2\). Its spectrum therefore
offers a way to investigate the range of possible multiplication
norms in this restricted family. The family and weight assumptions are part of the
spectral study; it does not give a general spectral classification of
classical Dirichlet cyclicity. Santana-Rivero and Seco \cite{SSR26} continue this study by determining
the essential spectrum for polynomial symbols in a class of weighted
Hardy spaces and analyzing the point spectrum for a restricted class
of quadratic polynomials in the Bergman space.
The survey of B\'en\'eteau and Centner
\cite{BC21} places such approximation questions beside the wider
literature on optimal polynomial approximants.

B\'en\'eteau, Cheng, Felder, Khavinson, Manolaki, and Maronikolakis \cite{BCFKMM26}
study metric projections, optimal approximants, and extremal problems in
\(H^p\) for \(1<p<\infty\). Its norm measures the \(L^p\) size of
the Hardy boundary values. Here the least-distance problem survives, but,
except at \(p=2\), its minimizer is a metric projection in a Banach space
rather than a Hilbert orthogonal projection. Banach-space orthogonality
and associated inequalities provide tools for studying these minimizers
and their zeros. This work shows how much of the optimal-approximation
program can be pursued without the Hilbert norm identities used in the
present paper. B\'en\'eteau and Centner's survey \cite{BC21} provides
broader orientation, while the specific metric projection questions here
belong to the \(H^p\) geometry.

A final comparison returns to Hardy factorization and the singular inner
factors used below. A singular inner factor records a positive singular
boundary measure in the inner part of the factorization. It has no zeros
inside the disk, despite its nontrivial boundary singularity. Prescribing
such a factor selects the Dirichlet functions whose Hardy factorizations
contain it; the norm on those functions remains the classical Dirichlet
norm. Hartz and Richter \cite[Appendix~B, proof of Theorem~B.1]{HR26} pass from Blaschke zeros to a singular measure on the circle using cutoff functions. Each cutoff is the average of the corresponding Blaschke factor and \(1\). This
explains why a boundary measure is a natural intermediate object after
finite interior zero sets. The measure specifies divisibility inside the
fixed space, not a new energy norm.

Dayan and Seco \cite[Theorems~1.1 and~1.2]{DS26} study the cyclic
behavior of singular inner functions in different ambient spaces. They
show how regularity of the singular measure can yield cyclicity in
specified Besov and analytic coefficient-sequence spaces, in a regime with
\(p>2\), and obtain a necessary condition involving Beurling--Carleson
sets. The measure's regularity is therefore part of a positive
approximation theory for singular inner functions in those norms.  This comparison shows why a singular factor's cyclic behavior
must be read in its own space. Hartz and Richter's theorem on the
existence of factors in the classical Dirichlet space \cite{HR26},
and the finite-atomic targets used below, serve a different purpose:
they organize subspaces inside the fixed classical norm. No cyclicity
of those atomic singular targets is supplied by the Besov or
sequence-space theorems. The next step chooses one finite
atomic measure and controls both the norm error of its canonical target
and the measure's total mass.

In the counterexample, the positive approximation error must persist
as the radial zero set is made polar. To retain outerness in the
limit, we also keep track of the boundary logarithms.

\subsection{Outerness, the capacity of radial zero sets, and noncyclicity}\label{subsec:01-mechanism}

We will take a limit of outer canonical functions with smaller and
smaller radial-zero capacity. At the same time, their polynomial
multiples must stay a positive distance from \(1\). To see how we can
keep this distance, let \(P_{[f]}\) denote orthogonal projection onto
\([f]\). For a canonical function \(f\), Section~\ref{subsec:02-certificate}
gives
\[
 P_{[f]}1=f(0)f,\qquad
 \|1-P_{[f]}1\|_{\D}^2=1-f(0)^2 .
\]
The value at zero varies continuously in the Dirichlet norm.
We can therefore start with a small ball about a proper outer canonical
function, so that every canonical function in that ball has value at
zero bounded away from one. Section~\ref{subsec:03-seed} supplies a
starting function. Staying in the ball preserves a uniform positive
approximation error.

To construct the sequence, we first ask a local question. Given a
proper outer canonical function \(f\), can we find another such
function close to it whose radial zero set has as small a capacity
as we wish? Two approximations answer this question together.

\textbf{Approximation A.}
Theorem~\ref{thm:A} starts with \(f\) and chooses a nonzero finite
positive measure \(\mu\) supported on finitely many points of the
circle. The measure determines a singular inner factor \(S_\mu\).
This factor has
boundary modulus one almost everywhere and value
\(S_\mu(0)=e^{-\mu(\T)}<1\). Thus a small total mass makes the
factor's value at zero close to one, while positive mass keeps the
factor nonconstant. We let \(G_\mu\) be the canonical
function of the space \(S_\mu H^2\cap\D\), whose members are
divisible by \(S_\mu\) with quotient in \(H^2\).
Given positive bounds on the Dirichlet norm error and the total mass,
the theorem chooses one \(\mu\) that meets both. These two bounds may be prescribed independently, and their simultaneous satisfaction
provides a single target for the next step. The mass estimate also
controls \(S_\mu(0)=e^{-\mu(\mathbb T)}\), in addition to the norm
approximation. The function \(G_\mu\) is the normalized extremal of the
whole space, with positive value at zero, and its inner factor is exactly
\(S_\mu\). The proof uses Richter and Sundberg's weighted representation
\cite[Theorem~3.2, p.~177]{RS92} to help characterize \([f]\) by a
weighted boundary integral. Hartz and Richter's cutoff construction
\cite[Appendix~B, proof of Theorem~B.1]{HR26} supplies a trial vector
in \(S_\mu H^2\cap\D\). For the reverse comparison,
Section~\ref{subsec:07-comparison} constructs a trial vector in the
finite-zero space from \(G_\mu/S_\mu\). The resulting estimates bound
the distance from the first trial vector to \(G_\mu\) in the
Dirichlet norm.

\textbf{Approximation B.}
Theorem~\ref{thm:B} keeps the measure \(\mu\) fixed, including the
location and mass of every atom. It constructs trace-zero spaces
whose outer canonical functions converge in Dirichlet norm to
\(G_\mu\). The defining compact sets satisfy explicit bounds on location
and capacity. Membership in a trace-zero space requires the boundary
trace to vanish on the defining compact set outside a polar set.
Each canonical function extends continuously to the closed disk
and its zero set there is exactly its defining
compact set. The projections onto the trace-zero spaces converge
strongly to the projection onto \(S_\mu H^2\cap\D\). We choose one
approximant \(g\) at an index that meets the prescribed norm and
capacity bounds. For a fixed compact set, Richter and Yilmaz use
analytic continuation of the canonical function across the
complementary arcs \cite[proof of Theorem~1.3]{RY19}.
Section~\ref{subsec:12-assembly} then extends the recovery to
nonzero finite positive singular measures \(\mu\) for which
\(S_\mu H^2\cap\D\ne\{0\}\).

Proposition~\ref{prop:retain-singular-factor} allows the approximants
to retain any prescribed singular inner divisor \(S_\nu\) with
\(0\le\nu\le\mu\). The remaining measure \(\mu-\nu\) is the
singular part of the weak-* limit of the measures given by the
negative boundary logarithms.

For the varying compact sets \(K_j\), put \(\eta_j=G_{\D_{K_j}}\).
To pass to the limit, the proof needs a lower bound on a fixed neighborhood
and an area bound for \((z\eta_j)'\) on a common neighborhood, uniform as the sets vary.
Sections~\ref{sec:08}--\ref{sec:09} establish these bounds
and identify the one-atom limit.

The Dirichlet triangle inequality through \(G_\mu\) makes \(g\)
close to \(f\). We then compare the boundary logarithms of \(f\) and \(g\) directly. For a fixed outer target \(f\), the
estimate in Section~\ref{subsec:04-continuity} turns sufficiently small Dirichlet error into small \(\dlog\)-distance. The estimate first controls the part of \(\log|f^*|\) where the boundary modulus
is very small, then uses the norm error. This determines the order of our choices. We fix the outer target and the desired metric accuracy, choose a small-modulus threshold that makes this logarithmic tail small, and only then choose an allowable norm error. This choice depends on the fixed target, including its value at zero and its logarithmic tail.
We allow half of the norm error for each approximation. Approximation~A selects the measure while meeting its separately prescribed mass bound;
with that measure held fixed, Approximation~B supplies one index meeting both the remaining norm bound and the desired capacity bound. The norm triangle inequality now gives the error required by the fixed-target estimate. The intermediate function \(G_\mu\) has a nonconstant inner factor, so the Dirichlet norm is the distance used in both approximation steps.
We apply the estimate for the boundary logarithms only after both approximations have been chosen. The analytic-logarithm criteria discussed earlier
\cite{BCLSS15,AR26} concern cyclicity; the estimate used here concerns continuity of real boundary logarithms at one fixed outer target.

A function with small radial-zero capacity gives a possible center
for the next ball. Its capacity bound must also hold at the eventual
limit. Section~\ref{subsec:03-neighborhood} shows how to choose a
radius for which the bound holds on the entire closed ball. We start
with a strict margin between the center's capacity and the desired
upper bound. A nearby function can vanish only where the center has
small trace or where the change in the trace is large, up to a polar
exception. As the threshold for a small trace decreases, the capacity
of the first set approaches the capacity of the center's zero set.
We choose this threshold so that the first contribution uses only
part of the available margin. Once the threshold is fixed, the norm
bound on the change of trace controls the capacity of the second set.
A sufficiently small radius keeps that contribution within the
remaining margin, including at the boundary of the ball. This is why
we choose the threshold first and the radius of the ball second.
Since the outer metric dominates Dirichlet distance, the same estimate
holds on its closed balls. The eventual limit may lie on a ball's
boundary, so this inclusion of the boundary is needed when passing to
the limit.

Section~\ref{sec:13} repeats this choice inside the preceding ball,
requiring radial-zero capacity less than \(2^{-n}\) at every point
of the \(n\)-th ball. Each new center is chosen strictly inside the
previous ball, leaving room for a smaller closed ball around it.
The balls are nested and their radii tend to zero, so their centers
form a Cauchy sequence for \(\dlog\). These balls are closed relative
to the complete space of outer canonical functions.
Completeness gives an outer canonical limit. Because every ball is
closed, the limit belongs to all of them. Its radial-zero capacity
is zero, and its polynomial multiples stay a positive distance
from \(1\) by the bound on the initial ball. That initial ball also
keeps the value at zero positive and bounded away from one.
The three conclusions concern this one common limit: completeness
preserves outerness and canonical normalization, membership in every
later ball forces the radial zero set to be polar, and membership in
the initial ball preserves the positive approximation error.

\subsection{Proof structure and further reading}\label{subsec:01-reading}

Sections~\ref{sec:02}--\ref{sec:04} prepare the shared facts on
canonical functions, boundary traces, and the complete metric.
The proof of Approximation~A, in
Sections~\ref{sec:05}--\ref{sec:07}, first approximates \(f\) by
the extremal of a space defined by finitely many interior zeros.
It uses Richter and Sundberg's weighted representation of
the space generated by \(f\)
\cite[Theorem~3.2]{RS92}.
Section~\ref{sec:05} combines this representation with a Hardy-quotient
estimate to identify a weighted boundary integral that is finite exactly on the
generated space. Bounded truncations of its density lead to the
finite-zero approximation in Section~\ref{sec:06}. The finite set must
give both a small norm error and a small sum of
\(\log(1/|\lambda|)\) over its points \(\lambda\). A determinantal
averaging identity helps select one set that meets both bounds
(Section~\ref{subsec:06-selection}). Its finite-subset weights are the
standard ones described by Kulesza and Taskar \cite{KT12}; the related
mean embedded-inverse identity of Mutn\'y, Derezi\'nski, and Krause
\cite{MDK20} and expected-projector identity of Fanuel, Schreurs, and
Suykens \cite{FSS21} are stated for real positive-definite matrices.
Section~\ref{subsec:06-determinant} gives the formulation and proof used
here. The second selection bound controls the mass when the zeros
are replaced by boundary atoms. We use Hartz and Richter's cutoff construction and apply the parallelogram law as they suggest \cite[Appendix~B, proof of Theorem~B.1]{HR26}. The estimates in
Section~\ref{sec:07} then bound the Dirichlet distance from the trial function to the canonical function of the resulting singular-factor space.  Thus the three tasks are to characterize
the target space by the finiteness of a weighted boundary integral, select
one finite set meeting both bounds, and transfer that selection to the
canonical atomic target.

For Approximation~B, we first recover the canonical function
associated with one atom. Sections~\ref{sec:08}--\ref{sec:09}
prove an analytic criterion for this recovery. Fix one atom of
positive mass, including its location. Suppose nonpolar Carleson sets
(Section~\ref{subsec:03-ideal}) shrink to that location in Hausdorff
distance and their canonical values at zero converge to the atomic
target's value at zero. Then the canonical functions converge to that
target in Dirichlet norm. Within this specified geometric family, convergence of the
values at zero to the target's value suffices. For a single compact set,
continuation of the extremal off that set is used by Richter and Yilmaz \cite[proof of Theorem~1.3]{RY19}. The distinction
between individual and common continuation was explained in the
stability discussion above. Section~\ref{sec:08} gives a small-disk energy bound and common
neighborhoods of continuation away from the atom. Section~\ref{sec:09}
uses these estimates to prove strong convergence and identify the limit.

Sections~\ref{sec:10}--\ref{sec:11} construct sets to which this
criterion applies. Kellay's inequality \cite[Theorem~2.2]{K11}
gives the subarc-length factor in the squared trace estimate for
analytic traces. Using Kellay's reflection argument \cite[Lemma~2.1]{K11},
Theorem~\ref{thm:packed-trace} gives the same factor, with an
explicit constant, for functions in \(H^{1/2}(\T)\) whose
quasi-continuous representatives vanish quasi-everywhere
on the Cantor copies.
We apply this estimate to the Cantor copies and use Jensen's
inequality, as in the proof of Kellay's Theorem~3.1 \cite{K11},
to obtain a lower bound on a boundary logarithmic integral.
After we apply suitable M\"obius transformations, this lower bound makes the canonical value at zero small at one end of a path of sets. At the other end the value tends
to one. Continuity gives sets with the prescribed value.

Section~\ref{sec:12} takes the union of the sets for the individual atoms. The trace-zero space of the union is the intersection of the
trace-zero spaces of the individual compact sets. We approximate target vectors by polynomial multiples of the products of the one-atom canonical functions. The limit of these products generates the target space because its outer factor is cyclic.
We use Hardy quotients to show that every subsequential weak limit is divisible by each prescribed singular inner factor and hence belongs to the target space. These two properties give Mosco convergence \cite{Mosco69}. The orthogonal projections therefore converge in norm on each fixed vector. Applying them to \(1\) and normalizing yields convergence to the canonical function for the fixed measure.

Section~\ref{sec:13} uses the two approximations to choose the nested balls. Figure~\ref{fig:whole-map} records the dependencies in that final step, including the separate roles of the initial bound on the values
at zero, the whole-ball capacity estimate, and completeness.

\begin{figure}[H]\centering
\begin{tikzpicture}[>=Latex,
 b/.style={draw=black,fill=white,rounded corners=2pt,align=center,
           font=\small,text width=3.5cm,minimum height=1.12cm},
 mid/.style={b,text width=4.3cm},
 arr/.style={->,semithick}]
\node[b] (f) at (-5.2,0) {Proper outer\\canonical target \(f\)};
\node[b] (s) at (0,0) {Canonical function \(G_\mu\)\\of one fixed atomic measure};
\node[b] (g) at (5.2,0) {Outer canonical \(g\)\\with small capacity\\of radial zeros};
\node[mid] (d) at (0,-1.9) {New center close to \(f\) in \(\dlog\)\\Section~\ref{subsec:04-continuity}};
\node[mid] (n) at (0,-3.65) {Nested closed balls\\Section~\ref{subsec:13-balls}};
\node[b] (c) at (5.2,-3.65) {Capacity bound\\on each closed ball\\Section~\ref{subsec:03-neighborhood}};
\node[b] (o) at (-5.2,-5.4) {Values at \(0\)\\bounded away from \(1\)\\Section~\ref{subsec:13-balls}};
\node[b] (m) at (5.2,-5.4) {Complete outer metric\\Section~\ref{subsec:04-complete}};
\node[mid,fill=gray!20] (l) at (0,-5.4) {Outer, noncyclic limit\\with polar radial zeros};
\draw[arr] (f)-- node[above,font=\small]{A} (s);
\draw[arr] (s)-- node[above,font=\small]{B} (g);
\draw[arr] (f.south) -- (d.west);
\draw[arr] (g.south) -- (d.east);
\draw[arr] (d)--(n);
\draw[arr] (c)--(n);
\draw[arr] (n)--(l);
\draw[arr] (o)--(l);
\draw[arr] (m)--(l);
\end{tikzpicture}
\caption{The Dirichlet-norm approximations A and B supply a new
center. The estimates at an outer target and on whole balls let the
construction repeat this step and pass to a limit.}
\label{fig:whole-map}\end{figure}
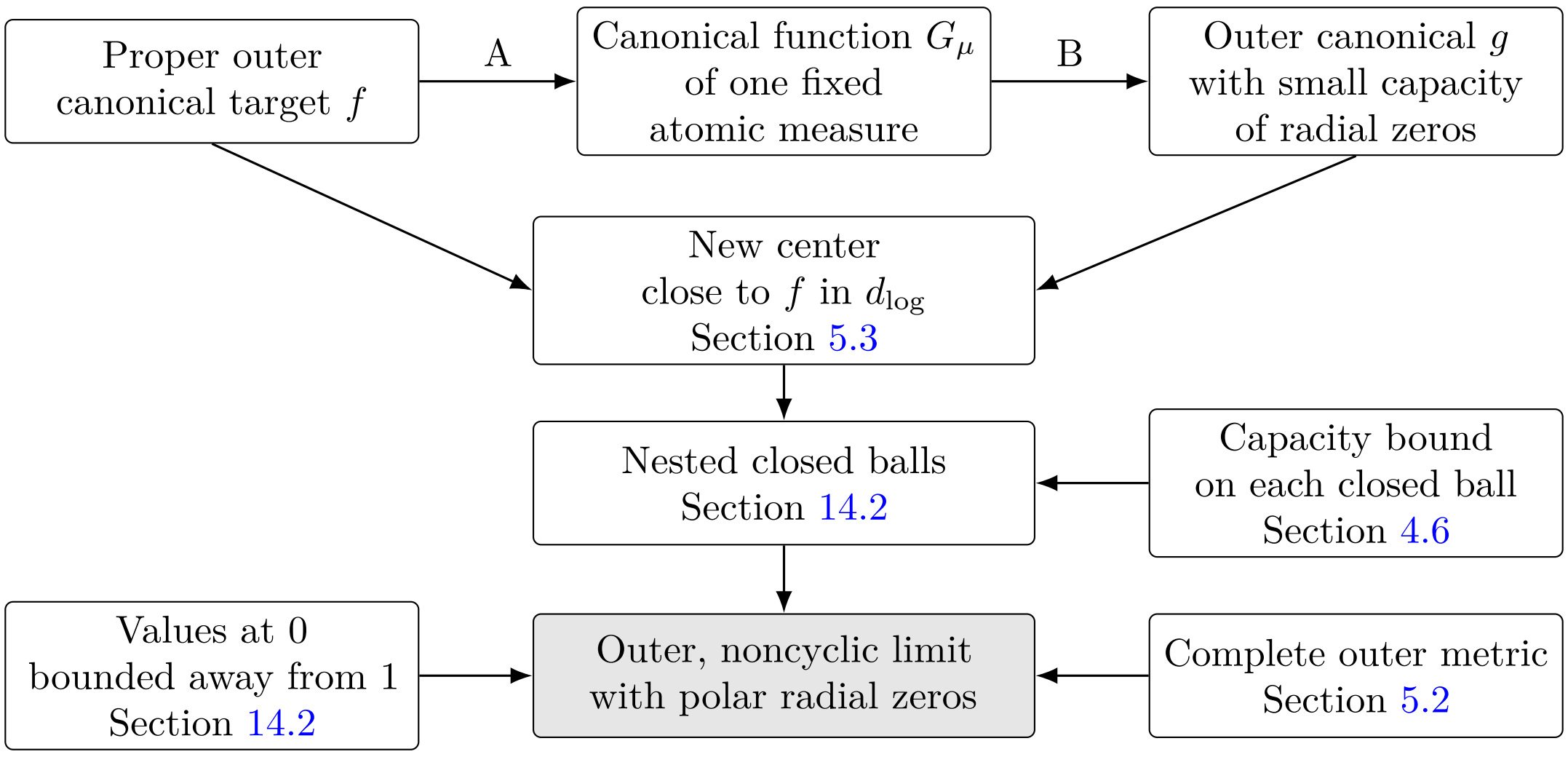\FloatBarrier

The paper can be read along several routes. For the conceptual
argument, continue to Section~\ref{sec:approximation-and-limiting-argument},
\emph{Overview of the proof}, which explains the two approximation theorems
and their use in the limiting argument. For the foundations, read
Section~\ref{sec:02}, \emph{Canonical functions in the Dirichlet space},
then Section~\ref{sec:03}, \emph{Boundary traces, capacity, and trace-zero
spaces}, and Section~\ref{sec:04}, \emph{A metric on outer canonical functions and two approximations}. These sections fix the normalizations,
boundary representatives, and topology used in every later step.
Readers following the technical construction of Approximation~A can
then take Sections~\ref{sec:05}--\ref{sec:07} together; readers following
Approximation~B can take Sections~\ref{sec:08}--\ref{sec:12} together,
starting with the boundary estimates and the one-atom criterion, then
proceeding to the geometric construction and recovery for singular measures.
Section~\ref{sec:13}, \emph{Constructing the counterexample}, is the common destination of these routes.
Figure~\ref{fig:whole-map} can be used to locate the role of a technical
estimate while reading its proof.

For systematic background on the Dirichlet space, the primer of
El-Fallah, Kellay, Mashreghi, and Ransford
\cite{Primer14} develops the foundational function theory.
Arcozzi, Rochberg, Sawyer, and Wick
provide an accessible survey \cite{ARSW11} and a separate, broader
monograph \cite{ARSW19}. The survey connects capacity, Carleson measures,
and multipliers; the monograph also develops the relations with
reproducing kernels, interpolation, tree models, and neighboring Besov
spaces. Ross's survey \cite{Ross06} offers a focused route
through classical boundary values, zero sets, and the forward and
backward shifts. Zhijian Wu \cite{Wu98} surveys the Dirichlet shift and its invariant
subspaces, as well as interpolation and Hankel operators.  For invariant subspaces and cyclicity, El-Fallah,
Kellay, and Ransford's survey \cite{EFKR10} brings the
principal questions and methods together. For optimal approximants,
B\'en\'eteau and Centner \cite{BC21} explain their
connections with reproducing kernels, orthogonal polynomials, zeros,
and convergence. 

The next section follows
the two approximations and the limiting argument that combines them.

\newpage 
\section{Overview of the proof}
\label{sec:approximation-and-limiting-argument}

Fix a proper outer canonical function \(e\). The first approximation
selects a subspace associated with a finite atomic singular inner function whose normalized generator is close to \(e\).
The second holds that subspace fixed and recovers its generator by outer functions. We then compare the boundary logarithms of \(e\) and the outer canonical function chosen in the second approximation and pass to a limit in a neighborhood where polynomial approximation of \(1\) has a positive error.

\subsection{Approximation A with small singular mass}\label{subsec:overview-approximation-a}

Theorem~\ref{thm:A} approximates any proper outer canonical function
in the Dirichlet norm by canonical functions with finite atomic
singular inner factors. The total mass of the singular measure can be
made arbitrarily small at the same time.

\medskip
\textbf{Precise statement.} Let \(e\) be a normalized outer function in \(\D\) satisfying
\[
 \|e\|_{\D}=1,\qquad e(0)>0,\qquad
 \langle z^ne,e\rangle_{\D}=0\quad(n\ge1).
\]
Assume \(e\ne1\). By the canonical characterization in
Section~\ref{subsec:02-certificate}, \(e\) is the canonical function of
its generated space \([e]\), and that space is proper.

Write a nonzero finite positive atomic measure on \(\T\) as
\[
 \mu=\sum_{i=1}^N\tau_i\delta_{\xi_i},
 \qquad \tau_i>0,\quad \xi_i\in\T,
\]
with distinct \(\xi_i\). Its singular inner function is
\[
 S_\mu(z)=
 \exp\left(-\sum_{i=1}^N\tau_i\frac{\xi_i+z}{\xi_i-z}\right).
\]
The notation in \eqref{eq:atomic-target} gives the invariant subspace
and its normalized canonical function as
\[
 M_\mu=S_\mu H^2\cap\D,\qquad G_\mu=G_{M_\mu}.
\]
Theorem~\ref{thm:A} allows separate positive parameters
\(\varepsilon_{\mathrm{norm}}\) and \(\varepsilon_{\mathrm{mass}}\).
In particular, for each \(\varepsilon>0\) it gives one nonzero
measure \(\mu\) for which
\[
 \|G_\mu-e\|_{\D}<\varepsilon,\qquad
        0<\mu(\T)<\varepsilon.
\]
The exact Hardy inner factor of \(G_\mu\) is \(S_\mu\), and \(G_\mu\)
generates \(M_\mu\).

\medskip
\textbf{Strong approximation by canonical functions with a singular inner factor.} The starting function \(e\) is outer. Every approximant \(G_\mu\)
has a nontrivial singular inner factor. Letting both parameters tend to zero gives strong convergence to \(e\) in \(\D\) and lets the total singular mass tend to zero along the same sequence. Each approximant
is the normalized canonical function of its own invariant subspace, so the convergence takes place within the class of canonical
functions.

\medskip
\textbf{Main steps of the proof.} The proof proceeds in three steps.

\begin{enumerate}
\item \textbf{A boundary integral condition for the generated space.}
Section~\ref{sec:05} represents \([e]\) as the finite domain of the
weighted boundary quadratic form \(Q_e\), defined in
\eqref{eq:cost-density}. The form makes the orthogonal
projection onto \([e]\) accessible by variational approximation.

To use a boundary condition to recover membership in the generated
space, we need to know which quotients by the generator belong to
the weighted space, as well as the product norm identity.
Richter and Sundberg's representation
\cite[Theorem~3.2]{RS92}
makes multiplication by the unit canonical generator an onto isometry
from the weighted space \(\mathcal H_w\) to the invariant subspace.
Functions in \(\mathcal H_w\) belong to \(H^2\) and have finite
weighted derivative energy; their norm combines the \(H^2\) norm
with that energy.
Section~\ref{sec:05} uses this exact model and a separate Hardy-quotient
estimate to identify the finite domain of the boundary form.

\item \textbf{Approximation by finitely many interior zeros.}
Positive quadrature and a finite determinantal selection argument in
Section~\ref{sec:06} choose finite sets \(A\subset\DD\setminus\{0\}\)
whose associated canonical functions approximate \(e\).
The same selection ensures that the corresponding sums satisfy
\[
 b(A)=\sum_{z\in A}-\log|z|\longrightarrow0.
\]
Controlling the Dirichlet error and the sum \(b(A)\) for the same choice
of \(A\) is the quantitative step used in the transfer to the boundary.

The selection uses an average of projections to find a suitable
subset of the evaluation vectors. Kulesza and Taskar describe the
standard determinant weights for finite subsets, called finite
\emph{L-ensemble} weights
\cite[Section~2.2, (2.6) and Theorem~2.1]{KT12}.
Mutn\'y, Derezi\'nski, and Krause give a related mean embedded-inverse
identity \cite[Theorem~1 and Remark~1]{MDK20}, and Fanuel, Schreurs,
and Suykens give the corresponding expected-projector identity
\cite[Corollary~2.3, p.~286]{FSS21}; both are
stated for real positive-definite matrices. The expected projector is
the average of the projections associated with the possible subsets.
Section~\ref{subsec:06-determinant} gives the complex-vector formulation
used here. Section~\ref{subsec:06-selection} then combines the error estimates and bounds on \(b(A)\) to choose one subset satisfying both.

\item \textbf{Transfer from interior zeros to boundary atoms.}
Section~\ref{sec:07} replaces the finite Blaschke factor for \(A\)
by a finite atomic singular inner factor. Comparison estimates
control the change in the canonical functions in the Dirichlet norm.
The resulting singular measure has total mass at most \(b(A)\).

We use Hartz and Richter's cutoff functions and apply the parallelogram law as they suggest \cite[Appendix~B, proof of Theorem~B.1]{HR26}. This gives a vector of controlled norm in \(M_\mu\). Their comparison starting from an extremal function
\cite[Sections~3 and~5]{HR26} suggests how to reach the
normalized generator itself. Starting with \(G_\mu\) for
the measure just selected, Section~\ref{subsec:07-comparison}
uses the quotient \(G_\mu/S_\mu\) to construct a trial vector
in the finite-zero space. The resulting bound on \(G_\mu(0)\) converts the trial-vector estimate into a Dirichlet
distance to \(G_\mu\). Both the distance and the mass concern
this one measure.
\end{enumerate}

Does \(b(A)\to0\) by itself control the corresponding projection error? The following example shows why that condition must be accompanied by a norm estimate.

\textbf{Projection error as \(b(A)\to0\).}
We use the constant-target case of Sundberg's construction,
recorded by Luo and Richter \cite[Theorem~6.5]{LR15}. The following
choice of radius lets us compute both \(b(A_N)\) and the
limiting projection error.
Fix \(L>0\). For each \(N\ge 1\), choose \(r_N\in(0,1)\) so that
\[
 1-r_N^{2N}=e^{-NL},
 \qquad
 A_N=\{r_N e^{2\pi i j/N}:0\le j<N\}.
\]
The corresponding sum is
\[
 b(A_N)=-N\log r_N=-\tfrac12\log(1-e^{-NL})
 \longrightarrow 0.
\]

Let \(M_{A_N}=\{h\in\D:h|_{A_N}=0\}\). Rotation by \(2\pi/N\) is
unitary on \(\D\), preserves \(M_{A_N}\), and fixes \(1\).
Consequently, the projections of \(1\) onto \(M_{A_N}\) and
\(M_{A_N}^{\perp}\) are both rotation invariant. In the rotation-invariant
subspace, the functions have the form
\[
 h(z)=\sum_{m\ge0}u_m z^{Nm},
 \qquad
 \|h\|_{\D}^2=\sum_{m\ge0}(Nm+1)|u_m|^2,
\]
and vanishing on \(A_N\) is equivalent to the single condition
\(h(r_N)=0\). The evaluation kernel at \(r_N\) in this subspace and its
squared norm are
\[
 k_N(z)=\sum_{m\ge0}\frac{r_N^{Nm}z^{Nm}}{Nm+1},
 \qquad
 K_N=\|k_N\|_{\D}^2
     =\sum_{m\ge0}\frac{r_N^{2Nm}}{Nm+1}.
\]
Since \(k_N(0)=1\), the projection of \(1\) onto \(M_{A_N}^{\perp}\) is \(k_N/K_N\).
The projection error and the normalized canonical function are therefore
\[
 \dist_{\D}(1,M_{A_N})^2=\frac1{K_N},
 \qquad
 e_{A_N}=\frac{1-k_N/K_N}{\sqrt{1-1/K_N}}.
\]

For \(m\ge1\), replacing \(1/(Nm+1)\) by \(1/(Nm)\) incurs an
error of at most \(1/(N^2m^2)\). Summing and using the logarithmic
series gives
\[
 K_N
 =1-\frac{\log(1-r_N^{2N})}{N}+O(N^{-2})
 =1+L+O(N^{-2}).
\]
Thus \(\dist_{\D}(1,M_{A_N})^2\to 1/(1+L)>0\), even though
\(b(A_N)\to0\). On every disk \(|z|\le R<1\), the bound
\[
 |k_N(z)-1|\le\frac{R^N}{1-R^N}
\]
shows that \(e_{A_N}\) converges locally uniformly to the constant
\(\sqrt{L/(1+L)}\). This constant has Dirichlet norm strictly less
than one, and all the functions \(e_{A_N}\) have norm one. Hence
the sequence has no strong Dirichlet limit. The example explains
why the finite selection in Approximation A controls the
projection error and \(b(A)\) for the same set of zeros.

\medskip
\textbf{Role in the counterexample.} Approximation A gives the intermediate canonical functions, whose
singular inner factors are determined by finitely many boundary
atoms. Approximation B then approximates each such function by
outer canonical functions whose radial zero sets have arbitrarily
small logarithmic capacity. The limiting construction combines these
approximations and produces a function that is outer and noncyclic
and whose radial zero set has capacity zero.

In particular, for an arbitrary proper outer canonical function,
the first theorem supplies the finite atomic target needed by the
second theorem, with arbitrarily small Dirichlet error and total
singular mass.

\subsection{Approximation B by outer canonical functions}

Theorem~\ref{thm:B} starts with the fixed subspace \(M_\mu\). It recovers its canonical function in the Dirichlet norm
by outer canonical functions whose radial zero sets have
arbitrarily small logarithmic capacity.

\medskip
\textbf{Finite atomic targets.} Fix a nonzero finite positive atomic measure
\[
 \mu=\sum_{i=1}^N\tau_i\delta_{\xi_i},
 \qquad \tau_i>0,\quad \xi_i\in\T,
\]
where the \(\xi_i\) are distinct, and let \(G_\mu\) be the canonical
function of
\[
 M_\mu=S_\mu H^2\cap\D.
\]
The theorem gives compact sets \(K_{\mu,\ell}\) and outer canonical
functions \(H_{\mu,\ell}\) with
\[
 [H_{\mu,\ell}]=\D_{K_{\mu,\ell}},\qquad
 0<\CR\bigl(\Zrad(H_{\mu,\ell})\bigr)\le\frac{5N}{\ell},
 \qquad \|H_{\mu,\ell}-G_\mu\|_{\D}\longrightarrow0.
\]
For each \(y\in\D\), the projections satisfy
\[
 \|P_{\D_{K_{\mu,\ell}}}y-P_{M_\mu}y\|_{\D}
 \longrightarrow0.
\]
Here \(\Zrad(H)\) is the set of boundary points at which the radial
limit of \(H\) is zero, as in \eqref{eq:radial-zero}.
The capacity \(\CR\) is the obstacle capacity defined in
Section~\ref{subsec:03-norm}. Its comparison with classical
logarithmic capacity in Section~\ref{subsec:03-pairing} gives
vanishing logarithmic capacity along this sequence.

Each \(H_{\mu,\ell}\) extends continuously to the closed disk. Its
zero set there is \(K_{\mu,\ell}\), so
\(\Zrad(H_{\mu,\ell})=K_{\mu,\ell}\). In particular, the radial zero
set has Hausdorff distance at most \(3\,2^{-\ell}\) from
\(\supp\mu=\{\xi_1,\ldots,\xi_N\}\).

\medskip
\textbf{Loss of outerness in the limit.} The exact inner factor of \(G_\mu\) is the nontrivial singular
inner function \(S_\mu\).
 This sequence has two properties
that matter for the construction.

\begin{itemize}
\item Each approximant is outer, and its radial zero set has positive
capacity that can be made arbitrarily small by taking \(\ell\) large.
\item The Dirichlet-norm limit has a nontrivial singular inner factor.
\end{itemize}

In particular, strong convergence in the Dirichlet norm can lose
outerness.

\medskip
\textbf{Main steps of the proof.} The proof also proceeds in three steps.

\begin{enumerate}
\item \textbf{Recovery of a single atom.}
For a prescribed boundary atom, Theorem~\ref{thm:G} supplies Carleson sets of small positive capacity. The analytic criterion in
Theorem~\ref{thm:C} gives strong convergence of their outer canonical
functions to the single-atom canonical function, as stated in
Corollary~\ref{cor:R}.

Fix the atom's positive mass and location first. Choose the prescribed
value at zero in Theorem~\ref{thm:G} to be the value at zero of this
fixed atomic canonical target. That value lies strictly between zero
and one, so the geometric theorem applies. Its compacts are nonpolar
Carleson sets shrinking to the prescribed location in Hausdorff distance, and their canonical functions have exactly the target's value at zero
at every index. To use these scalar equalities, fix an arc away from
the atom. Richter and Yilmaz obtain continuation and nonvanishing
for a fixed compact set \cite[proof of Theorem~1.3]{RY19}.
Propositions~\ref{prop:common-nonzero-collar}
and~\ref{prop:fixed-analytic-neighborhood} give a lower bound
and a continuation neighborhood that work for all sufficiently
large indices. These estimates control the derivatives away
from the atom, while the small-region estimate controls the remaining
energy. Together these bounds yield strongly convergent subsequences.
Cyclicity of the outer factor identifies the subspace generated by each limit.
The prescribed value at zero then identifies its atomic mass, so all
such limits equal the canonical generator of the given atom.

\item \textbf{Combination of finitely many atoms.}
Section~\ref{sec:12} uses finite products and the trace-zero spaces
of the unions to recover \(M_\mu\). Each union has \(\CR\)-capacity at most \(5N/\ell\).
Section~\ref{subsec:11-realization} gives the closed-disk zero set
of the canonical function of each union.

\item \textbf{Convergence of the canonical functions.}
The relevant orthogonal projections converge strongly on \(\D\).
Applying them to \(1\) and normalizing the projected vectors gives
the asserted convergence of the canonical functions.

The normalization passes to the limit because the limiting projection
is nonzero. Lemma~\ref{lem:atomic-witness} gives a function in \(M_\mu\) that does not vanish at the origin. Evaluation at the origin is represented by
the projection of the constant function onto that subspace, so the
projection cannot vanish. Strong convergence of the projections applied
to the constant therefore gives vector convergence to a nonzero vector.
Their norms converge to its strictly positive norm and are consequently
bounded away from zero eventually. Division by these norms preserves
convergence along the same recovery family and retains unit norm
and positivity at zero. Every approximant's canonical vector is already
defined; this positive lower bound explains the continuity of its
normalization at the limit.
\end{enumerate}

\medskip
\textbf{Role in the counterexample.} 
Starting with \(e\), the first theorem selects \(G_\mu\), which has
a finite atomic singular inner factor. With \(\mu\) fixed, the second
theorem supplies \(H_{\mu,\ell}\):
\[
 e\xrightarrow{\text{Theorem~\ref{thm:A}}}G_\mu
  \xrightarrow{\text{Theorem~\ref{thm:B}}}H_{\mu,\ell},
\]
where \(e\) is outer and \(G_\mu\) has a finite atomic singular
inner factor. The second theorem provides the recovery
\[
 G_\mu\xleftarrow{\text{Theorem~\ref{thm:B}}}H_{\mu,\ell},
\]
where \(H_{\mu,\ell}\) is outer and has radial zero set of arbitrarily
small capacity as \(\ell\) increases. The Dirichlet triangle
inequality gives
\[
 \|H_{\mu,\ell}-e\|_{\D}
 \le\|H_{\mu,\ell}-G_\mu\|_{\D}+\|G_\mu-e\|_{\D}.
\]
The second term is controlled by the choice of \(\mu\), and the first
by the subsequent choice of \(\ell\). This proves Dirichlet-norm
density of outer canonical functions with arbitrarily small
radial-zero capacity. The constant function \(1\) belongs to each
strict capacity sublevel because its radial zero set is empty.

\medskip
\textbf{The \(L^1\) distance between the boundary logarithms.} For this same fixed nonzero \(\mu\), the approximants also satisfy
\[
 \lim_{\ell\to\infty}
 \left\|\log|H_{\mu,\ell}^*|-\log|G_\mu^*|\right\|_{L^1(m)}
 =\mu(\T)>0.
\]
Thus the convergence in \(\D\) in this recovery is accompanied by
a positive limiting distance between the boundary logarithms.
The final construction uses additional control of the boundary
logarithms to keep its limit outer.

\textbf{Concentration of boundary logarithms.}
Keep \(\mu\) fixed and write \(G_\mu=S_\mu O_\mu\), with
\(O_\mu(0)>0\). The logarithmic measures converge as \(\ell\to\infty\):
\begin{equation}\label{eq:finite-atomic-log-measure-limit}
 \bigl(-\log|H_{\mu,\ell}^*|\bigr)\,dm
 \overset{*}{\rightharpoonup}
 \bigl(-\log|O_\mu^*|\bigr)\,dm+\mu.
\end{equation}
Proposition~\ref{prop:finite-recovery-logarithmic-defect} proves
\eqref{eq:finite-atomic-log-measure-limit}. The outer factor
\(O_\mu\) determines the absolutely continuous part of the limit.
The singular part is the prescribed measure \(\mu\).

In particular, the second theorem replaces the finite atomic
singular target by outer canonical functions of trace-zero spaces
with small positive radial-zero capacity. It retains both the
canonical form of the approximants and their Dirichlet-norm
convergence to the target, and supplies the outer approximants used
in the counterexample.

\subsection{From approximation to the counterexample}

The positive limiting \(L^1(m)\) distance between the boundary logarithm
of \(G_\mu\) and those of its outer approximants explains the order of
the final argument. We combine Approximations~A and~B in the Dirichlet
norm, then compare the boundary logarithms of the original outer target
and the outer canonical function selected in the second step. This gives
the density used below to construct nested closed balls.

\medskip
\textbf{Density from the two approximation theorems.} For each proper outer canonical target \(e\), the two theorems give
\[
 e\xrightarrow{\text{Theorem~\ref{thm:A}}}G_\mu
   \xrightarrow{\text{Theorem~\ref{thm:B}}}H_{\mu,\ell}.
\]
The outer canonical function \(H_{\mu,\ell}\) can be chosen
arbitrarily close to \(e\) in the Dirichlet norm, with
\[
 \CR\bigl(\Zrad(H_{\mu,\ell})\bigr)\le \frac{5N}{\ell}.
\]
Proposition~\ref{prop:outer-continuity} compares the boundary
logarithms of \(e\) and \(H_{\mu,\ell}\).
For a fixed \(e\), making the Dirichlet error sufficiently small
also makes their \(\dlog\)-distance small. Consequently, for every
\(\varepsilon>0\) the outer canonical functions with
\(\CR\)-capacity of their radial zero sets less than \(\varepsilon\)
are dense in the space \((\X,\dlog)\) of normalized outer canonical functions, as stated in
Proposition~\ref{prop:small-capacity-density}.

For this density step, fix the proper outer target, the desired metric
accuracy, and the capacity bound. The fixed-target logarithmic estimate
first supplies an allowable Dirichlet error, before either approximation
is applied. Allow half of this error bound for each approximation.
Approximation~A chooses the measure; then keep all its masses, locations,
and its number of atoms fixed. One sufficiently large index in that
measure's recovery family satisfies both the remaining bound on the norm error and the capacity bound. The required index may depend on the fixed measure
and the chosen family. The Dirichlet triangle inequality through the
singular intermediate target now gives the required norm error.
The logarithmic estimate is applied directly to \(e\) and \(H_{\mu,\ell}\),
which gives the desired metric accuracy.

\medskip
\textbf{Convergence through nested closed balls.} We work in the space \(\X\) of normalized outer canonical functions
with the metric
\[
 \dlog(f,g)=\|f-g\|_{\D}
   +\|\log|f^*|-\log|g^*|\|_{L^1(m)}.
\]
The logarithmic term permits reconstruction of the outer function
at the limit. Theorem~\ref{thm:outer-completeness} gives completeness
of \((\X,\dlog)\).

Using density to choose the centers and
Corollary~\ref{cor:capacity-ball} to choose their radii, we construct
nested closed balls
\[
 B_1\supset B_2\supset B_3\supset\cdots
\]
whose radii tend to zero and for which
\[
 g\in B_n\quad\Longrightarrow\quad
 \CR(\Zrad(g))<2^{-n}.
\]
Completeness gives a unique common point \(f\). Closedness puts
\(f\) in every \(B_n\), so all of the capacity bounds apply to the
same function. It follows that
\[
 \CR(\Zrad(f))=0.
\]
The bounds throughout the balls are essential to this argument.
The capacity of a radial zero set can increase at a limit, as
discussed in Section~\ref{subsec:03-neighborhood}. The closed-ball
bounds apply directly to the limit, including when it lies on the
boundary of one of the balls.

\medskip
\textbf{Preserving noncyclicity.} The initial ball is centered at a proper outer canonical function
\(f_0\). In the numerical specialization of
Corollary~\ref{cor:numerical}, the initial function is chosen so that
\[
 f_0(0)=\frac12.
\]
In either case the initial radius is chosen to keep the limit
sufficiently close to \(f_0\). The bounds on evaluation at zero
then give
\[
 0<f(0)<1.
\]
The limit is canonical. Writing \(a=f(0)\), the 
identity \eqref{eq:certificate} gives
\[
 P_{[f]}1=af,
\]
and therefore
\[
\dist_{\D}(1,[f])^2=1-a^2>0.
\]
The positive distance excludes \(1\) from \([f]\), which proves
noncyclicity.

\medskip
\textbf{The resulting counterexample.} The equivalence of polar sets in Section~\ref{subsec:03-pairing}
turns the zero \(\CR\)-capacity of the radial zero set into zero
classical logarithmic capacity. The resulting \(f\in\D\) satisfies
\[
f\text{ is outer};\qquad
f\text{ is noncyclic};\qquad
\operatorname{cap}_{\log}\bigl(\Zrad(f)\bigr)=0.
\]
This is the counterexample in Theorem~\ref{thm:main} to the
Brown--Shields question in the  radial-zero formulation stated
in the introduction.

The two approximation theorems supply density. The limiting argument
uses that density to choose centers, and the bounds on the nested
closed balls give a common limit with all three properties. In
particular, the capacity bounds at the centers permit the choice
of the balls, and the bounds throughout the balls give capacity
zero at their common point.

\newpage \section{Canonical functions in the Dirichlet space}\label{sec:02}

Projecting \(1\) onto a closed invariant subspace gives its best
approximation to \(1\). When the projection is nonzero, normalizing
it gives a canonical function that generates the subspace. We start
with the norm and evaluation formulas needed for that projection.
For a proper outer canonical function, the estimates below lead to
the characterization of its generated space by the finiteness of a
weighted boundary integral in Section~\ref{sec:05}.  We end by constructing the finite atomic
targets that connect Approximations~A and~B.

\subsection{The Dirichlet norm, point  evaluation, and invariant subspaces}\label{subsec:02-norm}

We will use the coefficient norm to control point evaluation and
multiplication by \(z\). Use the norm in \eqref{eq:full-norm}
and take all inner products linear in their first argument.
For $f,g\in\D$ with $f(z)=\sum_{n\ge0}a_nz^n$ and
$g(z)=\sum_{n\ge0}b_nz^n$,
\[
 \int_{\DD}f'\overline{g'}\,dA_0
   =\sum_{n\ge1}n a_n\overline{b_n},\qquad
 \langle f,g\rangle_{\D}=\sum_{n\ge0}(n+1)a_n\overline{b_n}.
\]
Indeed, orthogonality in the area integral for $D(f)$ gives
the factors $2n^2\int_0^1r^{2n-1}\,dr=n$. The squared identity
follows by Tonelli and the polarized one follows by polarization.
The coefficient norm identifies $\D$ with a weighted sequence Hilbert
space and gives polynomial density by truncation. For $0<r<1$, the
dilation $f(r\,\cdot)$ converges in $\D$ to $f$ as $r\uparrow1$,
by dominated convergence in that same coefficient norm.

In the coefficient inner product, evaluation at a point is an
inner product with an explicit vector. For $a\in\DD$ this vector is
\[
 k_a(z)=\sum_{n\ge0}\frac{\bar a^{\,n}z^n}{n+1},\qquad
 \langle f,k_a\rangle_{\D}=f(a),\qquad k_0=1.
\]
The same kernel normalization is used in \cite[p.~143]{GZ10}.
The coefficient formula gives
\[
 K(r):=\|k_{r\xi}\|_{\D}^2
   =-\frac{\log(1-r^2)}{r^2}
 \quad(0<r<1,\ \xi\in\T),\qquad K(0)=1.
\]
Cauchy--Schwarz gives $\sup_{|z|\le R}|f(z)|\le\sqrt{K(R)}\|f\|_{\D}$
for $0\le R<1$. Thus norm convergence implies locally uniform
convergence, and Cauchy's formula gives the same conclusion for
derivatives.

For each
$k\ge0$, direct calculation gives
\[
 \|z^kf\|_{\D}^2=\|f\|_{\D}^2+k\|f\|_{H^2}^2.
\]
Since $\|f\|_{H^2}\le\|f\|_{\D}$, the norm of multiplication by $z^k$
is $\sqrt{k+1}$, with equality on constants. In particular, all
polynomials are bounded multipliers, and the shift $M_zf=zf$ is bounded
below. Hence $zM$ is closed for closed $M\subset\D$. The space $[f]$ in
\eqref{eq:module} is the smallest closed shift-invariant subspace
containing $f$. If it contains $1$, it contains every polynomial, so
$[f]=\D$.

In the recovery argument of Section~\ref{sec:09}, we start with
weakly convergent sequences in \(\D\). Weak convergence gives
bounded norms and coefficient convergence.
The coefficient estimate
\[
 \sum_{n>N}|a_n|^2\le\frac{\|f\|_{\D}^2}{N+2}
\]
then gives strong convergence in $H^2$. The same estimate, together with
finite-dimensional approximation by coefficient truncations, proves
compactness of the inclusion $\D\hookrightarrow H^2$. The Hardy
evaluation bound
$\sup_{|z|\le R}|f(z)|\le(1-R^2)^{-1/2}\|f\|_{H^2}$ and Cauchy's
formula now give locally uniform convergence of the functions and their
derivatives for every weakly convergent sequence in $\D$. We also use
the Hilbert-space facts that bounded sequences have weakly convergent
subsequences and that weak convergence together with convergence of
norms implies strong convergence.

\subsection{Canonical projections and generation}\label{subsec:02-projection}

The function \(1\) represents evaluation at zero in \(\D\).
Projecting it onto a subspace gives the vector that represents the
same evaluation there. Let $M$ be a closed linear subspace of
$\D$, and let $P_M$ be the orthogonal projection onto $M$.
Evaluation at $0$ on $M$ is represented by $y=P_M1$, since
$\langle h,y\rangle_{\D}=\langle h,1\rangle_{\D}=h(0)$ for every
$h\in M$. Thus $P_M1\ne0$ exactly when some member of $M$ is nonzero
at $0$. In that case put $a=\|y\|_{\D}>0$. Taking $h=y$ gives
$y(0)=a^2$, so
\begin{equation}\label{eq:canonical-definition}
 G_M=\frac{P_M1}{\|P_M1\|_{\D}},\qquad
 G_M(0)=\|P_M1\|_{\D}=a.
\end{equation}
We call $G_M$ the unit canonical function of $M$, normalized to have positive value at $0$.
It is the extremal for evaluation at zero in the Dirichlet norm \cite[Section~2]{HR26}.
The identity $h(0)=a\langle h,G_M\rangle_{\D}$ shows that $a$ is the
maximum of $|h(0)|$ on the unit ball of $M$, and that $G_M$ is the unique
unit maximizer with positive value at $0$.

When $M$ is also shift-invariant, its shifted subspace \(zM\)
consists of functions vanishing at zero. Thus
$\langle zh,y\rangle_{\D}=(zh)(0)=0$
for $h\in M$, so $y$ belongs to $M\ominus zM$. The dimension of this
space follows from the codimension-one theorem of Richter and Shields
\cite[Theorem~2(c)]{RSh88}. Generation is due to Richter
\cite[Theorem~1]{Richter88}. The two statements are collected in
\cite[Theorem~2.1]{GZ10} and \cite[Theorem~3, p.~135]{EFKR10}.

\begin{published}[Wandering generation]\label{input:generation}
For every nonzero closed shift-invariant $M\subset\D$, the wandering
space $M\ominus zM$ is one-dimensional, and every nonzero member of it
generates $M$ under polynomial multiplication.
\end{published}

Consequently $[G_M]=M$ whenever $M$ is invariant and $P_M1\ne0$.
Moreover, $0<a\le1$ because projection is contractive and
$\|1\|_{\D}=1$. Equality holds exactly when $1\in M$, hence exactly
when $M=\D$ by invariance and polynomial density.

When \(P_{M_j}1\) converges in norm to a nonzero vector \(P_M1\),
the projections are eventually nonzero, and their normalized
canonical functions converge to \(G_M\).
Section~\ref{sec:12} uses this step for the projections onto the
trace-zero spaces.

\subsection{Moments and the distance from the constant}\label{subsec:02-certificate}

The projection formula gives norm and orthogonality identities
that pass to a strong limit. For $f=G_M$ in an invariant subspace,
$z^nf\in zM$ for $n\ge1$, so orthogonality to \(zM\) gives
\begin{equation}\label{eq:canonical-moments}
 \|f\|_{\D}=1,\qquad f(0)=a>0,\qquad
 \langle z^nf,f\rangle_{\D}=0\quad(n\ge1).
\end{equation}
Conversely, let $f\in\D$ satisfy \eqref{eq:canonical-moments}. Then
$1-af$ is orthogonal to $f$ because $\langle1,f\rangle_{\D}=a$.
It is orthogonal to every $z^nf$ with $n\ge1$ because
$z^nf$ vanishes at $0$ and the corresponding moment vanishes. Since
$af\in[f]$, this proves
\begin{equation}\label{eq:certificate}
 P_{[f]}1=af,\qquad
 \operatorname{dist}_{\D}(1,[f])^2=1-a^2.
\end{equation}
In particular $f=G_{[f]}$. The orthogonal decomposition
$1-pf=(1-af)+(a-p)f$ gives
\begin{equation}\label{eq:polynomial-error}
 \|1-pf\|_{\D}^2=1-a^2+\|(p-a)f\|_{\D}^2
 \qquad(p\in\mathbb C[z]).
\end{equation}
The constant polynomial $p=a$ attains the minimum. For a function
satisfying \eqref{eq:canonical-moments}, cyclicity is equivalent to
$a=1$, in which case the coefficient norm forces $f=1$.

In the language of optimal polynomial approximation from
Section~\ref{subsec:01-generation-history},
\eqref{eq:polynomial-error} identifies the answer at every finite
degree. Among polynomials of degree at most any prescribed
nonnegative integer, the unique minimizer is the constant
\(p=a\). Increasing the degree cannot improve its error.
These are the constant approximant and squared-distance identities of
B\'en\'eteau, Fleeman, Khavinson, Seco, and Sola
\cite[Theorem~2.3(ii)--(iii)]{BFKSS18}, in the phase convention
\(f(0)=a>0\). The moment
hypothesis is essential, as the following example shows.

The cyclic function \(1-z\) from Section~\ref{subsec:01-problem}
illustrates the role of the moments. Its normalization
$v(z)=(1-z)/\sqrt{3}$ is outer by the boundary factor calculation
below, has unit norm, and has $0<v(0)<1$. Here
$\langle zv,v\rangle_{\D}=-2/3$. The moment hypothesis singles
out the functions for which the value at zero determines the approximation error.

If \(f_j\to f\) in \(\D\), each fixed moment converges by boundedness of
multiplication by \(z^n\). Thus unit norm and the vanishing moments pass
to a strong limit. Once \(f(0)>0\) has been established,
\eqref{eq:canonical-moments} and \eqref{eq:certificate}
give the same distance formula for the limit. In
Section~\ref{sec:13} we keep the value at zero bounded away from one
on the first closed ball. This gives the same positive lower bound
for the squared error throughout the ball.

\subsection{Hardy factorization and quotients}\label{subsec:02-hardy}

We will need quotients whose Hardy norms can be read from boundary
values. The outer integral gives the regularization used in the
division lemma below. We use the standard Hardy boundary and
factorization theorems
\cite[Theorems~2.2, 2.6, and~2.8]{Duren70}. A nonzero $f\in H^2$ has a nontangential boundary
value $f^*$ almost everywhere, with $|f^*|>0$ almost everywhere and
$\log|f^*|\in L^1(m)$. Its Hardy norm is $\|f^*\|_{L^2(m)}$.
It factors as $f=IO$, where $I$ is inner and $O$ is outer. We fix
the phase by requiring $O(0)>0$, so $O$ is given by
\eqref{eq:outer-intro} with $\omega=1$ and boundary modulus $|f^*|$.
Conversely, for real $u\in L^1(m)$ with $e^u\in L^2(m)$, the
exponential of the integral in \eqref{eq:outer-intro} with $u$
in place of $\log|f^*|$ is an outer $H^2$ function. Its boundary
modulus is $e^u$ almost everywhere, and its value at zero is
$\exp(\int_{\T}u\,dm)$.
When $f(0)>0$ this also gives $I(0)>0$. The inner part consists of
a constant of modulus one, a Blaschke factor, and a singular
inner factor. Since $|I^*|=1$ almost everywhere, multiplication by
$I$ on $H^2$ is an isometry and $IH^2$ is closed. For boundary
integrals in this section, $f^*$ denotes the arclength equivalence
class. Section~\ref{sec:03} introduces the representative used in
capacity estimates.

The boundary factor $q_\xi(z)=1-\bar\xi z$ is outer for each $\xi\in\T$.
The mean of $\log|q_\xi(r\zeta)|$ over $\zeta\in\T$ is zero for every
$r<1$, by harmonicity. For $r\ge1/2$,
$|1-r\bar\xi\zeta|^2=(1-r)^2+r|1-\bar\xi\zeta|^2$ gives an integrable
lower bound on the logarithm, while its positive part is bounded by
$\log 2$. Dominated convergence gives
$\int_{\T}\log|q_\xi^*|\,dm=0=\log|q_\xi(0)|$. The inner factor
therefore has modulus one at $0$ and is constant.

Section~\ref{sec:05} estimates a quotient first on the circle.
The next lemma turns that estimate into membership in \(H^2\) inside
the disk. It is the bounded outer-denominator case of the Smirnov
maximum principle \cite[Theorem~2.11, p.~28]{Duren70}.

\begin{lemma}[Division by a bounded outer function]\label{lem:bounded-outer-division}
Let $q$ be nonzero, bounded, and outer, and let $h\in H^2$. If
$h^*/q^*\in L^2(m)$, then the analytic quotient $h/q$ belongs to $H^2$,
with
\[
 \|h/q\|_{H^2}^2=\int_{\T}|h^*/q^*|^2\,dm.
\]
\end{lemma}
\begin{proof}
We first use outer denominators bounded away from zero. For
$0<\varepsilon<1$, let $q_\varepsilon$ have phase $q(0)/|q(0)|$ and
the outer modulus $\max\{|q^*|,\varepsilon\}$. The clipped logarithm
converges in $L^1(m)$ to $\log|q^*|$, since their difference is bounded
by the integrable negative part of $\log|q^*|$. The outer integral
therefore gives $q_\varepsilon\to q$ locally uniformly. The Poisson
formula for its logarithmic modulus gives $|q_\varepsilon|\ge\varepsilon$
in the disk, so $h/q_\varepsilon\in H^2$. Its squared Hardy norm is
\[
 \int_\T\frac{|h^*|^2}{\max\{|q^*|,\varepsilon\}^2}\,dm
 \le \int_\T|h^*/q^*|^2\,dm.
\]
For each fixed radius $r<1$, local uniform convergence of the
quotients passes this upper bound to the circle mean of $|h/q|^2$.
Taking the supremum over $r$ proves $h/q\in H^2$. Its boundary values
are $h^*/q^*$ almost everywhere, which proves the norm identity.
\end{proof}

\subsection{Local Dirichlet integrals and contractive inner division}\label{subsec:02-energy}

The weighted boundary integral in Section~\ref{sec:05} uses the local energy of
the canonical generator at each point of the circle. The local
Douglas formula of Richter and Sundberg expresses this energy using
boundary values or the derivative in the disk \cite[Proposition~2.2, pp.~359--360]{RS91}. The
boundary formula also lets us compute the energy of a product with
an inner factor in Section~\ref{sec:07}. We use the normalization of
\cite[(2.1)--(2.2), pp.~2121--2122]{RY19}.

Richter and Sundberg's local theory isolates the
energy associated with one boundary point. Their restatement
\cite[p.~168]{RS92} gives two complementary descriptions: the
area formula weights the derivative energy by the Poisson kernel
of that point, while the boundary formula measures the change
from the specified boundary limit through a squared difference
quotient. Averaging these local energies over the circle recovers
the ordinary Dirichlet energy. A weighted average instead produces
the Poisson-weighted derivative energy used in the estimates
for the weighted boundary integral. These formulas concern the energy term; the Dirichlet
norm also retains the Hardy norm, with the normalization used by Richter and Yilmaz \cite{RY19}. Finite local energy
ensures the pointwise limit needed in the boundary formula.
The identities also allow infinite energy as a nonnegative
extended integral. The division argument following the local
Dirichlet integral formula separately establishes that the
difference quotient belongs to the Hardy space.

\begin{published}[Local Dirichlet integral]\label{input:local-douglas}
For $F\in H^2$ and $\xi\in\T$, define
\[
 D_\xi(F)=\int_{\DD}\frac{1-|z|^2}{|\xi-z|^2}|F'(z)|^2\,dA_0(z).
\]
If $F$ has nontangential limit $F(\xi)$ at $\xi$, then
\begin{equation}\label{eq:local-douglas}
 D_\xi(F)=\int_{\T}
 \frac{|F^*(\zeta)-F(\xi)|^2}{|\zeta-\xi|^2}\,dm(\zeta)
\end{equation}
in the nonnegative extended sense. Finite $D_\xi(F)$ ensures that
this limit exists.
\end{published}

For finite $D_\xi(F)$, the difference quotient
$(F-F(\xi))/(z-\xi)$ belongs to $H^2$ and has squared norm
$D_\xi(F)$. Apply Lemma~\ref{lem:bounded-outer-division} to the
numerator $F-F(\xi)$ and the bounded outer denominator $z-\xi$,
using \eqref{eq:local-douglas}. At a fixed boundary point the symbol
$F(\xi)$ refers to the specified limit. In an integral against
arclength it agrees with the Hardy boundary value almost everywhere.

For a nonnegative $v\in L^1(m)$, define
$P[v](z)=\int_{\T}(1-|z|^2)|\xi-z|^{-2}v(\xi)\,dm(\xi)$.
Tonelli gives the weighted identity
\begin{equation}\label{eq:weighted-local-douglas}
 \int_\T D_\xi(F)v(\xi)\,dm(\xi)
 =\int_{\DD}|F'(z)|^2P[v](z)\,dA_0(z),
\end{equation}
also with infinite values allowed. For $v=1$ it gives
$\int_\T D_\xi(F)\,dm(\xi)=D(F)$.

The Dirichlet norm can only decrease when we remove an inner
factor with a Hardy quotient. We will use this fact to show that
an extremal for an inner-factor subspace carries exactly the inner
factor prescribed by that subspace. The norm inequality also follows
from Carleson's energy formula \cite{Carleson60}, in the Dirichlet-norm
form of \cite[(3.2)]{HR26}. The backward-shift proof shows that the
Hardy adjoint of inner multiplication contracts the whole of $\D$.

\begin{lemma}[Contractive inner division]\label{lem:inner-division}
If $I$ is inner, $h\in H^2$, and $Ih\in\D$, then
$h\in\D$ and $\|h\|_{\D}\le\|Ih\|_{\D}$.
\end{lemma}
\begin{proof}
We use the Hardy backward shift $B$ to compare the Dirichlet
norms. Counting the occurrence of each Taylor coefficient in the
tails gives
\[
 \|v\|_{\D}^2=\sum_{k\ge0}\|B^kv\|_{H^2}^2\qquad(v\in H^2),
\]
with the value $+\infty$ allowed. Let $T_I$ denote multiplication by
$I$ on $H^2$ in this proof. Its Hardy adjoint $T_I^*$ is a
contraction and commutes with $B$, since $T_I$ commutes with the
forward shift. Consequently
\[
 \sum_{k\ge0}\|B^kT_I^*v\|_{H^2}^2
 \le\sum_{k\ge0}\|B^kv\|_{H^2}^2.
\]
We apply this to $v=Ih$, for which $T_I^*v=h$ by the Hardy isometry.
\end{proof}

This argument also shows that the Hardy adjoint of multiplication by
$I$ restricts to a contraction on $\D$. Later, when a product $Ih$ is constructed with $h\in H^2$, its Dirichlet membership will be proved by a separate energy formula \eqref{eq:finite-product} in
Section~\ref{subsec:07-product}.

\subsection{Canonical multipliers and Hardy quotients}\label{subsec:02-multiplier}

We now use the canonical moments to estimate products \(ep\) with a
polynomial \(p\). An upper bound will extend multiplication by \(e\)
to a contraction on \(\D\). A lower bound in the Hardy norm will
give a quotient for each member of \([e]\). We first use the moments
to link the boundary modulus of \(e\) to its derivative energy.
The polynomial norm identity is the restriction of the weighted isometry in
\cite[(3.1)]{RS92}. We give the coefficient
calculation to make the norm identity and the \(H^2\) estimate for the
quotient explicit. See also \cite[Theorem~0.2]{Wu98}.

For $f\in\D$, define the function \(\ell_f\) and measure \(\sigma_f\) by
\begin{equation}\label{eq:source}
 \ell_f(\xi)=\int_{\DD}P_z(\xi)|f'(z)|^2\,dA_0(z),
 \qquad \sigma_f=\ell_f\,m,
 \qquad P_z(\xi)=\frac{1-|z|^2}{|\xi-z|^2}.
\end{equation}
Thus \(\ell_f(\xi)\) is the local energy \(D_\xi(f)\) from the
preceding subsection. Tonelli gives $\int_\T\ell_f\,dm=D(f)$.
Writing $\|\cdot\|_{2,A}$
for the norm in $L^2(\DD,dA_0)$, integration and Cauchy--Schwarz give
\[
 \|\ell_f-\ell_g\|_{L^1(m)}
 \le (\|f'\|_{2,A}+\|g'\|_{2,A})\|f'-g'\|_{2,A}.
\]
Now let $e(z)=\sum_{j\ge0}c_jz^j$ have unit Dirichlet norm and
satisfy $\langle z^ne,e\rangle_{\D}=0$ for every $n\ge1$.
Put $w=|e^*|^2$.
For $v\in L^1(m)$ let $\widehat v(n)=\int_\T v(\xi)\bar\xi^n\,dm(\xi)$.
For $n\ge0$,
\[
 \widehat w(n)=\sum_{j\ge0}c_{j+n}\bar c_j,\qquad
 \widehat\ell_e(n)=\sum_{j\ge1}j c_{j+n}\bar c_j.
\]
For the second formula, the $n$th boundary coefficient of $P_z$
is $\bar z^n$. With polynomial $e$, angular integration of
$\bar z^n|e'(z)|^2$ leaves the terms with index difference $n$;
their radial integral gives the factor $j$. Taylor truncation passes
to the limit in both formulas, using $H^2$ convergence for $w$ and the
displayed $L^1$ estimate for $\ell_e$. The coefficient sums converge
absolutely by weighted Cauchy--Schwarz. With the Dirichlet inner product,
\[
 (n+1)\widehat w(n)+\widehat\ell_e(n)
 =\sum_{j\ge0}(j+n+1)c_{j+n}\bar c_j
 =\langle e,z^ne\rangle_{\D}=\delta_{n0},
\]
where $\delta_{n0}$ is $1$ for $n=0$ and $0$ otherwise. Negative
frequencies follow by conjugation. If $\Lambda$ denotes the Fourier
multiplier with symbol $|n|$, this is the distributional identity
\begin{equation}\label{eq:canonical-source-identity}
 (1+\Lambda)w+\ell_e=1.
\end{equation}
This is the Fourier form of the resolvent identity of Richter and
Sundberg \cite[Lemma~2.1]{RS94}, as stated in
\cite[Lemma~2.4, (2.1)]{GZ10}. In that formula the two resolvent
terms give, respectively, the coefficients of $w+\ell_e$ and
$\Lambda w$. For \(f=g=e\), \eqref{eq:canonical-moments} gives
\(\langle e/(1-\bar\lambda z),e\rangle_{\D}=1\)
for every \(\lambda\in\DD\). 
The inverse of $1+\Lambda$ on Fourier coefficients is the positive
operator $Rv=\int_0^1 P_r*v\,dr$, where $P_r$ is the circle Poisson
kernel at radius $r$ and $*$ denotes circle convolution. Its multiplier
is $\int_0^1r^{|n|}\,dr=(1+|n|)^{-1}$. It fixes constants and maps
nonnegative $L^1$ functions to nonnegative $L^1$ functions. Thus
$w+R\ell_e=1$ almost everywhere, and $0\le w\le1$.
The Hardy Poisson bound gives $|e(z)|^2\le P[w](z)\le1$.

For $F\in H^2$ write
$D_w(F)=\int_{\DD}|F'|^2P[w]\,dA_0$, allowing infinity, and use the
corresponding polarized form $D_w(p,q)$ on polynomials. For $j\ge k\ge0$,
with $d=j-k$, the coefficient inner product gives
\[
 \langle z^je,z^ke\rangle_{\D}
 =\sum_{m\ge0}(m+j+1)c_m\overline{c_{m+d}}
 =\delta_{jk}+k\widehat w(-d).
\]
The first part of the sum with weight $m+d+1$ is the moment for shift
$d$. When $d=0$, this part is $\|e\|_{\D}^2$. Direct polar integration gives
$D_w(z^j,z^k)=k\widehat w(-d)$ as well. For $j,k>0$ the angular
integral selects $r^d\widehat w(-d)$ from $P[w](re^{it})$, and the
remaining radial factor integrates to
$2jk\int_0^1r^{2j-1}\,dr=k$. If $k=0$, the derivative of the constant
is zero. Conjugation treats $j<k$, so for all polynomials $p,q$,
\begin{equation}\label{eq:canonical-gram}
 \langle ep,eq\rangle_{\D}
   =\langle p,q\rangle_{H^2}+D_w(p,q),
 \qquad \|ep\|_{\D}^2=\|p\|_{H^2}^2+D_w(p).
\end{equation}
Since $0\le P[w]\le1$, this gives $\|ep\|_{\D}\le\|p\|_{\D}$.
For $g\in\D$, approximate $g$ in $\D$ by polynomials. Their products
with $e$ are Cauchy in $\D$, and local uniform convergence identifies
the limit as the analytic product $eg$. Hence
\begin{equation}\label{eq:canonical-multiplier}
 eg\in\D,\qquad \|eg\|_{\D}\le\|g\|_{\D}\qquad(g\in\D).
\end{equation}

This is Richter and Sundberg's contractive multiplier
conclusion for unit Dirichlet-inner functions \cite{RS92}:
the norm bound holds when multiplying every vector in the
Dirichlet space. The calculation uses the Dirichlet norm in \eqref{eq:full-norm}. Seco's account also credits
Aleman with an alternative proof of the contractivity
theorem \cite[p.~1655]{Seco19}. Shimorin develops a further approach through reproducing kernels and
extremal identities, using an argument that can be extended to
vector-valued spaces \cite{Shimorin98}. These approaches help
explain the strength of the conclusion: the shifted orthogonality
controls multiplication on the whole space. The exact isometry
in the Richter--Sundberg representation uses the weighted norm of the quotient by the generator; multiplication on the original Dirichlet space
has the contractive bound just obtained.

Seco proved the converse in the classical Dirichlet space
\cite[Theorem~2.1, p.~1656]{Seco19}. A function of
unit Dirichlet norm whose multiplication operator is contractive
has vanishing inner products with all of its positive shifts.
Testing multiplication on the constant function already makes
its multiplier norm at least one, so contractivity here means
that the multiplier norm is exactly one. Seco's proof tests
products with a positive power of the variable plus an arbitrary
complex constant. Varying that constant isolates the mixed inner
product and forces the shifted moment to vanish. This recovers
the orthogonality condition; choosing a positive value at the
origin is a separate normalization when that value is nonzero.

Section~\ref{sec:08} uses this contraction to test the energy near a
boundary point with an evaluation vector \(k_a\). The lower bound in
\eqref{eq:canonical-gram} controls the quotient when we pass from
polynomial multiples of \(e\) to their closure.

\begin{lemma}[Hardy quotient for a canonical generator]\label{lem:completed-canonical-quotient}
If $e$ satisfies \eqref{eq:canonical-moments} and $h\in[e]$, there is
a unique $F\in H^2$ such that $h=eF$, with
$\|F\|_{H^2}\le\|h\|_{\D}$.
\end{lemma}
\begin{proof}
To pass from polynomial multiples to their closure, choose polynomials $p_j$
such that $ep_j\to h$ in $\D$.
The Gram identity applied to $p_j-p_k$ makes $(p_j)$ Cauchy in
$H^2$, since
$\|p_j-p_k\|_{H^2}\le\|e(p_j-p_k)\|_{\D}$. We write $F$ for
its Hardy limit. We can now identify the two limits: both convergences
are locally uniform, which gives
$h=eF$ throughout the disk, including at any zeros of $e$.
Uniqueness follows from the identity theorem, and
$\|p_j\|_{H^2}\le\|ep_j\|_{\D}$ gives the bound on $F$.
\end{proof}

The recovery argument in Section~\ref{sec:12} uses the Hardy bound
on these quotients to pass to a weak limit. The next elementary
estimate shows that the product of finitely many bounded Dirichlet
functions with prescribed inner factors still belongs to the Dirichlet space.

\begin{lemma}[Bounded products and smooth multipliers]\label{lem:bounded-products}
An analytic function $v$ with $v,v'\in H^\infty$ multiplies $\D$
boundedly. If $v\in A^\infty$, meaning analytic in the disk and smooth
on the closed disk, multiplication by $v$ preserves every closed
shift-invariant subspace of $\D$. Any finite product of functions
in $\D\cap H^\infty$ belongs to $\D\cap H^\infty$.
\end{lemma}
\begin{proof}
We first estimate both parts of the Dirichlet norm.
The product rule and the bound
$\int_{\DD}|g|^2\,dA_0\le\|g\|_{H^2}^2$ show that
multiplication by $v$ is bounded on $\D$. Its operator norm is at
most $2\|v\|_\infty+\|v'\|_\infty$, by separately estimating
the Hardy and derivative norms. Here the area bound follows by
integrating the Hardy circle means against $2r\,dr$.
This operator estimate will also let us pass from polynomial to smooth
multipliers in the second claim.
For $v\in A^\infty$, its Taylor coefficients decay faster than any
fixed power by integration by parts on the circle. Its Taylor
polynomials therefore converge uniformly along with their first
derivatives. They converge in multiplier norm by the estimate just
proved, so every closed shift-invariant subspace is invariant under
multiplication by $v$. The last claim needs no smoothness of the factors.
Finally, any finite product of functions in
$\D\cap H^\infty$ belongs to $\D\cap H^\infty$ by the product
rule, since each derivative term has one area-$L^2$ factor and all
other factors bounded.
\end{proof}

\subsection{Exact inner factors}\label{subsec:02-factor}

We will compare canonical functions of spaces specified by an inner
factor. Hartz and Richter use division by an inner factor, based on
Carleson's energy formula, in
\cite[Proposition~3.1, proof of (ii)$\Rightarrow$(i)]{HR26}.
We combine that division step with the extremal property for evaluation at zero: dividing out an additional nonconstant inner factor would
increase the value at zero without increasing the norm.
Luo and Richter make this exact-factor observation in the proof of
\cite[Theorem~4.7]{LR15}. We give the argument with positive value
at zero, since that value will determine the recovery parameter.

\begin{lemma}[The inner factor of the canonical function for \(IH^2\cap\D\)]\label{lem:exact-inner-factor}
Let $I$ be inner with $I(0)>0$, and set $M_I=IH^2\cap\D$.
Suppose $M_I$ contains a function nonzero at $0$. Then $M_I$ is
closed and invariant, $[G_{M_I}]=M_I$, and the Hardy inner factor
of $G_{M_I}$, normalized to have positive value at zero, is exactly $I$.
\end{lemma}
\begin{proof}
The inclusion into $H^2$ is continuous and $IH^2$ is closed there.
This proves closedness in $\D$. Invariance follows from the shift
invariance of $IH^2$ and boundedness of the Dirichlet shift.
The hypothesis at $0$ ensures $P_{M_I}1\ne0$, and
Theorem~\ref{input:generation} gives generation.

The extremal property now determines the whole inner factor.
Write $G_{M_I}=ICO$ by factoring the Hardy quotient by $I$, where
$C$ is inner, $O$ is outer, and $C(0),O(0)>0$. Since $IO\in H^2$
and $C(IO)=G_{M_I}\in\D$, Lemma~\ref{lem:inner-division} puts
$IO$ in $M_I$ with norm at most one. If $C$ is nonconstant, then
$C(0)<1$ by the maximum modulus principle and
$(IO)(0)=G_{M_I}(0)/C(0)>G_{M_I}(0)$, contrary to extremality.
The remaining constant $C$ equals $1$ by the phase convention.
\end{proof}

\subsection{The finite atomic targets}\label{subsec:02-atomic}

Approximation~A produces a measure with finitely many atoms, and
Approximation~B starts from the canonical function of its
inner-factor subspace. We construct a function in that subspace
with nonzero value at zero so that the canonical function is
defined. For a finite positive singular measure $\mu$ on $\T$, the
singular inner factor in our convention is the one in
\cite[Section~1.1, p.~134]{EFKR10}
\[
 S_\mu(z)=\exp\left(-\int_\T\frac{\xi+z}{\xi-z}\,d\mu(\xi)\right).
\]
It satisfies $S_\mu(0)=e^{-\mu(\T)}$, and multiplying factors
adds their measures. For the zero measure take $S_0=1$. If
$\mu\ne0$, then $1/S_\mu\notin H^2$. Indeed, it would have
Hardy norm one by its almost-everywhere boundary modulus but its
value at zero is greater than one.

For a nonzero finite positive measure of finite support, write
$\mu=\sum_{i=1}^N\tau_i\delta_{\xi_i}$ with $N\ge1$,
$\tau_i>0$, and distinct $\xi_i\in\T$ after merging repeated
directions. The sum in the exponent shows that $S_\mu$ is analytic
through each complementary arc of $\{\xi_1,\dots,\xi_N\}$. Put
\begin{equation}\label{eq:atomic-target}
 S_\mu(z)=\exp\left(-\sum_{i=1}^N\tau_i\frac{\xi_i+z}{\xi_i-z}\right),
 \qquad M_\mu=S_\mu H^2\cap\D,\qquad G_\mu=G_{M_\mu}.
\end{equation}
The proof below establishes the nonzero evaluation required for
$G_\mu$ to exist. Write $m_\mu=\mu(\T)$ and
$E_\mu=\{\xi_1,\ldots,\xi_N\}$. For one atom use
$S_{\tau,\xi}=S_{\tau\delta_\xi}$ and
$G_{\tau,\xi}=G_{\tau\delta_\xi}$.

Hartz and Richter use the function below for a single atom at $\xi=1$
\cite[paragraph before Theorem~5.2]{HR26}. This is the case \(\beta=1\) of Brown and Shields' multiplier
example \cite[Proposition~16]{BS84}.  Its squared Dirichlet norm
also follows from Carleson's formula \cite{Carleson60}: the polynomial
\(1-\bar\xi z\) contributes \(3\), and its squared boundary modulus
cancels the denominator in the atomic weight, leaving \(2\tau\).
The direct calculation below gives a function in \(M_\mu\)
with nonzero value at zero before any recovery sequence is constructed.

\begin{lemma}[A function in \(M_\mu\) nonzero at zero]\label{lem:atomic-witness}
For $\tau>0$ and $\xi\in\T$, the function
$F_{\tau,\xi}(z)=(1-\bar\xi z)S_{\tau,\xi}(z)$ belongs to
$\D\cap H^\infty$, with
$F_{\tau,\xi}(0)=e^{-\tau}$ and
$\|F_{\tau,\xi}\|_{\D}^2=3+2\tau$.
For every measure $\mu$ as in \eqref{eq:atomic-target}, $M_\mu$
has nonzero evaluation at zero and is proper. Its canonical
function generates $M_\mu$ and has exact inner factor $S_\mu$.
\end{lemma}
\begin{proof}
Rotation preserves the coefficient norm, so take $\xi=1$. The
Hardy norm squared of $F_{\tau,1}$ is $2$ by the inner modulus.
The Cayley coordinate makes the exponent of the singular factor
linear and makes the energy integral explicit. Use the coordinate
$u=(1+z)/(1-z)=x+iy$, which sends the disk to $x>0$ and writes
$F_{\tau,1}$ as $2e^{-\tau u}/(u+1)$. The conformal change of
variables cancels the derivative of the coordinate map in the
area integral. With $a=x+1$, differentiation gives
\[
 D(F_{\tau,1})=
 \frac4\pi\int_0^\infty e^{-2\tau x}
 \int_{\mathbb R}\frac{(\tau a+1)^2+\tau^2y^2}{(a^2+y^2)^2}
 \,dy\,dx.
\]
The exponential weight now depends only on the real coordinate, while
the remaining vertical integral is elementary. Tonelli applies before
finiteness is known. Substitution
$y=a\tan\theta$ gives the two inner integrals
$\pi/(2a^3)$ and $\pi/(2a)$ with numerators $1$ and $y^2$,
respectively. Hence the energy is
\[
 -\int_0^\infty\frac{d}{dx}
 \left[e^{-2\tau x}\left(\frac{2\tau}{x+1}
                         +\frac{1}{(x+1)^2}\right)\right]dx
 =1+2\tau.
\]
This proves the asserted norm formula and membership.

We now multiply the single-atom functions $F_{\tau_i,\xi_i}$,
$i=1,\ldots,N$. Their product belongs to $\D\cap H^\infty$
by Lemma~\ref{lem:bounded-products}.
It is
$S_\mu\prod_{i=1}^N(1-\bar\xi_i z)$, so it lies in $M_\mu$ and
has value $e^{-m_\mu}>0$ at zero. The condition $1/S_\mu\notin H^2$
gives $1\notin M_\mu$. Lemma~\ref{lem:exact-inner-factor} now
gives generation and the exact inner factor.
\end{proof}

The unitary rotation $f(z)\mapsto f(\bar\xi z)$ fixes $1$ and
takes $M_{\tau\delta_1}$ onto $M_{\tau\delta_\xi}$. It therefore carries the projection of $1$ onto $M_{\tau\delta_1}$ to the projection
onto $M_{\tau\delta_\xi}$, so $G_{\tau,\xi}(0)$ is independent
of $\xi$. Section~\ref{subsec:08-target} shows that the value is strictly decreasing in $\tau$. Once the direction of the atom is known, the value at zero of the limiting canonical function
determines the mass in Section~\ref{sec:09}.

\newpage \section{Boundary traces, capacity, and trace-zero spaces}\label{sec:03}

Approximation~B uses closed spaces of Dirichlet functions that vanish
on prescribed compact subsets of \(\T\) with zero arclength and positive logarithmic capacity. As an \(L^2(m)\) class, the Hardy trace \(f^*\)
is defined up to arclength-null sets. To impose vanishing on these sets, we use the quasi-continuous representative \(\widetilde f\)
of that class. It agrees with the  radial limit outside a
set of capacity zero.

The boundary norm gives quantitative control of exceptional sets.
We use this control to prove closedness of the trace-zero spaces
and to bound the capacity of radial zeros throughout the closed balls
in Section~\ref{sec:13}. The norm also controls integration against
finite-energy measures, which will test traces on the Cantor sets in
Section~\ref{sec:10}.

\subsection{The harmonic norm and obstacle capacity}\label{subsec:03-norm}

Taking the modulus of an analytic boundary trace can introduce
negative Fourier modes. To handle this operation we work in a space
that includes all Fourier modes. Let $m$ be normalized arclength on
$\T$. We use the Hilbert space $V=H^{1/2}(\T)$ with the norm
\begin{equation}\label{eq:V}
 \|u\|_V^2=\sum_{n\in\mathbb Z}(1+|n|)|\widehat u(n)|^2,
 \qquad \widehat u(n)=\int_{\T}\bar\zeta^{\,n}u(\zeta)\,dm(\zeta).
\end{equation}
The Hardy boundary trace $f^*$ of $f\in H^2$ has only nonnegative
Fourier modes, so $f^*\in V$ exactly when $f\in\D$, with
$\|f^*\|_V=\|f\|_{\D}$.

We use the following normalized $L^2$ version of Douglas's difference
formula \cite[Section~21, pp.~307--311]{Douglas31}. For every
$u\in L^2(m)$, with infinity allowed, it gives
\begin{equation}\label{eq:V-difference}
 \|u\|_V^2=\|u\|_2^2+
 \iint_{\T^2}\frac{|u(\zeta)-u(\eta)|^2}{|\zeta-\eta|^2}
                    \,dm(\zeta)\,dm(\eta).
\end{equation}
Indeed, for each angle $h$ the squared $L^2$ norm of
$u(e^{ih}\,\cdot)-u$ is
$\sum_n|e^{inh}-1|^2|\widehat u(n)|^2$ by Parseval. Divide by
$|e^{ih}-1|^2$ and integrate in $h$ with measure $dh/(2\pi)$.
Tonelli's theorem applies to the nonnegative summands. The integrated
quotient for mode $n$ is $|n|$, by the finite geometric sum for
$n>0$, conjugation for $n<0$, and direct evaluation for $n=0$.
This proves \eqref{eq:V-difference}. In particular,
$\||u|\|_V\le\|u\|_V$. The same argument applies to any real
$1$-Lipschitz contraction fixing zero when $u$ is real.

We measure the size of an open set by the infimum of the squared norm
among real functions that are at least one almost everywhere on it.
For an open set $U\subset\T$, define
\[
 C_R(U)=\inf\{\|v\|_V^2:v\in V\text{ real},\ v\ge1\text{ a.e. on }U\}.
\]
For an arbitrary set $E\subset\T$ put
$C_R(E)=\inf_{U\supset E,\ U\text{ open}}C_R(U)$. Formula
\eqref{eq:V-difference} identifies this norm and capacity with
$D^h(m)$ and its open-obstacle capacity in
\cite[pp.~2121--2122]{RY19} and \cite[Section~5]{RS14}. In that
specialization the local integral averaged over $m$ is precisely
the double integral in \eqref{eq:V-difference}. The $L^2$ term gives
$C_R(E)\ge m^*(E)$, where $m^*$ is arclength outer measure, and the
constant obstacle $1$ gives $C_R(\T)=1$.

Richter and Sundberg place these harmonic weighted
spaces in the abstract Dirichlet-space tradition of Beurling and Deny \cite{BD59} in their account of capacity
\cite[Section~5]{RS14}. Their discussion explains how the norm,
real obstacles, and exceptional sets belong to one potential-theoretic
framework. The harmonic space is useful here because it admits
moduli and real truncations; the analytic Dirichlet space remains
the space in which we study polynomial generation and cyclicity.
Richter and Sundberg's account, together with that of Richter and Yilmaz \cite{RY19}, supplies the modern setting for the particular norm and capacity fixed above.

A set of $C_R$-capacity zero is called \emph{polar}. A property
holds \emph{quasi-everywhere}, abbreviated q.e., if it holds
outside a polar set. The logarithmically polar sets are exactly
the polar sets for $C_R$. We give the comparison in
Section~\ref{subsec:03-pairing}.

The trace conditions below each hold outside a polar set. The union
estimate lets us combine countably many such exceptions.

\begin{lemma}[Capacity of unions]\label{lem:capacity-unions}
For any finite family of subsets of $\T$,
\begin{equation}\label{eq:capacity-finite-subadditivity}
 C_R\!\left(\bigcup_{j=1}^N E_j\right)\le\sum_{j=1}^N C_R(E_j).
\end{equation}
For a countable family,
\begin{equation}\label{eq:capacity-root-subadditivity}
 C_R\!\left(\bigcup_{j\ge1}E_j\right)^{1/2}
     \le \sum_{j\ge1} C_R(E_j)^{1/2}.
\end{equation}
In particular a countable union of polar sets is polar. Changing
a set on a polar set leaves its capacity unchanged.
\end{lemma}
\begin{proof}
For finite unions, choose real obstacles $v_j$ on open
neighborhoods of $E_j$ and take $v=\max(0,v_1,\dots,v_N)$.
Pointwise,
\[
 |v(\zeta)|^2\le\sum_j|v_j(\zeta)|^2,\qquad
 |v(\zeta)-v(\eta)|^2\le\sum_j|v_j(\zeta)-v_j(\eta)|^2.
\]
Thus $v$ is an obstacle on the union of the neighborhoods with
$\|v\|_V^2\le\sum_j\|v_j\|_V^2$. Taking the infimum proves
\eqref{eq:capacity-finite-subadditivity}.

For a countable family with summable right side in
\eqref{eq:capacity-root-subadditivity}, choose nonnegative obstacles
$v_j$ on open neighborhoods with
$\|v_j\|_V\le C_R(E_j)^{1/2}+\varepsilon 2^{-j}$.
Their series converges in $V$. The sum represents the nonnegative
pointwise sum almost everywhere, since the $L^1(m)$ norms are
summable and the series converges in $L^2(m)$. It is an obstacle
on the union of the neighborhoods. The triangle inequality and
$\varepsilon\downarrow0$ prove \eqref{eq:capacity-root-subadditivity}.
The polar-set statements follow by subadditivity and monotonicity.
\end{proof}

\subsection{Representatives and radial zeros}\label{subsec:03-trace}

To estimate the capacity of near-zero trace levels in
Section~\ref{subsec:03-neighborhood}, we will use continuity on
the complement of a small open set.
A representative $\widetilde u$ of $u\in V$ is quasi-continuous if
for every $\varepsilon>0$ it is continuous on the complement of an
open set of capacity less than $\varepsilon$. We use the following
boundary theorem in the normalization just fixed.

\begin{published}[Boundary representatives]\label{input:trace}
Every $u\in V$ has a quasi-continuous representative, unique q.e.
For $f\in\D$ the complex radial limit exists q.e. and equals
the quasi-continuous representative $\widetilde f$ of $f^*$:
\begin{equation}\label{eq:radial-trace-agreement}
 \lim_{r\uparrow1}f(r\zeta)=\widetilde f(\zeta)\qquad\text{q.e. on }\T.
\end{equation}
\end{published}

Beurling established radial existence outside a set of
logarithmic capacity zero for analytic Dirichlet functions
\cite{Beurling40}; Brown and Shields recall this
result in \cite[p.~295]{BS84}. The harmonic representative theory
is presented by Richter and Sundberg
\cite[Section~5]{RS14} and by Richter and Yilmaz
\cite[p.~2122]{RY19}. For the nontangential maximal estimate in Dirichlet spaces with
harmonic weights, see \cite[Theorem~1.3]{Chacon11}. Richter and Yilmaz also credit Chac\'on
with the nontangential limit agreement for harmonic extensions.
The boundary statement above rests on two facts:
existence of analytic boundary limits and their identification
with the quasi-continuous harmonic representative.
To identify the radial limit with the representative in our norm,
apply the harmonic result to the real and imaginary parts of $f^*$,
whose Poisson extensions are the real and imaginary parts of $f$.
The radial specialization gives \eqref{eq:radial-trace-agreement}.

Nagel, Rudin, and Shapiro developed a broader theory relating boundary approach regions to the capacity
of exceptional sets \cite{NRS82}. Their work highlights two
choices in a boundary-limit assertion: how the point is approached
from inside the disk, and how the set of possible failures is
measured. This distinction helps explain the  radial
formulation used here. The zero set in this paper specifies
approach along a radius, while the representative theorem controls
the exceptional set in capacity. The wider tangential theory
provides context for those choices.

We always choose Borel representatives. In fact, for a
quasi-continuous representative take open exceptional sets $U_j$ with
$C_R(U_j)<2^{-j}$. The representative restricted to each closed
$\T\setminus U_j$ is continuous, and their union covers $\T$
except for the polar Borel set $\bigcap_j U_j$. Defining the
representative to be zero on that intersection preserves
quasi-continuity and makes it Borel. A sum or difference of
quasi-continuous representatives is quasi-continuous by
\eqref{eq:capacity-finite-subadditivity}. Uniqueness identifies it
q.e. with the representative of the corresponding sum or difference
in $V$. The modulus has the same property by
\eqref{eq:V-difference} and continuity of the modulus map.

Representative uniqueness also holds locally. If $u=0$ a.e. on
an open set $U$, multiply by a smooth function $\chi$ supported in
$U$. Formula \eqref{eq:V-difference} gives $\chi u\in V$ since
$\chi$ is bounded and Lipschitz. The quasi-continuous representative
$\chi\widetilde u$ of the zero function is zero q.e. A countable
collection of such cutoffs, each equal to one on a compact subset
of $U$, yields $\widetilde u=0$ q.e. on $U$.

To pass a vanishing condition to a norm limit, we need convergence
of traces away from a polar set. Brown and Shields credit Shapiro with the general closedness argument for the trace-zero space \cite[p.~295]{BS84}. It uses the capacitary estimates of Nagel, Rudin, and Shapiro \cite{NRS82}.
The obstacle argument below gives the
coefficient one for the norm \eqref{eq:V}.

\begin{lemma}[Level bound and trace subsequence]\label{lem:trace-subsequence}
If $u\in V$ and $t>0$, then
\begin{equation}\label{eq:capacity-level}
 C_R(\{|\widetilde u|>t\})\le t^{-2}\|u\|_V^2.
\end{equation}
If $u_j\to u$ in $V$, a subsequence of $\widetilde u_j$ converges
to $\widetilde u$ q.e.
\end{lemma}
\begin{proof}
Choose an open exceptional set $G$ of arbitrarily small capacity
such that $\widetilde u$ is continuous on $\T\setminus G$. The
strict superlevel on $\T\setminus G$ is relatively open, so it
has the form $O\setminus G$ for some open $O\subset\T$. Let $b$
be a nonnegative obstacle on $G$. Outside the exception, the scaled
trace already exceeds one; the added obstacle covers the exception.
The maximum of $|u|/t$ and $b$ is an obstacle on $O\cup G$. The same maximum estimate as in
Lemma~\ref{lem:capacity-unions} bounds the squared norm by
$t^{-2}\|u\|_V^2+\|b\|_V^2$. Taking nearly minimizing $b$ and
then letting $C_R(G)$ decrease to zero proves \eqref{eq:capacity-level}.

For convergence, choose a subsequence whose level-set exceptions have
capacities with summable square roots. Take indices
$j_k$ with $\|u_{j_k}-u\|_V\le2^{-2k}$.
By \eqref{eq:capacity-level} and representative uniqueness the
set where $|\widetilde u_{j_k}-\widetilde u|>2^{-k}$ has capacity
at most $2^{-2k}$. The square root of the capacity of the union of these
exceptional sets for $k\ge n$ is at most $\sum_{k\ge n}2^{-k}$.
The tail bound tends to zero, so the limsup is polar. Outside it,
the trace differences eventually fall below the chosen thresholds.
\end{proof}

For the  radial zero set in \eqref{eq:radial-zero},
\[
 \Zrad(f)=
 \bigcap_{k\ge1}\ \bigcup_{j\ge2}\
 \bigcap_{\substack{r\in\mathbb Q\\1-1/j<r<1}}
 \{\zeta\in\T:|f(r\zeta)|\le 1/k\}.
\]
The sets in the innermost intersection are closed. Using all
rational radii in the indicated interval suffices by continuity
in $r$. Hence $\Zrad(f)$ is Borel. By
\eqref{eq:radial-trace-agreement}, it agrees with
$\{\widetilde f=0\}$ off a polar set, and therefore
\begin{equation}\label{eq:radial-zero-capacity}
 C_R(\Zrad(f))=C_R(\{\widetilde f=0\}).
\end{equation}

\subsection{Integration of boundary traces against finite-energy measures}\label{subsec:03-pairing}

In Section~\ref{sec:10}, probability measures carried by Cantor
sets will be used to test boundary traces. The Fourier weights
of the trace norm determine which measures give bounded integrals.
The corresponding energy uses the reciprocal weights.
For a finite measure $\sigma$ on $\T$, write
$\widehat\sigma(n)=\int_{\T}\bar\zeta^{\,n}\,d\sigma(\zeta)$ and
\[
 I_R(\sigma)=\sum_{n\in\mathbb Z}\frac{|\widehat\sigma(n)|^2}{1+|n|}.
\]
For a positive measure, the energy has a positive kernel formula.
With $P_r(\zeta)=(1-r^2)/|\zeta-r|^2$ and
$K_R(\zeta)=\int_0^1P_r(\zeta)\,dr$, we have
\begin{equation}\label{eq:reciprocal-energy-kernel}
 I_R(\sigma)=\iint_{\T^2}K_R(\zeta\bar\eta)\,
                      d\sigma(\zeta)\,d\sigma(\eta).
\end{equation}
For each $r<1$ the integral of $P_r(\zeta\bar\eta)$ against the
product measure is $\sum_n r^{|n|}|\widehat\sigma(n)|^2$ by its
absolutely convergent Fourier series. Integration in $r$ gives
\eqref{eq:reciprocal-energy-kernel} by Tonelli, including infinite
energy.

Richter and Yilmaz use the quasi-continuous representative when
integrating against a finite-energy measure \cite[p.~2123]{RY19}.
The Fourier weights above give the constant in the following bound.

\begin{proposition}[Integration of the boundary representative]\label{prop:trace-pairing}
Let $\sigma$ be a finite positive Borel measure with
$I_R(\sigma)<\infty$. Every polar set has $\sigma$-outer measure
zero, and
\begin{equation}\label{eq:pairing}
 \int_{\T}|\widetilde u|\,d\sigma
     \le I_R(\sigma)^{1/2}\|u\|_V\qquad (u\in V).
\end{equation}
The continuous extension to $V$ of the Fourier pairing on
trigonometric polynomials equals $\int_\T\widetilde u\,d\sigma$.
\end{proposition}
\begin{proof}
Weighted Cauchy--Schwarz applied to
$\int p\,d\sigma=\sum_n\widehat p(n)\widehat\sigma(-n)$ gives
$|\int p\,d\sigma|\le I_R(\sigma)^{1/2}\|p\|_V$.
For continuous $u\in V$ the Fej\'er means converge uniformly
and in $V$, so the same inequality holds for the ordinary
integral of $u$. Applying it to $|u|$ and using the modulus
contraction gives $\int |u|\,d\sigma\le I_R(\sigma)^{1/2}\|u\|_V$
for continuous $u\in V$. This is the absolute integral hypothesis
in the positive finite-energy measure theorem of
\cite[p.~2123]{RY19}. The theorem implies that $\sigma$ gives
zero mass to Borel polar sets. Every polar set lies in a Borel
polar set by outer regularization.

To integrate a general trace $u\in V$, we must identify the Fourier
pairing with the integral of its representative. Take trigonometric
polynomials $p_j$ tending to $u$ in $V$. The absolute bound makes
them Cauchy in $L^1(\sigma)$. Lemma~\ref{lem:trace-subsequence} gives a
subsequence converging to $\widetilde u$ q.e., hence
$\sigma$-a.e. A further subsequence converges to the $L^1(\sigma)$
limit almost everywhere, so this limit is $\widetilde u$.
Passing to the limit in the bound and in the integral identifies
both the $L^1(\sigma)$ norm and the extended Fourier pairing.
\end{proof}

If $\nu=\sigma_1-\sigma_2$ with both $\sigma_i$ positive and satisfying
$I_R(\sigma_i)<\infty$, then $|\nu|\le\sigma_1+\sigma_2$.
Proposition~\ref{prop:trace-pairing} gives absolute integrability
against $|\nu|$ and identifies the signed integral with the
difference of the two positive pairings. Fourier duality then gives
$|\int_\T\widetilde u\,d\nu|\le I_R(\nu)^{1/2}\|u\|_V$. A nonzero point mass
$c\delta_\xi$ has $I_R(c\delta_\xi)=
 |c|^2(1+2\sum_{n\ge1}(n+1)^{-1})=\infty$.

\paragraph{Sets of zero capacity.}
The logarithmic and obstacle capacities have the same polar sets
\cite[Section~5]{RS14}. The direct comparison below lets us use
the norm \eqref{eq:V} for quantitative estimates and interpret zero
capacity in the usual logarithmic convention. For the obstacle
capacity, the first step is to minimize the energy of a probability
supported on the compact set.
For a nonempty compact $K\subset\T$, let $\mathcal P(K)$ be
the positive probability measures supported on $K$. The exact
dual formula is
\begin{equation}\label{eq:capacity-energy-duality}
 C_R(K)=\left[\inf_{\sigma\in\mathcal P(K)}I_R(\sigma)\right]^{-1},
 \qquad \frac1\infty=0.
\end{equation}
To see this, the infimum in the obstacle definition for a compact set may be taken over smooth functions at least one on $K$. In one direction, multiply a smooth obstacle by $1+\varepsilon$
to obtain an obstacle on an open neighborhood. In the other,
smooth an open-neighborhood obstacle with a nonnegative probability
kernel supported in a sufficiently short arc. Smoothing is contractive
in $V$ and produces a function at least one on $K$. Testing with these
smooth obstacles and applying \eqref{eq:pairing} gives
$1\le C_R(K)I_R(\sigma)$ for every finite-energy $\sigma\in
\mathcal P(K)$. For the reverse bound when $c=C_R(K)>0$, the
equilibrium measure of $K$ from \cite[p.~2123]{RY19} has mass
$c$ and its representing vector in $V$ has squared norm $c$.
The squared dual norm of integration is $I_R$ by
\eqref{eq:V}, so normalizing the equilibrium measure by $c$
gives a probability of energy $1/c$. If $c=0$, the test with
smooth obstacles excludes every finite-energy probability.

For a probability $\sigma$ on $\T$ the logarithmic energy has the
Fourier form used in \cite[p.~294]{BS84}
\[
 I_{\log}(\sigma)=\iint_{\T^2}\log\frac1{|\zeta-\eta|}
            \,d\sigma(\zeta)\,d\sigma(\eta)
     =\sum_{n\ge1}\frac{|\widehat\sigma(n)|^2}{n}.
\]
For the equality, first integrate the Fourier expansion of
$-\log|1-r\zeta\bar\eta|$ at $r<1$. The resulting series
$\sum_{n\ge1}r^n|\widehat\sigma(n)|^2/n$ increases as $r\uparrow1$.
The regularized kernels are bounded below by $-\log 2$.
For $r\ge1/2$, the inequality
$|1-r\zeta\bar\eta|^2\ge r|\zeta-\eta|^2$ also bounds them
above by $\log(1/|\zeta-\eta|)+(\log 2)/2$. Dominated convergence
gives the identity when the logarithmic energy is finite, and
Fatou's lemma gives it when the energy is infinite.
Since $|\widehat\sigma(n)|\le1$,
\[
 2 I_{\log}(\sigma)-1\le I_R(\sigma)\le 2 I_{\log}(\sigma)+1.
\]
Indeed, the difference between the finite partial sums is
$1-2\sum_{n=1}^N|\widehat\sigma(n)|^2/[n(n+1)]$.
Using \eqref{eq:capacity-energy-duality}, the classical
capacity $\operatorname{cap}_{\log}(K)
=\exp[-\inf_{\sigma\in\mathcal P(K)} I_{\log}(\sigma)]$ satisfies
\[
 e^{-1/2}\exp\!\left(-\frac{1}{2C_R(K)}\right)
 \le \operatorname{cap}_{\log}(K)
 \le e^{1/2}\exp\!\left(-\frac{1}{2C_R(K)}\right),
\]
with the exponential interpreted as zero at $C_R(K)=0$.
On open sets both capacities are suprema over compact subsets,
using the compact inner regularity of $C_R$ from
\cite[p.~2122]{RY19}. Both are defined on arbitrary sets
by outer regularization. The two bounds pass through
these suprema and infima because their functions of $C_R$ are
increasing and continuous on $[0,1]$. In particular the two
conventions give the same polar sets.
\subsection{Trace-zero spaces and smooth outer functions}\label{subsec:03-ideal}
Approximation~B varies the boundary set and follows the canonical
functions of the resulting spaces. To apply the projection
construction from Section~\ref{subsec:02-projection}, we first verify
closedness and shift invariance. For a boundary set \(K\subset\T\),
define the trace-zero space
\[
 \D_K=\{f\in\D:\widetilde f=0\text{ quasi-everywhere on }K\}.
\]
By \eqref{eq:radial-trace-agreement}, membership is equivalent to
vanishing in the  radial sense quasi-everywhere on \(K\).
Brown and Shields discuss Carleson's closedness result for thin sets and credit Shapiro with the general trace-space argument
\cite[pp.~294--295]{BS84}. See also the closed invariant spaces in
\cite[p.~2119]{RY19}. We use the trace subsequence just proved.
Linearity follows by taking the union of the polar exceptions for the
two summands. If \(f_j\in\D_K\) converges to \(f\) in \(\D\),
Lemma~\ref{lem:trace-subsequence} gives a subsequence converging to
\(\widetilde f\) outside a polar set. Also exclude the countable union
of the vanishing exceptions of that subsequence. The limit is zero
on the remaining part of \(K\), so \(\D_K\) is closed. It is invariant
under multiplication by \(z\), since \(zf\in\D\) and its radial
limit at a point where \(f\) has limit zero is zero.

The projection of \(1\) onto the trace-zero space is nonzero as soon
as the space contains a function with nonzero value at zero.
The Carleson condition supplies a smooth outer function for this purpose.

A nonempty compact \(K\subset\T\) is a \emph{Carleson set}
\cite{Carleson52} if
\[
       \int_{\T}\log\frac{2}{\dist(\zeta,K)}\,dm(\zeta)<\infty ,
\]
where distance is Euclidean. Such a compact is arclength null.
Equivalently, the set is arclength null and the complementary arcs
\(I\) have finite gap entropy \(\sum_I m(I)\log(1/m(I))\)
\cite[Remark~1, (6)--(7), p.~113]{Kor72}.
Indeed, on an arc of angular length \(L\), the distance to \(K\) is
comparable, with absolute constants, to the shorter angular distance
to the two endpoints. Integrating the logarithm therefore gives
\(L\log(1/L)+O(L)\), with a uniform error bound. Rescaling \(L\) by
\(2\pi\) changes the sum by a finite constant.

The smooth outer function vanishes on the prescribed Carleson set.
Outerness keeps its value at the origin nonzero, while smoothness lets us use
it as a multiplier. Taylor and Williams construct outer
functions smooth to every order on the closed disk, with exactly
the prescribed Carleson set as their zero set and with every
derivative vanishing there \cite[Definition~3.1 and Theorem~3.3]{TW70}.
Wells records that the existence of such a flat outer function
was proved independently by Novinger and by Taylor and Williams \cite[Theorem~A]{Wells70}.
Novinger's published paper is \cite{Novinger71}.
This gives the smooth outer member needed to start the projection
construction. The canonical function is obtained by a separate extremal
construction; Section~\ref{subsec:03-seed} proves its outerness by
combining generation with the presence of the smooth outer function
in the trace-zero space.

\begin{published}[Smooth outer functions with prescribed Carleson zero sets]\label{input:flat-outer}
For a Carleson set \(K\), there is a Hardy outer function
\(F_K\in A^\infty\) whose closed-disk zero set is \(K\) and whose
derivatives of every order vanish there. Here \(A^\infty\) consists of
analytic functions whose derivatives of all orders extend continuously
to the closed disk. We use the formulation of Korenblum's result
\cite[Lemma~23, English translation p.~134]{Kor72} recorded by
Richter and Yilmaz \cite[p.~2119]{RY19}.
\end{published}

The function \(F_K\) belongs to \(\D_K\) and \(F_K(0)\ne0\). By
Lemma~\ref{lem:bounded-products}, \(F_Kh\in\D\) for every \(h\in\D\).
At quasi-every point of \(K\), \(h\) has a finite radial limit and
\(F_K\) has continuous trace zero. Thus \(F_Kh\in\D_K\).
The same lemma gives invariance of every closed shift-invariant
subspace under \(A^\infty\) multipliers.

Section~\ref{sec:12} combines sets associated with different atoms.
For finitely many compacts \(K_i\), the union of the individual
trace exceptions is polar, so
\[
       \D_{\bigcup_i K_i}=\bigcap_i\D_{K_i}.
\]
If each \(K_i\) is Carleson, then so is their union because \(\log(2/\dist(\zeta,\bigcup_i K_i))\) is the maximum of the nonnegative functions \(\log(2/\dist(\zeta,K_i))\) and is bounded by their sum. The capacity of the union is bounded by the sum in \eqref{eq:capacity-finite-subadditivity}.
If any \(K_i\) is nonpolar, monotonicity gives nonpolarity of the union.

\subsection{Outer canonical functions for nonpolar Carleson sets}\label{subsec:03-seed}
The nested balls of Section~\ref{sec:13} start from a proper outer
canonical function. A Carleson set gives a smooth outer member
of the trace-zero space, and positive capacity keeps \(1\) out of the
space. Nonpolar thin sets of this type go back to Carleson
\cite{Carleson52}. We verify the two conditions for a concrete compact.

Embed the middle-thirds Cantor set in a short circular arc. Its gap
entropy is finite, since at level \(j\) there are \(2^{j-1}\) gaps of
length proportional to \(3^{-j}\), and
\(\sum_{j\ge1}(1+j)(2/3)^j<\infty\). The remaining complementary
arc contributes a finite term. Give the Cantor set its fair
Bernoulli probability. Two independent points have their first
different digit at level \(j\) with probability \(2^{-j}\), and at
that level their chordal separation is at least \(c\,3^{-j}\)
for a fixed \(c>0\) from the arc embedding. Their logarithmic energy
is finite by \(\sum_{j\ge1}2^{-j}(1+j)<\infty\), so this compact is
nonpolar. This supplies a choice of \(K\) independent of the two
approximations.

For any nonpolar Carleson set \(K\), let \(M=\D_K\). The smooth
outer function \(F_K\) in \(M\) has nonzero value at zero, so \(P_M1\ne0\).
The constant \(1\) is outside \(M\) because \(K\) is nonpolar.
Section~\ref{subsec:02-projection} thus gives a unit canonical function
\(\eta_K=G_M\) with \(0<\eta_K(0)<1\) and \([\eta_K]=M\).
To prove outerness we use generation by \(\eta_K\) together with
the presence of \(F_K\) in \(M\). Factor \(\eta_K\) in \(H^2\) as \(IO\) with
\(I\) inner and \(O\) outer. Polynomial multiples of \(\eta_K\) lie in
the closed Hardy subspace \(IH^2\). Since \(\D\)-convergence implies
\(H^2\)-convergence, \(M\subset IH^2\). In particular \(F_K=IH\) for
some \(H\in H^2\), and outerness of \(F_K\) forces \(I\) to be constant
by uniqueness of inner--outer factorization. For the Cantor-type compactum
just constructed, set \(e_0=\eta_K\).

Since \(\eta_K\in\D_K\), the set \(K\setminus\Zrad(\eta_K)\) is
polar. By Lemma~\ref{lem:individual-continuation}, \(\eta_K\) extends holomorphically across each point of \(\T\setminus K\) and has a nonzero value at that point. Hence \(\Zrad(\eta_K)\subset K\).

The same closed-subspace argument explains Brown and Shields'
necessary conditions for cyclicity
\cite[Corollary~1 after Proposition~6 and Theorem~5]{BS84}.
If \(f\ne0\) is cyclic
and \(I\) is its Hardy inner factor, then \(\D=[f]\subset IH^2\);
density of polynomials in \(H^2\) gives \(IH^2=H^2\), hence \(I\) is
constant. Also, \([f]\subset\D_{\Zrad(f)}\) by the definition of the
trace-zero space and its closedness and shift invariance. Cyclicity
then forces \(1\in\D_{\Zrad(f)}\), so \(\Zrad(f)\) is polar.

\subsection{Capacity of radial zero sets throughout a closed ball}\label{subsec:03-neighborhood}
The limit in Section~\ref{sec:13} belongs to every chosen closed
ball. To control its radial zero set, we need a capacity bound
throughout each ball. We obtain it by estimating where the center
has small trace and where a perturbation can change the trace by
more than a fixed threshold. The estimates apply to every \(f\in\D\).
Fix a quasi-continuous trace
\(\widetilde f\) and set
\[
   \Phi(f)=C_R(\Zrad(f)),\qquad E_t(f)=\{|\widetilde f|\le t\}
              \quad(t\ge0).
\]
By \eqref{eq:radial-zero-capacity}, \(C_R(E_0(f))=\Phi(f)\).

\begin{lemma}[Near-zero levels and perturbations]\label{lem:capacity-perturbation}
For \(f,g\in\D\) and \(t>0\),
\begin{equation}\label{eq:capacity-perturbation}
 C_R(E_t(f))\downarrow\Phi(f)\quad(t\downarrow0),\qquad
 \Phi(g)\le C_R(E_t(f))+\frac{\|g-f\|_{\D}^2}{t^2}.
\end{equation}
\end{lemma}

\begin{proof}
For the near-zero limit, monotonicity gives the lower bound. Fix \(s>0\).
Take an open \(U\supset E_0(f)\) with \(C_R(U)<\Phi(f)+s\) and an
open \(G\) with \(C_R(G)<s\) such that \(\widetilde f\) is continuous
on \(\T\setminus G\). On the compact complement of \(U\cup G\),
\(|\widetilde f|\) is continuous and strictly positive. Its minimum
there is positive if the complement is nonempty. For all sufficiently
small \(t\), it follows that \(E_t(f)\subset U\cup G\), an inclusion
which is automatic if the complement is empty. This uses continuity
outside the exceptional set, rather than continuity from above of capacity. Finite subadditivity
bounds \(C_R(E_t(f))\) by \(\Phi(f)+2s\). Let \(s\downarrow0\).

For the perturbation bound, discard the union of the polar trace
exceptions for \(f\) and \(g\). At a point of \(\Zrad(g)\) outside \(E_t(f)\)
we then have \(|\widetilde g-\widetilde f|=|\widetilde f|>t\).
The difference is a quasi-continuous representative of \(g^*-f^*\).
Thus a new zero lies in one of those two level sets, up to polar
exceptions. Apply \eqref{eq:capacity-level}, finite subadditivity and the isometry
from \(\D\) to the analytic traces in \(V\).
\end{proof}

In particular \(\Phi\) is upper semicontinuous in \(\D\): for
\(g_j\to f\), take the upper limit in the second bound at fixed \(t\)
and then let \(t\downarrow0\). The size of the zero set can increase
at the limit, as the constants \(1/j\to0\) show.

Every nonzero constant in this example has an empty radial zero
set, whereas the zero limit vanishes on the whole circle. Thus
upper semicontinuity permits the capacity to increase at the
limit. This example takes place in the ambient Dirichlet space;
it is not a sequence in the unit canonical class. It illustrates why
bounds at successive centers alone do not supply a general
bound at their limit. The construction will instead impose a
capacity bound on every point of each chosen closed ball,
including its boundary. The next corollary fixes a small-trace
threshold with a strict capacity margin first, then chooses the
radius so that every perturbation in the closed ball stays within
that margin.

\begin{corollary}[A closed ball with a strict capacity margin]\label{cor:capacity-ball}
Let \(f\in\D\) and \(\varepsilon>\Phi(f)\), and put
\(b=\varepsilon-\Phi(f)>0\).
Choose \(t>0\) with \(C_R(E_t(f))\le\Phi(f)+b/4\). With
\(\delta=t\sqrt b/2\), every \(g\in\D\) satisfying
\(\|g-f\|_{\D}\le\delta\) obeys
\[
                 \Phi(g)\le\Phi(f)+b/2<\varepsilon .
\]
\end{corollary}
\begin{proof}
In \eqref{eq:capacity-perturbation} the second term is at most
\(\delta^2/t^2=b/4\), including at the boundary of the ball.
\end{proof}

Section~\ref{subsec:13-balls} applies this estimate after each
center has been selected. The metric in the next section dominates
Dirichlet distance, so its closed balls inherit the capacity bound
at the same radius.

\newpage \section{A metric on outer canonical functions and two approximations}\label{sec:04}

Section~\ref{subsec:03-neighborhood} gives a capacity bound throughout
a closed neighborhood of a function with small radial-zero capacity.
To use this bound repeatedly, we need to find such a function inside
each preceding neighborhood. We also need the centers to converge to
an outer canonical function. We address the limit by adding a boundary
logarithm term to the Dirichlet distance. We then show how a norm
approximation at a fixed outer function supplies the next center in
this metric.

\subsection{The outer canonical class}\label{subsec:04-class}

Let $\X$ be the class of outer $f\in\D$ satisfying
\eqref{eq:canonical-moments}. By \eqref{eq:certificate}, each such $f$
is the unit canonical function of $[f]$ with positive value at zero. The Dirichlet norm
gives $f(0)\le1$, with equality only for $f=1$. Write
$\Xpr=\X\setminus\{1\}$ for the proper part. The construction in
Section~\ref{subsec:03-seed} gives an element of $\Xpr$.

Every $f\in\X$ satisfies $|f(z)|\le1$ by the canonical bound in
Section~\ref{subsec:02-multiplier}. Its Hardy boundary values therefore
give a nonnegative function $u_f=-\log|f^*|$ in $L^1(m)$. For
$u\in L^1(m)$ define
\[
 \mathcal C_Hu(z)=\int_{\T}\frac{\zeta+z}{\zeta-z}\,u(\zeta)\,dm(\zeta).
\]
The normalization $f(0)>0$ forces $\omega=1$ in
\eqref{eq:outer-intro}, so
\begin{equation}\label{eq:outer-reconstruction}
 f=\exp(-\mathcal C_Hu_f),\qquad \int_{\T}u_f\,dm=-\log f(0).
\end{equation}
Define the metric
\begin{equation}\label{eq:dlog}
 \dlog(f,g)=\|f-g\|_{\D}
       +\|\log|f^*|-\log|g^*|\|_{L^1(m)}\qquad(f,g\in\X).
\end{equation}
The first term controls the canonical moments through strong Dirichlet
convergence. The second controls the boundary logarithms used to reconstruct the outer
function in \eqref{eq:outer-reconstruction}. We include $1$ when we
prove completeness and later choose a closed ball inside $\Xpr$ to
keep the distance from $1$ to the generated spaces positive.

\subsection{A complete metric retaining outerness}\label{subsec:04-complete}

The two terms of the metric give two limits. The proof identifies the
Dirichlet limit with the outer function reconstructed from the limit
of the boundary logarithms.

\begin{theorem}[Completeness]\label{thm:outer-completeness}
The metric space $(\X,\dlog)$ is complete.
\end{theorem}
\begin{proof}
Let $(f_j)$ be Cauchy for $\dlog$. It converges in $\D$ to a function
$f$, and $u_{f_j}$ converges in $L^1(m)$ to some $u\ge0$. The second
limit suggests an outer candidate; we will identify it with the
Dirichlet limit. For any
nonnegative $v\in L^1(m)$, put $F_v=\exp(-\mathcal C_Hv)$. The outer
representation gives a bounded outer function with
$-\log|F_v^*|=v$ almost everywhere. In fact
$|F_v|=\exp(-\operatorname{Re}\mathcal C_Hv)\le1$, and the radial limit
of the Poisson integral $\operatorname{Re}\mathcal C_Hv$ equals $v$
almost everywhere.

For $0\le r<1$, the integral kernel defining $\mathcal C_H$ has
modulus at most $(1+r)/(1-r)$ on $|z|\le r$. The map $A\mapsto e^{-A}$
is $1$-Lipschitz on the closed right half-plane, by integrating its
derivative along a line segment. Thus for $v,w\ge0$,
\[
 \sup_{|z|\le r}|F_v(z)-F_w(z)|
 \le \frac{1+r}{1-r}\|v-w\|_{L^1(m)}.
\]
By \eqref{eq:outer-reconstruction}, $f_j=F_{u_{f_j}}$. It follows that
$f_j\to F_u$ uniformly on each compact subdisk. Convergence in $\D$
also gives compact-subdisk convergence by the evaluation bound of
Section~\ref{subsec:02-norm}. Therefore
\[
 f=F_u,\qquad f(0)=\exp\!\left(-\int_{\T}u\,dm\right)>0,
 \qquad -\log|f^*|=u\quad\text{almost everywhere}.
\]

With the two limits identified, we check the canonical conditions.
Strong Dirichlet convergence gives $\|f\|_{\D}=1$. For every fixed
integer $k\ge1$, boundedness of multiplication by $z^k$ passes
$\langle z^kf_j,f_j\rangle_{\D}=0$ to the limit. Hence $f$ satisfies
\eqref{eq:canonical-moments} and belongs to $\X$. The two convergences
now give $\dlog(f_j,f)\to0$.
\end{proof}

\subsection{Continuity of boundary log moduli for outer canonical functions}\label{subsec:04-continuity}

Completeness gives the limit after the centers have been chosen.
To choose a new center, we will obtain an outer canonical function
with a bound on its Dirichlet error relative to the current outer
target. At a fixed target, that bound also controls the change in
the logarithms. The estimate isolates the region where the target
boundary modulus is small, since there the
logarithm is most sensitive to changes in modulus. The proof estimates
the positive-part integral in Scheff\'e's convergence argument
\cite{Scheffe47}.

\begin{proposition}[A logarithmic bound at a fixed outer function]
\label{prop:outer-continuity}
For $f,g\in\X$, put $a=f(0)$ and $\delta=\|g-f\|_{\D}$. For $0<t<1$
let
\[
 T_f(t)=\int_{\{|f^*|<t\}}-\log|f^*|\,dm.
\]
Then $T_f(t)\to0$ as $t\downarrow0$. If $\delta\le a/2$, then
\begin{equation}\label{eq:outer-continuity}
 \dlog(f,g)\le (1+2/a+2/t)\delta+2T_f(t).
\end{equation}
\end{proposition}
\begin{proof}
Put $u=-\log|f^*|$ and $v=-\log|g^*|$. Since $u\in L^1(m)$ and
$|f^*|>0$ almost everywhere, dominated convergence gives the claim
about $T_f(t)$. First we bound the part of the integral on which $u$
exceeds $v$. On $\{|f^*|<t\}$ we have $(u-v)_+\le u$, where
$s_+=\max(s,0)$. On the rest of the circle, the positive part can
be nonzero only where $|g^*|>|f^*|\ge t$. The mean-value theorem for
$\log$ then gives
\[
 (u-v)_+\le\frac{\big||g^*|-|f^*|\big|}{t}.
\]
Since $m(\T)=1$ and $\|g-f\|_{H^2}\le\delta$, integration yields
$\int(u-v)_+\,dm\le T_f(t)+\delta/t$.

The signed integral of $v-u$ is determined by the two values at zero.
Combining it with the positive-part bound gives the $L^1$
distance. By \eqref{eq:outer-reconstruction},
\[
 \|u-v\|_{L^1(m)}
 =\int(v-u)\,dm+2\int(u-v)_+\,dm
 =\log\frac{a}{g(0)}+2\int(u-v)_+\,dm.
\]
Evaluation at zero has norm one, so $|g(0)-a|\le\delta$. With
$\delta\le a/2$, the mean-value theorem gives
$|\log a-\log g(0)|\le 2\delta/a$. Adding the Dirichlet error proves
\eqref{eq:outer-continuity}.
\end{proof}

Given a desired radius $R>0$ about $f$, choose $t$ with $T_f(t)<R/8$.
With this tail contribution controlled, take
\[
 \eta=\min\!\left\{a/2,\ \frac{R}{2(1+2/a+2/t)}\right\}.
\]
Any $g\in\X$ with $\|g-f\|_{\D}<\eta$ satisfies
$\dlog(f,g)<3R/4$.
This choice of \(\eta\) works for all such $g$ once $f$ and $R$ are fixed. Its construction uses the value \(a=f(0)\) and the tail $T_f$ of that particular $f$.
The reverse continuity follows from
$\|g-f\|_{\D}\le\dlog(f,g)$, so the two metrics induce the same relative topology on $\X$. For the construction, we can therefore fix the allowable Dirichlet error before choosing either approximant.

Evaluation and \eqref{eq:outer-reconstruction} also give, for any
$f,g\in\X$,
\begin{equation}\label{eq:metric-origin}
\begin{aligned}
 |f(0)-g(0)|&\le\|f-g\|_{\D}\le\dlog(f,g),\\
 |\log f(0)-\log g(0)|&\le\|\log|f^*|-\log|g^*|\|_{L^1(m)}
                         \le\dlog(f,g).
\end{aligned}
\end{equation}
In particular, if $\dlog(f,f_0)\le r$ and $a_0=f_0(0)$, then
\[
 a_0-r\le f(0)\le a_0+r,
 \qquad a_0e^{-r}\le f(0)\le a_0e^r.
\]
For a proper center $f_0$, choosing $r$ small enough that $a_0+r<1$
gives a positive lower bound $1-(a_0+r)^2$ on the squared
approximation error in \eqref{eq:certificate} throughout the ball.

\subsection{Approximations A and B}\label{subsec:04-approximations}

To find a nearby outer canonical function with small radial-zero
capacity, we first approximate by the canonical function of a space
defined by a singular inner factor. Section~\ref{subsec:02-atomic}
defines this space $M_\mu$ and its canonical function $G_\mu$ for a
nonzero finite positive measure $\mu$ of finite support.
Approximation A chooses the measure using the starting function and the prescribed errors. Approximation B takes that fixed measure
and recovers its canonical function in norm using outer functions
whose radial-zero capacity tends to zero.

\begin{theorem}[Approximation A]\label{thm:A}
For every $e\in\Xpr$ and every $\varepsilon_{\rm norm},
\varepsilon_{\rm mass}>0$, there is a nonzero finite positive
measure $\mu$ of finite support such that
\[
 \|G_\mu-e\|_{\D}<\varepsilon_{\rm norm},\qquad
 0<\mu(\T)<\varepsilon_{\rm mass}.
\]
For this same measure, $[G_\mu]=M_\mu$ and $G_\mu=S_\mu O_\mu$ with
$O_\mu$ outer and $O_\mu(0)>0$.
\end{theorem}
Sections~\ref{sec:05}--\ref{sec:06} choose interior zeros while
controlling both the norm error and the sum that will bound the
atomic mass. Section~\ref{sec:07} transfers this approximation to the
canonical function of the resulting atomic space.

There is also a qualitative reason for small mass. Suppose \(g_j\to e\)
in \(H^2\), with \(e\) outer, and \(I_j\) is an inner divisor of \(g_j\).
For every polynomial \(p\), the distance from \(1\) to \(I_jH^2\) satisfies
\[
 \sqrt{1-|I_j(0)|^2}\le \|1-pe\|_{H^2}
                  +\|p\|_\infty\|e-g_j\|_{H^2}.
\]
Indeed \(pg_j\in I_jH^2\), and the projection of \(1\) onto
\(I_jH^2\) has norm \(|I_j(0)|\). First choose \(p\) by outerness,
then let \(j\) tend to infinity. This gives \(|I_j(0)|\to1\),
hence \(\mu_j(\T)\to0\) for \(I_j=S_{\mu_j}\) and \(b(A_j)\to0\) for \(I_j=B_{A_j}\). Sundberg's example, discussed in Section~\ref{subsec:overview-approximation-a},
shows why the condition \(b(A_j)\to0\) alone cannot give convergence
to \(e\).
Estimate~\eqref{eq:expected-blaschke-budget} controls the expected value of \(b(A)\) for a fixed total mass \(V_\nu\), while the subset is still being chosen.

\begin{theorem}[Approximation B]\label{thm:B}
Fix a nonzero finite positive measure $\mu$ of finite support as in
\eqref{eq:atomic-target}. Write $E_\mu=\supp\mu$, let $N$ be
the number of points in $E_\mu$, and put $m_\mu=\mu(\T)$.
There is a sequence of nonpolar Carleson sets
$K_{\mu,\ell}\subset\T$ whose canonical functions
$H_{\mu,\ell}=G_{\D_{K_{\mu,\ell}}}$ belong to $\Xpr$.
For every $\ell\ge1$, $H_{\mu,\ell}$ extends continuously to
$\overline{\DD}$ and satisfies
\[
 [H_{\mu,\ell}]=\D_{K_{\mu,\ell}},\qquad
 \{z\in\overline{\DD}:H_{\mu,\ell}(z)=0\}
   =\Zrad(H_{\mu,\ell})=K_{\mu,\ell}.
\]
The same sets satisfy
\[
 d_H(K_{\mu,\ell},E_\mu)\le3\,2^{-\ell},\qquad
 0<C_R(K_{\mu,\ell})\le5N/\ell,
\]
where $d_H$ is Hausdorff distance for the Euclidean chord metric on $\T$.
Along this sequence, as $\ell\to\infty$,
\[
 H_{\mu,\ell}\longrightarrow G_\mu\quad\text{in }\D,\qquad
 P_{\D_{K_{\mu,\ell}}}y\longrightarrow P_{M_\mu}y
 \quad\text{in }\D\quad(y\in\D).
\]
The target $G_\mu$ has exact Hardy inner factor $S_\mu$.
Writing $m$ for normalized arclength and $f^*$ for the Hardy boundary
values of $f$, we also have
\begin{align*}
 (-\log|H_{\mu,\ell}^*|)\,dm
 &\overset{*}{\rightharpoonup}(-\log|G_\mu^*|)\,dm+\mu,\\
 \|\log|H_{\mu,\ell}^*|-\log|G_\mu^*|\|_{L^1(m)}
 &\longrightarrow m_\mu>0 .
\end{align*}
The weak-* convergence is tested against $C(\T)$.
\end{theorem}
Theorem~\ref{thm:B} fixes the masses, directions, and number of atoms.
For any prescribed positive tolerances, the Dirichlet-norm error
and the capacity are below their respective tolerances for all
sufficiently large indices.
Qualitative decay of radial-zero capacity already follows from
strong convergence to this target. Put
\(Q(z)=\prod_{\xi\in E_\mu}(1-\bar\xi z)\).
The function \(S_\mu Q\) in Lemma~\ref{lem:atomic-witness}
belongs to \([G_\mu]\). Every polynomial multiple of \(G_\mu\)
vanishes quasi-everywhere on \(\Zrad(G_\mu)\), and this trace
condition is closed in \(\D\). It therefore holds for \(S_\mu Q\),
which extends with nonzero value off the finite set \(E_\mu\).
Consequently \(\Phi(G_\mu)=0\), and
Lemma~\ref{lem:capacity-perturbation} gives \(\Phi(H_{\mu,\ell})\to0\).
The theorem identifies the radial zero set of each approximant
with its defining compact set \(K_{\mu,\ell}\) and gives the
capacity bound \(5N/\ell\). For fixed \(\mu\), the index for a
prescribed norm error is obtained from strong convergence in \(\D\).

The positive limiting \(L^1\) distance between the boundary
logarithms explains why we use the Dirichlet distance for the two
steps through \(G_\mu\). For an outer approximant, the integral of
its negative boundary logarithm is determined by its value at zero.
At the limit, the singular inner factor has boundary modulus one
almost everywhere and multiplies the outer factor's value at zero
by \(e^{-m_\mu}\). The limiting integral of the approximants thus
exceeds the integral for \(G_\mu\) by the mass \(m_\mu\). The theorem
shows that the $L^1$ distance tends to this same positive gap.
By \eqref{eq:finite-atomic-log-measure-limit}, the singular part
of the weak-* limit of \((-\log|H_{\mu,\ell}^*|)\,dm\) is \(\mu\).

The proof of B starts with one atom. Theorem~\ref{thm:G} supplies
shrinking compacts whose canonical values at zero equal
$G_{\tau,\xi}(0)$. The analytic criterion in Theorem~\ref{thm:C} then
gives strong convergence along that sequence, as stated in
Corollary~\ref{cor:R}.
Section~\ref{sec:12} combines these families for the atoms of the
fixed measure. One subspace condition places weak limits of vectors
from the approximating subspaces in the target space. The other approximates
each vector of the target in norm by vectors in those subspaces.
Together they give convergence of the projections and the strong
limit of the canonical functions in Section~\ref{subsec:12-assembly}.

\subsection{How the approximations will be joined}\label{subsec:04-overview}

To choose a new center near $e\in\Xpr$, fix a metric radius $R>0$
and a capacity bound $\varepsilon>0$.
Section~\ref{subsec:04-continuity} supplies a bound $\eta$ on the
norm error for the radius $R$. Allow half of this bound for each step.
Apply Approximation A with a bound of $\eta/2$ on the norm error and any
positive bound on the total mass to choose $\mu$ with
$\|G_\mu-e\|_{\D}<\eta/2$. Hold that measure fixed and take an index
$\ell$ in Approximation B large enough that
\[
 \|H_{\mu,\ell}-G_\mu\|_{\D}<\eta/2,
 \qquad C_R(\Zrad(H_{\mu,\ell}))<\varepsilon.
\]
The norm triangle inequality gives $\|H_{\mu,\ell}-e\|_{\D}<\eta$.
Apply the fixed-target estimate directly to the two outer functions
$H_{\mu,\ell}$ and $e$ to get $\dlog(H_{\mu,\ell},e)<R$. This gives
the new center. Section~\ref{subsec:03-neighborhood} then supplies a
capacity bound throughout a sufficiently small closed ball about that
center, because $\dlog$ bounds the Dirichlet distance. Repeating these
two steps gives the nested closed balls of Section~\ref{sec:13}.

\newpage \section{A weighted boundary criterion for the generated space}\label{sec:05}

Approximation~A begins with the projection of \(1\) onto the space
generated by the outer canonical target \(e\). We characterize \([e]\)
as the set of functions in \(\D\) for which a weighted boundary integral
is finite. In Section~\ref{sec:06}, we truncate the density of this
integral and use the resulting bounded quadratic forms to approximate
the projection.

\subsection{The weighted representation and the target domain}\label{subsec:05-model}

Fix \(e\in\Xpr\) and set \(M=[e]\) and \(a=e(0)\). The identity \eqref{eq:certificate} gives \(P_M1=ae\) and
\(0<a<1\). Thus \(M\) is proper and \(e\) is nonconstant. Let
\(w=|e^*|^2\), and let \(\ell_e\) be the function defined in \eqref{eq:source}. Since \(e\) is outer and \(|e|\le 1\), we
have \(0<w\le 1\) almost everywhere. Also
\(\int_\T\ell_e\,dm=D(e)\), so \(\ell_e\) is finite almost everywhere.
Define the boundary density and the associated extended-valued functional by
\begin{equation}\label{eq:cost-density}
 \rho_e=\frac{\ell_e}{w}\quad\text{a.e.},\qquad
 Q_e(h)=\int_\T \rho_e |h^*|^2\,dm\quad(h\in\D).
\end{equation}
The density \(\rho_e\) weights \(|h^*|^2\) in \(Q_e\). For a
product \(h=eF\), cancellation of \(w\) rewrites \(Q_e(h)\) as the
integral of \(|F^*|^2\) against the local energy \(\ell_e\).
The weight \(w\) itself enters the derivative energy \(D_w(F)\)
defined below.

\begin{theorem}[Characterization by a weighted boundary integral]\label{thm:cost-domain}
Let \(e\in\Xpr\) and define \(Q_e\) by \eqref{eq:cost-density}. Then
\[
 [e]=\{h\in\D:Q_e(h)<\infty\},\qquad
 Q_e(h)\le 4\|h\|_{\D}^2\quad(h\in[e]).
\]
\end{theorem}

The onto weighted representation stated below will be used.
Identifying the finite domain of \(Q_e\) with \([e]\) additionally
requires the Hardy quotient and the positive lower bound for \(\ell_e\) proved here.
For the forward bound, we estimate \(Q_e\) on polynomial multiples
of \(e\) and pass to their closure. For the converse, we start with
\(Q_e(h)<\infty\) and prove that \(h/e\) is a Hardy function with finite
weighted energy. To conclude that \(h\in[e]\), we use the
representation below. In it \(w\) is the boundary weight for the
derivative energy of the quotient.

With \(P[w]\) the Poisson extension relative to \(m\), put
\[
 D_w(F)=\int_{\DD} P[w](z)|F'(z)|^2\,dA_0(z),\qquad
 \mathcal H_w=\{F\in H^2:D_w(F)<\infty\}.
\]

\begin{published}[Onto canonical weighted representation]\label{input:weighted}
For \(e\in\Xpr\) and \(w=|e^*|^2\), multiplication by \(e\) maps
\(\mathcal H_w\) onto \([e]\), with
\[
 \|eF\|_{\D}^2=\|F\|_{H^2}^2+D_w(F)\qquad(F\in\mathcal H_w).
\]
\end{published}

Apply Richter and Sundberg's weighted representation
\cite[(3.1) and Theorem~3.2]{RS92} with base measure \(m\).
Since \(e\) is a unit vector in \(M\ominus zM\), the representing
measure defined immediately before (3.1) on p.~175 of \cite{RS92}
is \(w\,dm\).
Their norm contains the \(H^2\) term, and
\eqref{eq:weighted-local-douglas} identifies its energy term with
\(D_w\) in the normalized area \(dA_0\). See also
\cite[Theorem~0.2 and proof on p.~188]{Wu98} for the weighted range.
The polynomial identity \eqref{eq:canonical-gram} will be enough
for the forward estimate. The onto representation is what allows
the converse to apply to a general \(F\in\mathcal H_w\). 

For related weighted spaces, Ghara and Patra
\cite[Theorem~1.2]{GP26} characterize \(\D(m_f)\) when
\(f\in\mathcal A_\epsilon\) for some \(\epsilon>0\), where
\[
\mathcal A_\epsilon
=
\left\{
f\in H^2:
\sup_{z\in\mathbb D}
|f'(z)|^2(1-|z|)^{1-\epsilon}<\infty
\right\}.
\]
Writing \(dm_f=|f^*|^2\,dm\), their theorem gives
\[
\D(m_f)=\{g\in H^2:fg\in\D\},
\]
with the \(\D(m_f)\)-norm equivalent to
\(\bigl(\|g\|_{H^2}^2+D(fg)\bigr)^{1/2}\).
The derivative regularity assumption is
not imposed on the canonical function \(e\) in
Theorem~\ref{thm:cost-domain}.

\paragraph*{Complete Pick property of the weighted space.}
The norm in Theorem~\ref{input:weighted} gives \(\mathcal H_w\)
the positive harmonic derivative weight \(P[w]\), so
Shimorin's theorem \cite[Theorem~1.1]{Shimorin02} applies with \(d\mu=w\,dm\).
The criterion of Luo and Rydhe \cite[Theorem~1.2]{LuoRydhe26}
gives another proof. Indeed, \(M=\mathbb Ce\oplus zM\) implies
that \(eF-F(0)e\in zM\) for \(F\in\mathcal H_w\), so
\((F-F(0))/z\in\mathcal H_w\).
The forward shift on \(\mathcal H_w\) is intertwined with
\(M_z|_M\) by multiplication by \(e\) and hence is a two-isometry.
Both shifts are bounded by the closed graph theorem.
Also \(\langle F,1\rangle_{\mathcal H_w}=F(0)\), so its kernel
is normalized at \(0\). If \(K_M\) and \(K_{\mathcal H_w}\)
denote the respective reproducing kernels, then
\[
 K_M(z,\lambda)
 =e(z)\overline{e(\lambda)}K_{\mathcal H_w}(z,\lambda).
\]
Since \(e\) has no zeros in \(\DD\), this identity also
gives the complete Pick property for \(K_M\).

For \(h\in M\) with \(\|h\|_{\D}\le1\), apply
\cite[Theorem~1.1]{AHMR17} to \(h/e\in\mathcal H_w\).
There are contractive multipliers \(\varphi,\psi\) of
\(\mathcal H_w\), with \(\psi(0)=0\), such that
\[
 h=e\,\frac{\varphi}{1-\psi},
 \qquad
 \overline{(1-\psi)\mathcal H_w}^{\,\mathcal H_w}
 =\mathcal H_w.
\]
Polynomial density in \(\mathcal H_w\) and the canonical
isometry then give
\[
 [e(1-\psi)]=M.
\]

\subsection{A positive lower bound for the local Dirichlet integral}\label{subsec:05-quotient}

The condition \(Q_e(h)<\infty\) means that \(h^*/e^*\) is square
integrable with weight \(\ell_e\). To obtain ordinary square integrability,
we use a positive lower bound on \(\ell_e\) for this fixed target. 
Indeed, \(P_z(\xi)\ge (1-|z|)/(1+|z|)\) for all \(\xi\in\T\).
Restrict the integral in \eqref{eq:source} to \(|z|\le 1/2\) to get
\begin{equation}\label{eq:cost-floor}
 \ell_e(\xi)\ge c_e,\qquad
 c_e=\frac13\int_{|z|\le 1/2}|e'(z)|^2\,dA_0(z)>0 .
\end{equation}
If the integral were zero, \(e'\) would vanish on the smaller disk
and hence on \(\DD\), contrary to nonconstancy. Once \(e\) is fixed,
the same positive \(c_e\) works for every \(\xi\in\T\).
Thus \(\rho_e\ge c_e\) almost everywhere. The
measure \(\rho_e\,dm\) is absolutely continuous and \(\sigma\)-finite.

Suppose \(h\in\D\) has \(Q_e(h)<\infty\). Then
\[
 \int_\T \left|\frac{h^*}{e^*}\right|^2\,dm
 \le \frac1{c_e}\int_\T\ell_e
                  \left|\frac{h^*}{e^*}\right|^2\,dm
 =\frac{Q_e(h)}{c_e}.
\]
The bounded-outer division lemma
(Lemma~\ref{lem:bounded-outer-division}) now gives
\(F=h/e\in H^2\). Its boundary function is \(h^*/e^*\)
almost everywhere, and
\begin{equation}\label{eq:cost-hardy-quotient}
 \|F\|_{H^2}^2\le Q_e(h)/c_e,\qquad
 \int_\T \ell_e|F^*|^2\,dm=Q_e(h).
\end{equation}

\subsection{The forward and reverse boundary estimates}\label{subsec:05-bounds}

The two estimates come from taking a boundary difference of the
product \(h=eF\). Split it into the change in \(e\) with \(F\) held
fixed and the change in \(F\) with \(e\) held fixed. For \(F\in H^2\),
put \(h=eF\). Boundedness of \(e\) gives \(h\in H^2\), so the
Hardy boundary functions satisfy
\begin{equation}\label{eq:cost-product-difference}
 h^*(\xi)-h^*(\eta)
 =F^*(\xi)(e^*(\xi)-e^*(\eta))
   +e^*(\eta)(F^*(\xi)-F^*(\eta))
 \quad\text{for a.e. }(\xi,\eta)\in\T^2.
\end{equation}
By \eqref{eq:local-douglas} and \eqref{eq:weighted-local-douglas},
the integrals of the squared moduli of the two terms on the right
against \(dm(\xi)dm(\eta)/|\xi-\eta|^2\) are, respectively,
\[
 \int_\T \ell_e(\xi)|F^*(\xi)|^2\,dm(\xi)
 \qquad\text{and}\qquad D_w(F).
\]
In the first term we integrate the difference of \(e\) to obtain
\(D_\xi(e)=\ell_e(\xi)\), which gives \(Q_e(h)\) when \(h\in\D\).
In the second term we integrate the difference of \(F\) to obtain
\(D_\eta(F)\) with outer weight \(w(\eta)\). For \(h\in\D\), the
integral for the left-hand difference is \(D(h)\).
These are identities of nonnegative integrals, including the case
of infinite energy.

For the forward bound, take \(F=p\) polynomial and isolate the first
term on the right of \eqref{eq:cost-product-difference}. Using
\(|u-v|^2\le 2|u|^2+2|v|^2\) and \eqref{eq:canonical-gram} gives
\begin{equation}\label{eq:cost-forward}
 Q_e(ep)\le 2D(ep)+2D_w(p)\le 4\|ep\|_{\D}^2 .
\end{equation}
The constant \(4\) holds for every \(e\in\Xpr\). If \(ep_j\to h\)
in \(\D\), the same estimate for \(p_j-p_k\)
makes \((ep_j)^*\) Cauchy in \(L^2(\rho_e\,dm)\). Since
\(\rho_e\ge c_e\), convergence in this weighted space implies
convergence in \(L^2(m)\), where the limit is \(h^*\) by Hardy
norm convergence. Thus the weighted limit is the ordinary boundary
function of \(h\), and the trace operator satisfies
\begin{equation}\label{eq:cost-trace}
 T_e:M\longrightarrow L^2(\rho_e\,dm),\qquad
 T_e h=h^*,\qquad \|T_e\|\le 2 .
\end{equation}

For the reverse bound, suppose \(Q_e(h)<\infty\) and take the
\(H^2\) quotient \(F=h/e\) from \eqref{eq:cost-hardy-quotient}.
To apply the weighted representation, we still need the energy
of \(F\). Isolate the second term on the right of
\eqref{eq:cost-product-difference} and use the same square inequality.
The weighted integral then gives
\begin{equation}\label{eq:cost-reverse}
 D_w(F)\le 2D(h)+2\int_\T \ell_e|F^*|^2\,dm
          =2D(h)+2Q_e(h).
\end{equation}
In particular \(F\in\mathcal H_w\).

\subsection{The domain where \texorpdfstring{\(Q_e\)}{Q\_e} is finite}\label{subsec:05-domain}

The two directions now meet at the domain of the weighted
representation. In the forward direction, the weighted trace
estimate gives \(Q_e(h)<\infty\) for every \(h\) in the closure
of the polynomial multiples of \(e\). In the reverse direction, the lower bound
on \(\ell_e\) and Lemma~\ref{lem:bounded-outer-division} give
\(F=h/e\in H^2\). The reverse boundary estimate then gives
\(D_w(F)<\infty\).
These two properties put \(F\) in \(\mathcal H_w\).
Richter and Sundberg's isometry onto \(M\) then gives \(h=eF\in M\)
\cite[Theorem~3.2, p.~177]{RS92}. The polynomial norm identity supplied
the initial forward estimate.

\begin{proof}[Proof of Theorem~\ref{thm:cost-domain}]
Equation \eqref{eq:cost-trace} gives the bound on \(M\), so
\(Q_e(h)<\infty\) for every \(h\in M\). Conversely, if \(h\in\D\)
and \(Q_e(h)<\infty\), then \(F=h/e\in H^2\) by
\eqref{eq:cost-hardy-quotient}, and
\eqref{eq:cost-reverse} gives \(D_w(F)<\infty\).
Theorem~\ref{input:weighted} now puts \(h=eF\) in \(M\).
\end{proof}

The extended functional \(Q_e\) is lower semicontinuous in the
\(\D\)-norm. It is the supremum of the continuous functionals
\(h\mapsto\int_\T\min(\rho_e,n)|h^*|^2\,dm\), by monotone
convergence and the continuous inclusion into \(H^2\). Its finite
domain \(M\) is a proper closed subspace. In particular \(Q_e(1)\)
is infinite, whereas \(Q_e(e)=D(e)\) by \eqref{eq:cost-density}.
In Section~\ref{sec:06} we minimize
\(J_t(f)=\|1-f\|_{\D}^2+tQ_e(f)\) for \(t>0\).
The functional \(J_t\) is finite exactly on \(M\), and its
minimizer converges in \(\D\) to \(P_M1\) as \(t\downarrow0\).
Once the coefficient is fixed, bounded density truncations let
us approximate that minimizer by problems on all of \(\D\).

\newpage \section{Approximation by finite zero sets}\label{sec:06}

To approximate the fixed outer canonical target in Approximation~A,
we first use the canonical function of a space defined by finitely
many interior zeros. The same zero set must give a small norm error
and a small value of \(b(A)\). The sum \(b(A)\) bounds the mass of the atoms
obtained from these zeros in Section~\ref{sec:07}.

Fix \(e\in\Xpr\), put \(a=e(0)\in(0,1)\), and write \(M=[e]\).
We work with the projection of \(1\) onto \(M^\perp\). We will
approximate it by projections of \(1\) onto spans of finitely many
evaluation kernels. The identity
\eqref{eq:certificate} gives
\begin{equation}\label{eq:finite-zero-target}
 P_M1=ae,\qquad q_*=P_{M^\perp}1=1-ae,\qquad
 q_*(0)=\|q_*\|_{\D}^2=1-a^2.
\end{equation}
Also \(\|1-2q_*\|_{\D}=1\), since \(1-2q_*=ae-q_*\) is the
difference of the two orthogonal components of the unit vector \(1\).

For a finite set \(A\subset\DD\setminus\{0\}\), let
\begin{equation}\label{eq:finite-zero-projection}
 M_A=\{h\in\D:h|_A=0\},\qquad
 q_A=P_{M_A^\perp}1,\qquad e_A=G_{M_A},\qquad
 b(A)=\sum_{z\in A}-\log |z|.
\end{equation}
The sum \(b(A)\) will be used to control the mass of the boundary atoms. The space \(M_A\) is closed by
boundedness of evaluation, and it is shift-invariant. The polynomial
\(\prod_{z\in A}(\,\cdot-z)\) is in \(M_A\) and is nonzero at the
origin, so \(e_A\) is defined by
\eqref{eq:canonical-definition}. Here and below an empty set gives
\(M_\varnothing=\D\), \(q_\varnothing=0\), and \(b(\varnothing)=0\).
Bounded evaluation gives
\[
 M_A^\perp=\operatorname{span}\{k_z:z\in A\}.
\]
These kernels are linearly independent: pairing a linear relation
with each Lagrange interpolation polynomial at the points of \(A\) shows that the corresponding coefficient vanishes.

\begin{theorem}[Finite-zero approximation]\label{thm:finite-zero-approximation}
Fix \(e\in\Xpr\), put \(a=e(0)\), and let \(q_*\) be as in
\eqref{eq:finite-zero-target}. For every \(s>0\), there is a nonempty
finite set \(A\subset\DD\setminus\{0\}\) of distinct points of one
modulus \(r\), with \(\max\{1/2,1-s\}<r<1\), such that
\[
 \|q_A-q_*\|_{\D}^2+b(A)<s,\qquad
 |e_A(0)-a|<\sqrt{s},\qquad
 \|e_A-e\|_{\D}<\frac{2\sqrt{s}}{a}.
\]
In particular, given any two positive bounds, the same set \(A\) can satisfy the first bound on \(\|e_A-e\|_{\D}\) and the second on
\(b(A)\).
\end{theorem}

We multiply the weighted boundary integral from Section~\ref{sec:05}
by a small positive factor and truncate its density. We add the
resulting bounded quadratic form to \(\|1-f\|_{\D}^2\). With the
factor chosen sufficiently small and the cutoff then chosen sufficiently
large, the minimizer of this functional is close to \(P_M1=ae\).
We then approximate the truncated integral by a weighted sum of
\(|f(z)|^2\) over points on a single circle and use the corresponding
evaluation kernels to select a subset of these points.
The operator error will control the projection error, and the total quadrature weight will control \(b(A)\). The determinantal mean identity
follows from known finite-dimensional projection formulas; here
we use it to select one radius and one subset satisfying both bounds.

\subsection{Truncation and convergence of the minimizers}\label{subsec:06-minimizers}

Let \(Q_e\) and \(\rho_e\) be the functional and density from Section~\ref{sec:05}.
By Theorem~\ref{thm:cost-domain}, the finite domain of \(Q_e\) is \(M\).
On that closed subspace let \(A_M=T_e^*T_e\), where \(T_e\) is the
weighted trace in \eqref{eq:cost-trace}. Then
\[
 \langle A_M f,g\rangle_{\D}
   =\int_\T \rho_e f^*\overline{g^*}\,dm \quad(f,g\in M),
 \qquad 0\le A_M\le 4I_M.
\]
The bound follows from \(\|T_e\|\le2\). In particular, the graph norm
\((\|f\|_{\D}^2+Q_e(f))^{1/2}\) lies between \(\|f\|_{\D}\) and
\(\sqrt{5}\|f\|_{\D}\) on \(M\). It is complete because \(M\) is closed.

For \(t>0\), we first minimize the extended functional
\(J_t(f)=\|1-f\|_{\D}^2+tQ_e(f)\), whose finite domain is
\(M\). The parameter \(t\) controls how far its minimizer is from
\(P_M1=ae\).
On \(M\), orthogonality gives
\[
 J_t(f)=\|q_*\|_{\D}^2+\|ae-f\|_{\D}^2
                     +t\langle A_M f,f\rangle_{\D}.
\]
Its unique minimizer is \(x_t=(I_M+tA_M)^{-1}ae\).
Expanding the quadratic at \(x_t\) cancels the linear terms by
\((I_M+tA_M)x_t=ae\) and leaves the difference
\(\|f-x_t\|_{\D}^2+tQ_e(f-x_t)\) for \(f\in M\).
The spectral bound \(0\le A_M\le4I_M\) gives
\begin{equation}\label{eq:finite-zero-penalty}
 \|x_t-ae\|_{\D}\le \frac{4ta}{1+4t}.
\end{equation}
Here the error operator on \(M\) is \(tA_M(I_M+tA_M)^{-1}\),
and \(t\lambda/(1+t\lambda)\) is increasing for \(\lambda\ge0\).

To discretize the boundary integral \(tQ_e(f)\), we need a bounded form on all of \(\D\).
We obtain it by cutting off the density at a finite height. Write
\(\iota:\D\to L^2(m)\) for the bounded trace map \(f\mapsto f^*\).
For any bounded nonnegative density \(\nu\), let
\(T_\nu=\iota^*M_\nu\iota\), where \(M_\nu\) is multiplication by
\(\nu\) on \(L^2(m)\). This operator acts on all of \(\D\) and satisfies
\(0\le T_\nu\le\|\nu\|_\infty I_{\D}\). For \(t>0\) and \(n\ge1\) put
\[
 \nu_{t,n}=\min\{t\rho_e,n\},\qquad
 J_{t,n}(f)=\|1-f\|_{\D}^2+\int_\T\nu_{t,n}|f^*|^2\,dm .
\]
The functional \(J_{t,n}\) has the  unique minimizer
\[
 x_{t,n}=(I_{\D}+T_{\nu_{t,n}})^{-1}1,\qquad
 \|x_{t,n}\|_{\D}\le1.
\]
If \(m_t=J_t(x_t)\) and \(m_{t,n}=J_{t,n}(x_{t,n})\), expansion at
this ambient minimizer gives, for all \(f\in\D\),
\begin{equation}\label{eq:finite-zero-ambient-defect}
 J_{t,n}(f)-m_{t,n}=\|f-x_{t,n}\|_{\D}^2+
        \int_\T \nu_{t,n}|(f-x_{t,n})^*|^2\,dm .
\end{equation}

The monotone convergence here is an instance of the strong-resolvent
theorem for increasing closed forms, with limiting domain \(M\)
\cite[Theorems~3.1 and 4.1]{Simon78}. The following argument
also gives a bound in terms of the difference of the minima.

\begin{lemma}[Strong convergence of the truncated minimizers]\label{lem:finite-zero-truncation}
For each fixed \(t>0\), \(m_{t,n}\uparrow m_t\) and
\begin{equation}\label{eq:finite-zero-truncation}
 \|x_t-x_{t,n}\|_{\D}^2\le m_t-m_{t,n}\longrightarrow0 .
\end{equation}
\end{lemma}

\begin{proof}
We first identify the limiting minimum. A fixed cutoff lets us test
a weak limit of the minimizers; removing the cutoff then shows that
\(Q_e\) is finite at this limit, so the limit belongs to \(M\).
The minima \(m_{t,n}\) increase by monotonicity of the densities.
Testing with \(x_t\) gives \(m_{t,n}\le m_t\), and testing \(J_t\)
with \(0\) gives \(m_t\le1\). The preceding norm bound makes the minimizers bounded in the
Dirichlet Hilbert space, so a weakly convergent subsequence is available.
Let \(L=\lim_n m_{t,n}\) and take a
weakly convergent subsequence \(x_{t,n_j}\rightharpoonup x\) in
\(\D\). For a fixed \(k\), the functional \(J_{t,k}\) is weakly
lower semicontinuous, since it is a sum of squared norms of bounded
linear images, with a fixed translation in the first term. For
\(n_j\ge k\), it follows that
\[
 J_{t,k}(x)\le\liminf_j J_{t,k}(x_{t,n_j})
                \le \lim_j m_{t,n_j}=L .
\]
Monotone convergence as \(k\to\infty\) gives \(J_t(x)\le L\).
Thus \(x\in M\) by the exact domain theorem. By minimality,
\(m_t\le J_t(x)\le L\le m_t\). Finally, \eqref{eq:finite-zero-ambient-defect}
at \(f=x_t\) bounds the squared error by \(m_t-m_{t,n}\) for every \(n\).
\end{proof}

For a prescribed error, choose \(t>0\) small enough in
\eqref{eq:finite-zero-penalty}, then keep it fixed while increasing
\(n\) in Lemma~\ref{lem:finite-zero-truncation}. Since \(t>0\),
removing the cutoff gives \(J_t\), whose finite domain is \(M\). The order matters. At fixed \(n\), \(\nu_{t,n}\to0\) almost
everywhere as \(t\downarrow0\). Comparing \(J_{t,n}(x_{t,n})\) with
\(J_{t,n}(1)=\int_\T\nu_{t,n}\,dm\) and using dominated convergence
gives \(x_{t,n}\to1\) in \(\D\). Removing the cutoff before letting
\(t\) tend to zero gives \(ae\). Letting \(t\) tend to zero at fixed
cutoff gives \(1\), so the subsequent limit in \(n\) stays at \(1\).

\subsection{Radial approximation and positive quadrature}\label{subsec:06-quadrature}

We approximate the boundary integral by integration over a circle inside the disk, where point evaluation is bounded. Fix the bounded density \(\nu=\nu_{t,n}\) and
its total weight \(V_\nu=\int_\T\nu\,dm>0\). The positivity follows
from that of \(\rho_e\) almost everywhere. For \(0<r<1\), define
\(\iota_r f(\xi)=f(r\xi)\).
The coefficient norm \eqref{eq:full-norm} gives
\[
 \|\iota_r-\iota\|_{\D\to L^2(m)}^2
   =\sup_{k\ge0}\frac{(1-r^k)^2}{k+1}\le1-r,\qquad
 \|\iota_r\|=\|\iota\|=1.
\]
In fact \(1-r^k\le\min\{1,k(1-r)\}\), whose square is at most
\(k(1-r)\). For \(T_{\nu,r}=\iota_r^*M_\nu\iota_r\), insert the mixed term
\(\iota^*M_\nu\iota_r\) in the difference to obtain
\begin{equation}\label{eq:finite-zero-radial}
 \|T_{\nu,r}-T_\nu\|\le 2\|\nu\|_\infty\sqrt{1-r}.
\end{equation}
These are operator norms on the indicated spaces. At a fixed
radius, we can approximate the radial operator by assigning all the
mass of each short arc to one point. This preserves the total mass
exactly.

\begin{proposition}[Positive quadrature]\label{prop:positive-quadrature}
Let \(\nu\in L^\infty(m)\) be nonnegative with total weight \(V_\nu>0\).
At each fixed \(r\in(0,1)\) and for any angular mesh bound \(\eta>0\),
there is a nonempty finite set of distinct points \(a_i\) of modulus
\(r\), with positive weights \(w_i\), such that the operator
\[
 T=\sum_i w_i(k_{a_i}\otimes k_{a_i}),\qquad
 (u\otimes u)f=\langle f,u\rangle_{\D}u,
\]
satisfies
\begin{equation}\label{eq:finite-zero-quadrature}
 \sum_i w_i=V_\nu,\qquad
 \operatorname{Tr}T=K(r)V_\nu,\qquad
 \|T-T_{\nu,r}\|\le 2V_\nu\sqrt{K(r)}\,\frac{r\eta}{1-r^2}.
\end{equation}
\end{proposition}

\begin{proof}
The evaluation identity from Section~\ref{subsec:02-norm} gives
\[
 T_{\nu,r}=\int_\T(k_{r\xi}\otimes k_{r\xi})\,\nu(\xi)\,dm(\xi),
\]
where the operator integral exists in norm by continuity of the
kernel and finiteness of the measure \(\nu\,dm\). For that continuity
and a quantitative modulus, the coefficient formula yields
\[
 \left\|\frac{d}{d\theta} k_{re^{i\theta}}\right\|_{\D}^2
    =\sum_{k\ge1}\frac{k^2r^{2k}}{k+1}
    \le\sum_{k\ge1}kr^{2k}=\frac{r^2}{(1-r^2)^2}.
\]
The derivative series converges uniformly in \(\theta\) in \(\D\)
at this fixed \(r\). Since
\(\|u\otimes u-v\otimes v\|\le(\|u\|+\|v\|)\|u-v\|\),
the rank-one kernel oscillates by at most
\(2\sqrt{K(r)}\,r\eta/(1-r^2)\) on any arc of angular length at most
\(\eta\).

Take a finite partition of \(\T\) into disjoint half-open arcs of
angular length at most \(\eta\). In each arc \(I_i\) of positive
\(\nu\,dm\)-mass, choose an interior point \(\xi_i\), set \(a_i=r\xi_i\)
and put \(w_i=\int_{I_i}\nu\,dm\). The points \(a_i\) are distinct, and the weights sum to \(V_\nu\). Integrating the oscillation bound over
the arcs gives the operator error in \eqref{eq:finite-zero-quadrature}.
The trace of each rank-one operator is \(K(r)\), which gives the
trace formula.
\end{proof}

We now carry the radial and quadrature errors over to the
minimizing vector. For bounded nonnegative self-adjoint operators
\(S_1,S_2\) on one Hilbert space,
\[
 (I+S_1)^{-1}-(I+S_2)^{-1}
      =(I+S_1)^{-1}(S_2-S_1)(I+S_2)^{-1}.
\]
The two inverse factors are contractions. Use this identity on \(\D\)
to compare \(T\) with \(T_{\nu,r}\) and \(T_{\nu,r}\) with \(T_\nu\).
Combining these comparisons with \eqref{eq:finite-zero-penalty} and
\eqref{eq:finite-zero-truncation} gives
\begin{equation}\label{eq:finite-zero-resolvent}
 \begin{split}
 \|(I+T)^{-1}1-ae\|_{\D}
 &\le \frac{4ta}{1+4t}+\sqrt{m_t-m_{t,n}}\\
 &\quad+2\|\nu\|_\infty\sqrt{1-r}
         +2V_\nu\sqrt{K(r)}\,\frac{r\eta}{1-r^2}.
 \end{split}
\end{equation}
The four terms come from adding \(tQ_e(f)\) to \(\|1-f\|_{\D}^2\),
truncating the weight in \(Q_e\), replacing \(f^*(\zeta)\) by
\(f(r\zeta)\), and approximating the resulting integral by a
weighted sum over finitely many points on the circle \(|z|=r\). We choose \(r\) after
fixing \(\nu\) and refine the mesh after fixing \(r\). In particular,
refinement leaves \(V_\nu\) unchanged. The estimate for \(b(A)\) below will use this fact.

\subsection{Determinantal average of projections}\label{subsec:06-determinant}

A Gram determinant records the geometry of a selected collection
of vectors: it is nonnegative and vanishes when the vectors are
linearly dependent. Normalized Gram determinants can therefore
assign probabilities to subsets of the fixed evaluation points. Lyons develops the theory of discrete determinantal measures
through Hilbert-space projections, positive contractions,
and combinatorial structure \cite{Lyons2003}. Hough, Krishnapur, Peres, and Vir\'ag explain these probability laws through mixtures of projection processes and the independent random variables underlying point counts
\cite[Theorem~7, p.~210]{HKPV2006}. This provides an operator and probabilistic
setting for the finite averaging argument below. 

The operator from positive quadrature includes every point with a weight.
To replace it by one zero set, we average the orthogonal projections
associated with subsets of those points. The identity gives the
average in terms of the quadrature operator. The next subsection uses
this formula to choose one subset. The determinant weights are the
standard \(L\)-ensemble weights \cite[Section~2.2, (2.6) and Theorem~2.1]{KT12}.
Mutn\'y, Derezi\'nski, and Krause give a formula for the mean
embedded inverse for real positive-definite matrices
\cite[Theorem~1 and Remark~1]{MDK20}. Fanuel, Schreurs, and Suykens
state the corresponding expected-projector identity for real
positive-definite kernels
\cite[Corollary~2.3, p.~286]{FSS21}.
For arbitrary real matrices, including rank-deficient ones,
Derezi\'nski, Khanna, and Mahoney give the expected-projector
identity explicitly
\cite[Lemma~5 in Appendix~B, with proof in Appendix~C]{DKM20}.
Lyons's conditional-projection formula
\cite[Lemma~7.17 and equation~(7.7)]{Lyons2003},
applied to the graph of \(X\), also yields the identity below.
We give a direct proof in the present notation, including the
complex-vector and linearly dependent cases.

\begin{lemma}[Determinantal average of projections]\label{lem:determinantal-projections}
Let \(v_1,\ldots,v_m\) be vectors in a complex Hilbert space \(H\)
whose inner product is linear in the first variable. Let
\(X:\mathbb C^m\to H\) be defined by \(Xc=\sum_i c_iv_i\), and put
\(L=X^*X\), \(T=XX^*\). For \(S\subseteq\{1,\ldots,m\}\), write
\(L_S\) for the principal submatrix and \(P_S\) for the orthogonal
projection onto the span of the selected vectors, with
\(P_\varnothing=0\). The numbers
\[
 p_S=\frac{\det L_S}{\det(I+L)},\qquad \det L_\varnothing=1
\]
are probabilities and satisfy
\begin{equation}\label{eq:determinantal-projector-mean}
 \sum_S p_S P_S=X(I+L)^{-1}X^*=T(I+T)^{-1}.
\end{equation}
\end{lemma}

\begin{proof}
The principal-minor expansion gives
\(\det(I+L)=\sum_S\det L_S\). Since \(L\) is a Gram matrix, all
these minors are nonnegative and the empty term is one.

We identify the average projection by differentiating the
principal-minor expansion. The derivative of each nonsingular minor
will be a multiple of the squared norm of a projection. Fix \(h\in H\) and set \(u=X^*h\).
Apply the same expansion to
\(L+suu^*\) and differentiate in the real parameter \(s\) at zero.
The derivative of
\(\det(I+L+suu^*)\) is
\(\det(I+L)u^*(I+L)^{-1}u\). For nonempty \(S\) with invertible
\(L_S\), the derivative of its term is
\[
 \det L_S\,u_S^*L_S^{-1}u_S
       =\det L_S\,\|P_Sh\|_H^2,
\]
using \(P_S=X_S L_S^{-1}X_S^*\), where \(X_S\) restricts \(X\) to
the selected columns. If \(L_S\) is singular, then
\(\ker L_S=\ker X_S\), so \(u_S=X_S^*h\) is orthogonal to that
kernel. The perturbation \(L_S+su_Su_S^*\) therefore annihilates
the same nonzero kernel vectors and its determinant is identically
zero. The empty term has derivative zero as well.

The singular and empty terms contribute zero, so the differentiated
expansion accounts for every subset. Divide by \(\det(I+L)\) and
polarize the equality of quadratic forms to obtain the first identity in
\eqref{eq:determinantal-projector-mean}. The relation
\((I+T)X=X(I+L)\) gives the second.
\end{proof}

\subsection{Joint bounds for projection error and the logarithmic sum over zeros}\label{subsec:06-selection}

We apply the projection identity to the weighted evaluation kernels.
Applying the average to \(1\) controls the mean squared projection
error. Taking the trace controls the expected number of selected points, hence the expected value of \(b(A_S)\) at the fixed radius. Both estimates use the same weights.

Use the points \(a_i\) and weights \(w_i\) of
Proposition~\ref{prop:positive-quadrature} and apply
Lemma~\ref{lem:determinantal-projections} to
\(v_i=\sqrt{w_i}\,k_{a_i}\). Its operator \(T\) is the quadrature
operator. For each subset \(S\), let \(A_S=\{a_i:i\in S\}\).
The selected span is \(M_{A_S}^\perp\), so \(P_S1=q_{A_S}\)
in the notation of \eqref{eq:finite-zero-projection}.

Write \(\mathbb E\) for the finite average with weights \(p_S\),
and let \(q=\mathbb E q_{A_S}\). Applying
\eqref{eq:determinantal-projector-mean} to the constant function gives
\[
 q=1-(I+T)^{-1}1,\qquad
 q_{A_S}(0)=\langle q_{A_S},1\rangle_\D=\|q_{A_S}\|_\D^2,
\]
where the last equality follows from \(q_{A_S}=P_S1\). The target
\(q_*=1-ae\) is also an orthogonal projection of \(1\).
These identities turn the mean squared error into an expression
linear in the mean projection vector. Expanding the square and using
\(q_*(0)=\|q_*\|_\D^2\) gives
\begin{align}
 \mathbb E\|q_{A_S}-q_*\|_\D^2
 &=\operatorname{Re}\langle q-q_*,1-2q_*\rangle_\D \notag\\
 &\le \|(I+T)^{-1}1-ae\|_\D . \label{eq:expected-finite-zero-error}
\end{align}
The inequality uses \(\|1-2q_*\|_\D=1\) and
\(q-q_*=ae-(I+T)^{-1}1\). 

The same weights control \(\mathbb E b(A_S)\). All the points \(a_i\) have
modulus \(r\), so \(b(A_S)=(-\log r)|S|\).
Each positive-weight subset has independent selected vectors and
\(\operatorname{Tr} P_S=|S|\). Taking traces in
\eqref{eq:determinantal-projector-mean} and using
\(0\le T(I+T)^{-1}\le T\) gives
\begin{align}
 \mathbb E b(A_S)
 &=(-\log r)\operatorname{Tr}\bigl(T(I+T)^{-1}\bigr)\notag\\
 &\le (-\log r)K(r)V_\nu
 \le 8V_\nu(1-r)\log\frac{1}{1-r}
 \qquad (r\ge1/2). \label{eq:expected-blaschke-budget}
\end{align}
Here \(\operatorname{Tr}T=\sum_i w_i\|k_{a_i}\|_\D^2=K(r)V_\nu\).
For the last bound use
\(-\log r\le 2(1-r)\), \(r^{-2}\le 4\), and
\(-\log(1-r^2)\le \log(1/(1-r))\) in the kernel diagonal
formula of Section~\ref{subsec:02-norm}. The bound tends to zero
as \(r\uparrow1\) at fixed \(V_\nu\), independently of the mesh.

Add the two nonnegative quantities before choosing the subset.
Some subset of positive weight has their sum at most its mean,
so one set gives the following bound.
\begin{equation}\label{eq:finite-zero-combined-cost}
 \|q_A-q_*\|_\D^2+b(A)
 \le \|(I+T)^{-1}1-ae\|_\D+(-\log r)K(r)V_\nu .
\end{equation}
For the empty subset, the sum is exactly
\(\|q_*\|_\D^2=1-a^2>0\). Thus an upper bound strictly below
\(1-a^2\) ensures a nonempty subset.

\subsection{The finite-zero canonical function}\label{subsec:06-normalization}

The bound on \(q_A\) gives the same bound on the error of its
complementary vector. We now normalize that vector and bound the approximation error
in terms of the fixed target value \(a\). For a nonempty finite set
\(A\subset\DD\setminus\{0\}\), let \(y_A=1-q_A=P_{M_A}1\) and
\(\gamma_A=\|y_A\|_\D>0\).
The strict positivity follows from the polynomial in \(M_A\)
used above to define \(e_A\). By
\eqref{eq:canonical-definition}, the unit canonical
function is \(e_A=y_A/\gamma_A\), with
\(e_A(0)=\gamma_A=\sqrt{1-q_A(0)}\).
The wandering-generation theorem, Theorem~\ref{input:generation},
gives \([e_A]=M_A\). The reverse triangle inequality gives
\[
 |e_A(0)-a|\le \|y_A-ae\|_\D=\|q_A-q_*\|_\D.
\]
Then \(a(e_A-e)=(y_A-ae)+(a-\gamma_A)e_A\) gives
\begin{equation}\label{eq:finite-zero-normalized-error}
 \|e_A-e\|_\D\le\frac{2}{a}\|q_A-q_*\|_\D .
\end{equation}

\begin{proof}[Completion of the proof of Theorem~\ref{thm:finite-zero-approximation}]
We choose the coefficient \(t\), the cutoff, the radius, and the angular mesh in that order. For the prescribed \(s>0\), let
\(\varepsilon=\min\{s,(1-a^2)/2\}\). The second bound in this minimum
excludes the empty subset. In
\eqref{eq:finite-zero-resolvent}, choose \(t>0\) so that the first
term is less than \(\varepsilon/8\). At this fixed \(t\), choose
\(n\) so that the second term is less than \(\varepsilon/8\).
Fix \(\nu=\min(t\rho_e,n)\) and its total \(V_\nu\in(0,n]\).
Next choose \(r>\max\{1/2,1-s\}\) close enough to \(1\) to
make the radial term less than \(\varepsilon/8\) and to have
\((-\log r)K(r)V_\nu<\varepsilon/2\).
At this fixed radius, refine the mesh to make the quadrature term
less than \(\varepsilon/8\). The total quadrature weight is still
\(V_\nu\), so the bound on \(\mathbb E b(A_S)\) holds for the resulting evaluation points.

The right side of \eqref{eq:finite-zero-combined-cost} is now
strictly below \(\varepsilon\). This gives a nonempty set \(A\) of
distinct points at the chosen radius with
\(\|q_A-q_*\|_\D^2+b(A)<s\). The estimates for \(|e_A(0)-a|\) and
\(\|e_A-e\|_\D\) follow from the preceding normalization. For separate positive bounds on \(\|e_A-e\|_\D\) and \(b(A)\), choose
\(s\) small enough that \(2\sqrt{s}/a\) and \(s\) are smaller than the respective bounds.
\end{proof}

For the fixed \(e\), choose \(s_j\downarrow0\) and use the theorem
at each \(s_j\) to select a nonempty finite set \(A_j\). The points
in \(A_j\) have a common radius \(r_j\). Along this sequence,
\[
 q_{A_j}\to q_*,\qquad e_{A_j}\to e\quad\text{in }\D,\qquad
 e_{A_j}(0)\to a,\qquad 0<b(A_j)\to0,\qquad r_j\to1.
\]

The sum \(b(A)\) also connects the construction to the
classical problem of infinite interior zero sequences. Caughran presents Carleson's example of a Blaschke
sequence with a single boundary accumulation point. No nonzero
analytic function of finite Dirichlet integral vanishes at all
points of this sequence \cite[Theorem~1]{Caughran69}. Caughran credits the modified
presentation to Shields, Duren, and himself.
Thus the Blaschke condition and a single boundary accumulation
point do not by themselves ensure a Dirichlet zero sequence.
His complementary positive result assumes, in addition to the
Blaschke condition, that the limit points lie in a Carleson set
and that the sequence has finite contact with that set
\cite[Theorem~2]{Caughran69}. The contact condition bounds the
distance of each zero from the set by a constant times a fixed
positive power of its distance from the circle. Under these
hypotheses he obtains a function with the prescribed zeros whose
derivative extends continuously to the closed disk.

Bogdan further organizes the interior-zero question
by geometry and randomness \cite{Bogdan96}. He characterizes
subsets of the disk on which every Blaschke sequence is a Dirichlet
zero sequence, and also studies the almost-sure zero-sequence
question for prescribed moduli with independent uniformly
distributed arguments. These questions concern the admissibility
of a prescribed infinite sequence. In the construction above,
each selected finite set is already the zero set of a polynomial;
the quantitative task is to choose successive finite sets whose
canonical functions converge in norm to the fixed target
while \(b(A_j)\to0\). The next subsection
applies that approximation to canonical functions of trace-zero
spaces.

\subsection{Approximation of trace-zero extremals}\label{subsec:06-trace-zero}

Luo and Richter ask whether the extremal of a nonzero trace-zero
space \(\D_E\), for a compact \(E\subset\T\) of positive logarithmic
capacity, is a locally uniform limit of interior zero-set extremals
\cite[Question~6.3]{LR15}. They prove finite-zero approximation for
spaces of the form \(SH^2\cap\D\) with \(S\) inner
\cite[Theorems~1.4 and 4.7]{LR15}. The finite-zero theorem above gives
the following consequence for the trace-zero class in their question.

\begin{corollary}[Trace-zero extremals]\label{cor:trace-zero-finite-approx}
Let \(E\subset\T\) be compact with positive logarithmic capacity and
\(\D_E\ne\{0\}\), where \(\D_E\) is the trace-zero space from
Section~\ref{subsec:03-ideal}. There are nonempty finite sets
\(A_j\subset\DD\setminus\{0\}\), each consisting of distinct points
of one radius \(r_j\to1\), such that
\[
 \|G_{M_{A_j}}-G_{\D_E}\|_{\D}\longrightarrow0,\qquad
 b(A_j)\longrightarrow0.
\]
Moreover, \(P_{M_{A_j}}\to P_{\D_E}\) in the strong operator
topology on \(\D\).
\end{corollary}

\begin{proof}
We first find an outer member of the trace-zero space. Its value at
zero ensures the canonical projection is nonzero, and its outerness
determines the inner factor of the extremal. Choose \(0\ne f\in\D_E\)
and write \(f=IO\) in its Hardy inner--outer factorization.
Lemma~\ref{lem:inner-division} gives
\(O\in\D\). Since \(|f^*|=|O^*|\) almost everywhere, their
quasi-continuous moduli agree quasi-everywhere by
Section~\ref{subsec:03-trace}. Hence \(O\in\D_E\) and \(O(0)\ne0\).
The space \(\D_E\) is closed and shift-invariant by
Section~\ref{subsec:03-ideal}. Positivity of capacity excludes
the constant \(1\), so \(G_{\D_E}\) has value at zero in \((0,1)\) by
\eqref{eq:canonical-definition}. If \(J\) is its Hardy inner factor,
generation by \(G_{\D_E}\) and the continuous inclusion
\(\D\subset H^2\) put \(\D_E\) in \(JH^2\). This Hardy subspace
contains \(O\), so \(J\) is constant. Thus \(G_{\D_E}\in\Xpr\).

Apply Theorem~\ref{thm:finite-zero-approximation} to \(G_{\D_E}\)
with \(s_j\downarrow0\). The resulting norm convergence implies local
uniform convergence. With normalized arclength on the unit circle as the defining measure,
Luo and Richter use the Dirichlet norm in \eqref{eq:full-norm}.
The canonical functions here are the corresponding unit extremals,
normalized to have positive value at zero. Their projection theorem now gives the stated
strong operator convergence, that is, norm convergence of the
projected vectors for each fixed vector in \(\D\)
\cite[Theorem~1.3 and Corollary~4.5]{LR15}.
\end{proof}

Luo and Richter deduced the finite-set reduction and the implication
for Question~6.2 immediately after \cite[Question~6.3]{LR15}.
The preceding corollary supplies the required convergence for this
trace-zero class. For a Carleson set of positive capacity, Yilmaz's dissertation places the closed
span of zero-based subspaces from Blaschke tails between the subspace generated by the smooth outer function of Theorem~\ref{input:flat-outer}
and the trace-zero space.
The inclusion of the canonical extremal in that span is the remaining step identified in the dissertation
\cite[Theorems~1.2.3 and 3.2.4]{Yilmaz17}.
The same partial inclusion for the classical Dirichlet space
is established in the published article
\cite[Theorem~4.2]{Yilmaz18}.
Yilmaz also gives finite-zero approximation of trace-zero spaces
in the harmonically weighted Dirichlet space defined by a finite
sum of point masses \cite[Theorem~3.5]{Yilmaz23}. Borichev,
El-Fallah, and Kellay characterize when a prescribed sequence is
a Dirichlet zero set by using its finite prefixes and successive
weighted kernel diagonals \cite[Lemma~2.2 and Theorem~3.1]{BEK26}.
Here Theorem~\ref{thm:finite-zero-approximation} selects a finite set of points for a given canonical target and both prescribed bounds.

\subsection{Approximation of invariant subspaces}

\begin{corollary}[Finite-codimensional approximation]
\label{cor:finite-codimensional-approximation}
Let \(M\) be a nonzero closed shift-invariant subspace of \(\D\).
There is a sequence of closed shift-invariant subspaces \(M_j\),
each obtained by requiring vanishing to prescribed orders at
finitely many points of \(\DD\), such that
\[
  \|P_{M_j}h-P_Mh\|_{\D}\longrightarrow 0
  \qquad (h\in\D).
\]
In particular, each \(M_j\) has finite codimension.
\end{corollary}

\begin{proof}
By Theorem~\ref{input:generation}, write \(M=[f]\) and take the
Hardy factorization \(f=IO\). The inner--outer identity
\cite[Theorem~5.3]{RS92} gives \(O\in\D\) and
\[
  M=N\cap L,\qquad N=[O],\qquad L=IH^2\cap\D.
\]
Since \(O(0)\ne0\), the canonical function \(e=G_N\) is defined.
If \(J_0\) is the Hardy inner factor of \(e\), then the continuous
inclusion \(\D\subset H^2\) gives \(N=[e]\subset J_0H^2\).
Since \(O\in N\) is outer, \(J_0\) is constant. Thus \(e\) is outer.

Let \(k\ge0\) be the least integer such that
\(h\mapsto h^{(k)}(0)\) is nonzero on \(L\), and let \(g\) be the
unit extremal maximizing \(\operatorname{Re}h^{(k)}(0)\) over
the unit ball of \(L\). The derivative functional vanishes on \(zL\),
so \(g\in L\ominus zL\) and \([g]=L\) by
Theorem~\ref{input:generation}. Write \(g=JV\) in its Hardy
factorization. The inner--outer identity gives \(V\in\D\) and
\([g]=[V]\cap JH^2\). Since \(e\) is a multiplier, both \(eV\)
and \(eg\) belong to \(\D\). The product \(eV\) is outer. By the
inner--outer identity and the outer product formula
\cite[Theorem~4.5]{RS92},
\[
  [eg]=[eV]\cap JH^2
      =[e]\cap[V]\cap JH^2
      =N\cap L=M.
\]

If \(e\ne1\), then \(e\in\Xpr\), and
Theorem~\ref{thm:finite-zero-approximation} gives finite-zero
spaces \(N_j\) whose unit extremals \(e_j\) converge to \(e\)
in \(\D\). If \(e=1\), take \(N_j=\D\) and \(e_j=1\).
With \(m\) as the defining measure, Luo and Richter's theorem
\cite[Theorem~4.7]{LR15} gives finite-zero spaces \(L_j\) whose
unit extremals \(g_j\) converge locally uniformly to \(g\).
Since point evaluations are continuous on \(\D\), every weak
subsequential limit of the bounded sequence \((g_j)\) agrees
pointwise with \(g\). Thus \(g_j\rightharpoonup g\) in \(\D\).
Since \(\|g_j\|_{\D}=\|g\|_{\D}=1\), the convergence is strong.
Moreover, \cite[Corollary~4.5]{LR15} gives \(P_{N_j}\to P_N\)
and \(P_{L_j}\to P_L\) in the strong operator topology.

Set \(M_j=N_j\cap L_j\). The vanishing conditions defining
\(M_j\) are those of both spaces, with the larger prescribed
order at any common zero. They are finitely many continuous
derivative evaluations, so \(M_j\) has finite codimension.
Let \((j_k)\) be strictly increasing and suppose that
\(u_k\in M_{j_k}\) converges weakly to \(u\). For every
\(v\in N^\perp\),
\[
  \langle u_k,v\rangle_{\D}
   =\langle u_k,P_{N_{j_k}}v\rangle_{\D}\longrightarrow0.
\]
Hence \(u\in N\), and the same argument with \(L_{j_k}\) gives
\(u\in L\). Thus \(u\in M\).

Finally, \(e_jg_j\in M_j\). The unit wandering vectors \(e_j\)
and \(g\) act by contractive multiplication by
\eqref{eq:canonical-multiplier}, so
\[
  \|e_jg_j-eg\|_{\D}
  \le \|g_j-g\|_{\D}+\|e_j-e\|_{\D}\longrightarrow0.
\]
For every \(h\in M\) and every polynomial \(p\), we have
\(pe_jg_j\in M_j\), and hence
\[
  \limsup_{j\to\infty}\operatorname{dist}_{\D}(h,M_j)
  \le \|h-peg\|_{\D}.
\]
Since \([eg]=M\), the right-hand side can be made arbitrarily small.
Thus \(P_{M_j}h\to h\) for every \(h\in M\).
Lemma~\ref{lem:finite-recovery-projections} now gives
\(P_{M_j}h\to P_Mh\) for every \(h\in\D\).
\end{proof}

\subsection{Locally uniform limits of Dirichlet extremals}

For a nonzero closed shift-invariant subspace \(M\) of \(\D\), the
extremal is the unit vector maximizing \(\operatorname{Re}h^{(k)}(0)\)
over the unit ball of \(M\), where \(k\) is the least order of a
derivative evaluation at zero that is nonzero on \(M\). In this
subsection we allow unimodular multiples of these extremals. We write
\(U_f\) for the function \(u_f\) defined by Luo and Richter
\cite[p.~433]{LR15}. Their Question~6.4 \cite[p.~436]{LR15}
asks whether the locally uniform limits of unimodular multiples of
Dirichlet extremals are exactly the functions \(f\in\D\) satisfying
\(u_f\le1\). In our notation,
\[
 U_f(\lambda)=\|f\|_{\D}^2+
 2\operatorname{Re}\sum_{n\ge1}\bar\lambda^{\,n}
                       \langle z^nf,f\rangle_{\D},
 \qquad \lambda\in\DD .
\]
The series converges locally uniformly by the shift bound in
Section~\ref{subsec:02-norm}.

\begin{theorem}[Locally uniform limits of extremals]\label{thm:locally-uniform-extremals}
The locally uniform closures of the unimodular multiples of all
Dirichlet extremals and of the unimodular multiples of the
extremals of spaces defined by finitely many interior vanishing
conditions, with multiplicities allowed, both equal
\[
 \{f\in\D: U_f(\lambda)\le1\text{ for every }\lambda\in\DD\}.
\]
If \(\|f\|_{\D}=1\), unimodular multiples of extremals of the
finite-zero spaces can be chosen to converge to \(f\) in \(\D\).
\end{theorem}

\begin{proof}
Luo and Richter's necessity theorem \cite[Theorem~5.2, p.~433]{LR15},
with normalized arclength as the defining measure, applies here. The quantity \(U_f\) is
invariant under unimodular multiplication.

For sufficiency, let \(f_s(z)=f(sz)\) for \(0<s<1\). We first show
that \(U_{f_s}(\lambda)\le U_f(s\lambda)\). For polynomials, the
matrix of the difference on coefficient vectors, up to a diagonal
unitary change of phase, has entries
\[
 (\max(i,j)+1)(s^{|i-j|}-s^{i+j})|\lambda|^{|i-j|}.
\]
The first two factors equal \(b_i b_j c_{\min(i,j)}\), where
\(b_j=(j+1)s^j\) and \(c_j=(s^{-2j}-1)/(j+1)\). The sequence
\(c_j\) is increasing: for \(j\ge1\), both
\((s^{-2j}-1)/j\) and \(j/(j+1)\) are increasing, and \(c_0=0\).
The matrix \((c_{\min(i,j)})_{i,j}\) is positive semidefinite, since
\(c_{\min(i,j)}=\sum_{1\le k\le \min(i,j)}(c_k-c_{k-1})\).
Multiplying the rows and columns by \(b_i\) and \(b_j\), respectively,
preserves positivity. The Schur product with the positive Poisson
matrix \((|\lambda|^{|i-j|})_{i,j}\) does so as well. Polynomial density
and \(\|M_z^n\|=\sqrt{n+1}\) give the inequality for every \(f\in\D\).

Suppose \(f\ne0\) and \(U_f\le1\). Choose \(s\uparrow1\) avoiding
the moduli of the zeros of \(f\), and \(0<t\uparrow1\). The functions
\(g=tf_s\) are analytic on a neighborhood of the closed disk and
have no boundary zeros. Their Taylor coefficients decay geometrically,
so the moments defining \(U_g\) are absolutely summable. The dilation
inequality therefore gives \(U_g\le t^2<1\) on \(\T\). Fix one such \(g\),
and let \(M\) be the space defined by its interior zeros with their
multiplicities. Let \(p_0\) be the monic polynomial whose zeros are exactly the
interior zeros of \(g\), counted with multiplicity. If \(g\) has
no interior zeros, set \(p_0=1\). Then \(M=p_0\D\). Indeed, for the contractive backward shift \(T_0\),
\[
 \frac{h(z)-h(a)}{z-a}=(I-aT_0)^{-1}T_0h,\qquad a\in\DD ,
\]
so division by each zero factor preserves membership in \(\D\).
Both \(g/p_0\) and its inverse are analytic near the closed disk
and are multipliers by Lemma~\ref{lem:bounded-products}. Polynomial
multiples of \(g\) are therefore dense in \(M\).

On \(\T\) set \(v=(1-U_g)/|g^*|^2\), put
\(\gamma=\int_\T v\,dm\), and set \(d\rho=v\,dm/\gamma\). The density
\(v\) is continuous and strictly positive. For boundary restriction
\(J:M\to L^2(\rho)\), use the equivalent norm
\(\|h\|_v^2=\|h\|_{\D}^2+\gamma\|Jh\|_{L^2(\rho)}^2\) on \(M\).
The Fourier coefficients of \(U_g\) give \(\|g\|_v=1\) and
\(\langle z^ng,g\rangle_v=0\) for every \(n\ge1\). Hence
\(g\perp_v zM\). The kernel in \(M\) of its first nonzero derivative
evaluation \(L\) at zero is \(zM\): division by \(z\) preserves the
other zero conditions and reduces the order at zero by one.
Thus \(g\) is, up to phase, the unit extremal for \(L\) in this norm.

We realize this modified norm by adding zeros near \(\T\). For \(R\)
sufficiently close to one, put \(N=\lfloor\gamma K(R)\rfloor\), where
\(K(R)=\|k_{R\xi}\|_{\D}^2\) is the diagonal kernel in
Section~\ref{subsec:02-norm}. Partition \(\T\) into \(N\) consecutive arcs of equal \(\rho\)-mass.
In the \(i\)-th arc, choose \(\xi_i\) so that the two subarcs on either
side of \(\xi_i\) each have \(\rho\)-mass \(1/(2N)\).
Let \(M_R=\{h\in M:h(R\xi_i)=0\text{ for }1\le i\le N\}\). 
These points are distinct and disjoint from the fixed zeros once
\(R\) is close enough to one. For fixed \(g\), the sum
\(N(-\log R)\) tends to zero, since \(K(R)=-\log(1-R^2)/R^2\).
On \(M\), we claim that
\[
 P_{M_R}|_M\longrightarrow (I+\gamma J^*J)^{-1}
 \quad\text{in the weak operator topology as }R\uparrow1 .
\]

Define \(S_Rh=(h(R\xi_i))_{i=1}^N\) in the space with normalized
counting measure, identified with the corresponding step functions
in \(L^2(\rho)\). The Gram operator \(G_R=S_RS_R^*\) has entries
\(k^M(R\xi_i,R\xi_j)/N\), where \(k^M\) is the kernel of \(M\).
The kernel of \(\D\) is \(k(z,w)=-\log(1-z\bar w)/(z\bar w)\),
with value \(1\) when \(z\bar w=0\). Its difference from \(k^M\) is a
finite-rank kernel continuous on the closed bidisk, since
\(M^\perp\) is spanned by kernels for interior evaluations and
their derivatives. Extend \(G_R\) to \(L^2(\rho)\) by the
step-function multiplier \(k^M(R\xi_i,R\xi_i)/N\) on arc \(i\)
plus the integral operator whose kernel equals
\(k^M(R\xi_i,R\xi_j)\) on different arcs and vanishes on each
same-arc rectangle. The step-function space reduces this extension.
Its multiplier part tends in norm to \(\gamma^{-1}I\).

The integral part tends in Hilbert--Schmidt norm to \(JJ^*\).
Indeed, for the boundary restriction \(J_0:\D\to L^2(\rho)\),
the orthonormal monomials give \(\sum_{n\ge0}(\zeta\bar\eta)^n/(n+1)\)
as the kernel of \(J_0J_0^*\). The series converges in \(L^2(m\times m)\)
and hence in \(L^2(\rho\times\rho)\), and its partial sums represent
the operators obtained by inserting the corresponding polynomial
projections between \(J_0\) and \(J_0^*\). Subtracting the finite-rank
correction gives the kernel of \(JJ^*=J_0P_MJ_0^*\). The step
kernels converge to it off the diagonal. The arc lengths and the distances between points \(\xi_i\) in
adjacent arcs are comparable to \(1/N\), because the density of
\(\rho\) is bounded above and below. If \(\zeta\) lies in arc \(i\)
and \(\eta\) lies in a different arc \(j\), then \(|\zeta-\eta|\)
is at most a constant times \(|\xi_i-\xi_j|\). The identity
\[
 |1-R^2\xi_i\bar\xi_j|^2
   =(1-R^2)^2+R^2|\xi_i-\xi_j|^2
\]
then bounds the absolute value of the step kernel at
\((\zeta,\eta)\) by \(C(1+\log(2/|\zeta-\eta|))\). This bound is
square-integrable against \(\rho\times\rho\), and dominated
convergence gives the asserted limit.

The extended Gram operators therefore converge in norm to
\(\gamma^{-1}I+JJ^*\), which is bounded below by \(\gamma^{-1}I\).
Their inverses converge in norm and restrict to \(G_R^{-1}\) on
the step-function spaces. Finite sums of interior kernels of \(M\)
are dense in \(M\) and extend continuously to \(\T\), so their
step samples converge to their boundary values in \(L^2(\rho)\).
Testing \(P_{M_R}|_M=I-S_R^*G_R^{-1}S_R\) on pairs of such functions
and using the uniform bound on the projections proves the claim,
since
\[
 I-J^*(\gamma^{-1}I+JJ^*)^{-1}J=(I+\gamma J^*J)^{-1}.
\]

Let \(h\) represent \(L\) in \(M\) for the original norm. Its
representer in the modified norm is
\(h_v=(I+\gamma J^*J)^{-1}h\). The claim gives
\(P_{M_R}h\rightharpoonup h_v\) and
\[
 \|P_{M_R}h\|_{\D}^2=L(P_{M_R}h)\longrightarrow
 L(h_v)=\|h_v\|_v^2>0.
\]
The first nonzero derivative evaluation in \(M_R\) is still \(L\).
Indeed, \(p_0(z)\prod_{i=1}^N(z-R\xi_i)\) belongs to \(M_R\)
and has the same order of vanishing at zero as \(p_0\). Thus normalizing \(P_{M_R}h\) and restoring the phase
of \(g\) gives unit extremals converging weakly to \(g\). Bounded
point evaluations and Montel's theorem give local uniform
convergence. A diagonal choice on expanding compact disks as
\(s,t\uparrow1\) gives \(f\). For \(f=0\), take
\(z^n/\sqrt{n+1}\), the unit extremals for \(z^n\D\). Finally,
a locally uniformly convergent sequence bounded in \(\D\) converges
weakly to its limit. If \(\|f\|_{\D}=1\), equality of norms gives
convergence in \(\D\).
\end{proof}

\textit{Relation to Theorem~\ref{thm:finite-zero-approximation}.}
Take \(e\in\Xpr\), \(a=e(0)\), and \(q_*=1-ae\). The regularized
targets \(g=te(s\,\cdot)\) have no zeros on the closed disk and
\(g(0)>0\). Thus the construction uses only the added simple
zeros \(A_R=\{R\xi_i\}\) and its extremals are \(e_{A_R}\).
For fixed \(g\) we have \(e_{A_R}\rightharpoonup g\), so
\[
 \|e_{A_R}-e\|_{\D}^2\longrightarrow
 2-2\operatorname{Re}\langle g,e\rangle_{\D},\qquad
 b(A_R)=\lfloor\gamma K(R)\rfloor(-\log R)\longrightarrow0.
\]
Given the tolerance \(s_0>0\) of
Theorem~\ref{thm:finite-zero-approximation}, first take \(g\)
sufficiently close to \(e\) in \(\D\), then take \(R\) sufficiently
close to one so that
\[
 R>\max\{1/2,1-s_0\},\qquad
 \|e_{A_R}-e\|_{\D}<\sqrt{s_0/8},\qquad b(A_R)<s_0/2 .
\]
Since \(q_{A_R}=1-e_{A_R}(0)e_{A_R}\) and evaluation at zero has
norm one, \(\|q_{A_R}-q_*\|_{\D}\le2\|e_{A_R}-e\|_{\D}\).
Hence \(\|q_{A_R}-q_*\|_{\D}^2+b(A_R)<s_0\).
The reverse triangle inequality for the norms of \(1-q_{A_R}\)
and \(1-q_*\) gives \(|e_{A_R}(0)-a|<\sqrt{s_0}\), and
\eqref{eq:finite-zero-normalized-error} gives
\(\|e_{A_R}-e\|_{\D}<2\sqrt{s_0}/a\). This gives another proof of
the finite-zero approximation used in Approximation~A.

We now return to Approximation~A. For the fixed outer canonical
target \(e\in\Xpr\), choose the sets \(A_j\) by applying
Theorem~\ref{thm:finite-zero-approximation} with \(s_j\downarrow0\).
A zero at \(r_j\xi\) will give an atom at \(\xi\) of mass \(1-r_j\).
Since \(1-r_j\le-\log r_j\), \(b(A_j)\) bounds the total
mass. The convergence of \(e_{A_j}\) gives the starting norm
approximation for the canonical comparison after the transfer.

\newpage \section{From finite zeros to singular atoms}\label{sec:07}

Theorem~\ref{thm:finite-zero-approximation} gives finite zero sets
\(A_j\) for which \(e_{A_j}\to e\) in the Dirichlet norm and
\(b(A_j)\to0\). We replace each interior zero by an atom at the
boundary point on the same ray. The proof uses the outer factor of
the finite-zero canonical function to construct a comparison vector
for the new atomic space. A comparison in both directions then gives
convergence of the canonical functions of those atomic spaces.

\subsection{Exact inner factors and convergent outer factors}\label{subsec:07-quotients}

Fix \(e\in\Xpr\) and choose \(A_j\) by applying
Theorem~\ref{thm:finite-zero-approximation} for \(s_j\downarrow0\).
To replace the factors
that record the zeros, we first identify  the exact inner factor
of each canonical function and control 	the outer factor.
For a finite set
\(A=\{r_k\xi_k:1\le k\le N\}\subset\DD\setminus\{0\}\) of distinct
points, with \(0<r_k<1\) and \(\xi_k\in\T\), put
\[
 b_k(z)=\frac{r_k-\overline{\xi_k}z}{1-r_k\overline{\xi_k}z},
 \qquad B_A=\prod_{k=1}^N b_k,\qquad
 B_A(0)=\prod_{k=1}^N r_k=e^{-b(A)}.
\]
The \(H^2\) factorization from Section~\ref{subsec:02-hardy} gives
\[
 M_A=B_A H^2\cap\D .
\]
Indeed, the Hardy Blaschke factor of any nonzero \(f\in M_A\)
contains each \(b_k\), so \(f/B_A\in H^2\). Conversely every
Dirichlet function in \(B_A H^2\) vanishes on \(A\). The polynomial
vanishing on \(A\) in Section~\ref{sec:06} gives nonzero evaluation
at \(0\) on \(M_A\). Since \(B_A(0)>0\),
Lemma~\ref{lem:exact-inner-factor} identifies \(B_A\) as the exact
inner factor of \(e_A\). Thus \(e_A=B_A h_A\) with \(h_A\in H^2\)
outer. Contractive inner division in Lemma~\ref{lem:inner-division}
then gives
\begin{equation}\label{eq:finite-factor-data}
 h_A\in\D,\qquad \|h_A\|_{\D}\le1,\qquad
 h_A(0)=e_A(0)e^{b(A)}>0 .
\end{equation}

The small value of \(b(A)\) makes the normalized inner factor close to
\(1\) at the origin. To pass from that value to a bound throughout
a compact subdisk, we use the following consequence of Schwarz--Pick.
If \(\varphi\) is analytic on \(\DD\),
\(|\varphi|\le1\), and \(\varphi(0)=u\in[0,1]\), then
\begin{equation}\label{eq:schur-origin-bound}
 \sup_{|z|\le R}|1-\varphi(z)|
       \le\frac{1+R}{1-R}(1-u),\qquad 0<R<1.
\end{equation}
For \(u<1\), write \(\varphi=(u+\psi)/(1+u\psi)\).
The Schwarz--Pick estimate \(|\psi(z)|\le |z|\) bounds
\((1-u)|(1-\psi)/(1+u\psi)|\) by the right side of
\eqref{eq:schur-origin-bound}. For \(u=1\), the maximum principle
gives \(\varphi\equiv1\).

Set \(B_j=B_{A_j}\) and \(h_j=h_{A_j}\). Since
\(0<b(A_j)\to0\), applying \eqref{eq:schur-origin-bound} to \(B_j\)
gives
\begin{equation}\label{eq:finite-factor-local}
 \sup_{|z|\le R}|1-B_j(z)|
   \le \frac{1+R}{1-R}\bigl(1-e^{-b(A_j)}\bigr)
   \le \frac{1+R}{1-R} b(A_j)\longrightarrow0 .
\end{equation}
Every zero in \(A_j\) has modulus at least \(e^{-b(A_j)}\), since
each radius is at least the product of the radii. Thus for each
\(R<1\), \(h_j=e_{A_j}/B_j\) on \(\{|z|\le R\}\) for all large \(j\).
The convergence \(e_{A_j}\to e\) in \(\D\) gives local uniform convergence
by the evaluation bound in Section~\ref{subsec:02-norm}.
Together with \eqref{eq:finite-factor-local}, it follows that
\(h_j\to e\) locally uniformly.

We will use the next lemma twice. It first upgrades local convergence
of the functions \(h_j\) to convergence in Dirichlet norm. Later it will
give the same upgrade for the comparison vectors in the atomic spaces.

\begin{lemma}[Strong convergence from a local limit and a norm bound]
\label{lem:local-norm-strong}
Suppose \(f_j\in\D\) converges locally uniformly to \(f\in\D\) and
\(\limsup_j\|f_j\|_{\D}\le\|f\|_{\D}\). Then \(f_j\to f\) in \(\D\).
\end{lemma}
\begin{proof}
Local uniform convergence gives convergence of every Taylor coefficient.
The coefficient norm therefore gives convergence of the inner products
against all polynomials. By boundedness of the norms and polynomial
density, \(f_j\) converges weakly to \(f\). Expanding
\(\|f_j-f\|_{\D}^2\) and taking the upper limit now gives zero.
\end{proof}

By \eqref{eq:finite-factor-data}, \(\|h_j\|_{\D}\le1=\|e\|_{\D}\),
so the lemma applies. Hence
\begin{equation}\label{eq:outer-quotients-strong}
 h_j\longrightarrow e\quad\text{in }\D,\qquad
 d_j:=1-\|h_j\|_{\D}^2\longrightarrow0,\qquad 0\le d_j<1 .
\end{equation}
The strict upper bound follows from \(h_j(0)>0\). \(d_j\)
is the contribution of the Blaschke product to the squared norm,
as the product formula below will show. In Section~\ref{subsec:07-competitor}, we use \(d_j\) to bound
the weighted boundary integral added to the squared Dirichlet norm
when we multiply by the singular inner factor. The sum \(b(A_j)\) bounds the total
mass of the new atoms.

\subsection{Multiplication by an inner factor}\label{subsec:07-product}

We need to test Dirichlet membership when we exchange an inner
factor for another. For a known Dirichlet product, the formula below
bounds the norm of the remaining factor and identifies the additional
integral contributed by the inner factor. Conversely, if that norm
and the integral for a new inner factor are finite, the new product
belongs to the Dirichlet space.
The inner--outer case appears in Richter and Sundberg's
factorization formula \cite[equation~(2.1), p.~171]{RS92}.
Luo and Richter also restate the broader local product identity
for arbitrary functions in $H^2$
\cite[p.~432, immediately before Theorem~4.7]{LR15},
and attribute it to Richter and Sundberg \cite[Lemma~3.4, pp.~365--366]{RS91}. We give a direct proof for an arbitrary $h\in H^2$ and a finite Blaschke or finite atomic singular factor. The integrated identity is the finite-factor case of Carleson's formula for the Dirichlet integral. The normalization and explicit weights also appear in \cite[(3.2)--(3.3)]{HR26},
which points to \cite[Theorems~7.1.3, 7.6.1 and 7.6.7]{Primer14}.

Let \(I\) be such an inner factor. Call \(\xi\in\T\) regular for \(I\)
when \(I\) extends analytically across a neighborhood of \(\xi\).
Every boundary point is regular for a finite Blaschke product, and
only the atomic directions are excluded for a finite atomic factor.
Continuity and the almost everywhere inner boundary modulus give
\(|I(\xi)|=1\) on each regular arc. At a point where \(f\in H^2\) has
a nontangential value \(f(\xi)\), put
\[
 q_{f,\xi}(z)=\frac{f(z)-f(\xi)}{z-\xi}.
\]

\begin{lemma}[The inner quotient at a regular boundary point]
\label{lem:regular-inner-quotient}
If \(\xi\) is regular for \(I\), then \(q_{I,\xi}\) lies in
\(K_I=H^2\ominus IH^2\) and
\begin{equation}\label{eq:regular-inner-weight}
 \|q_{I,\xi}\|_{H^2}^2=D_\xi(I)
   =\xi\overline{I(\xi)}I'(\xi)\ge 0.
\end{equation}
\end{lemma}
\begin{proof}
We use the evaluation kernel in \(K_I\) to recover the boundary
quotient and its squared norm. Convergence in norm retains orthogonality.
The Hardy evaluation kernel at \(\lambda\in\DD\) is
\(k_\lambda^{H^2}(z)=(1-\overline{\lambda}z)^{-1}\).
Write \(T_I\) for multiplication by \(I\) on \(H^2\).
It is an isometry, and its adjoint sends \(k_\lambda^{H^2}\) to
\(\overline{I(\lambda)}k_\lambda^{H^2}\). The projection
\(\operatorname{Id}_{H^2}-T_IT_I^*\) onto \(K_I\) therefore gives the kernel
\[
 k_\lambda^I(z)=
 \frac{1-\overline{I(\lambda)}I(z)}{1-\overline{\lambda}z},
 \qquad
 \|k_\lambda^I\|_{H^2}^2
   =\frac{1-|I(\lambda)|^2}{1-|\lambda|^2}.
\]
In particular, evaluation at \(\lambda\) on \(K_I\) is represented by
\(k_\lambda^I\).

Write \(d_I=\xi\overline{I(\xi)}I'(\xi)\). Differentiating the
constant boundary modulus along the arc gives \(d_I\in\mathbb R\).
Taylor expansion gives
\(I(r\xi)=I(\xi)(1-(1-r)d_I+O((1-r)^2))\). For \(r,s\uparrow1\),
\[
 \big\langle k_{r\xi}^I,k_{s\xi}^I\big\rangle_{H^2}
   =\frac{1-\overline{I(r\xi)}I(s\xi)}{1-rs}
   \longrightarrow d_I.
\]
Indeed, with \(u=1-r\) and \(v=1-s\), the numerator is
\(d_I(u+v)+O((u+v)^2)\) and the denominator is \(u+v-uv\).
The same limit with \(r=s\) gives \(d_I\ge0\), and the two-radius
limit gives \(\|k_{r\xi}^I-k_{s\xi}^I\|_{H^2}\to0\).
The norm limit belongs to \(K_I\). By point evaluation it is
\[
 k_\xi^I(z)=\frac{1-\overline{I(\xi)}I(z)}{1-\overline{\xi}z},
 \qquad \|k_\xi^I\|_{H^2}^2=d_I .
\]
Since \(q_{I,\xi}=(I(\xi)/\xi)k_\xi^I\), it has the same norm.
The boundary integral in \eqref{eq:local-douglas} is its squared
Hardy norm, proving \eqref{eq:regular-inner-weight}.
\end{proof}

\begin{proposition}[A Dirichlet norm identity for products with an inner factor]
\label{prop:finite-product-domain}
Let \(I\) be a finite Blaschke product or the singular inner function associated with a finite positive measure of finite support.
Put \(W_I(\xi)=D_\xi(I)\) at regular points, and set it to zero at
the finitely many remaining points. For every \(h\in H^2\),
\begin{equation}\label{eq:finite-product}
 \|Ih\|_{\D}^2
  =\|h\|_{\D}^2+\int_{\T}W_I|h^*|^2\,dm
\end{equation}
with nonnegative extended values. In particular,
\(Ih\in\D\) if and only if \(h\in\D\) and the displayed integral
is finite.
\end{proposition}
\begin{proof}
The product difference quotient splits into two terms. When both
terms belong to \(H^2\), their orthogonality gives the two
contributions to the energy.
For almost every \(\xi\), the point is regular for \(I\) and
\(h(\xi)=h^*(\xi)\) exists. The product \(Ih\) belongs to \(H^2\)
and has the value \(I(\xi)h(\xi)\) there. Its analytic quotient
satisfies
\[
 q_{Ih,\xi}=I q_{h,\xi}+h(\xi)q_{I,\xi}.
\]
If \(D_\xi(h)<\infty\), the local quotient result following
\eqref{eq:local-douglas} puts \(q_{h,\xi}\) in \(H^2\), with squared
norm \(D_\xi(h)\). The two terms in the decomposition lie in
\(IH^2\) and \(K_I\), respectively. Orthogonality and the Hardy
isometry give
\begin{equation}\label{eq:finite-local-product}
 D_\xi(Ih)=D_\xi(h)+|h(\xi)|^2D_\xi(I).
\end{equation}
If \(D_\xi(h)=\infty\), the boundary difference quotient of \(h\)
is outside \(L^2(m)\) by \eqref{eq:local-douglas}. Multiplying
that quotient by \(I^*\) preserves its modulus almost everywhere,
and the boundary function of \(h(\xi)q_{I,\xi}\) is in \(L^2(m)\).
The sum is thus outside \(L^2(m)\), so \(D_\xi(Ih)=\infty\)
by the same local formula. This proves \eqref{eq:finite-local-product}
also in this case.

Integrate \eqref{eq:finite-local-product} using
\eqref{eq:weighted-local-douglas} with \(v=1\), which applies to both
Hardy functions before either Dirichlet energy is known to be finite.
Adding the equal Hardy norms gives \eqref{eq:finite-product}
with all terms nonnegative. This also proves the stated criterion
for \(Ih\in\D\).
\end{proof}

For \(B=\prod_k b_k\) with \(b_k(z)=(r_k-\overline{\xi_k}z)/
(1-r_k\overline{\xi_k}z)\), \(0\le r_k<1\) and \(\xi_k\in\T\),
and for \(S_\mu\) with
\(\mu=\sum_k\tau_k\delta_{\xi_k}\), the weights are
\begin{equation}\label{eq:finite-product-weights}
 W_B(\zeta)=\sum_k\frac{1-r_k^2}{|\zeta-r_k\xi_k|^2},
 \qquad
 W_{S_\mu}(\zeta)=\sum_k\frac{2\tau_k}{|\zeta-\xi_k|^2}
 \quad\text{a.e. on }\T.
\end{equation}
At regular points \eqref{eq:regular-inner-weight} is
\(W_I(\zeta)=\zeta I'(\zeta)/I(\zeta)\), so the weights add under
finite products. Direct differentiation gives
\(\zeta b_k'(\zeta)/b_k(\zeta)=(1-r_k^2)/|\zeta-r_k\xi_k|^2\).
For the singular factor of mass \(\tau_k\) at \(\xi_k\), the
logarithmic derivative multiplied by \(\zeta\) is
\[
 -\frac{2\tau_k\zeta\xi_k}{(\xi_k-\zeta)^2}
   =\frac{2\tau_k}{|\zeta-\xi_k|^2}
   \qquad(\zeta\ne\xi_k).
\]
The equality uses
\((\xi_k-\zeta)^2=-\xi_k\zeta|\xi_k-\zeta|^2\) on the circle.
A zero at the origin contributes the constant weight \(1\), with its direction
arbitrary. Repeated atomic directions add their masses. The
singular weight for an atom has a squared reciprocal distance at
its boundary point. A cutoff with a simple zero at that point cancels
this singularity in the weighted integral.

\subsection{A trial function in the atomic subspace}\label{subsec:07-competitor}

The map from \(r_k\xi_k\) to an atom of mass \(1-r_k\), the
cutoffs below, and their boundary cancellation occur in Hartz and
Richter's factor-existence proof \cite[Appendix~B, Theorem~B.1]{HR26}.
They also point out that the parallelogram law can sharpen their norm
bound. We apply it at every replacement in the finite set to bound
the decrease in squared norm from below. The same integrals then
control the energy added by the singular factor.

Fix a nonempty finite set \(A=\{r_k\xi_k:1\le k\le N\}\) as in
Section~\ref{subsec:07-quotients}. Write \(B=B_A=\prod_k b_k\)
and \(e_A=Bh\), with \(b_k\) normalized as there. Put
\[
 b=b(A),\qquad d=1-\|h\|_{\D}^2,\qquad
 \gamma_B=e_A(0),\qquad \beta=B(0)=e^{-b}.
\]
The factorization \eqref{eq:finite-factor-data} gives
\(h\in\D\) and \(h(0)=\gamma_B/\beta>0\). Set
\[
 \mu=\sum_k(1-r_k)\delta_{\xi_k},\qquad S=S_\mu,\qquad
 s=S(0)=e^{-\mu(\T)},\qquad r_{\max}=\max_k r_k.
\]
Repeated directions add their masses in \(\mu\). For a set on one
circle, distinct zeros give distinct directions. We use a cutoff
that vanishes at \(\xi_k\) and is the midpoint of \(1\) and \(b_k\).
For the cutoffs put
\[
 q_k=\frac{1+b_k}{2},\qquad Q=\prod_k q_k,\qquad
 q_0=Q(0)=\prod_k\frac{1+r_k}{2}.
\]
Each \(b_k\) and \(q_k\) is rational and analytic through the closed
disk. They are Dirichlet multipliers by
Lemma~\ref{lem:bounded-products}. Their moduli on the closed disk
are at most one. In particular \(Qh\) and all intermediate products
of \(h\) with these factors belong to \(\D\). The product formula
\eqref{eq:finite-product} applied to \(Bh=e_A\) gives
\begin{equation}\label{eq:blaschke-row-charge}
 \int_{\T} W_B |h^*|^2\,dm=1-\|h\|_{\D}^2=d.
\end{equation}

The same stepwise weighted integrals give a lower bound for
the decrease in the squared norm and an upper bound for the
weighted boundary integral associated with the singular factor.
Keeping the earlier cutoffs in both calculations makes
these two bounds comparable. The sum of these integrals is bounded
by \(1-\|h\|_{\D}^{2}\).

Replace the Blaschke factors by the cutoffs in order, keeping the
earlier cutoffs in each new comparison. At the \(k\)-th step let
\[
 U_k=\Bigl(\prod_{j<k}q_j\Bigr)
     \Bigl(\prod_{j>k}b_j\Bigr)h,\qquad
 F_{k-1}=b_kU_k,\qquad F_k=q_kU_k.
\]
Thus \(F_0=Bh\) and \(F_N=Qh\). Since the remaining Blaschke factors
have boundary modulus one, \eqref{eq:finite-product} gives
\begin{equation}\label{eq:cutoff-charge}
 \Delta_k=\|b_kU_k\|_{\D}^2-\|U_k\|_{\D}^2
   =\int_{\T}W_{b_k}\prod_{j<k}|q_j|^2\,|h^*|^2\,dm,\qquad
 \Delta=\sum_k\Delta_k.
\end{equation}
The weights are nonnegative and \(W_B=\sum_k W_{b_k}\), so
\eqref{eq:blaschke-row-charge} gives \(0\le\Delta\le d\).
The midpoint identity expresses the decrease in the squared norm
at this step in terms of \(\Delta_k\). For \(F_k=(U_k+b_kU_k)/2\) it yields
\[
 \|F_k\|_{\D}^2=\|F_{k-1}\|_{\D}^2
   -\frac{\Delta_k}{2}-\frac14\|(1-b_k)U_k\|_{\D}^2.
\]
Summing over \(k\) gives \(\|Qh\|_{\D}^2\le1-\Delta/2\).

To bound the energy added by \(S\), we compare its atomic weights
after inserting the cutoff at each atom. Use the weights
\eqref{eq:finite-product-weights}. At boundary points
\(\zeta\ne\xi_k\), the formula
\[
 q_k(\zeta)=\frac{(1+r_k)(1-\overline{\xi_k}\zeta)}
                   {2(1-r_k\overline{\xi_k}\zeta)}
\]
gives
\[
 \frac{2(1-r_k)}{|\zeta-\xi_k|^2}|q_k(\zeta)|^2
       =\frac{1+r_k}{2}\,W_{b_k}(\zeta).
\]
To estimate the contribution of the \(k\)-th atom to the weighted
boundary integral, keep the cutoff factors with \(j<k\) and
bound the factors with \(j>k\) by one. The resulting
integrals are precisely those in \eqref{eq:cutoff-charge}:
\begin{align*}
 \int_{\T}W_S |(Qh)^*|^2\,dm
 &=\sum_k\frac{1+r_k}{2}\int_{\T}
             W_{b_k}\prod_{j\ne k}|q_j|^2\,|h^*|^2\,dm\\
 &\le\frac{1+r_{\max}}{2}\Delta .
\end{align*}
These are arclength integrals, so the finitely many atomic directions
do not affect the calculation.

We have \(Qh\in\D\), and the weighted boundary integral
\(\int_{\T}W_S|(Qh)^*|^2\,dm\) is finite. Applying
\eqref{eq:finite-product} to \(S(Qh)\) gives \(x=SQh\in\D\). Since \(Qh\in H^2\), this is a vector in
\(M_\mu=SH^2\cap\D\). Adding this integral to the estimate
for the squared norm of \(Qh\), and using the midpoint identity,
gives
\begin{equation}\label{eq:atomic-trial-bounds}
 \|x\|_{\D}^2\le 1+\frac{r_{\max}}2\Delta
       \le 1+\frac{r_{\max}}2d,\qquad
 x(0)=\frac{s q_0}{\beta}\gamma_B>0 .
\end{equation}
The next comparison uses the norm and the value at zero of \(x\)
to estimate its distance to the canonical function of \(M_\mu\).

\subsection{Distance from the trial function to the canonical function}\label{subsec:07-comparison}

The trial function \(x\) came from the finite-zero canonical function.
We now start with the atomic canonical function, divide out its
singular factor, and multiply the resulting quotient by the finite
Blaschke product. The extremal property in \(M_A\) will then give
an upper bound for the atomic canonical value at zero.
The reverse weight comparison and the construction of a Blaschke
trial vector from this quotient also occur in
\cite[Theorem~3.5, Corollary~3.6, and Lemma~5.1]{HR26}. We use
them with the very same finite set and measure as in the construction
of \(x\).

Set \(r_{\min}=\min_k r_k\) and \(c_{\mathrm{tr}}=2/(1+r_{\min})\).
For a boundary point \(\zeta\ne\xi_k\), put \(u=|\zeta-\xi_k|^2\le4\).
The explicit weights in \eqref{eq:finite-product-weights} and
\(|\zeta-r_k\xi_k|^2=(1-r_k)^2+r_ku\) give
\[
 \frac{W_{b_k}(\zeta)}{2(1-r_k)/|\zeta-\xi_k|^2}
 =\frac{(1+r_k)u}{2((1-r_k)^2+r_ku)}
 \le \frac{2}{1+r_k}.
\]
After clearing denominators, the last inequality reduces to
\((1-r_k)^2u\le4(1-r_k)^2\). Summing over the factors, with
repeated directions counted in the atomic mass, yields
\begin{equation}\label{eq:atomic-reverse-charge}
                  W_B\le c_{\mathrm{tr}} W_S
                  \qquad\text{almost everywhere}.
\end{equation}

Write \(\gamma_S=G_\mu(0)>0\) and \(G_\mu=Sv\). By definition of
\(M_\mu\), we have \(v\in H^2\).
Apply \eqref{eq:finite-product} to this known Dirichlet product:
\[
            \|v\|_\D^2+\int_\T W_S|v^*|^2\,dm=1 .
\]
The preceding identity gives \(v\in\D\), and
\eqref{eq:atomic-reverse-charge} implies
\(\int_{\T}W_B|v^*|^2\,dm<\infty\). Applying
\eqref{eq:finite-product} to \(Bv\) gives \(Bv\in\D\), hence
\(Bv\in M_A\). The same formula gives
\[
 \|Bv\|_\D^2\le \|v\|_\D^2+c_{\mathrm{tr}}\int_\T W_S|v^*|^2\,dm
      = c_{\mathrm{tr}}-(c_{\mathrm{tr}}-1)\|v\|_\D^2
      \le c_{\mathrm{tr}}.
\]
Its value at zero is \(\beta\gamma_S/s>0\). We can now apply the canonical
extremal property of \(e_A\) from Section~\ref{subsec:02-projection}:
\begin{equation}\label{eq:atomic-upper-origin}
 \frac{\beta\gamma_S}{s}\le \gamma_B\|Bv\|_\D,
 \qquad
 \gamma_S\le \sqrt{c_{\mathrm{tr}}}\,\frac{s}{\beta}\gamma_B .
\end{equation}

Section~\ref{sec:06} gives one common radius within each set \(A_j\),
and that radius can vary with \(j\). The scalar bounds below use only
\(0<r_k<1\) and the quantity \(b\), which bounds each radius from
below and controls the total mass. The inequalities \(1-r\le-\log r\) and
\((1+r)/2\ge\sqrt r\), for \(0<r<1\), give
\begin{equation}\label{eq:atomic-budget-scalars}
 \begin{gathered}
 0<\mu(\T)\le b,\qquad \beta=e^{-b}\le s\le1,\qquad
 e^{-b/2}\le q_0\le1,\qquad r_{\min}\ge e^{-b},\\
 1\le c_{\mathrm{tr}}\le r_{\min}^{-1/2}\le e^{b/2},
 \qquad q_0/\sqrt{c_{\mathrm{tr}}}\ge e^{-3b/4}.
 \end{gathered}
\end{equation}
Here \(r_{\min}\ge\prod_k r_k\) since all the radii are at most one.
These bounds hold for each finite set, regardless of its cardinality.

The canonical projection identity now converts the upper bound at
zero into a lower bound for the inner product of \(x\) with \(G_\mu\).
For \(x\in M_\mu\), that identity gives
\(\langle x,G_\mu\rangle_{\D}=x(0)/\gamma_S\). The value \(x(0)\)
in \eqref{eq:atomic-trial-bounds} and the upper bound
\eqref{eq:atomic-upper-origin} yield
\[
 \frac{x(0)}{\gamma_S}\ge \frac{q_0}{\sqrt{c_{\mathrm{tr}}}}
                         \ge e^{-3b/4}.
\]
Both \(x(0)\) and \(\gamma_S\) are positive. Using the norm estimate in
\eqref{eq:atomic-trial-bounds}, we get
\begin{align}
 \|x-G_\mu\|_\D^2
 &=\|x\|_\D^2+1-2x(0)/\gamma_S \notag\\
 &\le \frac{r_{\max}d}{2}+2(1-e^{-3b/4})
 \le \frac d2+\frac{3b}{2} . \label{eq:atomic-transfer-distance}
\end{align}
The coefficients \(1/2\) and \(3/2\) are independent of the number,
directions, and radii of the zeros. Along the sequence for the fixed
target \(e\), the quantity \(1-\|h\|_{\D}^{2}\)
and \(b(A_j)\) both tend to zero.
One can see exactly how much the cutoff estimate contributes.
The intermediate estimate in \cite[Appendix~B, proof of Theorem~B.1]{HR26}
bounds the squared trial norm by \(\|Bh\|_{\D}^{2}+
\int_{\T}W_B|h^*|^2\,dm=1+d\) in our notation.
With this bound the preceding calculation gives \(d+3b/2\);
using the decrease in the squared norm supplied by the midpoint
identity improves this to \(d/2+3b/2\).
Both calculations require the estimate for the normalized extremal
of the subspace \(M_\mu\) containing \(x\): its value
\(\gamma_S\) is the denominator in the projection identity.
Thus a controlled trial vector gives the asserted canonical
approximation only after the comparison in this subsection.

\subsection{Completion of Approximation A}\label{subsec:07-completion}

The functions \(h_j\) converge to \(e\), and the fact that \(b(A_j)\to0\) forces both replacement factors to approach \(1\) on compact subdisks. We combine the resulting local convergence with the norm
bound for \(x_j\), then pass to the canonical function of the
associated atomic space.

\begin{proof}[Proof of Theorem~\ref{thm:A}]
Fix the target \(e\in\Xpr\) and the two positive parameters
\(\varepsilon_{\rm norm},\varepsilon_{\rm mass}\).
Use the sequence \(A_j\) selected for this \(e\) in
Section~\ref{subsec:07-quotients}. Write \(e_{A_j}=B_jh_j\)
as in Section~\ref{subsec:07-quotients}. Then \(b_j=b(A_j)\to0\),
\(h_j\to e\) strongly in \(\D\), and
\(d_j=1-\|h_j\|_\D^2\to0\).
For each \(j\), let \(\mu_j\), \(S_j=S_{\mu_j}\), and \(Q_j\)
be the measure and factors of Section~\ref{subsec:07-competitor}.
Equations \eqref{eq:atomic-trial-bounds} and
\eqref{eq:atomic-budget-scalars} give
\[
  0<\mu_j(\T)\le b_j,\qquad S_j(0)\to1,\qquad Q_j(0)\to1 .
\]
Both \(S_j\) and \(Q_j\) are bounded by one in modulus.
The Schwarz--Pick estimate \eqref{eq:schur-origin-bound} gives
\(S_j,Q_j\to1\) locally uniformly. Thus \(x_j=S_jQ_jh_j\to e\)
locally uniformly, and \eqref{eq:atomic-trial-bounds} gives
\(\limsup_j\|x_j\|_\D^2\le1=\|e\|_\D^2\).
Lemma~\ref{lem:local-norm-strong} now gives \(x_j\to e\) in \(\D\).

We now pass to the canonical functions of the same measures. By
\eqref{eq:atomic-transfer-distance},
\begin{equation}\label{eq:atomic-transfer-convergence}
 \|G_{\mu_j}-e\|_\D
 \le \|x_j-e\|_\D+\sqrt{d_j/2+3b_j/2}\longrightarrow0 .
\end{equation}
Because \(0<\mu_j(\T)\le b_j\to0\) on the same sequence, choose a
single index \(j\) with \(\|G_{\mu_j}-e\|_\D<\varepsilon_{\rm norm}\)
and \(b_j<\varepsilon_{\rm mass}\). Set \(\mu=\mu_j\).
For this measure, the atomic target construction in
Section~\ref{subsec:02-atomic} gives \([G_\mu]=M_\mu\) and the exact
factorization \(G_\mu=S_\mu O_\mu\), with \(O_\mu\) outer and
\(O_\mu(0)>0\).
\end{proof}

Approximation~B starts from this fixed measure \(\mu\). Its norm
limit is the same \(G_\mu\) whose distance to the starting outer
function is controlled here.

\newpage \section{Boundary estimates for canonical functions of trace-zero spaces}\label{sec:08}

Approximation~B begins with a fixed atomic measure. We first
consider one atom at \(\xi\) and use canonical functions of
trace-zero spaces on nonpolar Carleson sets that shrink to \(\{\xi\}\).
Theorem~\ref{thm:C} proves convergence in Dirichlet norm when the
canonical values at zero converge to the target value.

We need to retain the Dirichlet energy and identify the mass of
the limit. We first show that the canonical value at zero
distinguishes different masses at \(\xi\). For norm convergence,
the estimates below serve two purposes. A kernel test bounds the
energy in a small disk about \(\xi\). On arcs away from \(\xi\),
we will use uniform analytic continuation to control convergence
of the derivatives.

\subsection{The one-atom target and its canonical value}\label{subsec:08-target}

We will compare the limit with the canonical function of the
one-atom space. Its value at zero is the largest possible modulus
at zero among vectors in the unit ball of that space.

Fix \(\xi\in\T\). For \(t>0\) write
\[
 S_{t,\xi}(z)=\exp\!\left(-t\frac{\xi+z}{\xi-z}\right),
 \qquad M_{t,\xi}=S_{t,\xi}H^2\cap\D,\qquad p_{t,\xi}=P_{M_{t,\xi}}1.
\]
These are the one-atom instances of \eqref{eq:atomic-target}.
The subspace \(M_{t,\xi}\) is closed because \(S_{t,\xi}H^2\) is
closed in \(H^2\) and the inclusion \(\D\hookrightarrow H^2\) is
continuous. It is invariant under multiplication by \(z\).
The function \((1-\overline{\xi}z)S_{t,\xi}\) in
Lemma~\ref{lem:atomic-witness} belongs to \(M_{t,\xi}\), has
Dirichlet norm \(\sqrt{3+2t}\), and takes the value \(e^{-t}\) at
zero. Thus \(p_{t,\xi}\ne0\). Put
\[
 G_{t,\xi}=\frac{p_{t,\xi}}{\|p_{t,\xi}\|_{\D}},\qquad
 a_{t,\xi}=G_{t,\xi}(0)=\|p_{t,\xi}\|_{\D}.
\]
By the extremal property in Section~\ref{subsec:02-projection},
the function \((1-\overline{\xi}z)S_{t,\xi}/\sqrt{3+2t}\)
gives the lower bound in
\begin{equation}\label{eq:one-atom-origin-bounds}
 0<\frac{e^{-t}}{\sqrt{3+2t}}\le a_{t,\xi}\le e^{-t}<1.
\end{equation}
For the upper bound write \(G_{t,\xi}=S_{t,\xi}v\) with \(v\in H^2\).
Since \(S_{t,\xi}\) is inner,
\(|v(0)|\le\|v\|_{H^2}=\|G_{t,\xi}\|_{H^2}\le1\).
Now \(a_{t,\xi}=e^{-t}|v(0)|\). Canonical generation gives
\([G_{t,\xi}]=M_{t,\xi}\). The rotation in
Section~\ref{subsec:02-atomic} shows that \(a_{t,\xi}\) is independent
of \(\xi\).

Increasing the mass makes the spaces smaller. To obtain strict
decrease of the canonical value, we use strict inclusion and the
fact that the normalized projection generates its space.

\begin{proposition}[Strict comparison of one-atom extremal values at zero]
\label{prop:atomic-origin-strict}
For \(0<s<t\) and \(\xi\in\T\),
\[
 M_{t,\xi}\subsetneq M_{s,\xi},\qquad
 a_{s,\xi}^{\,2}-a_{t,\xi}^{\,2}
       =\|p_{s,\xi}-p_{t,\xi}\|_{\D}^{2}>0 .
\]
In particular \(G_{t,\xi}(0)\) determines \(t\) for a fixed \(\xi\).
\end{proposition}

\begin{proof}
The factorization \(S_{t,\xi}=S_{s,\xi}S_{t-s,\xi}\) gives inclusion.
The function \((1-\overline{\xi}z)S_{s,\xi}\) belongs to \(M_{s,\xi}\).
If it belonged to \(M_{t,\xi}\),
its quotient \(F=(1-\overline{\xi}z)/S_{t-s,\xi}\) would belong to \(H^2\).
But for \(0<r<1\),
\[
 |F(r\xi)|=(1-r)\exp\!\left((t-s)\frac{1+r}{1-r}\right).
\]
The Hardy evaluation bound would give
\(|F(r\xi)|\sqrt{1-r^2}\le\|F\|_{H^2}\), whereas the displayed
expression multiplied by \(\sqrt{1-r^2}\) tends to infinity.
This proves strict inclusion.

We now compare the projections of \(1\) onto the two nested spaces.
Since \(1-p_{s,\xi}\perp M_{s,\xi}\) and
\(M_{t,\xi}\subset M_{s,\xi}\), we have
\(P_{M_{t,\xi}}p_{s,\xi}=p_{t,\xi}\). Hence
\(p_{s,\xi}-p_{t,\xi}\perp p_{t,\xi}\). Pythagoras gives
\(a_{s,\xi}^{\,2}-a_{t,\xi}^{\,2}=\|p_{s,\xi}-p_{t,\xi}\|_{\D}^{2}\).
If this difference were zero, the projections would coincide. Their nonzero norms and positive
normalizations would then give \(G_{s,\xi}=G_{t,\xi}\).
Canonical generation would force
\(M_{s,\xi}=[G_{s,\xi}]=[G_{t,\xi}]=M_{t,\xi}\), contrary to
strict inclusion.
\end{proof}

In their preprint, Hartz and Richter use the same function
\((1-z)S_{t,1}\) and compute the extremal for \(M_{t,1}\) in
the Dirichlet norm used here \cite[Section~5 and Theorem~5.2]{HR26}.
Their formula is
\[
 (zG_{t,1})'(z)
   =A_t^{-1/2}\frac{A_t-z}{1-z}\,S_{t,1}(z),\qquad
 A_t=\frac{\int_0^1 x\,u_t(x)\,dx}{\int_0^1 u_t(x)\,dx},
 \qquad u_t(x)=\frac{S_{t,1}(x)}{1-x}.
\]
In particular \(G_{t,1}(0)=A_t^{1/2}e^{-t}\).
Their proof uses the invariant subspaces associated with Blaschke products with an interior
zero at \(1-t/n\) of multiplicity \(n\). The norm convergence in
that proof uses convergence of the inner factors and domination
of their local Dirichlet integrals \cite[Lemma~5.1]{HR26}.
Theorem~\ref{thm:C} starts with trace-zero spaces on shrinking boundary
compact sets. In its proof, Proposition~\ref{prop:atomic-origin-strict}
identifies the mass from the value at zero of the limiting extremal function.

\subsection{Energy estimates and boundary multipliers}\label{subsec:08-source}

The identity \eqref{eq:canonical-source-identity} gives a bound on the potential of the measure \(\sigma_e\). We use it to integrate squared traces and then
to bound multiplication on the boundary test space \(V\). This
bound will let us multiply a distribution by the boundary trace
of \(e\) when constructing the analytic function \(\beta\) in the next subsection.
We also prove a separate energy estimate near the collapse point
using the analytic multiplier bound and a reproducing kernel.

Throughout this subsection let \(e\in\D\) have norm one and satisfy
\(\langle z^ne,e\rangle_{\D}=0\) for all \(n\ge1\). Set
\(w=|e^*|^2\) and use the measure \(\sigma_e=\ell_e\,m\) from
\eqref{eq:source}, where \(m\) is normalized arclength. The coefficient
calculation for \eqref{eq:canonical-source-identity} uses just the norm
and moments. Thus
\[
  (1+\Lambda)w+\ell_e=1,\qquad
  0\le w\le1,\qquad
  \sigma_e(\T)=D(e)\le1 .
\]
Here \(\Lambda\) has Fourier multiplier \(|n|\). This is the
boundary Fourier form of the extremal identity of Richter and
Sundberg \cite[(5.1)]{RS94}, reproduced in
\cite[(4.1)]{RY19}. It also follows from the Cauchy-transform
identity in \cite[Lemma~2.4, (2.1)]{GZ10}.
The measure \(\sigma_e\) is absolutely continuous with respect to
\(m\), so its integrals use ordinary almost everywhere boundary values.

The measure \(\sigma_e\) records the energy of the present canonical function; the later boundary-logarithm measures serve
the separate purpose of detecting singular mass. By absolute continuity, the integral
\(\int_{\T}|F|^2\,d\sigma_e\) depends only on the \(L^2(m)\) class
of \(F\). The estimate below bounds this integral by \(\|F\|_V^2\)
for every \(F\in V\). That estimate is
then used to prove bounded multiplication on the boundary
test space, including negative Fourier modes. Transposing this
multiplication lets us act on the conjugate derivative
distribution in the next subsection. The small-disk estimate
controls derivative energy near the collapse point uniformly,
complementing the continuation argument on arcs away from that
point.

\paragraph{The potential of the  measure \(\sigma_e\).}
Applying the inverse of \(1+\Lambda\) to the identity \eqref{eq:canonical-source-identity}
expresses the potential in terms of the squared modulus \(w\).
Write \(R=(1+\Lambda)^{-1}\) and use the nonnegative kernel
\(K_R=\int_0^1P_r\,dr\) from \eqref{eq:reciprocal-energy-kernel}.
For a finite positive measure \(\sigma\), define the potential by 
\[
R\sigma(\zeta)=\int_{\T}K_R(\zeta\bar\omega)\,d\sigma(\omega).
\]
The kernel has mass one, so the potential is finite almost everywhere
for \(m\) and belongs to \(L^1(m)\). Its Fourier coefficient at \(n\)
is \(\widehat\sigma(n)/(1+|n|)\), by Fubini and the Poisson Fourier
series. For \(\sigma=\sigma_e\), the identity \eqref{eq:canonical-source-identity} therefore gives
\begin{equation}\label{eq:source-potential}
 R\sigma_e=1-w\in[0,1]\qquad m\text{-a.e. and }\sigma_e\text{-a.e.}
\end{equation}
The second assertion uses \(\sigma_e\ll m\).
In \cite[(4.2)--(4.4)]{RY19}, the local energy \(D_\zeta(e)\)
is the density \(\ell_e(\zeta)\) by the local Douglas formula.
Their kernel \(r^{-1}\int_0^r P_s\,ds\) has Fourier coefficients
\(r^{|n|}/(1+|n|)\). Thus their identity is the Poisson
regularization of \eqref{eq:source-potential}. The positive
kernel formula \eqref{eq:reciprocal-energy-kernel} gives
\begin{equation}\label{eq:source-reciprocal-energy}
 I_R(\sigma_e)=\int_{\T}R\sigma_e\,d\sigma_e
       \le \sigma_e(\T)\le1 .
\end{equation}
For a logarithmic-potential bound for the same measure, see
\cite[Proposition~4.5]{RY19}. The Fourier coefficient form of \eqref{eq:source-potential} controls
the boundary norm of \(w\). In fact,
\(\widehat w(n)=\delta_{n0}-\widehat\sigma_e(n)/(1+|n|)\), where
\(\delta_{n0}\) is \(1\) for \(n=0\) and \(0\) otherwise.
Keeping the constant coefficient in the norm \eqref{eq:V} gives
\begin{equation}\label{eq:source-modulus-V}
 w\in V,\qquad
 \|w\|_V^2=1-2\sigma_e(\T)+I_R(\sigma_e)\le1 .
\end{equation}

\paragraph{Squared trace bound and multiplication.}
We need an estimate for the squared modulus of a test function
integrated against \(\sigma_e\). To obtain it, we first bound the
dual norm of the measures obtained by multiplying \(\sigma_e\) by
square-integrable functions.
Use the space \(V=H^{1/2}(\T)\) with norm \eqref{eq:V}, and
let \(V^*=H^{-1/2}(\T)\) denote its continuous dual. We use bilinear
distribution evaluation. For \(T\in V^*\), let
\(\widehat T(n)=T[\zeta\mapsto\bar\zeta^{\,n}]\). Then
\begin{equation}\label{eq:boundary-dual-norm}
 \|T\|_{V^*}^2=\sum_{n\in\mathbb Z}\frac{|\widehat T(n)|^2}{1+|n|}
\end{equation}
and
\begin{equation}\label{eq:boundary-dual-pairing}
 T[\phi]=\sum_{n\in\mathbb Z}\widehat T(n)\widehat\phi(-n)
 \qquad(\phi\in V).
\end{equation}
The pairing first holds for trigonometric polynomials and extends
by weighted Cauchy--Schwarz. Taking the supremum over polynomial tests
gives \eqref{eq:boundary-dual-norm}. Conjugation acts by
\(\bar T[\phi]=\overline{T[\bar\phi]}\), hence
\(\widehat{\bar T}(n)=\overline{\widehat T(-n)}\).

\begin{lemma}[Boundary multipliers]\label{lem:full-boundary-multipliers}
Let \(e\in\D\) have norm one and satisfy
\(\langle z^ne,e\rangle_\D=0\) for all \(n\ge1\).
Every \(F\in V\) satisfies
\begin{equation}\label{eq:source-square-trace}
 \int_{\T}|F|^2\,d\sigma_e\le\|F\|_V^2.
\end{equation}
Multiplication by either \(\psi=e^*\) or \(\psi=\overline{e^*}\)
is defined on all of \(V\) and satisfies
\begin{equation}\label{eq:full-boundary-multiplier-bound}
 \|\psi F\|_V\le2\|F\|_V \qquad(F\in V).
\end{equation}
\end{lemma}
\begin{proof}
We first control \(I_R(g\,\sigma_e)\).
Duality will then give the trace bound needed for multiplication.
For \(g\in L^2(\sigma_e)\), the complex measure \(g\,\sigma_e\)
is finite because \(\sigma_e\) is finite. The inequality
\(2|g(\zeta)g(\omega)|\le |g(\zeta)|^2+|g(\omega)|^2\) and
symmetry of \(K_R\) give the Schur estimate
\[
 \iint_{\T^2}K_R(\zeta\bar\omega)|g(\zeta)g(\omega)|
                           \,d\sigma_e(\zeta)d\sigma_e(\omega)
 \le\int_\T |g|^2 R\sigma_e\,d\sigma_e
 \le\|g\|_{L^2(\sigma_e)}^2 .
\]
At fixed \(r<1\), the Poisson Fourier series is absolutely convergent.
Integrating its terms against \(g(\zeta)\overline{g(\omega)}\)
gives \(\sum_n r^{|n|}|\widehat{g\,\sigma_e}(n)|^2\), using
symmetry of \(r^{|n|}\) in \(n\). The absolute bound above permits
integration in \(r\) on the kernel side, and Tonelli applies to
the nonnegative series. This gives
\[
 I_R(g\,\sigma_e)
 =\iint_{\T^2}K_R(\zeta\bar\omega)g(\zeta)\overline{g(\omega)}
                        \,d\sigma_e(\zeta)d\sigma_e(\omega)
 \le \|g\|_{L^2(\sigma_e)}^2 .
\]
The same bound applies to \(\bar g\,\sigma_e\).

We now turn this energy estimate into the squared trace bound,
first for polynomial tests and then for all boundary tests.
For a trigonometric polynomial \(F\), apply the dual pairing to
\(\bar g\,\sigma_e\) to obtain
\[
 \left|\int_\T F\bar g\,d\sigma_e\right|
   \le\|F\|_V I_R(\bar g\,\sigma_e)^{1/2}
   \le\|F\|_V\|g\|_{L^2(\sigma_e)}.
\]
Taking the supremum over \(\|g\|_{L^2(\sigma_e)}\le1\) proves
\eqref{eq:source-square-trace} for polynomials. If \(F_j\to F\) in
\(V\) through trigonometric polynomials, the inequality for differences
makes \(F_j\) Cauchy in \(L^2(\sigma_e)\). A subsequence converges to
\(F\) almost everywhere for \(m\), and hence for \(\sigma_e\).
A further subsequence identifies the \(L^2(\sigma_e)\) limit with
this representative. This proves \eqref{eq:source-square-trace}
on \(V\).

It remains to bound multiplication in the norm of \(V\).
The product difference has one term from the change in \(F\)
and one from the change in \(\psi\). The squared trace bound
controls the second term after we integrate the difference norm.
Let \(\mathcal E(F)\) denote the double integral in
\eqref{eq:V-difference}. The product \(\psi F\) belongs to \(L^2(m)\)
because \(|\psi|\le1\). For either choice of \(\psi\),
\[
 \begin{split}
 |\psi(\zeta)F(\zeta)-\psi(\omega)F(\omega)|^2
 &\le 2|F(\zeta)-F(\omega)|^2\\
 &\quad+2|F(\omega)|^2|e^*(\zeta)-e^*(\omega)|^2 .
 \end{split}
\]
By \eqref{eq:local-douglas} and \eqref{eq:source}, for \(m\)-almost
every \(\omega\in\T\) we have
\[
 \int_{\T}
 \frac{|e^*(\zeta)-e^*(\omega)|^2}{|\zeta-\omega|^2}\,dm(\zeta)
 = \ell_e(\omega).
\]
Thus, after division by \(|\zeta-\omega|^2\), the second term
contributes \(2\int_\T |F|^2\,d\sigma_e\) to the bound for
\(\mathcal E(\psi F)\). The difference formula for the norm of \(V\)
and \eqref{eq:source-square-trace} now give
\[
 \|\psi F\|_V^2\le \|F\|_2^2+2\mathcal E(F)+2\int_\T |F|^2\,d\sigma_e
             \le 3\|F\|_2^2+4\mathcal E(F)\le4\|F\|_V^2 .
\]
This controls the norm, including its constant mode.
Thus the ordinary product belongs to \(V\), with the stated bound.
\end{proof}

For \(T\in V^*\), Lemma~\ref{lem:full-boundary-multipliers}
defines the product by \((\psi T)[\phi]=T[\psi\phi]\) on \(V\).
The product belongs to \(V^*\) with norm at most \(2\|T\|_{V^*}\).
The test space includes both signs of Fourier mode, as needed to determine all Fourier coefficients of the product.

\paragraph{Uniform energy near the circle.}
Section~\ref{sec:09} needs a bound on the derivative energy in
a small disk, uniform over the canonical functions. We use the
analytic multiplier bound from Section~\ref{subsec:02-multiplier}
to control \(|e'|^2\,dA_0\) tested against a Dirichlet reproducing
kernel. On the small disk the kernel has size at least a constant
times \(\log(1/\delta)\), while its squared Dirichlet norm grows
at most on that logarithmic scale. This gives a bound tending to
zero with the radius.
This is a local version of the necessary single-box estimate for
Dirichlet Carleson measures
\cite[equation~(13) on p.~53 and Proposition~22 on p.~69]{ARSW11}.

\begin{proposition}[Uniform derivative energy near the circle]
\label{prop:small-region-energy}
Let \(e\in\D\) have norm one and satisfy
\(\langle z^ne,e\rangle_\D=0\) for all \(n\ge1\).
For every \(\xi\in\T\) and \(0<\delta\le 1/256\), with
\(B(\xi,\delta)=\{z:|z-\xi|<\delta\}\),
\begin{equation}\label{eq:small-region-energy}
 \int_{\DD\cap B(\xi,\delta)}|e'(z)|^2\,dA_0(z)
       \le\frac{4\pi}{\log(1/\delta)} .
\end{equation}
The bound is uniform over these \(e\) and all \(\xi\).
\end{proposition}
\begin{proof}
We bound multiplication by \(e'\) from \(\D\) to area
\(L^2\) using the product rule. A kernel test will convert this
multiplier bound into control of the derivative energy on a small
boundary disk.
The contractive multiplier bound \eqref{eq:canonical-multiplier}
and \(|e|\le1\) give
\[
 \|e'h\|_{2,A}
 \le \|(eh)'\|_{2,A}+\|eh'\|_{2,A}
 \le 2\|h\|_{\D}\qquad(h\in\D).
\]
We choose the kernel at an interior point whose distance from
\(\xi\) is the radius of the disk being tested.
Set \(a=(1-\delta)\xi\). The kernel for the norm from
Section~\ref{subsec:02-norm} satisfies
\[
 k_a(z)=\frac{1}{\bar a z}\log\frac{1}{1-\bar a z},
 \qquad
 \|k_a\|_{\D}^2=\frac{1}{|a|^2}\log\frac{1}{1-|a|^2}.
\]
The logarithm is the branch given by the power series, with a
removable value for the kernel at \(z=0\). For \(z\in\DD\cap
B(\xi,\delta)\) we have
\(|1-\bar a z|\le\delta+|a||z-\xi|\le2\delta\).
The real part of the logarithm and \(|\bar a z|\le1\) give
\(|k_a(z)|\ge\log(1/(2\delta))\). Test the derivative bound with
\(h=k_a\) to get
\[
 \log^2(1/(2\delta))\int_{\DD\cap B(\xi,\delta)} |e'|^2\,dA_0
 \le 4\|k_a\|_{\D}^2.
\]
Since \(\delta\le1/256\), we have \(|a|\ge255/256\),
\(1-|a|^2\ge\delta\), and
\(\log(1/(2\delta))\ge(7/8)\log(1/\delta)\).
The resulting upper bound is
\(4(256/255)^2(8/7)^2/\log(1/\delta)\), which is less than the
right side of \eqref{eq:small-region-energy}.
\end{proof}

\subsection{The function \texorpdfstring{\(\beta\)}{beta} and its boundary trace}\label{subsec:08-companion}

For an individual Dirichlet extremal \(e\), Richter and Sundberg's
pseudocontinuation argument gives an analytic function whose
boundary values agree almost everywhere with the product of
\(e^*\) and the conjugate of the nontangential boundary values
of \((ze)'\)
\cite[proof of Theorem~2.2(b), equation~(2.2), p.~869]{RS94}. 
For the estimates on a family, we express the product of \(e^*\)
with the conjugate trace of \((ze)'\) through an analytic function
with a uniform area bound. On a complementary arc where \(e\)
continues and is nonzero, dividing the trace of this function by \(e\)
gives the conjugate trace of \((ze)'\) used for reflection.

We first describe the distributional trace of a square-integrable
analytic function. The same Fourier weights appear in the area norm and in
the dual norm of \(V\).
Use \(dA_0=dx\,dy/\pi\) on the plane, including on neighborhoods
that cross the circle. For a domain \(\Omega\subset\mathbb C\),
write \(A^2(\Omega)\) for the holomorphic functions in
\(L^2(\Omega,dA_0)\). If
\(b(z)=\sum_{n\ge0}b_nz^n\in A^2(\DD)\), angular orthogonality and
radial integration give
\(\|b\|_{A^2(\DD)}^2=\sum_{n\ge0}|b_n|^2/(n+1)\).
Thus \(b\) has a boundary distribution \(b^*\in V^*\) with Fourier
coefficients \(b_n\) for \(n\ge0\) and zero for \(n<0\). The dual norm
\eqref{eq:boundary-dual-norm} gives
\begin{equation}\label{eq:bergman-radial-trace}
 \|b^*\|_{V^*}=\|b\|_{A^2(\DD)},\qquad
 \|b(r\,\cdot)-b^*\|_{V^*}^2
 =\sum_{n\ge0}\frac{|b_n|^2}{n+1}|r^n-1|^2\longrightarrow0
 \quad(r\uparrow1).
\end{equation}
The limit follows by dominated convergence in the series. The star on
a Bergman function denotes this distributional trace.

\begin{lemma}[The boundary-trace identity for \(\beta\)]\label{lem:analytic-companion}
Let \(e\) satisfy \eqref{eq:canonical-moments}, and put
\(k=(ze)'\) and \(w=|e^*|^2\). Then \(\|k\|_{A^2(\DD)}=1\).
The product \(e^*\overline{k^*}\) on the test space \(V\)
is the distributional trace of an analytic function
\(\beta\in A^2(\DD)\), with
\begin{equation}\label{eq:companion-identity}
 \beta(z)=1-\sum_{n\ge1}n\widehat w(n)z^n,\qquad
 \beta^*=e^*\overline{k^*}\ \text{in } V^*,\qquad
 \|\beta\|_{A^2(\DD)}^2\le2 .
\end{equation}
\end{lemma}

In this identity the star on \(e\) denotes its ordinary Hardy
boundary trace, whereas \(k^*\) and \(\beta^*\) are distributions
in \(V^*\). We read their identity by testing on every element of
\(V\), using the bilinear pairing
\eqref{eq:boundary-dual-pairing} and the multiplication bound
\eqref{eq:full-boundary-multiplier-bound}. Multiplication acts on
the test function before the distribution is evaluated. On an arc
where \(e\) continues holomorphically, this distribution agrees
with the ordinary smooth boundary product. The following
subsections have the further task of finding one neighborhood and
bounds that work uniformly for a family of canonical functions.

\begin{proof}
We define the product by transposition on the test space
\(V\), using the distributional derivative trace. The Fourier calculation
then shows that its negative frequencies vanish and gives an
area bound for the analytic function with the remaining coefficients.
Write \(e(z)=\sum_{m\ge0}c_mz^m\). The coefficient of \(z^m\)
in \(k\) is \((m+1)c_m\), so
\(\|k\|_{A^2(\DD)}^2=\sum_{m\ge0}(m+1)|c_m|^2=1\).
We use the bilinear evaluation \(T[\phi]\) of
\eqref{eq:boundary-dual-pairing}. Conjugation is isometric on \(V^*\)
by the symmetry of the weights \((1+|n|)^{-1}\).
The multiplier bound \eqref{eq:full-boundary-multiplier-bound}
ensures \(e^*\phi\in V\) for every \(\phi\in V\). Hence the dual product
\[
 (e^*\overline{k^*})[\phi]=\overline{k^*}[e^*\phi]
 \qquad(\phi\in V)
\]
gives a continuous functional, with norm at most \(2\|k^*\|_{V^*}\).
Its \(n\)-th coefficient is its value on \(\zeta\mapsto\zeta^{-n}\).

For each fixed \(n\in\mathbb Z\), multiplication by \(\zeta^{-n}\)
is bounded on \(V\), since \(1+|j-n|\le(1+|n|)(1+|j|)\).
The Taylor polynomials of \(e\) converge to \(e^*\) on the circle
in \(V\), so they can be used inside the pairing
with \(\overline{k^*}\) after this shift. The finite sums converge to
\begin{equation}\label{eq:companion-coefficients}
 \widehat{e^*\overline{k^*}}(n)
 =\sum_{m\ge\max(0,-n)}(m+1)c_{m+n}\overline{c_m}.
\end{equation}
The series is absolutely convergent by Cauchy--Schwarz, since
\(m+1\le(1+|n|)(m+n+1)\) on the indicated range.

The canonical moments cancel the negative frequencies; their
conjugates determine the positive ones.
If \(n=-q<0\), changing index to \(m=l+q\) gives the canonical
moment \(\sum_{l\ge0}(l+q+1)c_l\overline{c_{l+q}}
=\langle z^qe,e\rangle_{\D}=0\). For \(n=0\) the coefficient
is \(\|e\|_\D^2=1\). For \(n>0\), the conjugate moment
\(\langle e,z^ne\rangle_\D=0\) gives
\[
 \sum_{m\ge0}(m+1)c_{m+n}\overline{c_m}
 =-n\sum_{m\ge0}c_{m+n}\overline{c_m}
 =-n\widehat w(n).
\]
Here the last equality uses the Hardy coefficient product for \(w\)
from Section~\ref{subsec:02-multiplier}.

For the Bergman estimate, we compare the weighted coefficient sum
with the trace bound already available for \(w\in V\).
These coefficients define \(\beta\), since \eqref{eq:source-modulus-V}
gives
\[
 1+\sum_{n\ge1}\frac{n^2}{n+1}|\widehat w(n)|^2
 \le 1+\|w\|_V^2\le2 .
\]
In particular Cauchy--Schwarz gives local uniform convergence of the
analytic series, and the displayed sum is its squared \(A^2\) norm.
The boundary series in \eqref{eq:bergman-radial-trace} has the same
coefficients as the product at every integer frequency. Since
trigonometric polynomials are dense in \(V\), the two functionals agree.
\end{proof}

We use \(\beta\) in two ways. Its real part gives a
distributional identity on the circle. On arcs of holomorphic
continuation its trace agrees with the ordinary product.
Because \(w\) is real, the positive and negative coefficients in
\(\beta^*+\overline{\beta^*}\) combine to give
\begin{equation}\label{eq:companion-real-part}
 2\operatorname{Re}\beta^*=2-\Lambda w
 \qquad\text{in }V^* ,
\end{equation}
where \(\Lambda\) has Fourier symbol \(|n|\).
If \(e\) continues holomorphically across an open arc \(J\), then
\(k\) continues across the same arc. For a smooth test \(\phi\) compactly supported in \(J\), the functional in
\eqref{eq:companion-identity} is the limit
\[
 \beta^*[\phi]
 =\lim_{r\uparrow1}\int_{\T}
   e^*(\zeta)\overline{k(r\zeta)}\phi(\zeta)\,dm(\zeta)
 =\int_J e(\zeta)\overline{k(\zeta)}\phi(\zeta)\,dm(\zeta).
\]
The first equality uses \eqref{eq:bergman-radial-trace} paired with
\(e^*\phi\in V\). The second uses uniform convergence on the compact
support of \(\phi\) and agreement of the Hardy trace of \(e\) there with
its continuation. Thus
\begin{equation}\label{eq:companion-smooth-arc}
 \beta^*|_J=e\,\overline{k}\big|_J
 \qquad\text{as distributions on }J .
\end{equation}
The real-part identity will give a normal derivative formula on
complementary arcs. The local product identity will also identify
the reflected derivative on the common analytic neighborhood.

\subsection{Individual continuation and a uniform lower bound on the modulus}\label{subsec:08-individual}

To use the trace identity \eqref{eq:companion-smooth-arc} on an arc disjoint from a fixed compact
\(K\), we first continue the canonical function across the arc and
show that it is nonzero there. This gives smooth boundary values.
The sign of the outward normal derivative of the negative log modulus
will then give a lower bound for the family on a fixed region just
inside the circle.

For a nonpolar Carleson set \(K\), write \(M=\D_K\) and
\(\eta=\eta_K=G_M\). Section~\ref{subsec:03-seed} shows that
\(\eta\) is outer, that \([\eta]=M\), and that
\(0<\eta(0)<1\). The proof of Theorem~1.3 in \cite[p.~2126]{RY19}
applies Richter and Sundberg's boundary results
\cite[Lemma~3.1(a), Theorems~3.2 and~5.1]{RS94} to this fixed
trace-zero extremal. We give the local compressed-shift proof
for \(M\).
For the later sequence \(K_j\to\{\xi\}\), individual continuation
allows the neighborhoods and lower bounds to depend on \(j\).
Proposition~\ref{prop:common-nonzero-collar} instead fixes a
complementary arc first and uses a common positive lower bound at
the origin to obtain a uniform lower bound for \(|\eta_j|\)
on a fixed neighborhood of that arc.
Proposition~\ref{prop:fixed-analytic-neighborhood} then bounds the
area norm of the continued derivative \(k_j=(z\eta_j)'\)
on a common neighborhood of the arc.
In \eqref{eq:offpoint-derivative-convergence} the excluded disk is
held fixed as \(j\to\infty\). Only afterward does its radius decrease.

\begin{lemma}[Continuation and nonvanishing off the compact set]
\label{lem:individual-continuation}
Let \(K\subset\T\) be a nonpolar Carleson set, and let
\(\eta_K=G_{\D_K}\).
For each \(\zeta\in\T\setminus K\), the function \(\eta_K\)
extends holomorphically to a neighborhood of \(\zeta\).
The continued function satisfies \(\eta_K(\zeta)\ne0\).
\end{lemma}

\begin{proof}
The smooth outer function \(F_K\) will give an inverse for the shift compressed to
\(M^\perp\) at boundary points outside \(K\). Set \(N=M^\perp\).
If \(v\in A^\infty\), Lemma~\ref{lem:bounded-products}
gives a bounded multiplication operator \(M_v\) on \(\D\) with
\(M_vM\subset M\). Write \(\rho(v)=P_NM_v|_N\) and
\(T=\rho(z)\). For \(u,v\in A^\infty\), the term
\(P_NM_u P_M M_v|_N\) vanishes by invariance. Inserting
\(P_N+P_M\) between \(M_u\) and \(M_v\) therefore gives
\[
                 \rho(uv)=\rho(u)\rho(v),\qquad \rho(1)=I_N.
\]
The function \(F_K\) from Theorem~\ref{input:flat-outer}
satisfies \(F_K\D\subset M\) by Section~\ref{subsec:03-ideal},
so \(\rho(F_K)=0\).

Fix \(\zeta\in\T\setminus K\). Since \(F_K(\zeta)\ne0\), we can use a
divided difference whose numerator vanishes at \(\zeta\) and
compresses to \(I_N\). Define
\[
 A_\zeta(z)=\frac{1-F_K(z)/F_K(\zeta)}{z-\zeta}
 =-\frac{1}{F_K(\zeta)}
     \int_0^1 F_K'\bigl(\zeta+s(z-\zeta)\bigr)\,ds .
\]
The segment lies in the closed disk. The integral gives the
closed-disk extension at \(z=\zeta\). Its \(m\)th analytic
derivative is obtained by replacing the integrand by
\(s^m F_K^{(m+1)}(\zeta+s(z-\zeta))\), which is continuous
and bounded on the closed disk for each \(m\). Thus
\(A_\zeta\in A^\infty\), so its compression is bounded.
Compressing \((z-\zeta)A_\zeta=1-F_K/F_K(\zeta)\) gives
\[
 (T-\zeta I_N)\rho(A_\zeta)=I_N
                 =\rho(A_\zeta)(T-\zeta I_N).
\]

We now express the derivative of the projected vector by a resolvent.
The inverse at \(\zeta\) will extend the formula across the circle.
Put \(a=\eta(0)\) and \(q=P_N1=1-a\eta\).
Multiplicativity gives \(T^n=P_NM_{z^n}|_N\), so
\(\|T^n\|\le\sqrt{n+1}\) by Section~\ref{subsec:02-norm}.
Since \(M_{z^n}P_M1\in M\), it also gives \(T^nq=P_Nz^n\).
For \(q(z)=\sum_{n\ge0}q_nz^n\), the coefficient norm and
the linear-first convention yield
\(\langle (T^*)^n q,q\rangle_{\D}=\langle q,T^nq\rangle_{\D}
=\langle q,z^n\rangle_{\D}=(n+1)q_n\).
The series for \((I_N-zT^*)^{-1}\) converges in operator norm
when \(|z|<1\). Hence
\begin{equation}\label{eq:individual-resolvent}
                 (zq)'(z)
      =\big\langle (I_N-zT^*)^{-1}q,q\big\rangle_{\D}.
\end{equation}
At \(z=\zeta\), the operator \(I_N-\zeta T^*\) is invertible
because it equals \(-\zeta(T-\zeta I_N)^*\).
The right side of \eqref{eq:individual-resolvent} is therefore
holomorphic on a small disk about \(\zeta\). The disk may depend on
the compact set and the chosen boundary point. Choose this disk to
avoid \(0\), and take a primitive with value \(z_0q(z_0)\) at an
interior point \(z_0\) in its intersection with \(\DD\).
That intersection is connected. The primitive agrees with \(zq\)
there, so division by \(z\) continues \(q\) and thus \(\eta\).

We test a possible zero by dividing \(\eta\) by \(z-\zeta\) and using
the Hardy quotient property of \(M\). Suppose the continued function
vanishes at \(\zeta\). The quotient
\(h(z)=\eta(z)/(z-\zeta)\) then extends holomorphically through
\(\zeta\). On a smaller neighborhood of \(\zeta\), both \(h\)
and \(h'\) are bounded. On the rest of the disk the denominator
is bounded away from zero, so the quotient rule bounds \(|h'|\)
by a constant times \(|\eta'|+|\eta|\). Since \(\eta\in\D\)
and \(\int_\DD|\eta|^2\,dA_0\le\|\eta\|_{H^2}^2\), the derivative
of \(h\) is square integrable. The coefficient formula for the
Dirichlet integral then gives \(h\in\D\).

For quasi-every \(\omega\in K\), the radial limit of \(\eta\) is
zero and the denominator of \(h(r\omega)\) tends to
\(\omega-\zeta\ne0\). By \eqref{eq:radial-trace-agreement},
this puts \(h\) in \(M=[\eta]\). The conclusion of
Lemma~\ref{lem:completed-canonical-quotient} then gives
\(h=\eta H\) for some \(H\in H^2\). Because \(\eta\) is
outer, \(H(z)=1/(z-\zeta)\) in \(\DD\).
But this reciprocal has the expansion
\[
 \frac{1}{z-\zeta}=-\bar\zeta\sum_{n\ge0}(\bar\zeta z)^n ,
\]
whose coefficients all have modulus one. This contradicts
\(H\in H^2\) and proves nonvanishing.
The common interior and boundary lower bound below uses a
separate estimate for the reflected logarithm.
\end{proof}

By continuity, \(|\eta_K|\) is bounded below on a sufficiently
small disk about each point of \(\T\setminus K\). In particular
\(\Zrad(\eta_K)\subset K\), and
\(C_R(K\setminus\Zrad(\eta_K))=0\) by the trace condition on \(K\).
The same local lower bound gives
\(\liminf_{\DD\ni z\to\omega}|\eta_K(z)|>0\) at every
\(\omega\in\T\setminus K\).
The size of the continuation disk can depend on \(K\) and on the
chosen boundary point. The tangential modulus limits in
\cite[Theorem~1.4]{RY19} concern one extremal at a time. For a family
of shrinking compact sets, we will use the lower bound at the origin
to obtain a fixed region next to a complementary arc on which
division by the canonical functions is uniformly bounded.

\paragraph{The normal derivative on a complementary arc.}
We first compute the normal derivative of the squared modulus.  The trace identity \eqref{eq:companion-smooth-arc} gives one expression for this derivative on the arc, and the identity \eqref{eq:canonical-source-identity} gives its sign after taking the logarithm.
By Lemma~\ref{lem:individual-continuation}, \(\eta\) extends
holomorphically across every point of \(\T\setminus K\) and is
nonzero there. Write \(w=|\eta^*|^2\) and \(k=(z\eta)'\). At a point \(\zeta\) of an
open arc \(J\subset\T\setminus K\), write \(\partial_rw(\zeta)\)
for the derivative of \(r\mapsto|\eta(r\zeta)|^2\) at \(r=1\)
in the direction of increasing radius. This derivative is
\(2\operatorname{Re}(\zeta\eta'(\zeta)\overline{\eta(\zeta)})\) there.
The Cauchy-transform identity in \cite[Lemma~2.4, (2.2) and
the calculation at (2.8)]{GZ10} gives a parallel calculation
on such an arc.

The distributional trace \(\beta^*\) agrees with the smooth
function \(\eta\overline{k}\) on \(J\), by
\eqref{eq:companion-smooth-arc}. Since \(k=\eta+z\eta'\), this gives
\(2\operatorname{Re}\beta^*=2w+\partial_rw\) on \(J\). On the whole
circle, \eqref{eq:companion-real-part} gives
\(2\operatorname{Re}\beta^*=2-\Lambda w\) in distributions.
Combining these identities
with \eqref{eq:canonical-source-identity} gives
\begin{equation}\label{eq:arc-normal-source}
 \partial_rw=2-2w-\Lambda w=1-w+\ell_\eta\qquad\text{on }J.
\end{equation}
The corresponding almost-everywhere identity for extremal functions
appears in \cite[equation~(3.2), p.~873]{RS94}. The identity on \(J\)
initially holds in the sense of distributions.
It identifies \(\ell_\eta\) there with the smooth function
\(\partial_rw-1+w\). This representative is nonnegative everywhere
on the arc, since \(\ell_\eta\) is nonnegative almost everywhere.

Let \(U=-\log|\eta|\) in \(\DD\). The function \(\eta\) is outer and
has modulus at most one, so \(U\) is nonnegative and harmonic.
It is smooth near each point of \(J\) by individual nonvanishing and
continuation. On \(J\), \(0<w\le1\) and
\(\partial_rw=-2w\,\partial_rU\). Hence
\begin{equation}\label{eq:arc-log-normal}
  2\partial_rU=1-\frac{1+\ell_\eta}{w}\le0 .
\end{equation}

For the constant canonical function \(e=1\), we have
\(w=1\), \(\ell_e=0\), and \(\beta=1\).
Thus \(\partial_r w=0\), as expected.

\paragraph{A reflected logarithm.}
We use the sign in \eqref{eq:arc-log-normal} to reflect the logarithm
across the arc. Subharmonicity across the circle will let us estimate
its value from an area mean on a disk that crosses the arc.
Set \(\mathcal U(z)=U(z)\) for \(|z|<1\) and
\(\mathcal U(z)=U(1/\bar z)\) for \(|z|>1\), with its continuous
value on any complementary arc. The inversion preserves harmonicity
on the exterior side. Across the arc, the outer piece has
increasing-radius derivative \(-\partial_rU\), and the inner piece
has derivative \(\partial_rU\). Integration by parts on the two
sides against a smooth test supported near the arc cancels the terms
involving the common boundary value. The remaining term in the
distributional Laplacian is
\[
            (-2\partial_r U)\,ds,
\]
where \(ds\) is ordinary arclength on the circle. It is nonnegative
by \eqref{eq:arc-log-normal}. Thus \(\mathcal U\) is subharmonic on
each small disk that meets the circle only off \(K\). On either
side of the circle, its angular mean on a centered circle is
\(U(0)=-\log\eta(0)\). On the exterior side this follows by taking
the angular mean at the reciprocal radius inside the disk.

\begin{proposition}[A common lower bound near a complementary arc]
\label{prop:common-nonzero-collar}
Fix \(\xi\in\T\) and nonpolar Carleson sets \(K_j\) with
\(d_H(K_j,\{\xi\})\to0\) in the chord metric. Put
\(\eta_j=G_{\D_{K_j}}\). Suppose that \(0<a_*<1\) and that
\(\eta_j(0)\ge a_*\) for all sufficiently large \(j\).
Fix a compact arc \(A\subset\T\setminus\{\xi\}\), let
\(d=\dist(A,\xi)>0\), and choose
\(0<\rho\le\min(d/16,1/16)\). Define
\[
                    M_*=-\log a_*,\qquad b_\rho=e^{-3M_*/\rho}.
\]
Then for all sufficiently large \(j\),
\begin{equation}\label{eq:common-nonzero-collar}
 -\log|\eta_j(z)|\le 3M_*/\rho,\qquad |\eta_j(z)|\ge b_\rho
 \quad\text{if } z\in\DD\ \text{and }\dist(z,A)\le\rho .
\end{equation}
The same bounds hold at the continued boundary values where
\(z\in\T\) and \(\dist(z,A)\le\rho\).
\end{proposition}

\begin{proof}
We use the separation from \(K_j\) to keep the averaging disk in the
subharmonic region. The means on centered circles then bound the
integral over the surrounding annulus.
Choose \(J_0\) so that \(\eta_j(0)\ge a_*\) and
\(K_j\subset B(\xi,d/2)\) for \(j\ge J_0\).  For a point \(z\) at distance
at most \(\rho\) from \(A\), every point of \(B(z,2\rho)\) lies
within \(3\rho\) of \(A\). Since \(3\rho<d/2\), the disk avoids \(K_j\)
and lies in the annulus \(1-3\rho<|u|<1+3\rho\).
Thus the reflected logarithm \(\mathcal U_j\) is subharmonic on
\(B(z,2\rho)\). It is nonnegative. Its ordinary area integral
over the annulus equals
\[
 2\pi\bigl(-\log\eta_j(0)\bigr)\int_{1-3\rho}^{1+3\rho} r\,dr
    =12\pi\rho\bigl(-\log\eta_j(0)\bigr)\le12\pi\rho M_*.
\]
Here the possible values on \(K_j\) do not affect an area integral.
The subharmonic area mean inequality therefore bounds
\(\mathcal U_j(z)\) by \(12\pi\rho M_*/(4\pi\rho^2)=3M_*/\rho\).
This proves both assertions.
\end{proof}

The bound \(b_\rho\) depends only on \(a_*\) and \(\rho\). The
starting index \(J_0\) is chosen after fixing the arc and these
parameters. It uses the two thresholds \(\eta_j(0)\ge a_*\)
and \(d_H(K_j,\{\xi\})<d/2\). The positive limiting value at zero in
Theorem~\ref{thm:C} gives a uniform lower bound \(a_*\)
for all sufficiently large \(j\). The lower bound \(b_\rho\) allows us to bound
\(\beta_j/\eta_j\) in \(L^2\) inside the disk. Reflecting that quotient will give an
exterior formula on a region fixed by \(A\) and \(\rho\).

\subsection{Analytic continuation to a common neighborhood}\label{subsec:08-gluing}

Individual continuation supplies the local boundary identities we
need to reflect. The sizes of those individual neighborhoods may
depend on \(j\). For instance, if \(0<r_j\uparrow1\), then each
\(1/(1-r_jz)\) extends across \(\T\), while its pole at \(1/r_j\)
precludes a common neighborhood of \(1\) for the whole family.
The lower bound for \(|\eta_j|\) on the fixed neighborhood
allows us to define an exterior formula on a fixed region. To join it to the interior derivative, we use the boundary
identity \eqref{eq:companion-smooth-arc}. The next lemma turns the distributional
identity there into equality of analytic functions on the inner
side of the arc.

\begin{lemma}[Local uniqueness of an analytic boundary trace]
\label{lem:local-trace-uniqueness}
Let \(g\in A^2(\DD)\), and let \(g_0\) be holomorphic on a
neighborhood of an open arc \(J\subset\T\). If the distributional
trace of \(g\) agrees with the smooth trace of \(g_0\) on \(J\),
then \(g=g_0\) on the inner side of some neighborhood of each
point of \(J\).
\end{lemma}

\begin{proof}
We extend the difference by zero across the arc and show that the
extension is holomorphic. The trace equality will make the boundary
term vanish when testing the distributional \(\bar\partial\) of the
extension. To localize the boundary pairing to the arc where the
traces agree, we first record a cutoff estimate on the trace space. Let \(\chi\)
be a smooth function on \(\T\), with chordal Lipschitz constant
\(L_\chi\). Write \(\mathcal E(F)\) for the double integral term in
\eqref{eq:V-difference}. Expanding the difference of the product
at \(\zeta,\omega\in\T\) and using
\(|a+b|^2\le2|a|^2+2|b|^2\) gives
\begin{equation}\label{eq:trace-cutoff-localization}
 \mathcal E(\chi F)\le
     2\|\chi\|_\infty^2\mathcal E(F)+2L_\chi^2\|F\|_{L^2(m)}^2 .
\end{equation}
Together with the \(L^2\) estimate, this makes multiplication by
\(\chi\) bounded on \(V\), and duality gives the same property
on \(V^*\).

Take a disk \(U_0\) centered at a point of \(J\), small enough
that its closure lies in the domain of \(g_0\) and
\(\overline U_0\cap\T\subset J\). On \(U_0\), define \(H=g-g_0\)
on the inner side and \(H=0\) on the exterior side. This function
is locally square integrable. Let \(\phi\in C_c^\infty(U_0)\), with
its smooth extension by zero outside \(U_0\). For \(r<1\) close to
one, Cauchy--Green with \(\bar\partial=(\partial_x+i\partial_y)/2\)
gives
\begin{equation}\label{eq:trace-zero-extension-stokes}
 \int_{U_0\cap\{|z|<r\}}(g-g_0)\bar\partial\phi\,dx\,dy
  =\frac r2\int_0^{2\pi}
    (g-g_0)(re^{it})\,\phi(re^{it})e^{it}\,dt .
\end{equation}
The integrand on the right is defined to be zero where the test
vanishes. The boundary term on \(\partial U_0\) vanishes because
\(\phi\) has compact support. On the counterclockwise arc of
\(|z|=r\), \(dz=ir e^{it}\,dt\), giving the factor \(r/2\).

For all \(r\) sufficiently close to one, the nonzero tests on the
right are supported in a fixed compact subarc of \(J\). Choose
\(\chi\) supported in \(J\), equal to one on that subarc, with \(g_0\)
defined on a fixed annulus over its support. By
\eqref{eq:bergman-radial-trace} and \eqref{eq:trace-cutoff-localization},
the localized radial traces of \(g\) converge in \(V^*\). The
localized traces of \(g_0\) converge smoothly. Their difference
tends to zero in \(V^*\), by the hypothesis on \(J\).
The tests \(e^{it}\phi(re^{it})\) converge smoothly and hence in \(V\),
by the difference norm. The boundary pairing in
\eqref{eq:trace-zero-extension-stokes} tends to zero, since the
localized difference tends to zero and the test norms stay bounded.
The area integrals converge by local integrability on the support of
\(\bar\partial\phi\). Thus \(\bar\partial H=0\) in distributions.

For completeness, convolve \(H\) with smooth approximate identities
inside disks with closure in \(U_0\). Each convolution is holomorphic
on a slightly smaller disk, since \(\bar\partial H=0\).
The convolutions converge in local \(L^2\). On disks farther
inside, the holomorphic area mean inequality applied to their
differences makes them uniformly Cauchy. Their limits give a
holomorphic representative of \(H\). These representatives agree
where they overlap because they equal the same \(H\) almost
everywhere. The representative is zero on the exterior open part of
the connected disk \(U_0\), hence zero throughout \(U_0\) by the
identity theorem. This proves the asserted equality on the inner
part of \(U_0\).
\end{proof}

We now apply the lemma to the Bergman functions \(\beta_j\). The uniform lower
bound from Proposition~\ref{prop:common-nonzero-collar} fixes the
region where the quotient is controlled. The individual
continuations supply the local equalities needed for gluing.

\begin{proposition}[Fixed continuation and convergence away from the collapse point]
\label{prop:fixed-analytic-neighborhood}
Assume the hypotheses of Proposition~\ref{prop:common-nonzero-collar},
and fix \(A,\rho,b_\rho\) as there. For all sufficiently large \(j\),
\(k_j=(z\eta_j)'\) extends holomorphically to
\(W=\{z:\dist(z,A)<\rho/4\}\), with
\begin{equation}\label{eq:common-continuation-bounds}
 \begin{split}
 \|k_j\|_{A^2(W)}^2&\le B_\rho:=1+32b_\rho^{-2},\\
 |k_j(z)|^2&\le64B_\rho/\rho^2
             \qquad(\dist(z,A)<\rho/8).
 \end{split}
\end{equation}
If \(\eta_j\rightharpoonup\eta\) in \(\D\), then \(\eta\) continues
holomorphically and is nonzero across each point of
\(\T\setminus\{\xi\}\). The \(\eta_j\) and their derivatives
converge locally uniformly on fixed smaller neighborhoods of
those points. In particular, for every \(\varepsilon>0\),
\begin{equation}\label{eq:offpoint-derivative-convergence}
 \|\eta_j'-\eta'\|_{L^2(\DD\setminus B(\xi,\varepsilon),dA_0)}
                         \longrightarrow 0 .
\end{equation}
\end{proposition}

\(B_\rho\) depends only on \(a_*\) and \(\rho\). For the
convergence assertion, each smaller neighborhood is fixed before
the sequence limit is taken.

\begin{proof}
\emph{Define the exterior expression on the fixed region.}
On the neighborhood inside the disk \(E_\rho=\{z\in\DD:\dist(z,A)<\rho\}\), \(\beta_j\) from Lemma~\ref{lem:analytic-companion} satisfies
\[
        \|\beta_j/\eta_j\|_{L^2(E_\rho,dA_0)}^2
                        \le2b_\rho^{-2}.
\]
We reflect and conjugate the holomorphic quotient to obtain a
function on the fixed exterior region. For
\(z\in W\setminus\overline\DD\), define
\begin{equation}\label{eq:reflected-companion-quotient}
          \widetilde k_j(z)=
                 \overline{(\beta_j/\eta_j)(1/\bar z)}.
\end{equation}
For every \(\zeta\in A\),
\(|1/\bar z-\zeta|=|z-\zeta|/|z|\), so the reflected point lies
in \(E_\rho\). The inversion is antiholomorphic, and the conjugate
in \eqref{eq:reflected-companion-quotient} makes the expression
holomorphic. In the area change of variables \(u=1/\bar z\),
the inverse area Jacobian is \(|u|^{-4}\). Here \(|u|>1/2\),
since \(|z|<1+\rho/4\). It follows that
\[
 \|k_j\|_{L^2(W\cap\DD,dA_0)}^2\le1,\qquad
 \|\widetilde k_j\|_{L^2(W\setminus\overline\DD,dA_0)}^2
                    \le 32 b_\rho^{-2}.
\]

\emph{Agreement of \(\widetilde k_j\) with the continuation of \(k_j\).}
Fix one sufficiently large \(j\) and one \(\zeta\in W\cap\T\).
This point lies off \(K_j\). Individual continuation supplies
a small neighborhood on which \(\eta_j\) and \(k_j\) are
holomorphic and \(\eta_j\) is nonzero. Shrink the neighborhood
so that the reflected points are in the domain of the continued
\(k_j\). The function
\[
                  g_0(z)=\eta_j(z)\overline{k_j(1/\bar z)}
\]
is holomorphic there. On the corresponding arc its trace
is the smooth trace of the product \(\eta_j\overline{k_j}\).
By \eqref{eq:companion-smooth-arc} this is the distributional
trace of \(\beta_j\). Lemma~\ref{lem:local-trace-uniqueness}
gives \(\beta_j=g_0\) in the adjoining inner neighborhood.
Divide there by \(\eta_j\), then reflect and conjugate.
The individual continuation of \(k_j\) agrees with
\(\widetilde k_j\) on the adjoining exterior neighborhood.
These local agreements show that the definitions on the two
sides glue holomorphically throughout \(W\). Each local equality uses
an individual continuation disk. The exterior formula and its area
bound already hold on the exterior part of the fixed region \(W\),
so the extension of \(k_j\) has the same domain for all sufficiently large \(j\).

Adding the area bounds gives the first line of
\eqref{eq:common-continuation-bounds}. For a point at distance less
than \(\rho/8\) from \(A\), its disk of radius \(\rho/8\) is
contained in \(W\). The area mean inequality for \(|k_j|^2\)
with normalized area \(dA_0=dx\,dy/\pi\) gives the second line.
The same argument on disks in any compact subset of \(W\),
followed by the Cauchy estimates, bounds all fixed derivatives
there uniformly for sufficiently large \(j\).

\emph{Pass to a weak limit.}
Weak Dirichlet convergence gives locally uniform convergence
inside \(\DD\), by the evaluation and compactness estimates
of Section~\ref{subsec:02-norm}. In particular
\(k_j\to(z\eta)'\) there. Choose a disk \(W_0\) centered at a
point of \(A\), with \(\overline{W_0}\subset W\) and
\(0\notin\overline{W_0}\). The uniform bounds just proved give
normal-family compactness on \(W_0\). Every locally uniform
subsequence limit agrees with \((z\eta)'\) on the inner part,
where weak Dirichlet convergence has already determined the limit.
The identity theorem then fixes the continuation throughout \(W_0\).
Consequently the whole weakly convergent sequence of \(k_j\)
converges locally uniformly there.

Fix \(z_0\in W_0\cap\DD\). Integration of the extended derivatives
along paths in the disk gives
\[
        z\eta_j(z)=z_0\eta_j(z_0)+\int_{z_0}^z k_j(s)\,ds .
\]
The right side defines the continuation because it agrees with
\(z\eta_j\) on the connected intersection \(W_0\cap\DD\).
The primitives converge locally uniformly, and \(z\ne0\) on
\(W_0\) permits division. This gives the claimed convergence and,
by Cauchy's formula, convergence of every fixed derivative on
smaller disks. The limiting continued boundary value has modulus
at least \(b_\rho\) on \(A\), by
\eqref{eq:common-nonzero-collar}. Each boundary point off \(\xi\)
admits such an arc \(A\), so the limit is nonzero there.

For a fixed \(\varepsilon>0\), finitely many smaller disks cover
the portion of \(\T\) outside \(B(\xi,\varepsilon)\). The remaining
part of \(\overline\DD\setminus B(\xi,\varepsilon)\) is a compact
subset of \(\DD\) after deleting these neighborhoods. Convergence
of the derivatives on the two kinds of sets gives uniform convergence
on the entire complement, with the continued values on the circle.
Integrating over its portion inside \(\DD\) proves
\eqref{eq:offpoint-derivative-convergence}.
\end{proof}

Section~\ref{sec:09} will first fix a small disk about \(\xi\).
Equation \eqref{eq:offpoint-derivative-convergence} gives convergence
of the derivatives on its complement. The uniform small-region
estimate \eqref{eq:small-region-energy} makes the energy on the disk
small as its radius decreases. Together these two estimates control
the derivative term of the Dirichlet norm. The Hardy term is controlled
by the compact inclusion from Section~\ref{subsec:02-norm}.

\newpage \section{Strong convergence to the canonical function of one atom}\label{sec:09}

We recover the canonical function of one atom from trace-zero
spaces whose defining compact sets shrink to a single point. First we
establish strong compactness by controlling the derivative energy
near that point and the derivatives away from it.
Next we follow the boundary logarithms to locate the singular part
of any limit. To identify its mass, we prove that the outer factor
is cyclic and use the canonical value at zero.

\subsection{An analytic criterion for strong convergence}\label{subsec:09-criterion}

\begin{theorem}[Analytic recovery criterion C]\label{thm:C}
Fix \(\tau>0\) and \(\xi\in\T\). Let \(K_j\subset\T\) be nonpolar
Carleson sets with \(d_H(K_j,\{\xi\})\to0\), and let
\(\eta_j=G_{\D_{K_j}}\). If
\[
                  \eta_j(0)\longrightarrow G_{\tau,\xi}(0),
\]
then
\[
                  \|\eta_j-G_{\tau,\xi}\|_{\D}\longrightarrow0 .
\]
\end{theorem}

Convergence of values at zero alone would not control the Dirichlet
norm of arbitrary functions. For the canonical family in this theorem,
the small-disk energy estimate \eqref{eq:small-region-energy}
and the \(L^2\) convergence of the derivatives away from \(\xi\)
in \eqref{eq:offpoint-derivative-convergence} provide
the further control needed for strong convergence.
Hartz and Richter \cite[Lemma~5.1]{HR26} compare invariant subspaces defined by inner functions under local uniform convergence of the inner
functions and pointwise convergence of their boundary local energies.
They require each local energy to be bounded by a fixed multiple
of the local energy of the limiting inner function, uniformly in the
index and boundary point. Here each \(\eta_j\) is outer:
its Hardy inner factor is constant, although the limiting canonical
function is to have a positive singular mass. The lower bound near each compact arc away from \(\xi\) gives a
positive lower bound there for the modulus of the outer factor of
any subsequential limit. Once this outer factor is shown to be cyclic,
Theorem~\ref{input:module} identifies the invariant subspace generated
by the limit.
The functions \(\eta_j\) are the unit canonical functions of the trace-zero spaces, normalized to have positive value at zero. Section~\ref{subsec:03-seed}
gives their existence and outerness. For nonempty \(K_j\), the
Hausdorff distance to \(\{\xi\}\) is
\(\max_{\zeta\in K_j}|\zeta-\xi|\). Thus the location hypothesis
says that the entire set eventually lies in every neighborhood
of \(\xi\). A family with \(\eta_j(0)=G_{\tau,\xi}(0)\) at each stage
satisfies the convergence hypothesis.

Write \(a=G_{\tau,\xi}(0)\in(0,1)\). In the compactness step, take
\(a_*=a/2\) for the eventual lower bound \(\eta_j(0)\ge a_*\)
in the estimates on this neighborhood. The starting index is chosen using the convergence
\(\eta_j(0)\to a\). In the identification step we use the
exact limit \(a\). The positive gap \(1-a^2>0\) will exclude a cyclic
limit, and Proposition~\ref{prop:atomic-origin-strict} will determine the mass once we know the subspace generated by the limit.

\subsection{Retaining the Dirichlet norm}\label{subsec:09-norm}

Section~\ref{sec:08} gives two estimates for the derivative term
in the Dirichlet norm. Continuation gives convergence away from
\(\xi\), and the small-disk energy bound makes the contribution near
\(\xi\) uniformly small. The Hardy term is handled by the compact inclusion into \(H^2\).

Keep a small neighborhood of the collapse point fixed while
passing to the sequence limit on its complement. The uniform
energy bound controls the contribution inside that neighborhood
for both the sequence and its weak limit; Fatou's lemma gives the
latter bound before the limit is known to have unit norm. Only
after taking the sequence limit do we shrink the neighborhood.
Compact inclusion supplies convergence of the separate Hardy
term. With both parts of the norm controlled, the unit norm
and canonical moments can then pass to the limit.

\begin{proposition}[Strong subsequential compactness]
\label{prop:recovery-norm-retention}
Fix \(\xi\in\T\). Let \(K_j\) be nonpolar Carleson sets in \(\T\) with
\(d_H(K_j,\{\xi\})\to0\), and put \(\eta_j=G_{\D_{K_j}}\).
Assume \(0<a_*<1\) and \(\eta_j(0)\ge a_*\) for all sufficiently
large \(j\).
Then every weakly convergent subsequence of \((\eta_j)\) converges
strongly in \(\D\), and every subsequence has a further subsequence
converging strongly in \(\D\). Each such limit \(\eta\) is zero-free
in \(\DD\), satisfies \(|\eta|\le1\), and is the unit canonical function of \([\eta]\) with positive value at zero. If \(\eta_j(0)\to a\), then every
such limit satisfies \(\eta(0)=a\).
\end{proposition}

\begin{proof}
Take a weakly convergent subsequence, keep the index \(j\) for that
subsequence, and write \(\eta_j\rightharpoonup\eta\) in \(\D\).
Weak lower semicontinuity gives \(\|\eta\|_\D\le1\).
Section~\ref{subsec:02-norm} gives locally uniform convergence of
the functions and their derivatives. In particular
\(\eta(0)=\lim_j\eta_j(0)\ge a_*\), and the canonical bound
\(|\eta_j|\le1\) passes to \(\eta\). Each \(\eta_j\) is outer
by Section~\ref{subsec:03-seed} and therefore has no zeros in
\(\DD\). The positive value at zero ensures a nonzero limit, so
Hurwitz's theorem gives that \(\eta\) is zero-free in \(\DD\).

Fix \(0<\delta\le1/256\) and put \(E_\delta=\DD\cap B(\xi,\delta)\).
For every \(j\), the uniform energy bound
\eqref{eq:small-region-energy} gives
\[
 \int_{E_\delta}|\eta_j'|^2\,dA_0
       \le \frac{4\pi}{\log(1/\delta)} .
\]
Since \(\eta_j'\to\eta'\) pointwise in \(\DD\), Fatou's lemma gives
\begin{equation}\label{eq:limit-small-region-energy}
 \int_{E_\delta}|\eta'|^2\,dA_0
 \le \liminf_j\int_{E_\delta}|\eta_j'|^2\,dA_0
 \le \frac{4\pi}{\log(1/\delta)} .
\end{equation}
Fatou's lemma gives the bound for \(\eta\) before its canonical
normalization is established. For the fixed \(\delta\), use
the common neighborhoods away from
\(B(\xi,\delta)\) for sufficiently large \(j\). The resulting derivative convergence in
\eqref{eq:offpoint-derivative-convergence} gives
\(\|\eta_j'-\eta'\|_{L^2(\DD\setminus E_\delta,dA_0)}\to0\).
On \(E_\delta\), bound the difference using the energies of the
sequence and its limit and the inequality
\(|\eta_j'-\eta'|^2\le2|\eta_j'|^2+2|\eta'|^2\).
It follows that
\[
 \limsup_j D(\eta_j-\eta)\le
       \frac{16\pi}{\log(1/\delta)} .
\]
For the derivative energy, the radius \(\delta\) was fixed during
the sequence limit. Only then let
\(\delta\downarrow0\) to obtain \(D(\eta_j-\eta)\to0\).

For the Hardy part of the Dirichlet norm, we separate a finite coefficient
sum from the tail. We use the compact inclusion
\(\D\hookrightarrow H^2\) from Section~\ref{subsec:02-norm}.
In detail, if \(c_{j,n}\) and \(c_n\) are the Taylor coefficients of
\(\eta_j\) and \(\eta\), weak convergence gives
\(c_{j,n}\to c_n\) for each fixed \(n\), while
\[
 \sum_{n>N}|c_{j,n}-c_n|^2
 \le\frac{\|\eta_j-\eta\|_\D^2}{N+2}
 \le\frac4{N+2}.
\]
Taking the limit on the first \(N+1\) coefficients and then sending
\(N\) to infinity gives \(\|\eta_j-\eta\|_{H^2}\to0\). Together with
the derivative estimate this gives \(\|\eta_j-\eta\|_\D\to0\) in
the norm \eqref{eq:full-norm}.

Every subsequence has a weakly convergent further subsequence because
all \(\eta_j\) have unit norm. The two norm estimates show that
convergence is strong along each such further subsequence. We can
now pass the canonical conditions to the limit. Its norm is one,
and boundedness of multiplication by \(z^n\) on \(\D\) for each
fixed \(n\ge1\) gives convergence of the moments
\(\langle z^n\eta_j,\eta_j\rangle_\D=0\) to the corresponding
moments of \(\eta\).
With \(a=\eta(0)>0\), \eqref{eq:canonical-moments} and \eqref{eq:certificate} give
\begin{equation}\label{eq:recovery-limit-certificate}
 P_{[\eta]}1=a\eta,\qquad
 \operatorname{dist}_\D(1,[\eta])^2=1-a^2 .
\end{equation}
Thus \(\eta=G_{[\eta]}\), and if the values \(\eta_j(0)\) converge for the original
sequence then \(a\) equals that common limit.
\end{proof}

Under the hypothesis \(\eta_j(0)\to G_{\tau,\xi}(0)\) in Theorem~\ref{thm:C}, the squared
distance in \eqref{eq:recovery-limit-certificate} is strictly
positive. To decide which subspace the limit generates, we first recover
the inner and outer factors of each subsequential limit.

\subsection{Boundary logarithms and factorization of the limit}\label{subsec:09-logarithm}

Under Proposition~\ref{prop:recovery-norm-retention}, work along a
subsequence for which \(\eta_j\to\eta\) in \(\D\), keeping the index
\(j\) for that subsequence. The outer representation expresses
each function through its boundary logarithm. We keep track of these
logarithms as measures so that a possible concentration at \(\xi\)
can appear in the limit. Put
\(a=\eta(0)=\lim_j\eta_j(0)\ge a_*\).
Each \(\eta_j\) is outer and bounded by one in modulus. The normalized
outer representation \eqref{eq:outer-intro} gives
\begin{equation}\label{eq:outer-logarithmic-measures}
 u_j=-\log|\eta_j^*|\in L^1(m),\qquad
 \lambda_j=u_j\,m,\qquad
 u_j\ge0,\qquad \lambda_j(\T)=-\log\eta_j(0).
\end{equation}
For any finite measure \(\lambda\) on \(\T\), write
\[
 \mathcal H\lambda(z)=\int_\T\frac{\zeta+z}{\zeta-z}\,d\lambda(\zeta).
\]
Thus \(\eta_j=\exp(-\mathcal H\lambda_j)\), with phase fixed by
\(\eta_j(0)>0\). The measures \(\lambda_j\) record the outer logarithms.
By \eqref{eq:source}, the measure \(\sigma_{\eta_j}\) on \(\T\)
has density \(\ell_{\eta_j}\) with respect to \(m\).

For all sufficiently large \(j\), \(\lambda_j(\T)\le-\log a_*\). A further
subsequence of these positive measures therefore converges weak-*
on \(C(\T)\) to a positive finite measure \(\lambda\). Testing against
the constant function gives \(\lambda(\T)=-\log a\). For a fixed
\(z\in\DD\), the kernel in \(\mathcal H\lambda(z)\) is a continuous
test function of \(\zeta\). These kernels are bounded and uniformly
continuous on \(\T\) times any compact subdisk. A finite net in the
subdisk and the uniform mass bound thus give local uniform
convergence of the transforms. Strong convergence in \(\D\) already
gives local uniform convergence of \(\eta_j\), so
\begin{equation}\label{eq:logarithmic-limit-representation}
                  \eta(z)=\exp(-\mathcal H\lambda(z)).
\end{equation}
Thus \(\mathcal H\lambda\) is the negative analytic logarithm of
\(\eta\) with value \(-\log a\) at zero. This normalization determines it
uniquely. Expanding the kernel as
\(1+2\sum_{n\ge1}(z\bar\zeta)^n\) shows that this logarithm
determines the mass and all positive Fourier coefficients of
\(\lambda\). Reality determines the negative coefficients, so
Fourier uniqueness determines \(\lambda\). Every weak-* convergent
further subsequence of the logarithmic measures gives the same
representation of \(\eta\) and hence the same \(\lambda\). Thus the
measures converge weak-* to \(\lambda\) along the chosen strongly
convergent analytic subsequence.

To locate the singular part of \(\lambda\), we look on arcs away
from \(\xi\). Uniform convergence there can be passed through the
logarithm once the moduli are bounded below. Let \(J\subset\T\) be
an open arc whose closure avoids \(\xi\). Proposition~\ref{prop:common-nonzero-collar}
applies with the fixed arc \(A=\overline J\) and \(\rho\) chosen
from its separation from \(\xi\). It gives
\(|\eta_j(\zeta)|\ge b_\rho>0\) for \(\zeta\in\overline J\)
and all sufficiently large \(j\).
By Proposition~\ref{prop:fixed-analytic-neighborhood}, \(\eta_j\)
and \(\eta\) are holomorphic on fixed neighborhoods of the points
of \(\overline J\), with local uniform convergence there. A finite
cover gives uniform convergence on \(\overline J\). Their Hardy
traces agree almost everywhere there with the continued values, so
\(u_j\to-\log|\eta|\) uniformly for these representatives on
\(\overline J\). In particular, for any continuous \(\phi\) compactly
supported in \(J\),
\[
 \int_\T\phi\,d\lambda
     =\lim_j\int_\T\phi u_j\,dm
     =\int_J\phi(-\log|\eta|)\,dm .
\]
Thus \(\lambda|_J=(-\log|\eta|)\,m|_J\). A countable cover of
\(\T\setminus\{\xi\}\) by such arcs yields
\begin{equation}\label{eq:logarithmic-measure-decomposition}
 \begin{gathered}
 \lambda=u\,m+T\delta_\xi,\qquad T=\lambda(\{\xi\})\ge0,\\
 u=-\log|\eta^*|\quad m\text{-a.e.},\qquad
 u\ge0,\quad u\in L^1(m).
 \end{gathered}
\end{equation}
The density is locally bounded off \(\xi\) and integrable because
\(\lambda\) has finite mass. Since \(m(\{\xi\})=0\), the mass
\(T\) satisfies \(0\le T\le\lambda(\T)=-\log a\).

Use \(S_{T,\xi}=S_{T\delta_\xi}\) with \(S_{0,\xi}=1\), as in
Section~\ref{subsec:02-atomic}. Inserting
\eqref{eq:logarithmic-measure-decomposition} into
\eqref{eq:logarithmic-limit-representation} gives the analytic
factorization
\begin{equation}\label{eq:logarithmic-limit-factorization}
         \eta=S_{T,\xi}O,\qquad O(z)=\exp(-\mathcal H(u\,m)(z)).
\end{equation}
The real part of \(\mathcal H(u\,m)\) is the nonnegative Poisson
integral of \(u\). Hence \(|O|\le1\), and the outer representation in
Section~\ref{subsec:02-hardy} gives that \(O\) is outer with
\(O(0)=\exp(-\int_\T u\,dm)=e^Ta>0\).
In particular \(O\in H^\infty\subset H^2\).
Because the product in \eqref{eq:logarithmic-limit-factorization}
belongs to \(\D\), Lemma~\ref{lem:inner-division} gives
\(O\in\D\) and \(\|O\|_{\D}\le\|\eta\|_{\D}=1\).

The lower bound on the fixed neighborhood in Proposition~\ref{prop:common-nonzero-collar}
also passes to \(\eta\) at each fixed disk point of that neighborhood.
Since \(|O|\ge|\eta|\) in the disk, this yields
\begin{equation}\label{eq:cofactor-local-lower-bound}
   \liminf_{\DD\ni z\to\zeta_0}|O(z)|>0
          \qquad(\zeta_0\in\T\setminus\{\xi\}).
\end{equation}
Indeed, an arc through \(\zeta_0\) supplies one such neighborhood that contains
an entire disk neighborhood of \(\zeta_0\) intersected with \(\DD\).
The factorization above locates any singular mass at \(\xi\).
To identify the mass, we need the subspace generated by \(\eta\). We next show that polynomial multiples of \(O\) are dense
in \(\D\).

\subsection{Cyclicity of the outer factor}\label{subsec:09-cyclicity}

Fix a strong subsequential limit \(\eta\) as above. By
\eqref{eq:logarithmic-limit-factorization}, write
\(\eta=S_{T,\xi}O\) with \(T\ge0\) and \(O\) bounded, outer,
and in \(\D\). The cyclicity criterion below uses approach to a
boundary point from anywhere in \(\DD\). The lower bound from the common neighborhood gives exactly this control away from \(\xi\).

For \(h\) holomorphic on \(\DD\), define its \emph{unrestricted
approach-zero set} by
\begin{equation}\label{eq:unrestricted-approach-zeros}
 \Zapp(h)=\{\zeta\in\T:\liminf_{\DD\ni z\to\zeta}|h(z)|=0\}.
\end{equation}
The lower bound in \eqref{eq:cofactor-local-lower-bound} holds
throughout an inner neighborhood of each boundary point other
than \(\xi\), so \(\Zapp(O)\subset\{\xi\}\).
This is already enough for cyclicity. Hedenmalm and Shields proved
that every countable closed set is Bergman--Smirnov exceptional
\cite[Theorem~3, p.~102]{HS90}. In the terminology used here, Richter and Sundberg's
cyclicity criterion \cite[Corollary~3.3, p.~874]{RS94}
is stated in \cite[Theorem~2.2, p.~65]{EFKR06}:
an outer Dirichlet function is cyclic if its unrestricted
approach-zero set is Bergman--Smirnov exceptional.
Their combined countable-set theorem is stated in
\cite[Theorem~1.2]{EFKR06}. Apply it to \(O\): we proved \(O\in\D\)
and outerness before locating its approach zeros, and
\(\Zapp(O)\) is either empty or a singleton. Hence \([O]=\D\).
The limiting mass \(T\) is still unknown.

For comparison, the neighborhood criterion below gives a second
proof. Put \(F(z)=(1-\bar\xi z)O(z)\). The polynomial and \(O\)
are outer, and \(F\in\D\) by polynomial multiplication. Since
\(|O|\le1\), \(F\) tends to zero on every approach to \(\xi\).
At every other boundary point the polynomial is nonzero and the lower bound on this neighborhood bounds \(|O|\) away from zero. Thus
\begin{equation}\label{eq:exact-singleton-approach-zero}
                  \Zapp(F)=\{\xi\}.
\end{equation}

For this second argument we recall a sufficient condition of
El-Fallah, Kellay, and Ransford.
For a set \(E\subset\T\), let \(E^{(t)}\) be its \(t\)-neighborhood
in shortest arclength distance and let \(|E^{(t)}|\) be its
unnormalized arclength.

\begin{published}[Approach-zero cyclicity]\label{input:approach-cyclicity}
Let \(H\in\D\) be outer and put \(E=\Zapp(H)\). If
\[
 |E^{(t)}|=O(t^\sigma)\ \text{as }t\downarrow0
 \ \text{for some }\sigma>0,\qquad
 \int_0^\pi\frac{dt}{|E^{(t)}|}=\infty,
\]
then \(H\) is cyclic in \(\D\).
\end{published}

These are the hypotheses of \cite[Theorem~1.1, p.~2197]{EFKR09}, also
stated in \cite[Theorem~9.4.1]{Primer14}. For \(E=\{\xi\}\),
\(|E^{(t)}|=2t\) for \(0<t<\pi\). Thus the power bound holds with
\(\sigma=1\), and the integral diverges at zero.
Applying Theorem~\ref{input:approach-cyclicity} to \(F\)
gives \([F]=\D\). Every polynomial multiple of \(F\) is a
polynomial multiple of \(O\). Taking Dirichlet closures gives
\([F]\subset[O]\), so \(O\) is cyclic.

For the relation of this neighborhood condition to the earlier
Bergman--Smirnov exceptional-set method, see
\cite[Section~7, pp.~2212--2214]{EFKR09}.

\subsection{The limit subspace and the singular mass}\label{subsec:09-identification}

Cyclicity of the outer factor makes the inner factor determine the subspace generated by the limit by the following identity. The canonical value at zero can then be compared with the value for the prescribed mass.

\begin{published}[Invariant subspace generated by an inner--outer product]\label{input:module}
Suppose \(f\in\D\setminus\{0\}\) has Hardy factorization \(f=IO\)
with \(I\) inner and \(O\) outer. Then \(O\in\D\), and
\[
                         [f]=[O]\cap IH^2 .
\]
\end{published}

Richter and Sundberg give the invariant-subspace identity
in \cite[Theorem~5.3, p.~184]{RS92}; see also
\cite[Theorem~4, p.~135]{EFKR10}. The membership \(O\in\D\) also follows
directly from Lemma~\ref{lem:inner-division}.

Since $O$ is cyclic, $[O]=\D$, so the intersection identity gives
$[\eta]=\D\cap S_{T,\xi}H^2$. The positive approximation error excludes zero singular mass. Because the limit already has unit norm, canonical moments, and positive value at zero, it is the normalized canonical function of this subspace. Only at
that point can strict monotonicity of the extremal value at zero
identify the mass. Applying this identification to every strongly
convergent further subsequence then yields the conclusion for the original
sequence.

\begin{proof}[Completion of the proof of Theorem~\ref{thm:C}]
Let \(a=G_{\tau,\xi}(0)\in(0,1)\). From any subsequence take a
further subsequence converging strongly in \(\D\), as in
Proposition~\ref{prop:recovery-norm-retention}.
The strong limit \(\eta\) belongs to \(\D\), has norm one,
satisfies \(\eta(0)=a\), and obeys the canonical moments. Its
product \(\eta=S_{T,\xi}O\) in
\eqref{eq:logarithmic-limit-factorization} meets the hypothesis
of Theorem~\ref{input:module}. Since \(O\) is cyclic, we obtain
\begin{equation}\label{eq:atomic-limit-module}
 [\eta]=[O]\cap S_{T,\xi}H^2=\D\cap S_{T,\xi}H^2 .
\end{equation}
If \(T=0\), then \(S_{0,\xi}=1\) and \eqref{eq:atomic-limit-module}
would give \([\eta]=\D\). By \eqref{eq:certificate}, we have instead
\(\operatorname{dist}_{\D}(1,[\eta])^2=1-a^2>0\).
Therefore \(T>0\) and \([\eta]=M_{T,\xi}\).

Applying \eqref{eq:certificate} to \(M_{T,\xi}=[\eta]\) gives
\(P_{M_{T,\xi}}1=a\eta\).
Its norm is \(a>0\), so normalizing this projection gives
\(\eta=G_{T,\xi}\). Thus
\(G_{T,\xi}(0)=a=G_{\tau,\xi}(0)\).
Proposition~\ref{prop:atomic-origin-strict} forces \(T=\tau\).
The choice of the starting subsequence was arbitrary. Each one
therefore has a further subsequence converging strongly to the
same \(G_{\tau,\xi}\). A subsequence at distance at least some
\(\varepsilon>0\) from \(G_{\tau,\xi}\) would contradict that
conclusion. This proves convergence of the original sequence.
\end{proof}

\begin{corollary}[Outer factor and boundary zeros of the recovered target]
\label{cor:one-atom-cofactor}
Under the hypotheses of Theorem~\ref{thm:C}, write
\(G_{\tau,\xi}=S_{\tau,\xi}O\) with \(O(0)>0\). Then \(O\) is bounded
by one in modulus, is outer and cyclic in \(\D\), and multiplies
\(\D\) contractively. The recovered target satisfies
\begin{equation}\label{eq:one-atom-limit-zeros}
              \Zrad(G_{\tau,\xi})=\Zapp(G_{\tau,\xi})=\{\xi\}.
\end{equation}
\end{corollary}

\begin{proof}
The factorization and cyclicity were proved above. To obtain
the multiplier bound, first multiply by the canonical function
and then divide out its inner factor. For \(h\in\D\), the bound
on \(O\) gives \(Oh\in H^2\). By
\eqref{eq:canonical-multiplier}, \(G_{\tau,\xi}h\in\D\) and
\(\|G_{\tau,\xi}h\|_{\D}\le\|h\|_{\D}\).
Now \(G_{\tau,\xi}h=S_{\tau,\xi}(Oh)\), so
Lemma~\ref{lem:inner-division} gives
\begin{equation}\label{eq:one-atom-cofactor-contraction}
                  Oh\in\D,\qquad \|Oh\|_{\D}\le\|h\|_{\D}.
\end{equation}

Along the radius to \(\xi\), the factorization gives
\[
 |G_{\tau,\xi}(r\xi)|\le
           \exp\!\left(-\tau\frac{1+r}{1-r}\right)\longrightarrow0.
\]
Thus \(\xi\) is a radial zero and an unrestricted approach
zero. Apply Proposition~\ref{prop:fixed-analytic-neighborhood}
to the family now known to converge strongly to \(G_{\tau,\xi}\).
At each other boundary point it gives a holomorphic continuation
with nonzero boundary value.
Continuity gives a positive lower bound on a small disk about that
point, excluding every approach to zero from inside \(\DD\).
\end{proof}

The identity \eqref{eq:one-atom-limit-zeros} also follows from
\cite[Theorem~3.2]{RS94}. The function
\((1-\bar\xi z)S_{\tau,\xi}\in M_{\tau,\xi}\) from
Lemma~\ref{lem:atomic-witness} excludes approach zeros away
from \(\xi\), while the radial estimate in the proof above shows
that \(\xi\in\Zrad(G_{\tau,\xi})\).

The cyclic multiplier in Corollary~\ref{cor:one-atom-cofactor}
will be used in Section~\ref{sec:12} to form products of the outer factors associated with the individual atoms.

\paragraph{The two boundary measures.}
The \(L^1(m)\) estimate after \eqref{eq:source} gives
\(\ell_{\eta_j}\to\ell_{G_{\tau,\xi}}\) in \(L^1(m)\) by strong
Dirichlet convergence. For the logarithmic measures, the argument
of Section~\ref{subsec:09-logarithm} now applies to the original
strongly convergent sequence. It gives the weak-* limit
\(-\log|G_{\tau,\xi}^*|\,m+\tau\delta_\xi\).
Fix an open arc \(J\) around \(\xi\). A continuous cutoff supported
in \(J\), between zero and one and equal to one at \(\xi\), gives
\(\liminf_{j\to\infty}\int_J u_j\,dm\ge\tau\).
Arcs of arbitrarily small arclength retain this positive lower
bound, so the logarithmic densities fail uniform integrability.

The remaining task for one-atom recovery is geometric. We need
nonpolar Carleson sets shrinking to \(\xi\) whose canonical
values at zero converge to \(G_{\tau,\xi}(0)\).
Sections~\ref{sec:10}--\ref{sec:11} construct these sets with the
prescribed canonical value at zero and with capacities tending to
zero. Section~\ref{subsec:11-realization} shows that their canonical
functions extend continuously to the closed disk and vanish exactly
on the defining sets. Figure~\ref{fig:recovery-map} shows how the
geometric construction fits into the recovery proof.

\begin{figure}[H]\centering
\begin{tikzpicture}[>=Latex,
 b/.style={draw=black,fill=white,rounded corners=2pt,align=center,font=\small,text width=5.25cm,minimum height=1.1cm},
 arr/.style={->,semithick}]
\node[b] (n) at (-3.15,0) {Small-disk energy bound\\Common continuation};
\node[b] (g) at (3.15,0) {Trace estimate for Cantor copies\\and logarithmic mass};
\node[b] (l) at (-3.15,-1.65) {Strong compactness\\and factorization of the limit function};
\node[b] (t) at (3.15,-1.65) {Transport of the subspaces\\Continuity of the extremal value\\at zero for the Dirichlet norm};
\node[b] (an) at (-3.15,-3.4) {Cyclicity and mass identification\\Analytic criterion C\\
 Theorem~\ref{thm:C}};
\node[b] (ge) at (3.15,-3.4)
 {Prescribed value at zero\\
 Collapse and capacity decay\\
 Theorem~\ref{thm:G}};
\node[b] (r) at (0,-5.25) {One-atom recovery \\Corollary~\ref{cor:R}};
\node[b] (b) at (0,-6.95) {Norm approximation and\\weak limits in \(M_\mu\)\\Finite-atomic approximation B};
\draw[arr] (n)--(l);\draw[arr] (l)--(an);
\draw[arr] (g)--(t);\draw[arr] (t)--(ge);
\draw[arr] (an)--(r);\draw[arr] (ge)--(r);\draw[arr] (r)--(b);
\end{tikzpicture}
\caption{Criterion C applies when nonpolar Carleson sets satisfy
the stated collapse condition and the values at zero of their
canonical functions converge as required. The geometric construction
supplies such a family.}
\label{fig:recovery-map}\end{figure}
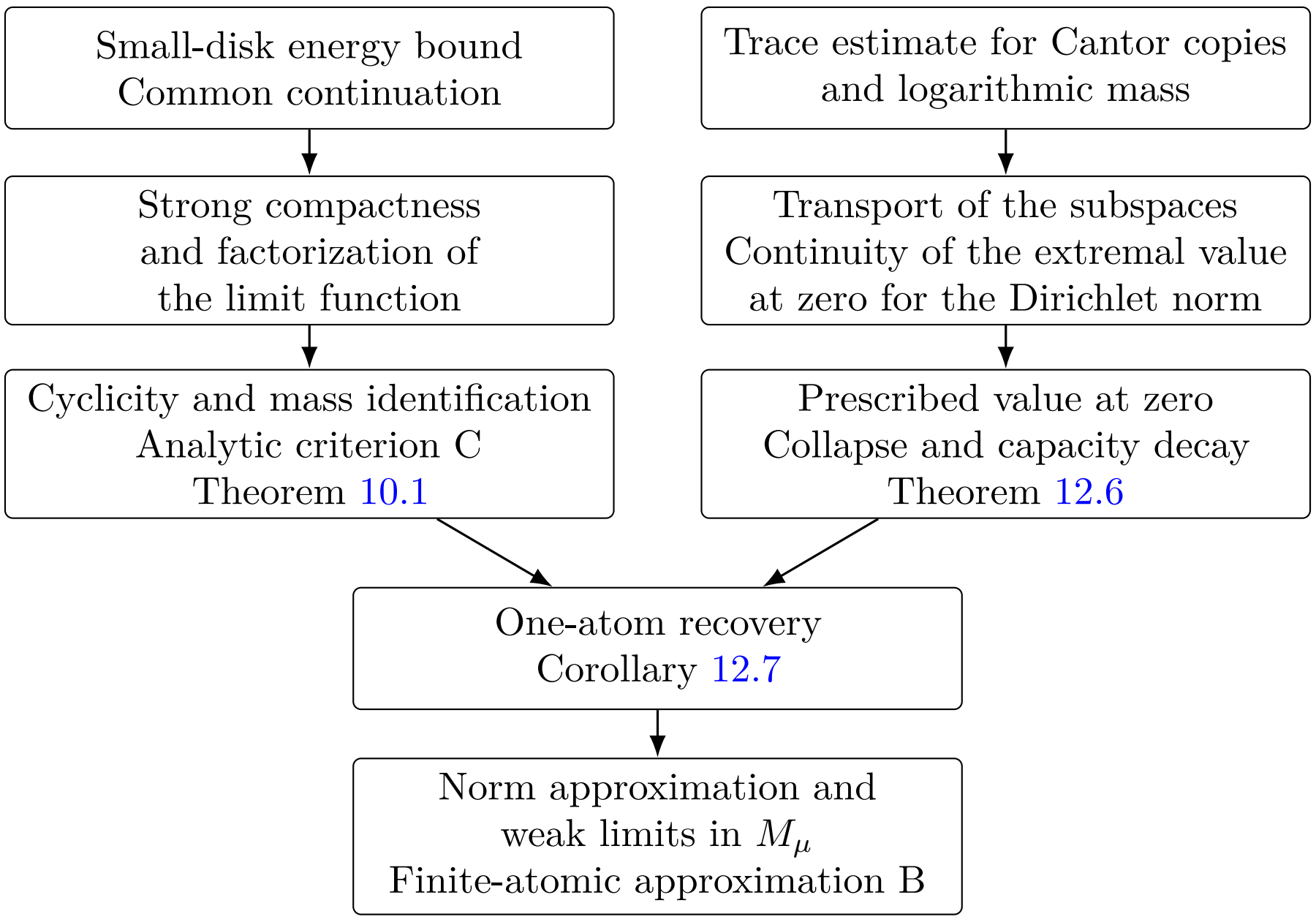\FloatBarrier

\newpage \section{Cantor copies on equal subarcs and logarithmic mass}\label{sec:10}

We construct a fixed compact set whose canonical function has enough
logarithmic mass in each of a sequence of short arcs. Section~\ref{sec:11}
will use that mass to make a canonical value at zero for a trace-zero space associated with an image of \(K\) under a Möbius transformation small. In each arc we place a Cantor copy in every equal
subarc. Vanishing quasi-everywhere on each copy gives zero integral
against its Cantor probability. A local energy estimate then controls
the squared trace averaged over the whole arc. Jensen's inequality
turns this average into the required logarithmic mass.
We adapt the reflection argument used in Kellay's proof of the
capacitary Poincar\'e inequality \cite{K11}. Summing the resulting estimates over all
subarcs gives the quantitative bound needed for the subsequent
selection.

\subsection{Cantor probabilities in the cells}\label{subsec:10-kernel}

A trace that vanishes quasi-everywhere on a Cantor copy vanishes
outside a polar exceptional set. To use that condition in an
integral, we choose finite-energy probabilities, which give zero
mass to polar sets. We need one on the original cell and one on
the semicircle used in the reflection argument.

Let \(C\subset[0,1]\) be the closed middle-third Cantor set,
including both endpoints. Let \(\nu_C\) be the probability assigning
mass \(2^{-n}\) to each level-\(n\) interval. Put
\(\alpha=\log 2/\log 3>1/2\). For \(0<r\le1\), choose \(n\ge0\)
with \(3^{-(n+1)}<r\le3^{-n}\). An interval of length \(2r\)
meets at most three level-\(n\) intervals, because their gaps
have length at least \(3^{-n}\). Hence
\begin{equation}\label{eq:packed-cell-growth}
 \nu_C(B(x,r))\le6r^\alpha.
\end{equation}
In particular \(\nu_C\) has no atoms. Integration of this ball bound
over logarithmic scales gives
\begin{equation}\label{eq:packed-cell-log-integrability}
 \iint_{[0,1]^2}\log\frac1{|x-y|}\,d\nu_C(x)d\nu_C(y)
 \le \int_0^\infty 6e^{-\alpha s}\,ds=\frac6\alpha.
\end{equation}

For a cell \(J=\{e^{i(t+lx)}:0\le x\le1\}\) of angular length
\(0<l\le1/4\), let \(\nu_J\) be the image of \(\nu_C\) under
\(x\mapsto e^{i(t+lx)}\). Let \(\nu_\pi\) be the image of \(\nu_C\) under
\(x\mapsto e^{i\pi x}\). The chord bounds
\[
 |1-e^{ilv}|\ge (2/\pi)l|v|,\qquad
 |1-e^{i\pi v}|\ge 2|v|\qquad(0<|v|\le1)
\]
and \eqref{eq:packed-cell-log-integrability} give finite
logarithmic energy for both probabilities. The comparison
\(I_R\le 2I_{\log}+1\) from Section~\ref{subsec:03-pairing}
therefore gives \(I_R(\nu_J)<\infty\) and
\(I_R(\nu_\pi)\le12/\alpha+1<27\).
Proposition~\ref{prop:trace-pairing} identifies their continuous
pairings on \(V\) with integration of the quasi-continuous trace.
The bound on \(I_R(\nu_\pi)\) is independent of the original
cell and will set the constant in the local estimate. The pairing
on the original cell uses only finiteness of \(I_R(\nu_J)\).

\subsection{A trace estimate for Cantor copies on equal subarcs}\label{subsec:10-cancellation}

For analytic traces, Kellay's capacitary Poincar\'e inequality
\cite[Theorem~2.2]{K11} already gives the factor \(L/N\), up to a
constant. Take his parameters \(\alpha=1,\ \beta=1/2\) on an open
cell \(J\) of length \(l=L/N\le1/4\). The Cantor probability on \(J\)
has \(1/2\)-Riesz energy at most \(C l^{-1/2}\) by
\eqref{eq:packed-cell-growth}, since \(\log 2/\log 3>1/2\).
It gives zero mass to the endpoints and to the polar exceptional set.
Thus the \(L^2\)-capacity \(C_{1/2,2}\) of the zero set in the open
cell is at least \(c l^{1/2}\). His inequality bounds the square of the cell mean of \(|h|\) by a constant times the local difference
energy. The variance of \(|h|\) has the same bound: in the double-integral formula for that variance use
\(\big||h(x)|-|h(y)|\big|\le|h(x)-h(y)|\) and \(|x-y|\le l\).
Multiplying the resulting mean-square bound by the cell length and summing gives the analytic estimate. Kellay also uses Jensen's inequality to pass to logarithmic mass \cite[proof of Theorem~3.1]{K11}.
The reflection proof below, following \cite[Lemma~2.1]{K11}, treats
arbitrary two-sided traces in \(V\) and gives the constant \(300\)
in the present normalization.

\begin{theorem}[A trace estimate for functions vanishing on Cantor copies]\label{thm:packed-trace}
Let \(I\) be a closed angular interval of length \(0<L\le1/4\), and
divide it into \(N\) equal closed cells \(J_k\) of length \(l=L/N\),
where \(N\) is a positive integer.
For the left endpoint \(t_k\) in angular coordinates, put the copy
\(C_k=\{e^{i(t_k+lx)}:x\in C\}\) in the image of \(J_k\)
under \(t\mapsto e^{it}\), including both endpoints of that arc.
If \(h\in V\) and its quasi-continuous representative vanishes
quasi-everywhere on
\(\bigcup_{k=1}^N C_k\), then
\begin{equation}\label{eq:packed-trace}
 \int_I|h|^2\,dm\le 300\,\frac LN\,\|h\|_V^2.
\end{equation}
Here the integrals on angular intervals use their images in \(\T\),
with \(dm=d\theta/(2\pi)\). The constant is uniform in the location
and length of \(I\) and in \(N\).
\end{theorem}

\begin{proof}
Let \(\mathcal E(h)\) be the double integral in the norm formula
\eqref{eq:V-difference}, and let \(\mathcal E_J(h)\) be its restriction
to \(J^2\). To use the fixed semicircle probability, rescale a cell \(J\)
to \([0,1]\), and extend its
trace by even reflection on \([-1,1]\) with the endpoints
identified. Choose an \(L^2(m)\) representative of \(h\), put
\(f(x)=h(e^{i(t+lx)})\) on \([0,1]\), and set
\[
 \mathcal J(f)=\iint_{[0,1]^2}
       \frac{|f(x)-f(y)|^2}{|x-y|^2}\,dx\,dy,\qquad
 U(e^{i\pi x})=f(|x|)\quad(-1\le x\le1).
\]
The bound \(|e^{ilx}-e^{ily}|\le l|x-y|\) and the two factors
\(l/(2\pi)\) from arclength give
\(\mathcal J(f)\le4\pi^2\mathcal E_J(h)\).

For \(x,y\in[0,1]\), the chord in either same-sign quadrant of the
reflected integral is at least \(2|x-y|\). In either opposite-sign
quadrant it is
\[
 2\sin(\pi(x+y)/2)\ge 2\min\{x+y,2-x-y\}\ge 2|x-y|.
\]
Each quadrant has product arclength measure \(dx\,dy/4\).
Consequently
\[
 \mathcal E(U)\le \mathcal J(f)/4,\qquad
 \|U\|_{L^2(m)}^2=\int_0^1|f(x)|^2\,dx
       \le (2\pi/l)\|h\|_{L^2(m)}^2.
\]
In particular restriction and reflection give a bounded map
\(V\to V\) for the fixed cell.

The definition of \(U\) depends only on the \(L^2(m)\) class of \(h\).
To transfer the vanishing condition for the quasi-continuous
trace, we check the two integrals on continuous approximants and
pass to the limit in \(V\). Choose trigonometric polynomials
\(h_j\to h\) in \(V\), and reflect their restrictions to obtain
\(U_j\to U\) in \(V\). Each \(U_j\) is continuous on the circle,
including at the seams of the reflection. Its integral against
\(\nu_\pi\) equals the integral of \(h_j\) against \(\nu_J\) by
the definitions of the two probabilities. Continuity of both
pairings gives
\[
 \int_\T\widetilde U\,d\nu_\pi
    =\int_\T\widetilde h\,d\nu_J=0.
\]
The last equality uses the vanishing condition on the Cantor
copy and the fact that \(\nu_J\) gives zero mass to polar sets.

The difference energy is invariant under adding a constant. We
use the zero pairing to control that constant part and the
energy to control the variation about it. Let
\(c=\int_0^1 f(x)\,dx=\int_\T U\,dm\).
The nonconstant Fourier modes give
\(\|U-c\|_V^2\le 2\mathcal E(U)\). Since \(\nu_\pi\) has mass
one and pairs to zero with \(U\), Proposition~\ref{prop:trace-pairing}
gives
\[
 |c|^2
 \le I_R(\nu_\pi)\|U-c\|_V^2
 \le54\mathcal E(U)\le\frac{27}{2}\mathcal J(f).
\]
The variance identity gives
\[
 \int_0^1|f-c|^2\,dx
   =\frac12\iint_{[0,1]^2}|f(x)-f(y)|^2\,dx\,dy
   \le\frac12\mathcal J(f).
\]
Thus \(\int_0^1|f|^2\,dx\le14\mathcal J(f)\). Rescaling to
\(J\) yields \(\int_J|h|^2\,dm\le28\pi l\,\mathcal E_J(h)\).
Adjacent cells meet in endpoints of arclength measure zero.
The cell squares therefore have total energy at most
\(\mathcal E(h)\), since the integrand defining \(\mathcal E(h)\) is
nonnegative. Summing the cell estimates leaves the common
factor \(l=L/N\) and gives \eqref{eq:packed-trace}, since
\(28\pi<300\).
\end{proof}

\subsection{One fixed nonpolar Carleson set}\label{subsec:10-compact}

Kellay uses Cantor copies on shrinking arcs in a construction of
uniqueness sets in weighted Dirichlet spaces \cite[Remark~3.3(2)]{K11}. Here the compactum we construct must satisfy
Carleson's integrability condition \cite[Section~I, no.~2]{Carleson52}.
The gaps left inside an arc of length \(L\) by Cantor copies
in \(N\) equal cells have entropy
of order \(L\log(N/L)\). Increasing \(N\) improves the trace
estimate. The entropy grows only logarithmically in \(N\), so
we can use very many cells in short arcs while keeping the sum
of the gap entropies finite. We choose the scales so that \(L_j\log N_j\) is of
order \(j^{-2}\). This scale will also set the lower bound
for the logarithmic mass.

For \(j\ge6\), set
\begin{equation}\label{eq:fixed-packed-scales}
 \begin{gathered}
 \theta_j=2^{-j},\qquad L_j=2^{-4j-10},\qquad
 I_j=[\theta_j-L_j/2,\theta_j+L_j/2],\\
 \omega_j=\frac{1}{j^2L_j},\qquad
 N_j=\lceil\exp(\omega_j)\rceil,\qquad l_j=\frac{L_j}{N_j}.
 \end{gathered}
\end{equation}
Divide the angular interval \(I_j\) into \(N_j\) equal closed cells.
In the \(k\)-th cell, for \(1\le k\le N_j\), put the image of the
Cantor set \(C\subset[0,1]\) under
\[
  x\longmapsto
  \exp\!\big(i[\theta_j-L_j/2+(k-1)l_j+l_j x]\big).
\]
Let \(E_j\) be the union of these images, and define
\begin{equation}\label{eq:fixed-packed-compact}
                     K=\{1\}\cup\bigcup_{j\ge6}E_j .
\end{equation}
When integrating against \(m\) on the circle, we also use \(I_j\) to denote its image under \(t\mapsto e^{it}\), so \(m(I_j)=L_j/(2\pi)\).
The set is specified by \eqref{eq:fixed-packed-scales} and
\eqref{eq:fixed-packed-compact} before choosing its canonical
function or prescribing a value at zero for the canonical functions
of its conformal images.

Here \(\theta_j\) locates the arc \(I_j\) and \(L_j\) gives its width;
the much smaller widths keep the arcs disjoint as they approach
\(1\). Inside \(I_j\), \(N_j\) is the number of cells and \(l_j\)
is the cell length appearing in the trace estimate. Each count is
finite, even though it grows very rapidly. The parameter
\(\omega_j\) prescribes the scale of \(\log N_j\): multiplying
by \(L_j\) places this contribution to entropy on the summable
\(j^{-2}\) scale. The same choice supplies the later
logarithmic-mass lower bound through Jensen's inequality. These parameters define the fixed set; the subsequent gap and radial-depth choices will be used when applying disk automorphisms to this set.

\begin{lemma}[Compactness, nonpolarity, and entropy]
\label{lem:fixed-packed-compact}
The set \(K\) in \eqref{eq:fixed-packed-compact} is a nonpolar
Carleson set with \(m(K)=0\).
\end{lemma}

\begin{proof}
We have \(L_j\le\theta_j/4\) and \(L_{j+1}=L_j/16\).
The gap between \(I_{j+1}\) and \(I_j\) has angular length
\[
 \gamma_j=\frac{\theta_j}{2}-\frac{L_j+L_{j+1}}{2},
 \qquad \frac{\theta_j}{4}\le \gamma_j\le\frac{\theta_j}{2}.
\]
Thus the intervals \(I_j\) are pairwise disjoint and approach angle zero.
Each \(E_j\) is compact and null. A convergent sequence of points
from infinitely many distinct arcs \(I_j\) can only accumulate at \(1\),
which belongs to \(K\). This proves compactness, and countable
subadditivity gives \(m(K)=0\). One of the Cantor copies in \(K\) supports a probability measure
of finite \(I_R\)-energy by Section~\ref{subsec:10-kernel}. Proposition~\ref{prop:trace-pairing}
says that this probability gives zero mass to polar sets. Since it
is supported on \(K\) and has mass one, \(K\) is nonpolar.

We check the Carleson condition by separating gaps inside the arcs \(I_j\),
gaps between these arcs, and the remaining gap around the circle.
In a cell of angular length \(l_j\), the complementary gaps of
angular length \(l_j3^{-q}\) have multiplicity \(2^{q-1}\), for
\(q\ge1\). The copies retain the endpoints of each cell, so
all gaps inside \(I_j\) are of this form. The geometric sums
\[
 \sum_{q\ge1}2^{q-1}3^{-q}=1,\qquad
 \sum_{q\ge1}q\,2^{q-1}3^{-q}=3
\]
therefore give the exact normalized entropy of the gaps inside \(I_j\):
\begin{equation}\label{eq:fixed-packed-entropy}
 \sum_{\substack{J\text{ complementary}\\J\subset I_j}}
   m(J)\log\frac1{m(J)}
    =\frac{L_j}{2\pi}
       \left(\log\frac{2\pi N_j}{L_j}+3\log3\right).
\end{equation}
Here the factor \(2\pi\) inside the logarithm comes from
\(m(J)=|J|/(2\pi)\), with \(|J|\) its angular length.

The ceiling gives \(\omega_j\le\log N_j\le\omega_j+\log2\).
The contribution of \(\log N_j\) to the entropy is summable
because \(L_j\omega_j=j^{-2}\). The other terms depend only on
the lengths of the short arcs. Thus the total entropy of the gaps inside
these arcs is finite, bounded by a constant times
\(\sum_{j\ge6}[j^{-2}+(j+1)2^{-4j}]\).
The bounds for \(\gamma_j\) give
\(\gamma_j\log(2\pi/\gamma_j)\le C(j+1)2^{-j}\), which is also summable.
The remaining gap goes from the top of \(I_6\) to angle zero
around the other side of the circle and contributes one finite
term. The gap-entropy characterization in Section~\ref{subsec:03-ideal}
now gives
\(\log(2/\dist(\,\cdot\,,K))\in L^1(m)\).
\end{proof}

Section~\ref{sec:11} compares canonical values along a radius ending
in the gap next to each arc. At the starting point the Poisson
kernel must be large enough on the arc for the logarithmic mass
to give an interior bound. The depth of the starting point also
sets the maximum depth for the selected parameter, which controls the location of the resulting images of \(K\). For
\(j\ge6\) put
\begin{equation}\label{eq:fixed-gap-geometry}
 \begin{gathered}
 g_j=\sqrt{L_j}=2^{-2j-5},\qquad
 \alpha_j=\theta_j+\frac32g_j,\qquad \zeta_j=e^{i\alpha_j},\\
 s_j=g_j^{3/2},\qquad b_j^0=(1-s_j)\zeta_j .
 \end{gathered}
\end{equation}
The point \(\zeta_j\) lies in the open gap immediately above \(I_j\).
Indeed, \(L_j/2<3g_j/2<\theta_j/2\). If \(j>6\), the lower
endpoint of \(I_{j-1}\) is at least \(3\theta_j/2\). For \(j=6\)
the gap above \(I_6\) is the remaining large gap. The choice $s_j=g_j^{3/2}$ gives $s_j/g_j^2=g_j^{-1/2}\to\infty$ and $s_j/g_j=g_j^{1/2}\to0$.
The first gives a large lower bound for the Poisson kernel on $I_j$ at $b_j^0$; the second will be used to control the images of $K$ under disk automorphisms.

\subsection{Logarithmic mass in the chosen arcs}\label{subsec:10-mass}

The fixed set now determines a single canonical function.
Its outer representation expresses its interior
modulus through the boundary logarithm. We estimate the mass of
that logarithm on each arc \(I_j\).

Let \(M=\D_K\) and \(e=G_M\). By Section~\ref{subsec:03-ideal},
\(M\) is closed and invariant, and it contains a smooth outer function
with a nonzero value at the origin. Nonpolarity excludes the constant
\(1\) from \(M\), so the canonical function has the properties proved in
Section~\ref{subsec:03-seed}:
\begin{equation}\label{eq:fixed-canonical-origin}
 [e]=M,\qquad \|e\|_{\D}=1,\qquad e\text{ is outer},\qquad
 a_0=e(0)\in(0,1),\qquad \kappa=1-a_0^2>0 .
\end{equation}
The canonical bound from Section~\ref{subsec:02-multiplier} gives
\(|e|\le1\) in the disk. Thus the outer representation
\eqref{eq:outer-intro} supplies a nonnegative integrable density
\(u=-\log|e^*|\) with
\begin{equation}\label{eq:fixed-logarithmic-density}
 \int_\T u\,dm=-\log a_0,\qquad
 -\log|e(b)|=\int_\T P_b u\,dm\qquad(b\in\DD).
\end{equation}
The function \(e\), its value \(a_0=e(0)\), the gap \(\kappa\),
and the density \(u\) are fixed by \(K\) and do not depend on \(j\).
Section~\ref{sec:11} will prescribe the value at zero of the canonical function for a suitable image of \(K\) separately.

Kellay uses Jensen's inequality to pass from a local trace estimate
to a logarithmic integral \cite[proof of Theorem~3.1]{K11}. The choice
of \(N_j\) gives a quantitative lower bound on the mass of \(u\) in
each arc \(I_j\).

\begin{lemma}[Logarithmic mass on each arc]
\label{lem:fixed-block-logarithmic-mass}
For \(e=G_{\D_K}\) and \(u=-\log|e^*|\) as above, every \(j\ge6\)
satisfies
\begin{equation}\label{eq:fixed-block-logarithmic-mass}
 \int_{I_j}|e^*|^2\,dm\le300\,\frac{L_j}{N_j},
 \qquad
 \int_{I_j}u\,dm\ge\frac{1}{8\pi j^2}.
\end{equation}
\end{lemma}

\begin{proof}
The \(V\)-norm of the trace \(e^*\) equals \(\|e\|_{\D}=1\).
Its quasi-continuous representative vanishes on \(E_j\) outside a
polar set, because \(e\in\D_K\). Since \(L_j\le1/4\), the trace estimate \eqref{eq:packed-trace} gives the first inequality in
\eqref{eq:fixed-block-logarithmic-mass}. The normalized average over \(I_j\)
is therefore at most \(600\pi/N_j\). To turn this upper bound into a lower bound on logarithmic mass,
apply Jensen's inequality for the logarithm to \(|e^*|^2\) and the probability \(dm/m(I_j)\).
Using \(u=-\tfrac12\log|e^*|^2\) gives
\[
 \int_{I_j}u\,dm
 \ge -\frac{m(I_j)}{2}
       \log\!\left(\frac1{m(I_j)}\int_{I_j}|e^*|^2\,dm\right)
 \ge\frac{L_j}{4\pi}\bigl(\log N_j-\log(600\pi)\bigr).
\]
For \(j\ge6\) the ratio \(\omega_{j+1}/\omega_j=16(j/(j+1))^2\)
is at least \(4\). Also \(\omega_6=2^{34}/36>2\log(600\pi)\).
Hence \(\log N_j-\log(600\pi)\ge\omega_j/2\) for every \(j\ge6\).
Substituting \(L_j\omega_j=j^{-2}\) proves the second inequality.
\end{proof}

Section~\ref{sec:11} combines the mass on \(I_j\) with a
Poisson-kernel lower bound at \(b_j^0\). This makes
\(|e(b_j^0)|\) small, and a comparison of evaluation norms then gives a small canonical value at zero for the resulting trace-zero space.

\newpage \section{Prescribing extremal values at zero via Möbius maps}
\label{sec:11}

To apply Theorem~\ref{thm:C}, we need nonpolar Carleson sets
that shrink to the chosen point and whose canonical values at zero
approach the prescribed one-atom value. Section~\ref{sec:10} supplies one compact set with logarithmic mass in every selected arc. We use automorphic images of this fixed set to attain the prescribed value exactly.

First we express the squared canonical value at zero for each image
by an extremal problem on the original trace-zero space. The mass on \(I_j\)
makes this squared value small at the starting points as the
ray index grows. On each fixed ray it tends to one as the parameter moves outward toward the circle. Continuity gives a point with the prescribed value, and the last step controls the location and capacity of the entire image at that point.

\subsection{Trace-zero spaces under Möbius transformations}
\label{subsec:11-ideals}

Let \(K\) be the set defined in
\eqref{eq:fixed-packed-compact}. To
describe the trace-zero space of an image of \(K\), we track two
effects of composition. The change of variables gives the Dirichlet norm
of the composed function. The capacity comparison shows that images of polar exceptional sets are again polar.

\paragraph{Boundary change of variables.}
For \(b\in\DD\), put
\begin{equation}\label{eq:automorphism-definition}
 \varphi_b(z)=\frac{b-z}{1-\overline b z},\qquad
 J_b(\zeta)=P_b(\zeta)=\frac{1-|b|^2}{|1-\overline b\zeta|^2},
 \qquad L_b=\frac{1+|b|}{1-|b|}\quad(\zeta\in\T).
\end{equation}
The identity
\[
 1-|\varphi_b(z)|^2
   =\frac{(1-|b|^2)(1-|z|^2)}{|1-\overline b z|^2}
\]
shows that \(\varphi_b\) maps the disk into itself and its boundary
to itself. Substitution gives \(\varphi_b\circ\varphi_b=\mathrm{id}\)
and \(\varphi_b(0)=b\). The boundary map is a smooth diffeomorphism.
Its arclength Jacobian is \(J_b=|\varphi_b'|\) on \(\T\), so
\(dm(\varphi_b(\zeta))=J_b(\zeta)\,dm(\zeta)\). The triangle
inequality in the denominator gives \(L_b^{-1}\le J_b\le L_b\).
Subtracting the two fractions in \eqref{eq:automorphism-definition}
gives the chord identity
\begin{equation}\label{eq:automorphism-chord}
 |\varphi_b(\zeta)-\varphi_b(\omega)|^2
      =J_b(\zeta)J_b(\omega)|\zeta-\omega|^2 .
\end{equation}
Thus distances change by factors between \(L_b^{-1}\) and \(L_b\),
and the boundary map preserves arclength-null sets in both directions.

The chord identity accounts for the difference term of the boundary
norm. Let \(\mathcal E(u)\) be the double integral in
\eqref{eq:V-difference}. In that integral for \(u\circ\varphi_b\),
set \(\zeta'=\varphi_b(\zeta)\) and \(\omega'=\varphi_b(\omega)\).
The two Jacobians cancel the factors from
\eqref{eq:automorphism-chord}. This gives
\(\mathcal E(u\circ\varphi_b)=\mathcal E(u)\), with the nonnegative
integrals well defined also before finiteness is known. For the
\(L^2\) term, the inverse map is \(\varphi_b\), so
\[
 \|u\circ\varphi_b\|_{L^2(m)}^2=\int_\T |u|^2 J_b\,dm .
\]
Together these formulas give the bounds on the whole space \(V\)
\begin{equation}\label{eq:transport-boundary-norm}
 L_b^{-1}\|u\|_V^2\le\|u\circ\varphi_b\|_V^2\le L_b\|u\|_V^2
 \qquad(u\in V).
\end{equation}

\paragraph{Capacity and boundary representatives.}
If \(v\in V\) is a real admissible obstacle for an open set \(U\),
then \(v\circ\varphi_b\) is at least one arclength-almost everywhere
on \(\varphi_b(U)\). The set \(\varphi_b(U)\) is open.
Taking the infimum of \eqref{eq:transport-boundary-norm} over
obstacles for \(U\) gives
\(C_R(\varphi_b(U))\le L_b C_R(U)\). Taking the infimum over open
\(U\supset E\) gives the same upper bound for any set \(E\subset\T\).
Applying that upper bound to \(\varphi_b(E)\) gives the reverse
comparison:
\begin{equation}\label{eq:transport-capacity}
 L_b^{-1}C_R(E)\le C_R(\varphi_b(E))\le L_b C_R(E).
\end{equation}
In particular, \(\varphi_b\) maps polar sets to polar sets. This lets us
identify the representative used in the trace-zero condition.

For \(u\in V\), arclength-null preservation shows that
\(\widetilde u\circ\varphi_b\) represents the class
\(u\circ\varphi_b\) almost everywhere. The capacity comparison
gives quasi-continuity as well. Given \(\varepsilon>0\), choose an open
\(U\) of capacity less than \(\varepsilon/L_b\) on whose complement
\(\widetilde u\) is continuous. Then \(\widetilde u\circ\varphi_b\) is continuous on the complement
of \(\varphi_b(U)\). By \eqref{eq:transport-capacity},
\(C_R(\varphi_b(U))\le L_b C_R(U)<\varepsilon\). 
It is therefore quasi-continuous. Uniqueness of the representative
from Theorem~\ref{input:trace} gives
\begin{equation}\label{eq:transport-qc-trace}
       \widetilde{u\circ\varphi_b}
            =\widetilde u\circ\varphi_b\qquad\text{q.e. on }\T .
\end{equation}

\paragraph{Analytic composition and its Dirichlet norm.}
Composition by a disk automorphism preserves Dirichlet membership
\cite[Lemma~10, p.~296]{BS84}. For the extremal problem we need the identity
for the norm that includes the Hardy term. For analytic \(f\), the
chain rule and area change of variables
give \(D(f\circ\varphi_b)=D(f)\). For a polynomial \(p\),
the composition extends through the closed disk. Boundary change
of variables gives
\[
 \|p\circ\varphi_b\|_{H^2}^2
      =\int_\T P_b|p|^2\,dm .
\]
To pass to \(f\in\D\), take polynomials \(p_n\to f\) in \(\D\).
Applied to \(p_n-p_m\), the two identities above and \(P_b\le L_b\)
show that \(p_n\circ\varphi_b\) is Cauchy in \(\D\).
The norm limit agrees with \(f\circ\varphi_b\) on compact subdisks,
since evaluation is bounded and \(\varphi_b\) maps compact subdisks
to compact subsets of \(\DD\). On the circle, the \(L^2(m)\) limit
of the composed polynomial traces is \(f^*\circ\varphi_b\) by the
same boundary change of variables. Thus the composition
\(f\circ\varphi_b\) belongs to \(\D\), and its Hardy trace satisfies
\begin{equation}\label{eq:transport-hardy}
 (f\circ\varphi_b)^*=f^*\circ\varphi_b\quad\text{in }L^2(m),
 \qquad
 \|f\circ\varphi_b\|_{H^2}^2=\int_\T P_b|f^*|^2\,dm .
\end{equation}
In particular,
\begin{equation}\label{eq:transport-full-norm}
 \begin{gathered}
 \|f\circ\varphi_b\|_\D^2=D(f)+\int_\T P_b|f^*|^2\,dm ,\\
 L_b^{-1}\|f\|_\D^2\le\|f\circ\varphi_b\|_\D^2
                              \le L_b\|f\|_\D^2 .
 \end{gathered}
\end{equation}
Thus the energy remains fixed under composition, and the Hardy term
carries the dependence on \(b\). For instance, take \(f(z)=1+z\) and
real \(0<b<1\). The Poisson
integral of \(|1+\zeta|^2=2+2\operatorname{Re}\zeta\) at \(b\)
is \(2+2b\), while \(D(f)=1\). Hence
\begin{equation}\label{eq:transport-norm-example}
 \|f\circ\varphi_b\|_\D^2=3+2b,\qquad \|f\|_\D^2=3.
\end{equation}
The extremal problem in the next subsection keeps this Hardy term.

\paragraph{Trace-zero spaces under disk automorphisms.}
We can now impose the boundary condition on the composed analytic
function. The Hardy trace formula \eqref{eq:transport-hardy} and the
representative identity \eqref{eq:transport-qc-trace} give the
quasi-continuous trace \(\widetilde f\circ\varphi_b\) q.e. for
\(f\circ\varphi_b\). Applying \eqref{eq:radial-trace-agreement} to
this composed function identifies its  complex radial limit
with that trace q.e.

Let \(K_0\subset\T\) be compact. If \(f\in\D_{K_0}\), the
trace of \(f\) vanishes on \(K_0\) outside a polar set. Its image
under \(\varphi_b\) is polar by \eqref{eq:transport-capacity}.
The trace identity therefore gives \(f\circ\varphi_b\in
\D_{\varphi_b(K_0)}\). Conversely, if \(g\in\D_{\varphi_b(K_0)}\),
the same implication applied to \(\varphi_b(K_0)\) gives
\(g\circ\varphi_b\in\D_{K_0}\). Composing once more recovers \(g\).
Hence
\begin{equation}\label{eq:transport-ideals}
 \D_{\varphi_b(K_0)}
       =\{f\circ\varphi_b:f\in\D_{K_0}\}.
\end{equation}

If \(K_0\) is a nonpolar Carleson set, its image is compact and arclength-null, and is nonpolar by
\eqref{eq:transport-capacity}. The lower chord bound gives
\(\dist(\varphi_b(\zeta),\varphi_b(K_0))\ge
L_b^{-1}\dist(\zeta,K_0)\). Changing variables on the circle
then gives
\[
 \int_\T\log\frac2{\dist(\eta,\varphi_b(K_0))}\,dm(\eta)
 \le L_b\int_\T
       \left(\log\frac2{\dist(\zeta,K_0)}+\log L_b\right)dm(\zeta)<\infty .
\]
Thus the image is Carleson. The bounds above hold with a common
constant \((1+r)/(1-r)\) on each parameter disk \(|b|\le r<1\).
The constant depends on the chosen disk. For parameters approaching
the circle, Section~\ref{subsec:11-location} will estimate the
capacity of the image by an explicit obstacle.

By Section~\ref{subsec:03-ideal}, \(\D_{\varphi_b(K)}\) is closed and shift-invariant. To use the canonical projection, we also
need nonzero evaluation at zero. For the fixed compact set \(K\),
\eqref{eq:transport-ideals} puts \(F_K\circ\varphi_b\) in
\(\D_{\varphi_b(K)}\).
Here \(F_K\) is the smooth outer function from
Section~\ref{subsec:03-ideal}, so the value \(F_K(b)\) at zero
after composition is nonzero. Evaluation at zero on \(\D_{\varphi_b(K)}\) is therefore nonzero, and we can put
\begin{equation}\label{eq:transport-canonical-level}
 M_b=\D_{\varphi_b(K)},\qquad
 \eta_b=G_{M_b}=\frac{P_{M_b}1}{\|P_{M_b}1\|_\D},\qquad
 c(b)=\eta_b(0)^2=\|P_{M_b}1\|_\D^2.
\end{equation}
The image is a nonpolar Carleson set, so
Section~\ref{subsec:03-seed} gives
\([\eta_b]=M_b\), \(\eta_b\) outer, and \(0<c(b)<1\).
The value \(c(b)\) is therefore obtained from the orthogonal
projection in the Dirichlet norm.

\subsection{Continuity and bounds for the extremal value at zero}\label{subsec:11-continuity}

Composition with \(\varphi_b\) allows us to formulate all the extremal problems on the fixed space
\(M=\D_K\). Both evaluation and the Hardy weight vary with \(b\).
To compare the resulting suprema, we first choose a common norm
on \(M\) using its derivative energy.

Keep \(e\), \(a_0\), and \(\kappa=1-a_0^2>0\) from
\eqref{eq:fixed-canonical-origin}. The projection formula for \(e\)
gives \(|f(0)|\le a_0\|f\|_{\D}\) on \(M\). The strict bound
\(a_0<1\) lets us absorb the constant term when estimating the
Dirichlet norm. The coefficient
formula for the energy gives \(\|f\|_{H^2}^2\le |f(0)|^2+D(f)\).
Adding \(D(f)\) and applying the evaluation bound at zero yields
\begin{equation}\label{eq:transport-energy-equivalence}
 \frac{\kappa}{2}\|f\|_{\D}^2\le D(f)\le\|f\|_{\D}^2
 \qquad(f\in M).
\end{equation}
Since \(M\) is closed in \(\D\), it is a Hilbert space for the
equivalent norm \(D(f)^{1/2}\). Write \(\mathcal H\) for this
space, with inner product
\(\langle f,g\rangle_D=\int_\DD f'\overline{g'}\,dA_0\), linear in
the first argument.

\begin{proposition}[Continuity of the extremal value at zero]
\label{prop:transport-origin-continuity}
The function \(c(b)\) in
\eqref{eq:transport-canonical-level} is continuous on \(\DD\).
The value \(\eta_b(0)=\sqrt{c(b)}\) is
therefore continuous there as well.
\end{proposition}

\begin{proof}
For \(b\in\DD\), define
\[
 q_b(f,g)=\langle f,g\rangle_D+\int_\T P_b f^*\overline{g^*}\,dm.
\]
By \eqref{eq:transport-full-norm},
\(q_b(f,f)=\|f\circ\varphi_b\|_\D^2\). The squared norm of evaluation
at zero on \(M_b\) is \(\|P_{M_b}1\|_{\D}^2=c(b)\), since
\(1\) represents evaluation at zero on \(\D\). Under the
bijection \eqref{eq:transport-ideals}, the value of \(f\circ\varphi_b\) at zero is \(f(b)\). Thus the exact extremal problem is
\begin{equation}\label{eq:transport-c-variational}
 c(b)=\sup_{0\ne f\in M}\frac{|f(b)|^2}{q_b(f,f)}
     =\sup_{0\ne f\in M}
       \frac{|f(b)|^2}{D(f)+\int_\T P_b|f^*|^2\,dm}.
\end{equation}

We control the numerator on parameter disks strictly inside \(\DD\).
For \(b\in\DD\), evaluation \(f\mapsto f(b)\) is bounded on
\(\mathcal H\) by \eqref{eq:transport-energy-equivalence} and
the evaluation bound of Section~\ref{subsec:02-norm}.
Fix \(0<r<1\). For \(f(z)=\sum_{n\ge0}f_n z^n\) and \(|b|,|b'|\le r\),
\[
 |f(b)-f(b')|
 \le D(f)^{1/2}
       \left(\sum_{n\ge1}\frac{|b^n-(b')^n|^2}{n}\right)^{1/2}
 \le \frac{|b-b'|}{1-r^2}\,D(f)^{1/2}.
\]
The last inequality follows from
\(|b^n-(b')^n|\le n r^{n-1}|b-b'|\) and
\(\sum_{n\ge1}n r^{2n-2}=(1-r^2)^{-2}\). Thus evaluation is
continuous in the dual norm. By \eqref{eq:transport-energy-equivalence}
and the kernel diagonal of Section~\ref{subsec:02-norm},
\(C_r=\sup_{|b|\le r}\|\operatorname{ev}_b\|_{\mathcal H^*}
\le(2K(r)/\kappa)^{1/2}\), where \(\operatorname{ev}_b(f)=f(b)\).

The quotient in \eqref{eq:transport-c-variational} is homogeneous,
so every parameter uses the same sphere \(D(f)=1\). The denominator \(q_b(f,f)\) includes the Hardy term and is at least one on this sphere because \(P_b\ge0\). The energy equivalence bounds the difference of
denominators by
\[
 |q_b(f,f)-q_{b'}(f,f)|
 \le \frac{2}{\kappa}\|P_b-P_{b'}\|_\infty .
\]
The Poisson kernels depend continuously on \(b\) in the uniform norm
on \(\T\), by joint continuity on \(\{|b|\le r\}\times\T\).
The two numerators are at most \(C_r^2\) and their difference is
at most \(2C_r|b-b'|/(1-r^2)\). Since both denominators are at
least one, the difference of quotients satisfies
\[
 \left|\frac{|f(b)|^2}{q_b(f,f)}
      -\frac{|f(b')|^2}{q_{b'}(f,f)}\right|
 \le \frac{2C_r|b-b'|}{1-r^2}
       +\frac{2C_r^2}{\kappa}\|P_b-P_{b'}\|_\infty .
\]
This bound holds for every \(f\) on the fixed sphere and all
\(b,b'\) on the fixed parameter disk. Since the bound is independent of the function on the sphere,
taking suprema proves continuity of \(c\) on every such parameter disk, hence on \(\DD\).
\end{proof}

\paragraph{Evaluation bounds for the endpoints.}
The quantity \(c(b)\) is the squared canonical value at zero in
the Dirichlet norm, including its Hardy term. For the endpoint estimates
we first bound an auxiliary quantity \(B(b)\), defined below as
the squared norm of evaluation in the derivative-energy Hilbert
space. The two-sided comparison that follows transfers those
bounds to the extremal problem defining \(c(b)\) while retaining its
Hardy contribution. Let \(v_b\) be the Riesz vector for
evaluation on \(\mathcal H\), so that \(f(b)=\langle f,v_b\rangle_D\).
Set
\[
 B(b)=\sup_{0\ne f\in M}\frac{|f(b)|^2}{D(f)}
      =\|\operatorname{ev}_b\|_{\mathcal H^*}^2=\|v_b\|_{\mathcal H}^2 .
\]
The proof of Proposition~\ref{prop:transport-origin-continuity}
shows that the evaluation functionals are bounded and depend
continuously on \(b\) in \(\mathcal H^*\). Thus \(B\) is finite
and continuous. It is positive because \(e\) is nonvanishing
in \(\DD\).
For \(f\in M\)
with \(f(b)=1\), Cauchy--Schwarz gives \(D(f)\ge B(b)^{-1}\),
with equality at \(f=v_b/B(b)\). The identity \eqref{eq:transport-ideals} and invariance of
\(D\) give the same minimum of \(D(g)\) over \(g\in M_b\) with
\(g(0)=1\). This is the homogeneous extremal problem studied by
Carleson
\cite[Section~II, nos.~6--7]{Carleson52}. In the Dirichlet norm, the Hardy
contribution includes the same value \(f(b)\) as the numerator of the
extremal quotient. Separate that constant coefficient from the others.
Writing
\[
 Q_f(b)=\|f\circ\varphi_b\|_{H^2}^2-|f(b)|^2,
\]
the nonconstant coefficients give
\(0\le Q_f(b)\le D(f\circ\varphi_b)=D(f)\). Hence
\[
 c(b)=\sup_{0\ne f\in M}
          \frac{|f(b)|^2}{|f(b)|^2+D(f)+Q_f(b)} .
\]
For each such \(f\), put \(t=|f(b)|^2/D(f)\).
The last quotient lies between \(t/(t+2)\) and \(t/(t+1)\).
Both expressions are increasing and continuous in \(t\), so
taking suprema gives
\begin{equation}\label{eq:transport-cB}
 \frac{B(b)}{B(b)+2}\le c(b)\le \frac{B(b)}{B(b)+1}.
\end{equation}

For a lower bound on \(B(b)\) we test evaluation with a kernel times
the canonical generator of \(M\). Let \(k_b\) be the Dirichlet
evaluation kernel of
Section~\ref{subsec:02-norm}. Its squared norm is
\(k_b(b)=K(|b|)\), where \(K(0)=1\) and
\(K(r)=-\log(1-r^2)/r^2\) for \(0<r<1\).
Contractive multiplication by \(e\) in \eqref{eq:canonical-multiplier}
gives \(ek_b\in\D\) and
\(\|ek_b\|_{\D}^2\le k_b(b)\). Approximating \(k_b\) by polynomials
puts \(ek_b\) in \([e]=M\). The product is nonzero, so its energy is
positive by \eqref{eq:transport-energy-equivalence}. Using it as a
test function in \(B(b)\) gives
\[
 B(b)\ge\frac{|e(b)k_b(b)|^2}{D(ek_b)}
            \ge |e(b)|^2 k_b(b).
\]

For an upper bound we control evaluation of every function in \(M\)
through its Hardy quotient by the generator. Take \(f\in M\). By
Lemma~\ref{lem:completed-canonical-quotient}, the quotient
\(F=f/e\) belongs to \(H^2\), with
\(\|F\|_{H^2}\le\|f\|_\D\).
The Hardy evaluation bound then gives
\[
 |f(b)|^2\le
     \frac{|e(b)|^2}{1-|b|^2}\|f\|_\D^2
     \le \frac{2|e(b)|^2}{\kappa(1-|b|^2)}D(f).
\]
Combining the estimates yields
\begin{equation}\label{eq:transport-B-bounds}
 |e(b)|^2 k_b(b)\le B(b)\le
          \frac{2|e(b)|^2}{\kappa(1-|b|^2)} .
\end{equation}
At the chosen interior starting points, the logarithmic mass on \(I_j\) makes
\(|e(b)|^2\) small enough to overcome the factor
\((1-|b|^2)^{-1}\) in the upper bound. Along each fixed ray toward a
nonzero boundary value of \(e\), the lower bound uses
\(K(r)\to\infty\) as \(r\uparrow1\). The comparison with \(c\)
in \eqref{eq:transport-cB} transfers these estimates to the
canonical value.

\subsection{Two endpoint limits and exact selection}\label{subsec:11-endpoints}

Use the points \(\zeta_j\) and depths \(s_j\) from
\eqref{eq:fixed-gap-geometry}. We seek two values of the function \(c(b)\) that bracket the prescribed level \(a^2\).
The logarithmic mass in \(I_j\) makes \(e\) small at \(b_j^0\).
The upper evaluation estimate then gives \(c(b_j^0)\to0\) as
\(j\to\infty\). With a sufficiently large index \(j\) fixed, letting \(r\uparrow1\) along that ray gives \(c(r\zeta_j)\to1\). This second limit uses
continuation and nonvanishing of \(e\) at \(\zeta_j\).

\begin{lemma}[Separation and the two endpoints]\label{lem:transport-endpoints}
Use the fixed compact set \(K\), the function \(e=G_{\D_K}\), and its
gap \(\kappa=1-e(0)^2>0\) from Section~\ref{sec:10}.
For every \(j\ge6\), each angular representative
\(t\in\{0\}\cup\bigcup_{k\ge6} I_k\) satisfies
\begin{equation}\label{eq:fixed-gap-separation}
 |t-\alpha_j|\ge g_j,\qquad
 \dist(\zeta_j,K)\ge\frac{2g_j}{\pi}.
\end{equation}
The starting points \(b_j^0=(1-s_j)\zeta_j\) satisfy
\begin{align}
 -\log|e(b_j^0)|
 &\ge \frac{1}{128\pi j^2\sqrt{g_j}},
 \label{eq:transport-low-log}\\
 c(b_j^0)
 &\le \frac{2}{\kappa}g_j^{-3/2}
       \exp\!\left(-\frac{1}{64\pi j^2\sqrt{g_j}}\right)
       \longrightarrow0\qquad (j\to\infty).
 \label{eq:transport-low-endpoint}
\end{align}
For each fixed \(j\ge6\),
\begin{equation}\label{eq:transport-high-endpoint}
                 \lim_{r\uparrow1}c(r\zeta_j)=1.
\end{equation}
\end{lemma}

\begin{proof}
\emph{Separation from all the arcs.}
For \(t\in I_j\) we have
\(\alpha_j-t\ge3g_j/2-g_j^2/2\ge g_j\), since \(L_j=g_j^2\)
and \(g_j\le1\). For every \(k\ge6\), \(L_k/2\le\theta_k/4\).
Thus \(t\in I_k\) with \(k>j\) gives \(t\le5\theta_j/8\),
while \(k<j\) gives \(t\ge3\theta_j/2\).
Since \(\theta_j/g_j=2^{j+5}\), both comparisons give
\(|t-\alpha_j|\ge g_j\). The angle zero also satisfies this
inequality. All these angles lie in an interval of length less
than one, so the chord estimate
\[
 |e^{it}-e^{i\alpha_j}|=2\sin\!\left(\frac{|t-\alpha_j|}{2}\right)
            \ge \frac2\pi |t-\alpha_j|
\]
gives the asserted distance from \(K\). In particular \(\zeta_j\)
lies in a complementary arc of \(K\). Lemma~\ref{lem:individual-continuation}
gives analytic continuation of \(e\) through that point and
\(e(\zeta_j)\ne0\).

\emph{The starting point on each ray.}
For \(t\in I_j\), \(|t-\alpha_j|\le3g_j/2+g_j^2/2\le2g_j\).
The radial depth \(s_j=g_j^{3/2}\) is at most \(g_j\), so
\(|e^{it}-b_j^0|\le 4g_j\). Also
\(1-|b_j^0|^2=2s_j-s_j^2\ge s_j\). With the normalized Poisson
kernel from \eqref{eq:source}, this gives
\[
             P_{b_j^0}(e^{it})\ge\frac{s_j}{16g_j^2}.
\]
Since \(u=-\log|e^*|\) is nonnegative, we can retain only the
contribution from \(I_j\) to its Poisson integral at \(b_j^0\).
The kernel lower bound weights the logarithmic mass from
Section~\ref{sec:10}. By \eqref{eq:fixed-logarithmic-density} and
\eqref{eq:fixed-block-logarithmic-mass},
\[
 -\log|e(b_j^0)|
 \ge \frac{s_j}{16g_j^2}\int_{I_j}u\,dm
 \ge \frac{1}{128\pi j^2\sqrt{g_j}},
\]
which proves \eqref{eq:transport-low-log}. Exponentiating bounds
the modulus of \(e\), which appears squared in the upper estimate
\eqref{eq:transport-B-bounds} for \(B\). By \eqref{eq:transport-cB},
\(c\le B/(B+1)\le B\).
Because \(1-|b_j^0|^2\ge s_j\), this yields
\[
 c(b_j^0)\le
 \frac{2|e(b_j^0)|^2}{\kappa s_j}
 \le\frac{2}{\kappa}g_j^{-3/2}
      \exp\!\left(-\frac{1}{64\pi j^2\sqrt{g_j}}\right).
\]
Here \(g_j^{-1/2}=2^{j+5/2}\). The logarithm of the prefactor grows
linearly in \(j\), whereas \(2^j/j^2\) grows faster than \(j\).
Thus the bound tends to zero.

\emph{The boundary end of a fixed ray.}
Fix \(j\ge6\). Continuation at the gap point gives
\(|e(r\zeta_j)|\to |e(\zeta_j)|>0\) as \(r\uparrow1\).
The kernel diagonal from
Section~\ref{subsec:02-norm} satisfies
\[
 k_{r\zeta_j}(r\zeta_j)=\frac{-\log(1-r^2)}{r^2}
                                      \longrightarrow\infty.
\]
The positive limit of the modulus therefore gives
\(B(r\zeta_j)\ge|e(r\zeta_j)|^2 k_{r\zeta_j}(r\zeta_j)\to\infty\)
by the lower evaluation bound in \eqref{eq:transport-B-bounds}.
The lower bound \(B/(B+2)\le c\) in \eqref{eq:transport-cB} and
\(c\le1\) give \eqref{eq:transport-high-endpoint}. For the boundary
values of \(e\), nonvanishing at this fixed gap point suffices.
\end{proof}

\begin{proposition}[Selection of a prescribed value along each ray]
\label{prop:transport-exact-selection}
For every \(a\in(0,1)\) there is an integer \(J(a)\ge6\) such that,
for each \(j\ge J(a)\), there is a least radius \(r_j\in[1-s_j,1)\)
satisfying \(c(r_j\zeta_j)=a^2\). Write \(b_j=b_j(a)=r_j\zeta_j\)
for the first crossing point and set \(d_j=1-r_j\). Then
\begin{equation}\label{eq:exact-origin-depth}
                  0<d_j<s_j,\qquad c(b_j)=a^2 .
\end{equation}
For \(K^{(j)}=\varphi_{b_j}(K)\), the outer canonical function
\(G_{\D_{K^{(j)}}}\) satisfies \(G_{\D_{K^{(j)}}}(0)=a\), and \(K^{(j)}\)
is a nonpolar Carleson set.
\end{proposition}

\begin{proof}
Prescribe \(a\) after the compact set \(K\) has been fixed. Choose
\(J(a)\ge6\) from
\eqref{eq:transport-low-endpoint}, so that \(c(b_j^0)<a^2\) for
all \(j\ge J(a)\). For each such \(j\), hold the index fixed
and take \(r\uparrow1\) in \eqref{eq:transport-high-endpoint}.
There is an interior radius \(R_j\in(1-s_j,1)\) with
\(c(R_j\zeta_j)>a^2\), where \(R_j\) can depend on \(a\) and \(j\).
We can therefore choose the first crossing
on a compact interval of radii whose depths are at most \(s_j\).
The continuity of \(c\) from
Proposition~\ref{prop:transport-origin-continuity} allows the choice
\[
 r_j=\inf\{r\in[1-s_j,1):c(r\zeta_j)\ge a^2\},
 \qquad b_j=r_j\zeta_j,\qquad d_j=1-r_j.
\]
The infimum is attained because it equals the minimum of the
nonempty closed set of such radii in \([1-s_j,R_j]\).
The strict starting inequality and continuity give \(r_j>1-s_j\),
and \(r_j\le R_j<1\). If \(c(r_j\zeta_j)>a^2\), continuity gives
a smaller radius in that closed set, contrary to minimality.
Hence \(c(r_j\zeta_j)=a^2\). The argument gives the first
crossing even if the function crosses the level several times.

Section~\ref{subsec:11-ideals} shows that \(K^{(j)}\) is a nonpolar
Carleson set and that \(\D_{K^{(j)}}\) has a canonical function.
Since that function is normalized to have positive value at zero, \(c(b_j)=a^2\) gives \(G_{\D_{K^{(j)}}}(0)=a\).
\end{proof}

\subsection{Location of the entire compact set and its capacity}
\label{subsec:11-location}

The first crossing fixes the image \(K^{(j)}\), whose canonical
function has the prescribed value \(a\) at zero. We now estimate
the location and capacity of \(K^{(j)}\) using
that same parameter. The location estimate holds for every \(b=(1-d)\zeta_j\)
with \(0<d\le s_j\). It uses separation of \(\zeta_j\)
from every point of the original compactum.

\begin{proposition}[Hausdorff distance estimate]
\label{prop:transport-collapse}
For \(j\ge6\) and every \(0<d\le s_j\), set \(r=1-d\) and
\(b=r\zeta_j\). Then the chordal Hausdorff distance satisfies
\begin{equation}\label{eq:transport-collapse}
 d_H(\varphi_b(K),\{1\})<2^{1-j}.
\end{equation}
\end{proposition}

\begin{proof}
For every \(\zeta\in K\), subtracting \(\zeta_j\) from
\(\varphi_b(\zeta)\) gives
\begin{equation}\label{eq:transport-displacement-identity}
 |\varphi_b(\zeta)-\zeta_j|
   =\frac{d|1+\overline{\zeta_j}\zeta|}
           {|1-r\overline{\zeta_j}\zeta|}.
\end{equation}
On the circle the square of the denominator is
\(d^2+r|\zeta-\zeta_j|^2\). The separation
\eqref{eq:fixed-gap-separation} holds for all \(\zeta\in K\),
including \(1\). As \(r\ge1/2\) and \(|\zeta-\zeta_j|\ge2g_j/\pi\),
the denominator is at least \(\sqrt{2}g_j/\pi\). The numerator is
at most \(2d\). Thus
\begin{equation}\label{eq:transport-displacement-bound}
 \sup_{\zeta\in K}|\varphi_b(\zeta)-\zeta_j|
 \le\pi\sqrt{2}\,d/g_j
 \le\pi\sqrt{2}\,\sqrt{g_j}
 =\frac{\pi}{4}\,2^{-j},
\end{equation}
where \(s_j=g_j^{3/2}\) and \(g_j=2^{-2j-5}\).
To compare \(\zeta_j\) with the limiting point, we use
\[
 |\zeta_j-1|\le\theta_j+\frac32g_j
 =\bigl(1+3\cdot2^{-j-6}\bigr)2^{-j}
 \le\bigl(1+3/1024\bigr)2^{-j}\qquad(j\ge6).
\]
Since \(\pi/4+1+3/1024<2\), the triangle inequality puts every
point of the image within \(2^{1-j}\) of \(1\). For a nonempty
compact set and a singleton, the Hausdorff distance equals this
supremum of distances to the singleton. This proves
\eqref{eq:transport-collapse}.
\end{proof}

To bound the capacity of the same images, we put each one in a
small arc around \(1\). For normalized arclength as defining
measure, the arc-capacity estimate of El-Fallah, Elmadani, and
Kellay gives an estimate of order \(1/\log(1/\delta)\) for arcs of length \(\delta\)
\cite[Theorem~2, p.~870]{EEK19Kernel}. The obstacle
below is greater than one on an open neighborhood of the arc and
gives an explicit upper bound for the capacity defined using \(\|\cdot\|_V\).

\begin{lemma}[An open-neighborhood obstacle for a small arc]
\label{lem:transport-arc-obstacle}
For \(0<\delta<1/4\), let
\(A_\delta=\{e^{it}:|t|\le\delta\}\) and
\(\Lambda_\delta=\log(1/\delta)\). Then
\begin{equation}\label{eq:transport-arc-capacity}
 C_R(A_\delta)\le
 \frac{1}{2\Lambda_\delta}
 +\frac{(\log 2)^2+\pi^2/12}{\Lambda_\delta^2}.
\end{equation}
\end{lemma}

\begin{proof}
Put \(\rho=1-\delta\) and define the smooth real function
\[
 H_\delta(e^{it})=\frac{\log(2/|1-\rho e^{it}|)}{\Lambda_\delta}.
\]
It is positive on the circle because \(|1-\rho e^{it}|<2\).
When \(|t|\le\delta\), we have
\[
 |1-\rho e^{it}|\le(1-\rho)+\rho|1-e^{it}|
 \le\delta+\rho\delta<2\delta .
\]
Hence \(H_\delta>1\) on the closed arc. The open set
\(\{H_\delta>1\}\) contains the arc, so \(H_\delta\) is admissible
in the obstacle definition of Section~\ref{subsec:03-norm}.

We compute \(\|H_\delta\|_V^2\) to obtain the capacity bound.
Taking real parts of
\(-\log(1-\rho e^{it})=\sum_{n\ge1}\rho^ne^{int}/n\)
gives
\[
 \widehat H_\delta(0)=\frac{\log2}{\Lambda_\delta},
 \qquad
 \widehat H_\delta(\pm n)=\frac{\rho^n}{2n\Lambda_\delta}\quad(n\ge1).
\]
In the \(V\)-norm, the constant coefficient and both
nonzero frequency signs give
\begin{align}
 \|H_\delta\|_V^2
 &=\frac{(\log 2)^2+\tfrac12\sum_{n\ge1}
                 (1+n)\rho^{2n}/n^2}{\Lambda_\delta^2}
                       \label{eq:transport-obstacle-norm}\\
 &\le\frac{(\log 2)^2+\pi^2/12
                   -\tfrac12\log(1-\rho^2)}{\Lambda_\delta^2}. \notag
\end{align}
Here we bound the terms with \(1/n^2\) by their sum at \(\rho=1\)
and sum the terms with \(1/n\) exactly.
Since \(1-\rho^2=2\delta-\delta^2\ge\delta\), the last display
is bounded by the right side of \eqref{eq:transport-arc-capacity}.
\end{proof}

The location estimate holds for every depth \(0<d\le s_j\).
For any such image, use angular representatives \(|t|\le\pi\). The chord
comparison
\(|t|\le(\pi/2)|e^{it}-1|\) and \eqref{eq:transport-collapse}
give the containing arc below. In particular this containment
holds for \(K^{(j)}=\varphi_{b_j}(K)\).
\begin{equation}\label{eq:transport-containing-arc}
 K^{(j)}\subset A_{\delta_j},\qquad \delta_j=6\cdot2^{-j}.
\end{equation}
For \(j\ge6\), we have \(\delta_j<1/4\) and
\(\log(1/\delta_j)=j\log 2-\log6\ge j/3\), using
\(\log2>2/3\) and \(\log6<2\). Also
\((\log2)^2+\pi^2/12<11/6\). Applying the arc estimate and then
using \(j\ge6\) bounds the capacity at every such depth by
\(17/(4j)\). For \(K^{(j)}\) this gives
\begin{equation}\label{eq:transport-selected-capacity}
 0<C_R(K^{(j)})\le\frac{3}{2j}+\frac{33}{2j^2}
        \le\frac{17}{4j}<\frac5j\qquad(j\ge J(a)).
\end{equation}
Positivity follows from the nonpolarity of \(K^{(j)}\),
as proved in Section~\ref{subsec:11-ideals}. Thus the same sets \(K^{(j)}\), whose canonical functions satisfy
\(G_{\D_{K^{(j)}}}(0)=a\), shrink to \(1\) and have capacities tending
to zero.

\subsection{The geometric family and one-atom recovery}
\label{subsec:11-realization}

For each prescribed value, we reindex the exact-level images
and rotate their limiting point from \(1\) to the desired
direction. The location bound and the prescribed value at zero will
give the hypotheses of Criterion C. The capacity bound will also
be used when we take unions in Section~\ref{sec:12}.

\paragraph*{Continuous extension and the zero set on the closed disk.}
Let $E$ be a nonempty finite union of disk automorphism images of the
compact set $K$ in \eqref{eq:fixed-packed-compact}. The canonical function
$u=G_{\D_E}$ extends continuously to $\overline{\DD}$ and satisfies
\[
 \{z\in\overline{\DD}:u(z)=0\}=\Zrad(u)=E.
\]

\begin{proof}
Sections~\ref{subsec:11-ideals} and~\ref{subsec:03-ideal} show that
$E$ is a nonpolar Carleson set. Its trace-zero space contains
a smooth outer function, so $u$ is outer and $|u|\le1$ in $\DD$.
Fix an automorphism $\phi$ with $\phi(K)\subset E$ and put $h=u\circ\phi$.
By \eqref{eq:transport-ideals} and rotation invariance, we have $h\in\D_K$. Since $h$ is nonzero and $|h|\le1$, Hardy factorization gives $v=-\log|h^*|\in L^1(m)$, $v\ge0$,
and
\[
 -\log|h(z)|\ge\int_{\T}P_zv\,dm\qquad (z\in\DD).
\]

Take $x\in K\setminus\{1\}$ and choose a Cantor copy in $K$ containing
$x$. At each stage of its construction choose a closed arc
$J$ containing $x$ among the arcs retained at that stage, nested
inside the arc at the previous stage. This includes their endpoints. Write $l$ for
the angular length of $J$ and $\mathcal E_J(h^*)$ for the double integral
in \eqref{eq:V-difference} restricted to $J\times J$. The copy in each
such $J$ is a rescaled middle-thirds Cantor set on which the trace of $h$
vanishes quasi-everywhere. The estimate on one cell in the proof of
Theorem~\ref{thm:packed-trace}, using $m(J)=l/(2\pi)$, gives
\[
 \frac{1}{m(J)}\int_J |h^*|^2\,dm
 \le 56\pi^2\mathcal E_J(h^*)\longrightarrow0.
\]
The convergence follows from integrability of the full difference
kernel and $l\to0$. Jensen's inequality shows that
$m(J)^{-1}\int_J v\,dm\to+\infty$. For $l<1$ and $\zeta\in J$,
$|\zeta-(1-l)x|\le |\zeta-x|+l\le2l$. Since $1-(1-l)^2\ge l$,
we have $P_{(1-l)x}(\zeta)\ge1/(4l)$ on $J$. The Poisson integral
of $v$ at $(1-l)x$ is at least $(8\pi m(J))^{-1}\int_J v\,dm$,
which tends to $+\infty$. Hence $h((1-l)x)\to0$.

For $x=1$, use the arcs $I_j$ and integers $N_j$ in
\eqref{eq:fixed-packed-scales}. Theorem~\ref{thm:packed-trace},
$\|h^*\|_V=\|h\|_\D$, and Jensen's inequality yield
\[
 \int_{I_j}v\,dm
 \ge \frac{L_j}{4\pi}\log\frac{N_j}{600\pi\|h\|_\D^2}
 \ge \frac{1}{8\pi j^2}
\]
for all sufficiently large $j$. The last inequality uses
$\log N_j\ge 1/(j^2L_j)\to+\infty$ with $h$ fixed. Every $\zeta\in I_j$
satisfies $|\zeta-1|\le2\cdot 2^{-j}$, so
$|\zeta-(1-2^{-j})|\le3\cdot2^{-j}$ and
$P_{1-2^{-j}}(\zeta)\ge 2^j/9$. Therefore
$-\log|h(1-2^{-j})|\ge 2^j/(72\pi j^2)\to+\infty$.

The two cases give a sequence on the radius to each $x\in K$ along
which $h$ tends to zero. For a fixed $\phi$ and $x$, smoothness of $\phi$
on the closed disk gives $|\phi(rx)-\phi(x)|=O(1-r)$ as $r\uparrow1$.
Writing $\phi$ as a rotation composed with $\varphi_b$ in
\eqref{eq:automorphism-definition} gives
$1-|\phi(rx)|\asymp 1-r$ by the identity immediately following that
formula. The images of these radial points thus lie in one fixed
Stolz region at $\phi(x)$ as $r\uparrow 1$.
The function $u$ has norm one and satisfies the extremal moments
\eqref{eq:canonical-moments}. Richter and Yilmaz
\cite[Theorem~1.4, p.~2120]{RY19} show that the limit of $|u|$ in this Stolz region
equals its unrestricted upper limit
from within $\DD$ at $\phi(x)$. Since $u(\phi(rx))=h(rx)$ tends to zero
along the sequence chosen above, the unrestricted limit of $u$ at
$\phi(x)$ is zero. This holds at every point of $E$.

On $\T\setminus E$, Lemma~\ref{lem:individual-continuation} gives a
holomorphic extension with nonzero value. Define $u$ on $E$ by zero.
For any $q_n\in\overline{\DD}$ tending to $q\in E$, choose interior
points $z_n$ with $|z_n-q_n|<1/n$ and $|u(z_n)-u(q_n)|<1/n$, using the
interior limits at boundary points. Then $z_n\to q$ and $u(z_n)\to0$,
so $u(q_n)\to0$. This proves continuity at every point of $E$.
Outerness excludes zeros in the disk, and continuity identifies the
boundary values with the radial limits at every point of $\T$.
\end{proof}

\begin{theorem}[Carleson sets with a prescribed extremal value at zero]\label{thm:G}
For every \(a\in(0,1)\) and \(\xi\in\T\), there are nonpolar
Carleson sets \(K_n\) and their normalized outer extremal functions with positive value at zero \(\eta_n=G_{\D_{K_n}}\in\Xpr\) such that, for every
\(n\ge1\),
\[
 \eta_n(0)=a,\qquad [\eta_n]=\D_{K_n},\qquad
 d_H(K_n,\{\xi\})\le3\cdot2^{-n},\qquad 0<C_R(K_n)\le5/n.
\]
Each $\eta_n$ extends continuously to $\overline{\DD}$ and satisfies
\[
 \{z\in\overline{\DD}:\eta_n(z)=0\}=\Zrad(\eta_n)=K_n.
\]
Here \(d_H\) is for chord distance. The fixed compact set \(K\) in
\eqref{eq:fixed-packed-compact} supplies all these families by applying Möbius transformations and then rotations.
\end{theorem}

\begin{proof}
First take \(\xi=1\). For the prescribed \(a\), let
\(J(a)\ge6\) and the selected parameters \(b_j\) be those of
Proposition~\ref{prop:transport-exact-selection}.
Set \(K_n=K^{(J(a)+n-1)}\) and
\(\eta_n=G_{\D_{K_n}}\). The construction gives nonpolar Carleson sets and normalized outer extremal functions satisfying \(\eta_n(0)=a\) and generated space \(\D_{K_n}\).
Since \(a<1\), each canonical function is in \(\Xpr\).
Since \(J(a)+n-1\ge n\), the location bound is at most \(2^{1-n}\)
and the capacity bound at most \(5/n\). They hold from the first
index on the same sets.

For general \(\xi\), rotate the sets by \(\xi\) and replace each function \(f(z)\) by \(f(\bar\xi z)\). Rotation preserves
radial paths, normalized arclength and chord distance. It preserves
the boundary norm, and rotating admissible open obstacles in both
directions preserves \(C_R\). Therefore it maps \(\D_{K_n}\) onto \(\D_{\xi K_n}\),
preserves nonpolarity and the Carleson condition, and carries the
geometric bounds to \(\{\xi\}\). It is unitary in the Dirichlet
norm and fixes \(1\), so it carries the projection of \(1\) and its positive unit normalization to those of the rotated space.

Each $K_n$ is the image of $K$ under a disk automorphism. The preceding
argument therefore gives the continuous extension and the zero-set
identity for $\eta_n=G_{\D_{K_n}}$.
\end{proof}

For a fixed atom of mass \(\tau>0\) at \(\xi\), apply
Theorem~\ref{thm:G} with \(a=G_{\tau,\xi}(0)\).
It gives nonpolar Carleson sets \(K_n\) with
\(d_H(K_n,\{\xi\})\to0\) and canonical functions \(\eta_n\)
with \(\eta_n(0)=a\) for every \(n\).
Theorem~\ref{thm:C} therefore applies to \((K_n,\eta_n)_{n\ge1}\). The capacity
bound remains available on this same sequence for the finite
unions used in the next section.

\begin{corollary}[One-atom recovery]\label{cor:R}
Fix \(\tau>0\) and \(\xi\in\T\). There is one family as in
Theorem~\ref{thm:G} with
\begin{equation}\label{eq:one-atom-recovery}
 \eta_n(0)=G_{\tau,\xi}(0)\quad(n\ge1),\qquad
 \|\eta_n-G_{\tau,\xi}\|_{\D}\longrightarrow0 .
\end{equation}
The target has the exact factorization
\(G_{\tau,\xi}=S_{\tau,\xi}O_{\tau,\xi}\), where
\(O_{\tau,\xi}\in\D\) is bounded by one, outer, positive at zero,
and cyclic in \(\D\). It is a contractive Dirichlet multiplier.
Both the radial zero set and the unrestricted approach-zero set of \(G_{\tau,\xi}\) equal \(\{\xi\}\).
\end{corollary}

\begin{proof}
Fix \(\tau,\xi\) and put \(a=G_{\tau,\xi}(0)\in(0,1)\)
by \eqref{eq:one-atom-origin-bounds}. Apply Theorem~\ref{thm:G}
at this \(a\) and direction \(\xi\). The equalities \(\eta_n(0)=a\)
and the Hausdorff bound verify the hypotheses of Theorem~\ref{thm:C}
for \((K_n,\eta_n)_{n\ge1}\), giving
\eqref{eq:one-atom-recovery}.
The same approximants satisfy the capacity bound and the closed-disk
zero-set conclusion of Theorem~\ref{thm:G}.
Corollary~\ref{cor:one-atom-cofactor} gives the properties of
$O_{\tau,\xi}$ and identifies the radial and unrestricted approach-zero
sets of the target $G_{\tau,\xi}$.
\end{proof}

Section~\ref{sec:12} will use these families for the atoms of a fixed finite atomic measure. Their capacity bounds control the capacities of the unions. The outer factors are cyclic Dirichlet multipliers, allowing the product of the one-atom limits to generate \(M_\mu\). The projection argument then gives convergence of the canonical functions of the corresponding trace-zero spaces.

\newpage \section{Recovery of canonical functions for singular measures}
\label{sec:12}

We now recover the canonical function for a fixed finite atomic
measure. At each index we take the union of the compact sets supplied
by the one-atom families. We then compare the resulting trace-zero
spaces with the target space in two ways. Every vector of the target
space must be approximated in norm by vectors from these spaces,
and every subsequential weak limit from the spaces must belong to the
target. Together these properties give norm convergence of the
projections of each fixed vector. Applying this to \(1\) and
normalizing gives the canonical functions in Approximation~B.
Section~\ref{subsec:12-assembly} then uses the finite atomic result
to treat a general singular measure \(\mu\) with
\(S_\mu H^2\cap\D\ne\{0\}\).

\subsection{Unions of Carleson sets and products of canonical functions}\label{subsec:12-products}

Fix a nonzero finite positive measure \(\mu\) of finite support as in
\eqref{eq:atomic-target}. Add the masses at each repeated direction
before writing
\[
 \mu=\sum_{i=1}^N\tau_i\delta_{\xi_i},\qquad
 \tau_i>0,\quad \xi_i\in\T,\quad \xi_i\ne\xi_k\quad(i\ne k).
\]
The number \(N\ge1\), the masses, and the directions are fixed throughout
the proof. The one-atom target at \(\xi_i\) uses the mass
\(\tau_i\) at that direction. Write
\begin{equation}\label{eq:finite-atomic-components}
 S_i=S_{\tau_i,\xi_i},\qquad
 G_i=G_{S_iH^2\cap\D},\qquad G_i(0)>0.
\end{equation}
In particular \(S_\mu=\prod_{i=1}^N S_i\). Section~\ref{subsec:02-atomic}
constructs \(M_\mu=S_\mu H^2\cap\D\) and its canonical function
\(G_\mu\) independently of the recovery argument.

\paragraph{The unions of Carleson sets.}
For each \(i\), fix a one-atom family from Corollary~\ref{cor:R}
and index all \(N\) families by \(\ell\ge1\). We use the
\(\ell\)-th compactum of each family to form the \(\ell\)-th union.
Thus \(K_{i,\ell}\) is a nonpolar Carleson set and
\(H_{i,\ell}=G_{\D_{K_{i,\ell}}}\) is its unit outer canonical function with positive value at zero. Along these families,
\begin{equation}\label{eq:finite-recovery-single-input}
 [H_{i,\ell}]=\D_{K_{i,\ell}},\qquad H_{i,\ell}(0)=G_i(0),\qquad
 H_{i,\ell}\longrightarrow G_i\quad\text{in }\D.
\end{equation}
For each family and every \(\ell\ge1\), the same compactums
satisfy the geometric bounds, with \(d_H\) in the chord metric:
\begin{equation}\label{eq:finite-recovery-single-geometry}
 d_H(K_{i,\ell},\{\xi_i\})\le 3\cdot2^{-\ell},
 \qquad 0<C_R(K_{i,\ell})\le 5/\ell .
\end{equation}
The corollary also gives \(G_i=S_iO_i\), with \(O_i\) outer,
positive at zero, cyclic in \(\D\), and a contractive Dirichlet
multiplier bounded by one on \(\DD\).

For each \(\ell\), put
\begin{equation}\label{eq:finite-recovery-union}
 K_\ell=\bigcup_{i=1}^N K_{i,\ell},\qquad M_\ell=\D_{K_\ell}.
\end{equation}
This finite union is compact and arclength-null. Its capacity is
positive by monotonicity, since it contains \(K_{1,\ell}\). The
capacity bounds \eqref{eq:finite-recovery-single-geometry} and
finite subadditivity \eqref{eq:capacity-finite-subadditivity} give
\begin{equation}\label{eq:finite-recovery-union-capacity}
 0<C_R(K_\ell)\le\sum_{i=1}^N C_R(K_{i,\ell})\le5N/\ell .
\end{equation}
The bound \(5N/\ell\) depends on the number of atoms and the
index. It tends to zero along this sequence with \(N\) fixed.
To verify the Carleson condition, observe that for
\(\zeta\in\T\),
\[
 \log\frac{2}{\dist(\zeta,K_\ell)}
   =\max_i\log\frac{2}{\dist(\zeta,K_{i,\ell})}
   \le \sum_{i=1}^N\log\frac{2}{\dist(\zeta,K_{i,\ell})}.
\]
The summands are nonnegative and integrable. By
Section~\ref{subsec:03-ideal}, \(M_\ell\) is closed, shift-invariant,
and proper. Its exact intersection description is
\begin{equation}\label{eq:finite-recovery-ideal-intersection}
 M_\ell=\bigcap_{i=1}^N \D_{K_{i,\ell}} .
\end{equation}
Indeed, vanishing quasi-everywhere on the union implies the same on
each constituent. Conversely, for a fixed function in the
intersection, the union of its finitely many polar exceptions
is polar. Its trace thus vanishes quasi-everywhere on \(K_\ell\).
This argument allows overlaps among the compactums.

\paragraph{Products and their norm bounds.}
To build vectors vanishing on the union, we will multiply the
one-atom approximants. The multiplier estimates make these
products legitimate in the Dirichlet norm and let us compare
them with the corresponding product of the limits. Each
\(H_{i,\ell}\) and \(G_i\) multiplies \(\D\) contractively by
\eqref{eq:canonical-multiplier}. For \(h\in\D\), boundedness gives \(O_i h\in H^2\), and the
canonical multiplier property gives \(G_i h\in\D\).
Since \(G_i h=S_i(O_i h)\), contractive inner division gives
\(\|O_i h\|_\D\le\|G_i h\|_\D\le\|h\|_\D\).
We can therefore compose these multiplication operators and define
\begin{equation}\label{eq:finite-recovery-products}
 P_\ell=\prod_{i=1}^N H_{i,\ell},\qquad
 P=\prod_{i=1}^N G_i=S_\mu O,\qquad
 O=\prod_{i=1}^N O_i .
\end{equation}
Each finite product acts contractively on \(\D\). Applied to \(1\),
this proves Dirichlet membership and norm bounds for the products
themselves. In particular, the bounded product \(O\) belongs to
\(H^2\), so the product \(P=S_\mu O\) belongs to
\(M_\mu\). The finite products \(O\) and \(P_\ell\) are outer:
their integrable boundary logarithms add in the outer representation
\eqref{eq:outer-intro}. The factors are positive at zero, and
\eqref{eq:finite-recovery-single-input} gives
\begin{equation}\label{eq:finite-recovery-product-origin}
 \|P_\ell\|_\D\le1,\qquad\|P\|_\D\le1,\qquad
 P_\ell(0)=P(0)=c_\mu:=\prod_{i=1}^N G_i(0)>0 .
\end{equation}

Replacing the factors one at a time gives
\[
 P_\ell-P=\sum_{i=1}^N
 \left(\prod_{k<i}H_{k,\ell}\right)
 (H_{i,\ell}-G_i)
 \left(\prod_{k>i}G_k\right).
\]
In each term the surrounding multipliers contract the Dirichlet norm
of the difference. Thus
\begin{equation}\label{eq:finite-recovery-product-convergence}
 \|P_\ell-P\|_\D
   \le\Delta_\ell:=\sum_{i=1}^N\|H_{i,\ell}-G_i\|_\D
   \longrightarrow0 .
\end{equation}
Here the number of terms is fixed and each factor converges as
\(\ell\to\infty\). The error \(\Delta_\ell\) depends on the
one-atom families chosen for this \(\mu\).

To place \(P_\ell\) in \(M_\ell\), take a point
\(\zeta\in K_{i,\ell}\) outside the polar exceptional set.
At this point, the trace of \(H_{i,\ell}\) agrees with its
radial limit and vanishes. The canonical bound
\(|H_{k,\ell}|\le1\) for all \(k\) then gives
\[
 |P_\ell(r\zeta)|\le |H_{i,\ell}(r\zeta)|\longrightarrow0
      \qquad(r\uparrow1).
\]
Combining the finitely many exceptions proves radial vanishing
quasi-everywhere on \(K_\ell\), which is the trace condition of
Section~\ref{subsec:03-ideal}. Shift invariance now gives
\begin{equation}\label{eq:finite-recovery-trial-membership}
   P_\ell\in M_\ell,\qquad pP_\ell\in M_\ell
       \quad\text{for every polynomial }p,\qquad P\in M_\mu.
\end{equation}
In particular, \eqref{eq:finite-recovery-product-origin}
gives a function in each space \(M_\ell\) with nonzero value at zero. These products
serve as trial vectors. The projection step below determines the
canonical functions. For strong recovery we first show that
polynomial multiples of the limit \(P\) are dense in \(M_\mu\).

\subsection{Strong recovery of every target vector}\label{subsec:12-strong}

The inner factor of \(P\) is \(S_\mu\). To show that \(P\) generates the whole target space, we prove that its outer factor
\(O=\prod_{i=1}^N O_i\) is cyclic.

The outer-product identity \([vw]=[v]\cap[w]\) for outer factors
\(v,w\in\D\) with \(vw\in\D\) is stated in
\cite[Theorem~2.4, p.~2199]{EFKR09}, with attribution to
\cite[Lemma~2.2 and Theorem~4.5]{RS92}. In our multiplier setting
the needed cyclicity implication is already Brown and Shields'
product result \cite[Proposition~8, pp.~275--276]{BS84}. We recall its short proof.
Approximate \(1\) by polynomial
multiples of the first factor and multiply by the second factor.
Suppose \(v\) and \(w\) are cyclic Dirichlet multipliers.
Multiplication by \(w\) gives \(vw\in\D\).
Choose polynomials \(q_j\) with
\(q_jv\to1\) in \(\D\). Since \(w\) multiplies all of \(\D\)
boundedly,
\[
                  q_jvw\longrightarrow w\quad\text{in }\D .
\]
Therefore \(w\in[vw]\). The space \([vw]\) is closed and invariant
under multiplication by \(z\), so it contains \([w]=\D\).
Applying this argument finitely many times to the contractive
multipliers \(O_i\) proves \([O]=\D\).

By \eqref{eq:finite-recovery-products} and
\eqref{eq:finite-recovery-trial-membership}, \(P=S_\mu O\) is a
nonzero function in \(\D\), with \(S_\mu\) inner and \(O\)
outer. Theorem~\ref{input:module} now gives
\begin{equation}\label{eq:finite-recovery-target-identification}
 [P]=[O]\cap S_\mu H^2=\D\cap S_\mu H^2=M_\mu .
\end{equation}
The identity above allows us to approximate any given \(h\in M_\mu\) by a polynomial multiple \(pP\). 
We now need approximants from the spaces \(M_\ell\). Once \(p\) is fixed, \(pP_\ell\) belongs to \(M_\ell\) and
converges to \(pP\) in \(\D\). This is the polynomial approximation argument in \cite[Proposition~4.3(a)]{LR15}.
The next lemma uses this observation to choose \(h_\ell\in M_\ell\) for every \(\ell\) so that \(h_\ell\to h\).

\begin{lemma}[Strong recovery along all indices]\label{lem:finite-recovery-strong}
Use the fixed \(\mu\) and the sequence \(M_\ell\) in
\eqref{eq:finite-recovery-union}.
For every \(h\in M_\mu\) there are \(h_\ell\in M_\ell\) with
\(\|h_\ell-h\|_\D\to0\) as \(\ell\to\infty\).
\end{lemma}
\begin{proof}
To control the error for one fixed polynomial
\(p(z)=\sum_{k=0}^d p_k z^k\), set
\(C_p=\sum_{k=0}^d|p_k|\sqrt{k+1}\).
The shift bounds in Section~\ref{subsec:02-norm} give
\(\|pf\|_\D\le C_p\|f\|_\D\) for every \(f\in\D\).
Using \eqref{eq:finite-recovery-product-convergence}, we obtain
\begin{equation}\label{eq:finite-recovery-polynomial-error}
 \|pP_\ell-h\|_\D\le\|pP-h\|_\D+C_p\Delta_\ell .
\end{equation}
The term involving \(\Delta_\ell\) tends to zero for this fixed polynomial. When we choose a new polynomial to meet a smaller bound on the approximation error, we allow its multiplier bound to change and choose a new index threshold afterward. This gives approximants for every sufficiently large index.

Fix \(h\in M_\mu\). By \eqref{eq:finite-recovery-target-identification},
for each \(j\ge1\) choose \(p_j\) with
\(\|p_jP-h\|_\D<2^{-j}\). After choosing \(p_j\), take a threshold
\(n_j\) such that \(C_{p_j}\Delta_\ell\le2^{-j}\) for every
\(\ell\ge n_j\). Choose these integers strictly increasing with
\(n_j\ge j\). Set \(h_\ell=0\) for \(\ell<n_1\), and use
\(h_\ell=p_jP_\ell\) on each block \(n_j\le\ell<n_{j+1}\).
These vectors lie in \(M_\ell\) by
\eqref{eq:finite-recovery-trial-membership} and invariance.
By \eqref{eq:finite-recovery-polynomial-error},
\(\|h_\ell-h\|_\D\le2^{1-j}\) throughout the \(j\)-th block.
As \(\ell\to\infty\) the block index tends to infinity, which
proves strong recovery for \((M_\ell)_{\ell\ge1}\).
\end{proof}

The polynomials and thresholds can depend on \(h\), while the
spaces \(M_\ell\) are fixed by \(\mu\). We next locate weak limits of
vectors chosen from these same spaces.

\subsection{Every subsequential weak limit lies in the target}
\label{subsec:12-weak}
We have shown that every function in \(M_\mu\) can be approximated
by functions in \(M_\ell\). We now consider a weakly convergent
sequence \(u_j\in M_{\ell_j}\) and show that its limit belongs to
\(M_\mu\).
Luo and Richter's Corollary~4.5 \cite{LR15} gives strong convergence of the projections for each one-atom family. Testing the weakly convergent sequence against the orthogonal complement of each limiting one-atom space also shows that the weak limit of \(u_j\) belongs to \(M_\mu\). The proof below uses Hardy quotients.

\begin{lemma}[Limits of Hardy quotients]
\label{lem:finite-recovery-quotient-limit}
Let \(g_j\) be unit canonical functions with \(g_j(0)>0\), and suppose
\(g_j\to g\ne0\) in \(\D\). If \(u_j\in[g_j]\) and
\(u_j\rightharpoonup u\) in \(\D\), then there is a Hardy function
\(F\in H^2\) with
\begin{equation}\label{eq:finite-recovery-quotient-limit}
               u=gF,\qquad
               \|F\|_{H^2}\le\liminf_j\|u_j\|_{\D}.
\end{equation}
\end{lemma}

\begin{proof}
Lemma~\ref{lem:completed-canonical-quotient} gives the factorizations
\(u_j=g_jF_j\), where \(F_j\in H^2\) and
\(\|F_j\|_{H^2}\le\|u_j\|_{\D}\). Weak convergence bounds the
Dirichlet norms, so the quotients are bounded in \(H^2\).
To retain the stated norm bound, choose a subsequence along which
\(\|u_j\|_{\D}\) tends to its liminf. Weak compactness gives a further
subsequence with \(F_j\rightharpoonup F\) in \(H^2\).

Evaluation at each fixed \(z\in\DD\) is a bounded functional on
\(\D\) and \(H^2\), as in Section~\ref{subsec:02-norm}.
Along the extracted indices \(u_j(z)\to u(z)\), \(g_j(z)\to g(z)\),
and \(F_j(z)\to F(z)\). Passing to the limit in \(u_j(z)=g_j(z)F_j(z)\)
therefore gives \(u=gF\) throughout the disk. Along this further
subsequence, the Dirichlet norms still converge to the liminf of the
original sequence, \(\liminf_{j\to\infty}\|u_j\|_{\D}\).
Weak lower semicontinuity of the Hardy norm gives the bound in
\eqref{eq:finite-recovery-quotient-limit}.
\end{proof}

\begin{lemma}[Intersection of invariant subspaces of \(H^2\)]
\label{lem:finite-recovery-hardy-intersection}
For the fixed measure \(\mu\) with distinct directions, let \(S_i\)
be as in \eqref{eq:finite-atomic-components}. Then
\begin{equation}\label{eq:finite-recovery-hardy-intersection}
                     \bigcap_{i=1}^N S_i H^2=S_\mu H^2 .
\end{equation}
\end{lemma}

\begin{proof}
The identity \(S_\mu=\prod_i S_i\) gives the inclusion of \(S_\mu H^2\)
in every \(S_i H^2\), since a product of the remaining inner factors
multiplies \(H^2\) into itself.

For the reverse inclusion, let \(v\ne0\) belong to every \(S_iH^2\).
The singular measure of \(v\) records divisibility at each direction.
Use the Hardy factorization from Section~\ref{subsec:02-hardy} to write
\[
                     v=\omega B_v S_{\nu_v}O_v ,
\]
where \(|\omega|=1\), \(B_v\) is the Blaschke factor, \(\nu_v\) is
its positive singular inner measure, and \(O_v\in H^2\) is outer
with positive value at zero. For each \(i\), write \(v=S_i v_i\) with
\(v_i\in H^2\). Factoring \(v_i\) in \(H^2\) and multiplying by
\(S_i=S_{\tau_i\delta_{\xi_i}}\) adds \(\tau_i\delta_{\xi_i}\) to
its singular measure. Uniqueness of the Hardy factorization therefore gives \(\nu_v(\{\xi_i\})\ge\tau_i\). As the \(\xi_i\)
are distinct, for every Borel set \(E\subset\T\) we have
\(\nu_v(E)\ge\sum_{\xi_i\in E}\tau_i=\mu(E)\).
The function \(\omega B_vS_{\nu_v-\mu}O_v\)
belongs to \(H^2\), since it is the product of an inner factor
and the original Hardy outer factor. Its product with \(S_\mu\)
is \(v\), so \(v\in S_\mu H^2\). The zero function belongs to both
sides of \eqref{eq:finite-recovery-hardy-intersection}.
\end{proof}

At a repeated direction with \(s\le t\),
\(S_{t,\xi}=S_{s,\xi}S_{t-s,\xi}\) gives
\(S_{t,\xi}H^2\subseteq S_{s,\xi}H^2\), so their intersection corresponds to the larger mass \(t\).
We merged repeated directions before defining the spaces \(S_iH^2\), so that each direction carries the sum of the
masses assigned to it.
Equation \eqref{eq:finite-recovery-hardy-intersection} then uses one factor per distinct direction.

\begin{proposition}[Every subsequential weak limit lies in the target]
\label{prop:finite-recovery-weak-containment}
Use the same fixed \(\mu\) and the union sequence \((M_\ell)_{\ell\ge1}\)
in \eqref{eq:finite-recovery-union}. Let \((\ell_j)\) be any
strictly increasing sequence of indices.
If \(u_j\in M_{\ell_j}\) and \(u_j\rightharpoonup u\) in \(\D\),
then \(u\in M_\mu\).
\end{proposition}

\begin{proof}
For each \(i\), the identity
\eqref{eq:finite-recovery-ideal-intersection} and the canonical
generation in \eqref{eq:finite-recovery-single-input} give
\(u_j\in[H_{i,\ell_j}]\). The convergence in
\eqref{eq:finite-recovery-single-input} also gives
\(H_{i,\ell_j}\to G_i\) in \(\D\). Apply
Lemma~\ref{lem:finite-recovery-quotient-limit} to obtain
\[
 u=G_i F_i=S_i(O_iF_i),\qquad F_i\in H^2,\qquad
 \|F_i\|_{H^2}\le\liminf_j\|u_j\|_{\D}.
\]
The extraction in the quotient lemma may depend on \(i\), since
every subsequence retains the original weak limit \(u\).
The outer factor \(O_i\) is bounded, so
\(\|O_iF_i\|_{H^2}\le\|O_i\|_\infty\|F_i\|_{H^2}\) and
\(u\in S_iH^2\) for every \(i\).
By \eqref{eq:finite-recovery-hardy-intersection},
\(u\in S_\mu H^2\). Since the weak limit lies in \(\D\),
this gives \(u\in S_\mu H^2\cap\D=M_\mu\).
\end{proof}

Thus every target vector has a norm approximation along the sequence
of spaces \((M_\ell)_{\ell\ge1}\), and every subsequential weak limit from those spaces lies in the target. We use the two statements together in
the projection argument.

\subsection{Projection convergence and canonical normalization}\label{subsec:12-projections}

The preceding results show that every \(v\in M_\mu\) is the norm
limit of a sequence \(v_\ell\in M_\ell\) and that every subsequential
weak limit of vectors from \(M_\ell\) belongs to \(M_\mu\).
These are the two conditions for Mosco convergence of \(M_\ell\) to
\(M_\mu\) \cite[Definition~1.1, pp.~519--520]{Mosco69}.

For fixed \(y\in\D\), the second condition puts each subsequential
weak limit of \(P_{M_\ell}y\) in \(M_\mu\). The first condition supplies
the test vectors \(v_\ell\) needed to pass to the limit in the
projection orthogonality relation. This identifies each such weak
limit as \(P_{M_\mu}y\).
The projection consequence below follows from
\cite[Corollary of Theorem~A, p.~542]{Mosco69}
in the inner-product setting; see also
\cite[Definition~2.2 and Lemma~2.5]{CHMB21}.

\begin{lemma}[Strong convergence of subspace projections]
\label{lem:finite-recovery-projections}
Let \(L_\ell\) for \(\ell\ge1\) and \(L\) be closed subspaces of a
Hilbert space. Assume that for every \(v\in L\) there are
\(v_\ell\in L_\ell\) for all \(\ell\) with \(v_\ell\to v\) in norm.
Assume also that for every strictly increasing sequence
\((\ell_j)\), every weakly convergent sequence \(u_j\in L_{\ell_j}\)
has its weak limit in \(L\). Then for every fixed vector \(y\) in
the Hilbert space, \(P_{L_\ell}y\to P_Ly\) in norm.
\end{lemma}

\begin{proof}
Fix \(y\) and put \(z_\ell=P_{L_\ell}y\). The projections are
contractive, so \(\|z_\ell\|\le\|y\|\) for all \(\ell\).
Every subsequence therefore has a further weakly convergent
subsequence \(z_{\ell_j}\rightharpoonup z\). By the second assumption,
\(z\in L\).

We now test orthogonality to the target space. For any \(v\in L\),
choose the recovery sequence \(v_\ell\) of the first assumption.
Projection orthogonality gives
\(\langle y-z_{\ell_j},v_{\ell_j}\rangle=0\).
Replacing \(v_{\ell_j}\) by \(v\) changes this pairing by at most
\(\|y-z_{\ell_j}\|\|v_{\ell_j}-v\|\), which tends to zero.
Taking the weak limit in the pairing with \(v\) shows
\(\langle y-z,v\rangle=0\). As \(v\) was arbitrary, \(z=P_Ly\).

It follows that \(z_\ell\) converges weakly to \(P_Ly\) as \(\ell\to\infty\).
Indeed, failure of convergence against any fixed test vector would
give a subsequence whose pairings stay a positive distance from
the asserted limit, contrary to the weak convergence of a further
subsequence. The projection identities express the squared norms as
pairings against the fixed vector \(y\).
\[
 \|z_\ell\|^2=\langle z_\ell,y\rangle
   \longrightarrow\langle P_Ly,y\rangle=\|P_Ly\|^2 .
\]
Weak convergence and convergence of norms give the norm convergence.
\end{proof}

Hartz and Richter prove norm convergence of the canonical functions
of \(B_nH^2\cap\D\) to that of \(SH^2\cap\D\) when the inner
functions \(B_n\) converge locally uniformly to \(S\), with
\(S(0)\ne0\) and \(SH^2\cap\D\ne\{0\}\). Their hypotheses require
pointwise convergence of the local Dirichlet integrals on \(\T\)
and a uniform bound \(D_\zeta(B_n)\le c D_\zeta(S)\)
\cite[Lemma~5.1]{HR26}.

For the spaces \(M_\ell\) in this proof, closedness follows from
Section~\ref{subsec:03-ideal}, and \(M_\mu\) is closed by
Lemma~\ref{lem:exact-inner-factor}.
Lemma~\ref{lem:finite-recovery-strong} gives recovery along
\((M_\ell)_{\ell\ge1}\) for each target vector.
Proposition~\ref{prop:finite-recovery-weak-containment} shows that every
subsequential weak limit from these spaces belongs to \(M_\mu\). Both apply to the same union
sequence, so
Lemma~\ref{lem:finite-recovery-projections} yields
\begin{equation}\label{eq:finite-recovery-projections}
 \|P_{M_\ell}y-P_{M_\mu}y\|_{\D}\longrightarrow0
       \qquad\text{for every fixed }y\in\D.
\end{equation}
This is strong operator convergence. The recovery sequence used to
pass projection orthogonality to the limit can depend on the test
vector in \(M_\mu\). In particular, writing
\(p_\ell=P_{M_\ell}1\) and \(p=P_{M_\mu}1\)
gives \(p_\ell\to p\) in \(\D\).

The functions \(P_\ell\in M_\ell\) and \(P\in M_\mu\) from
\eqref{eq:finite-recovery-trial-membership} have norm at most one
and the same positive value \(c_\mu\) at zero by
\eqref{eq:finite-recovery-product-origin}. Evaluation at zero is
therefore nonzero on both spaces. Since \(1\) represents evaluation
at zero in the Dirichlet inner product, \(p_\ell\) and \(p\) are
nonzero. The normalization in \eqref{eq:canonical-definition}
gives the canonical functions and their values at zero.
\begin{equation}\label{eq:finite-recovery-canonical-projections}
 \begin{aligned}
 H_\ell&=G_{M_\ell}=\frac{p_\ell}{\|p_\ell\|_{\D}},
    &\gamma_\ell&=H_\ell(0)=\|p_\ell\|_{\D},\\
 G_\mu&=G_{M_\mu}=\frac{p}{\|p\|_{\D}},
    &\gamma_\mu&=G_\mu(0)=\|p\|_{\D}.
 \end{aligned}
\end{equation}
For \(f\) in either subspace, pairing against \(1\) equals pairing
against the corresponding projection of \(1\). Since
\(\langle f,1\rangle_{\D}=f(0)\), division by the projection norm gives
\[
 \langle f,H_\ell\rangle_{\D}=\frac{f(0)}{\gamma_\ell}\quad(f\in M_\ell),
 \qquad
 \langle f,G_\mu\rangle_{\D}=\frac{f(0)}{\gamma_\mu}\quad(f\in M_\mu).
\]
The target \(G_\mu\) is the function defined independently in
Section~\ref{subsec:02-atomic}. As \(K_\ell\) is a nonpolar
Carleson set, Section~\ref{subsec:03-seed} gives that \(H_\ell\)
is outer and \([H_\ell]=M_\ell\). The functions
\(H_\ell\) and \(G_\mu\), each of unit norm and positive at zero,
are contractive multipliers by
\eqref{eq:canonical-multiplier} and have modulus at most one in
the disk by Section~\ref{subsec:02-multiplier}.

Testing the norm of evaluation at zero on the two trial vectors
gives \(\gamma_\ell,\gamma_\mu\ge c_\mu\).
For the finite-stage upper bound,
\(M_\ell\subset\D_{K_{i,\ell}}\) by
\eqref{eq:finite-recovery-ideal-intersection}. The norm of evaluation
on the containing space is \(G_i(0)\) by
\eqref{eq:finite-recovery-single-input}, so
\(\gamma_\ell\le G_i(0)\le e^{-\tau_i}<1\) for each \(i\), by
\eqref{eq:one-atom-origin-bounds}. For the target,
Lemma~\ref{lem:atomic-witness} gives
\([G_\mu]=M_\mu\) and the exact factorization
\(G_\mu=S_\mu O_\mu\), where \(O_\mu\in H^2\) is outer with
\(O_\mu(0)>0\). Since \(|S_\mu^*|=1\) almost everywhere, we have
\(|O_\mu^*|=|G_\mu^*|\) almost everywhere and
\(\|O_\mu\|_{H^2}=\|G_\mu\|_{H^2}\le1\). Thus
\(\gamma_\mu=e^{-m_\mu}O_\mu(0)\le e^{-m_\mu}\), where
\(m_\mu=\mu(\T)>0\).
Multiplying the individual lower bounds in
\eqref{eq:one-atom-origin-bounds} now gives
\begin{equation}\label{eq:finite-recovery-origin-bounds}
 \begin{gathered}
  0<\kappa_\mu:=\frac{e^{-m_\mu}}{\prod_{i=1}^N\sqrt{3+2\tau_i}}
      \le c_\mu,\qquad
  c_\mu\le\gamma_\ell\le\min_i G_i(0)<1\quad(\ell\ge1),\\
  c_\mu\le\gamma_\mu\le e^{-m_\mu}<1.
 \end{gathered}
\end{equation}

The lower bound on \(\gamma_\mu\) controls the division needed
to compare the normalized projections. By the reverse triangle
inequality for \(p_\ell,p\), the difference of their norms is bounded
by \(\|p_\ell-p\|_{\D}\). Also
\[
 \gamma_\mu(H_\ell-G_\mu)
   =(p_\ell-p)+(\gamma_\mu-\gamma_\ell)H_\ell .
\]
Since \(\|H_\ell\|_{\D}=1\), these two observations give
\begin{equation}\label{eq:finite-recovery-normalized-convergence}
 \begin{aligned}
  |\gamma_\ell-\gamma_\mu|
    &\le\|p_\ell-p\|_{\D}\longrightarrow0,\\
  \|H_\ell-G_\mu\|_{\D}
    &\le \frac{2}{\gamma_\mu}\|p_\ell-p\|_{\D}\longrightarrow0.
 \end{aligned}
\end{equation}
The coefficient in the second line is bounded by
\(2/c_\mu\le2/\kappa_\mu\) uniformly in \(\ell\). These upper
bounds depend only on the list of masses, because each
\(G_{\tau_i,\xi_i}(0)\) is independent of \(\xi_i\) by the
rotation invariance established in Section~\ref{subsec:02-atomic}.

The residuals \(1-p_\ell=1-\gamma_\ell H_\ell\) belong to
\(M_\ell^\perp\) and converge in \(\D\) to
\(1-p=1-\gamma_\mu G_\mu\in M_\mu^\perp\).
By Pythagoras and \(\|1\|_{\D}=1\), we obtain
\begin{equation}\label{eq:finite-recovery-origin-gap}
 \begin{aligned}
  \dist_{\D}(1,M_\ell)^2
    &=\|1-p_\ell\|_{\D}^2=1-\gamma_\ell^2
      \ge1-(\min_i G_i(0))^2>0,\\
  \dist_{\D}(1,M_\ell)^2
    &\longrightarrow \dist_{\D}(1,M_\mu)^2
      =\|1-p\|_{\D}^2=1-\gamma_\mu^2
      \ge1-e^{-2m_\mu}>0 .
 \end{aligned}
\end{equation}

\subsection{Boundary zeros and boundary logarithms}\label{subsec:12-boundary}

\paragraph{Location and capacity of the  zero sets.}
Let \(E_\mu=\{\xi_1,\ldots,\xi_N\}\) be the support of the fixed
atomic measure. We use \(K_\ell\) from
\eqref{eq:finite-recovery-union} and the  canonical function
\(H_\ell=G_{\D_{K_\ell}}\) from
\eqref{eq:finite-recovery-canonical-projections}. The generic result in
Section~\ref{subsec:03-seed} gives outerness and
\([H_\ell]=\D_{K_\ell}\) because \(K_\ell\) is a nonpolar Carleson set. The canonical moments, unit norm, and positive value at zero of \(H_\ell\), together with the properness of \(\D_{K_\ell}\), give \(H_\ell\in\Xpr\).

We use chordal Hausdorff distance to compare \(K_\ell\) and
\(E_\mu\). The estimate has two directions. First,
every point of the union lies close to one of the atom directions.
For each \(\zeta\in K_\ell\), choose \(i\) with
\(\zeta\in K_{i,\ell}\). The one-atom bound
\eqref{eq:finite-recovery-single-geometry} then gives
\[
 \dist(\zeta,E_\mu)\le |\zeta-\xi_i|\le 3\cdot2^{-\ell}.
\]
Conversely, each direction \(\xi_i\) has a nonempty compact
\(K_{i,\ell}\subset K_\ell\) at distance at most \(3\cdot2^{-\ell}\).
Thus \(\dist(\xi_i,K_\ell)\le\dist(\xi_i,K_{i,\ell})
\le3\cdot2^{-\ell}\). The two directed bounds yield
\begin{equation}\label{eq:finite-recovery-union-location}
                   d_H(K_\ell,E_\mu)\le 3\cdot2^{-\ell}.
\end{equation}

Put $Z_\ell=\Zrad(H_\ell)$. By Theorem~\ref{thm:G} and
Corollary~\ref{cor:R}, each constituent $K_{i,\ell}$ is the image of $K$
under a disk automorphism. Hence $K_\ell$ in
\eqref{eq:finite-recovery-union} is a nonempty finite union of such
images. Applying the result on exact boundary zeros in
Section~\ref{subsec:11-realization} to $K_\ell$ gives a continuous
extension of $H_\ell$ to $\overline{\DD}$ and
\begin{equation}\label{eq:finite-recovery-radial-zeros}
 \{z\in\overline{\DD}:H_\ell(z)=0\}=Z_\ell=K_\ell.
\end{equation}
By \eqref{eq:finite-recovery-union-capacity},
\begin{equation}\label{eq:finite-recovery-zero-capacity}
 0<C_R(Z_\ell)=C_R(K_\ell)\le5N/\ell.
\end{equation}
The bound tends to zero because $N$ is fixed.

\paragraph{Boundary logarithms at the fixed atomic target.}
The outer mean formula gives the integrals of the two negative
logarithms. The singular inner factor affects the target's value at
zero and has boundary modulus one almost everywhere. An adaptation
of Scheff\'e's argument
\cite[pp.~435--437, especially footnote~2 on p.~436]{Scheffe47}
uses the pointwise minimum to compute the limiting \(L^1\) distance.
The Poisson formula gives the weak-* limit of the logarithmic measures.

\begin{proposition}[Boundary logarithms in finite atomic recovery]
\label{prop:finite-recovery-logarithmic-defect}
For the fixed \(\mu\), let \(H_\ell\) be the canonical functions of
the trace-zero spaces associated with the unions in
\eqref{eq:finite-recovery-union}. Then
\begin{equation}\label{eq:finite-recovery-logarithmic-defect}
 \bigl\|\log|H_\ell^*|-\log|G_\mu^*|\bigr\|_{L^1(m)}
                       \longrightarrow m_\mu=\mu(\T)>0 .
\end{equation}
The measures \((-\log|H_\ell^*|)\,dm\) converge weak-* against
\(C(\T)\) to \((-\log|G_\mu^*|)\,dm+\mu\), as in
\eqref{eq:finite-atomic-log-measure-limit}.
\end{proposition}

\begin{proof}
Set \(u_\ell=-\log|H_\ell^*|\) and \(u=-\log|G_\mu^*|\), using Hardy
boundary values almost everywhere for \(m\). The canonical bound
from Section~\ref{subsec:02-multiplier} gives \(u_\ell,u\ge0\).
The Hardy boundary and factorization facts in
Section~\ref{subsec:02-hardy} give nonzero boundary moduli almost
everywhere and integrable logarithms.

Write \(\gamma_\ell=H_\ell(0)\) and \(\gamma_\mu=G_\mu(0)>0\),
with \(\gamma_\ell\to\gamma_\mu\) by
\eqref{eq:finite-recovery-normalized-convergence}. Each \(H_\ell\)
is outer and positive at zero. For the target,
Lemma~\ref{lem:atomic-witness} gives the exact factorization
\(G_\mu=S_\mu O_\mu\) with \(O_\mu\) outer and positive at zero.
Since \(|S_\mu^*|=1\) almost everywhere and
\(S_\mu(0)=e^{-m_\mu}\), we have
\(|G_\mu^*|=|O_\mu^*|\) almost everywhere and
\(O_\mu(0)=e^{m_\mu}\gamma_\mu\). The outer mean formula
\eqref{eq:outer-intro} therefore gives
\begin{equation}\label{eq:finite-recovery-log-means}
 \begin{aligned}
 \int_\T u_\ell\,dm &= -\log\gamma_\ell
                            \longrightarrow -\log\gamma_\mu,\\
 \int_\T u\,dm &= -\log O_\mu(0)
                            =-\log\gamma_\mu-m_\mu .
 \end{aligned}
\end{equation}

The strong Dirichlet convergence in
\eqref{eq:finite-recovery-normalized-convergence} gives strong
\(L^2(m)\) convergence of the Hardy boundary values. The logarithm
can change rapidly when the moduli are small. We first set aside the
points where the target modulus is small and use \(L^2(m)\)
convergence on the rest. Fix \(\varepsilon>0\) and
\(0<\delta<1/2\). On the set where
\(|G_\mu^*|\ge2\delta\) and
\(|H_\ell^*-G_\mu^*|\le\delta\min\{1,\varepsilon\}\), both moduli
are at least \(\delta\). The logarithm is \(1/\delta\)-Lipschitz
there, so \(|u_\ell-u|\le\varepsilon\). Consequently
\[
 m\{|u_\ell-u|>\varepsilon\}
 \le m\{|G_\mu^*|<2\delta\}
    +m\{|H_\ell^*-G_\mu^*|>\delta\min\{1,\varepsilon\}\}.
\]
The second term tends to zero for fixed \(\delta\).
As \(\delta\downarrow0\), the first tends to zero because
\(|G_\mu^*|>0\) almost everywhere. Thus \(u_\ell\to u\) in measure.

To compute the limiting \(L^1\) distance, we integrate the smaller
of the two densities at each point. We have
\(\min\{u_\ell,u\}\to u\) in measure because \(u_\ell\to u\)
in measure, and \(0\le\min\{u_\ell,u\}\le u\in L^1(m)\).
From every subsequence
one can take a further subsequence on which the minima converge
almost everywhere. Dominated convergence gives convergence of
their integrals along that further subsequence, so the integrals
converge to \(\int_\T u\,dm\) as \(\ell\to\infty\). Finally,
\[
 \|u_\ell-u\|_{L^1(m)}
    =\int_\T u_\ell\,dm+\int_\T u\,dm
         -2\int_\T\min\{u_\ell,u\}\,dm
    \longrightarrow m_\mu
\]
by \eqref{eq:finite-recovery-log-means}.
For the weak-* limit, \eqref{eq:finite-recovery-log-means} bounds the
total masses of the positive measures \(u_\ell\,dm\). Since \(C(\T)\)
is separable, every subsequence has a further weak-* convergent
subsequence. Write \(u_{\ell_j}\,dm\overset{*}{\rightharpoonup}\lambda\)
for any such subsequence. For each \(z\in\DD\), the Poisson kernel
\(P_z\) is continuous on \(\T\), so the outer representation gives
\[
 \int_\T P_z\,d\lambda
   =\lim_{j\to\infty}\int_\T P_z u_{\ell_j}\,dm
   =\lim_{j\to\infty}-\log|H_{\ell_j}(z)|
   =-\log|G_\mu(z)|.
\]
The last equality follows from
\eqref{eq:finite-recovery-normalized-convergence} and the fact that
\(G_\mu=S_\mu O_\mu\) is zero-free in \(\DD\). This factorization also
gives
\[
 -\log|G_\mu(z)|=\int_\T P_z\,d(u\,m+\mu),
\]
since \(|O_\mu^*|=|G_\mu^*|\) almost everywhere. The Poisson transform
determines the Fourier coefficients of its measure. Hence each
weak-* subsequential limit is \(u\,m+\mu\), and the full sequence
converges weak-* to this measure.
\end{proof}

Section~\ref{sec:13} will use the atomic target as an intermediate
point for the Dirichlet norm. To control the logarithms, it will apply
Proposition~\ref{prop:outer-continuity} directly to the original outer
canonical function and its outer approximant.

\subsection{Proof and extensions of Approximation B}\label{subsec:12-assembly}

For the fixed measure \(\mu\), the projections of \(1\) converge
in \(\D\), and the limiting projection has positive norm.
Normalizing gives convergence of the canonical functions in \(\D\).
Each \(H_\ell\) extends continuously to the closed disk with zero
set \(K_\ell\). These sets converge in Hausdorff distance to the
atomic support, and their capacities tend to zero. The two limits
for the boundary logarithms are proved in
Section~\ref{subsec:12-boundary}.

\begin{proof}[Proof of Theorem~\ref{thm:B}]
For the fixed measure \(\mu\), take \(K_{\mu,\ell}=K_\ell\) from
\eqref{eq:finite-recovery-union} and \(H_{\mu,\ell}=H_\ell\).
Each union uses the same index \(\ell\) in every one-atom family
from Corollary~\ref{cor:R}. The bounds for those families hold from
the first index, and the estimates above keep the union sequence
fixed throughout. Each union is nonpolar and Carleson, and its
canonical function belongs to \(\Xpr\) and generates \(\D_{K_\ell}\).
For the spaces \(M_\ell\), the two subspace conditions give the
projection convergence
\eqref{eq:finite-recovery-projections} on this sequence, for every
fixed \(y\). Normalization gives
\eqref{eq:finite-recovery-normalized-convergence}.

Section~\ref{subsec:12-boundary} gives the continuous extension of
each \(H_\ell\) to the closed disk. The location, zero-set, and
capacity conclusions are \eqref{eq:finite-recovery-union-location},
\eqref{eq:finite-recovery-radial-zeros}, and
\eqref{eq:finite-recovery-zero-capacity}, respectively.

Lemma~\ref{lem:atomic-witness} gives the exact inner factor of
\(G_\mu\). Proposition~\ref{prop:finite-recovery-logarithmic-defect}
gives the \(L^1\) distance and the weak-* limit along this same
sequence.
\end{proof}

For fixed \(\mu\), its chosen families, and \(\delta,\varepsilon>0\),
both inequalities below hold for every sufficiently large integer \(\ell\).
The first follows from \eqref{eq:finite-recovery-normalized-convergence}.
For the second, use \(C_R(K_\ell)\le 5N/\ell\) and take
\(\ell>5N/\varepsilon\). Then
\begin{equation}\label{eq:finite-recovery-simultaneous-tolerances}
 \|H_\ell-G_\mu\|_{\D}<\delta,\qquad
 C_R(\Zrad(H_\ell))=C_R(K_\ell)<\varepsilon .
\end{equation}
The joining argument in Section~\ref{sec:13} will use these estimates.

\paragraph*{Recovery for a general singular inner factor.}
Let $\mu$ be a nonzero finite positive singular measure on $\T$ with
$M_\mu=\D\cap S_\mu H^2\ne\{0\}$, and put $E_\mu=\supp\mu$.
Then the canonical function $G_\mu=G_{M_\mu}$ is defined,
$[G_\mu]=M_\mu$, and its exact inner factor is $S_\mu$.
Given positive sequences $\epsilon_j,\delta_j,\eta_j$ tending to zero,
there are nonpolar Carleson sets $K_j$, each a finite union of disk
automorphism images of the compact set $K$ from
\eqref{eq:fixed-packed-compact}, whose canonical functions
$H_j=G_{\D_{K_j}}$ are outer and continuous on $\overline{\DD}$.
On this sequence,
\[
 [H_j]=\D_{K_j},\qquad
 \{z\in\overline{\DD}:H_j(z)=0\}=\Zrad(H_j)=K_j ,
\]
\[
 \|H_j-G_\mu\|_{\D}<\eta_j,\qquad
 d_H(K_j,E_\mu)<\delta_j,\qquad 0<C_R(K_j)<\epsilon_j,
\]
and $P_{\D_{K_j}}y\to P_{M_\mu}y$ in $\D$ for every $y\in\D$.
With $m$ the normalized arclength measure and the stars denoting
Hardy boundary values,
\begin{align*}
 (-\log|H_j^*|)\,dm
 &\overset{*}{\rightharpoonup}(-\log|G_\mu^*|)\,dm+\mu,\\
 \big\|\log|H_j^*|-\log|G_\mu^*|\big\|_{L^1(m)}
 &\longrightarrow\mu(\T).
\end{align*}
Weak-* convergence is tested against $C(\T)$.

\begin{proof}
Take a nonzero $f\in M_\mu$. If $f$ vanishes to order $k$ at zero,
then so does $f/S_\mu$, since $S_\mu(0)>0$. The coefficient
norms show that $f/z^k\in\D$ and $(f/S_\mu)/z^k\in H^2$.
Thus $f/z^k$ belongs to $M_\mu$ and is nonzero at zero. Lemma~\ref{lem:exact-inner-factor}
gives $G_\mu=S_\mu h$ with $h$ outer, $h(0)>0$, and
$[G_\mu]=M_\mu$. Lemma~\ref{lem:inner-division} gives $h\in\D$.
The boundary modulus of $h$ equals that of $G_\mu$, which is at
most one by \eqref{eq:canonical-multiplier}. Hence $\|h\|_\infty\le1$.

For $\xi\in\T$ put
\[
 W_h(\xi)=2\int_\T\frac{|h^*(\zeta)|^2}{|\zeta-\xi|^2}\,dm(\zeta).
\]
For every finite positive singular measure $\sigma$, apply Carleson's
formula \cite{Carleson60} in the form given by Hartz and Richter
\cite[(3.2)--(3.3)]{HR26}. Tonelli's theorem gives
\[
 \|S_\sigma h\|_{\D}^2
   =\|h\|_{\D}^2+\int_\T W_h\,d\sigma,
\]
with the left side interpreted as $+\infty$ when $S_\sigma h\notin\D$.
With $\sigma=\mu$ the left side is one, so $W_h\in L^1(\mu)$.
For each $n$ partition $\T$ into finitely many half-open arcs of
chord diameter at most $2^{-n}$. In each arc $I$ of positive
$\mu$-mass, choose $\xi_I\in I\cap E_\mu$ with
\[
 W_h(\xi_I)\le \frac{1}{\mu(I)}\int_I W_h\,d\mu.
\]
Such a point exists because $\mu$ is supported on $E_\mu$ and the
integral of $W_h$ would exceed its own average times $\mu(I)$ if
$W_h$ were strictly above the average almost everywhere on $I$.
Set $\mu_n=\sum_I\mu(I)\delta_{\xi_I}$, summing over the arcs of
positive mass. Then $\mu_n(\T)=\mu(\T)$ and
$\int W_h\,d\mu_n\le\int W_h\,d\mu$. For each continuous test
function the integral error is at most $\mu(\T)$ times its
modulus of continuity at $2^{-n}$, so $\mu_n\to\mu$ weak-*.
If $\xi\in E_\mu$, every chordal ball about $\xi$ of radius $r>0$
meets some positive-mass arc in a set of positive $\mu$-mass.
The atom selected in that arc is within $r+2^{-n}$ of $\xi$.
Since $r$ is arbitrary and $\supp\mu_n\subset E_\mu$, this gives
$d_H(\supp\mu_n,E_\mu)\le2^{-n}$.

The functions $S_{\mu_n}h$ belong to the unit ball of $M_{\mu_n}$
and satisfy $(S_{\mu_n}h)(0)=G_\mu(0)$ since all the masses
agree. Extremality yields $G_{\mu_n}(0)\ge G_\mu(0)$.
We identify the weak limits as in \cite[proof of Lemma~5.1]{HR26}.
Take a weakly convergent Dirichlet subsequence with limit $g$ and write
$G_{\mu_n}=S_{\mu_n}h_n$ along it. The Hardy norms of $h_n$ are
at most one, by the isometry of multiplication by an inner function
in $H^2$. Passing to a further weak Hardy subsequence, let $h_0$ be
its limit. For every $z\in\DD$, the singular factor formula and
weak-* convergence of $\mu_n$ give $S_{\mu_n}(z)\to S_\mu(z)$.
Continuity of point evaluations in both spaces yields
$g=S_\mu h_0\in M_\mu$. The bounds $\|g\|_{\D}\le1$ and
$g(0)\ge G_\mu(0)>0$ identify $g$ with the unique canonical
extremal. Every subsequence has a further weakly convergent
subsequence, so $G_{\mu_n}\rightharpoonup G_\mu$ in $\D$. Since
their norms and the norm of $G_\mu$ are one, convergence is strong.

For each $j$ first choose $n$ such that
$\|G_{\mu_n}-G_\mu\|_{\D}<\eta_j/2$ and
$d_H(\supp\mu_n,E_\mu)<\delta_j/2$. For this fixed $\mu_n$,
Theorem~\ref{thm:B} gives all sufficiently large $\ell$ with
norm error below $\eta_j/2$, Hausdorff distance below $\delta_j/2$,
and capacity below $\epsilon_j$. Set $H_j=H_{\mu_n,\ell}$ and
$K_j=K_{\mu_n,\ell}$ for one such $\ell$. The triangle
inequalities give the claimed errors. Theorem~\ref{thm:B} gives
generation, continuity, and the exact radial zero set, and
\eqref{eq:finite-recovery-union} gives the description of $K_j$.

For each polynomial $p$ we have $pH_j\to pG_\mu$ in $\D$.
For $v\in M_\mu=[G_\mu]$, it follows that
\[
 \limsup_{j\to\infty}\dist(v,\D_{K_j})
 \le \inf_{p\,\mathrm{polynomial}}\|v-pG_\mu\|_{\D}=0.
\]
Thus $P_{\D_{K_j}}v\to v$ in $\D$ for every $v\in M_\mu$.
To identify weak limits, suppose
$v_k\in\D_{K_{j_k}}$ and $v_k\rightharpoonup v$ in $\D$ for
a strictly increasing sequence $j_k$. By
Lemma~\ref{lem:completed-canonical-quotient},
$v_k=H_{j_k}q_k$ with $\|q_k\|_{H^2}\le\|v_k\|_{\D}$.
Take a weak Hardy subsequence with limit $q$. Point evaluations
and strong convergence of $H_j$ identify
$v=G_\mu q=S_\mu(hq)\in M_\mu$, since $h\in H^\infty$.
Lemma~\ref{lem:finite-recovery-projections} then gives
the convergence of the projections.

Finally, put $u_j=-\log|H_j^*|$ and $u=-\log|G_\mu^*|$.
They are nonnegative and finite almost everywhere. Hardy
boundary convergence and $|G_\mu^*|>0$ almost everywhere
give $u_j\to u$ in measure. The outer mean formula and the
exact factorization give
\[
 \int_\T u_j\,dm=-\log H_j(0)\longrightarrow
 -\log G_\mu(0)=\int_\T u\,dm+\mu(\T).
\]
Since $0\le\min(u_j,u)\le u\in L^1(m)$, dominated convergence on
almost-everywhere convergent further subsequences gives
$\int\min(u_j,u)\,dm\to\int u\,dm$. The identity
$|u_j-u|=u_j+u-2\min(u_j,u)$ now gives the stated $L^1$ limit.
The positive measures $u_j\,dm$ have bounded mass, and $C(\T)$
is separable. Thus each subsequence has a further weak-* convergent
subsequence. For any such limit $\lambda$ and any $z\in\DD$, the
outer formula and convergence at $z$ give
\[
 \int_\T P_z\,d\lambda
   =\lim -\log|H_j(z)|
   =-\log|G_\mu(z)|
   =\int_\T P_z\,d(u\,dm+\mu),
\]
where the limit is along the chosen subsequence.
Uniqueness of the Poisson integral identifies every such
$\lambda$ with $u\,dm+\mu$, giving convergence of the full sequence.
\end{proof}

\begin{corollary}[Cyclicity of the outer factor for finitely supported measures]
\label{cor:finite-recovery-outer-cofactor}
Let \(\mu\) be a nonzero finite positive measure of finite support.
In the exact factorization \(G_\mu=S_\mu O_\mu\) from
Lemma~\ref{lem:atomic-witness}, the outer factor \(O_\mu\) belongs
to \(\D\) and is cyclic there.
\end{corollary}

\begin{proof}
Lemma~\ref{lem:atomic-witness} identifies \([G_\mu]=M_\mu\)
and supplies \(S_\mu Q\in M_\mu\), where
\(Q(z)=\prod_{\xi\in\supp\mu}(1-\bar\xi z)\).
Choose polynomials \(p_j\) such that \(p_jG_\mu\to S_\mu Q\) in \(\D\).
The Hardy factorization gives \(O_\mu\in H^2\), so
Lemma~\ref{lem:inner-division} applied to \(G_\mu\) gives \(O_\mu\in\D\).
Since \(p_j O_\mu-Q\in H^2\) and its product with \(S_\mu\)
is \(p_jG_\mu-S_\mu Q\in\D\), the same lemma gives
\[
 \|p_jO_\mu-Q\|_{\D}\le\|p_jG_\mu-S_\mu Q\|_{\D}\longrightarrow0.
\]
Thus \(Q\in[O_\mu]\). All zeros of \(Q\) lie on the circle, so
\(Q\) is cyclic \cite[corollary to Lemma~8, p.~289]{BS84}.
Consequently \(\D=[Q]\subset[O_\mu]\).
\end{proof}

\begin{proposition}[Recovery with a prescribed singular divisor]
\label{prop:retain-singular-factor}
Let $\mu$ be a nonzero finite positive singular measure on $\T$ with
$M_\mu=\D\cap S_\mu H^2\ne\{0\}$, and let $0\le\nu\le\mu$. Put
$\rho=\mu-\nu$, and use $M_\nu=\D\cap S_\nu H^2$, with $M_0=\D$.
If $\rho\ne0$, then for any positive sequences $\delta_j,\epsilon_j$
tending to zero there are nonpolar Carleson sets $K_j$, each a finite union of disk automorphism images of the compact set $K$ in
\eqref{eq:fixed-packed-compact}, such that
\[
 d_H(K_j,\supp\rho)<\delta_j,\qquad 0<C_R(K_j)<\epsilon_j.
\]
For $L_j=M_\nu\cap\D_{K_j}$, the canonical function $g_j=G_{L_j}$
is defined, generates $L_j$, and has exact Hardy inner factor $S_\nu$
under the normalization $S_\nu(0)>0$. On this sequence,
\[
 \|g_j-G_\mu\|_{\D}\longrightarrow0,
 \qquad \|P_{L_j}y-P_{M_\mu}y\|_{\D}\longrightarrow0
 \quad (y\in\D),
\]
and the Hardy boundary logarithms satisfy
\begin{align*}
 (-\log|g_j^*|)\,dm
 &\overset{*}{\rightharpoonup}(-\log|G_\mu^*|)\,dm+\rho,\\
 \big\|\log|g_j^*|-\log|G_\mu^*|\big\|_{L^1(m)}
 &\longrightarrow\rho(\T).
\end{align*}
Weak-* convergence is tested against $C(\T)$. If $\rho=0$, the same
function, projection, and logarithm conclusions hold with
$L_j=M_\mu$ and $g_j=G_\mu$ for all $j$.
\end{proposition}

\begin{proof}
The case $\rho=0$ follows by taking the constant sequence, so suppose
$\rho\ne0$. The preceding general singular recovery gives
$G_\mu(0)>0$. Since $S_\mu=S_\nu S_\rho$, the function $G_\mu$ belongs
to both $M_\nu$ and $M_\rho$. In particular, evaluation at zero is
nonzero on each space. By Lemma~\ref{lem:exact-inner-factor}, write
\[
 F=G_{M_\nu}=S_\nu U,\qquad G_\rho=G_{M_\rho}=S_\rho V,
\]
with $U,V$ outer. This includes $F=U=1$ when $\nu=0$.
Contractive inner division gives $U,V\in\D$. Their boundary moduli
are those of $F$ and $G_\rho$, so $\|U\|_\infty,\|V\|_\infty\le1$
by the canonical multiplier bound and the Hardy maximum principle.
Apply the preceding general singular recovery to $\rho$ to obtain
$K_j$ with the stated location and capacity bounds and outer canonical functions
$H_j=G_{\D_{K_j}}$ satisfying
\[
 [H_j]=\D_{K_j},\qquad H_j\longrightarrow G_\rho\text{ in }\D,
 \qquad P_{\D_{K_j}}y\longrightarrow P_{M_\rho}y\text{ in }\D
 \quad(y\in\D).
\]

All products below belong to $\D$ by Lemma~\ref{lem:bounded-products}.
The outer product identity \cite[Theorem~4.5]{RS92} and
Theorem~\ref{input:module} give
\begin{align*}
 [FH_j]
   &=[UH_j]\cap S_\nu H^2
     =[U]\cap[H_j]\cap S_\nu H^2
     =M_\nu\cap\D_{K_j}=L_j,\\
 [FG_\rho]
   &=[UV]\cap S_\mu H^2
     =[U]\cap[V]\cap S_\mu H^2=M_\mu.
\end{align*}
Indeed, the inner--outer identity gives $M_\nu=[U]\cap S_\nu H^2$ and
$M_\rho=[V]\cap S_\rho H^2$. Thus $M_\mu\subset[U]\cap[V]$, and the
reverse inclusion in the second identity follows from the factor
$S_\mu$ and membership in $\D$.

Since $(FH_j)(0)>0$, projection of $1$ onto $L_j$ is nonzero and its
unit normalization $g_j$ generates $L_j$. Every member of $L_j$ lies
in $S_\nu H^2$, so the inner factor of $g_j$ is divisible by $S_\nu$.
Conversely, the closed subspace of $H^2$ defined by that inner factor
contains $[g_j]=L_j$ by the continuous inclusion $\D\subset H^2$.
It contains $FH_j$, whose exact inner factor is $S_\nu$ because $U$ and
$H_j$ are outer. This proves that the inner factor of $g_j$ is $S_\nu$
with the stated normalization.

Multiplication by the fixed canonical function $F$ is bounded on
$\D$, so $FH_j\to FG_\rho$ in $\D$. For each polynomial $p$, the
vectors $pFH_j$ belong to $L_j$ and converge to $pFG_\rho$ in $\D$.
Since those limits are dense in $M_\mu$, we get
$P_{L_j}v\to v$ in $\D$ for every $v\in M_\mu$.

To identify weak limits, suppose that $x_k\in L_{j_k}$ and
$x_k\rightharpoonup x$ in $\D$ for strictly increasing $j_k$.
The closed subspace $M_\nu$ contains $x$. For each $k$, approximate
$x_k$ in $\D$ by polynomial multiples of $FH_{j_k}$. Dividing by
$S_\nu$ and using Lemma~\ref{lem:inner-division} shows that
$x_k/S_\nu$ is a Dirichlet norm limit of polynomial multiples of
$UH_{j_k}$. Hence $x_k/S_\nu\in[H_{j_k}]=\D_{K_{j_k}}$.
Division by $S_\nu$ is bounded on $M_\nu$, so
$x_k/S_\nu\rightharpoonup x/S_\nu$ weakly in $\D$.
For $y\perp M_\rho$, we have
\[
 \left\langle x_k/S_\nu,y\right\rangle_{\D}
 =\left\langle x_k/S_\nu,P_{\D_{K_{j_k}}}y\right\rangle_{\D}
 \longrightarrow0.
\]
Thus $x/S_\nu\in M_\rho$, and $x\in M_\mu$ because
$S_\nu S_\rho=S_\mu$. This argument also covers measures $\nu$ and
$\rho$ with overlapping supports. Lemma~\ref{lem:finite-recovery-projections}
now gives $P_{L_j}y\to P_{M_\mu}y$ for every $y\in\D$. Apply this to
$1$ and normalize the projections to conclude $g_j\to G_\mu$.

For the logarithms, set $u_j=-\log|g_j^*|$ and
$u=-\log|G_\mu^*|$. These functions are nonnegative and integrable.
For $z\in\DD$, the inner factors and the outer representation give
\[
 -\log|g_j(z)|=\int_{\T}P_z\,d(u_jm+\nu),
 \qquad -\log|G_\mu(z)|=\int_{\T}P_z\,d(um+\mu).
\]
At zero, $g_j(0)\to G_\mu(0)>0$, so
\[
 \int_{\T}u_j\,dm-\int_{\T}u\,dm
   =\rho(\T)+\log\frac{G_\mu(0)}{g_j(0)}
   \longrightarrow\rho(\T).
\]
The positive measures $u_j\,dm$ have bounded mass. For any weak-*
subsequential limit $\lambda$, convergence at each $z\in\DD$ and
continuity of $P_z$ on $\T$ show that the Poisson integrals of
$\lambda+\nu$ and $u\,dm+\mu$ are equal. Uniqueness of the Poisson
integral gives $\lambda=u\,dm+\rho$. Since $C(\T)$ is separable, this
identifies the weak-* limit of the full sequence.

Dirichlet convergence also gives $H^2$ convergence, and therefore
$u_j\to u$ in measure since $|G_\mu^*|>0$ almost everywhere.
The bound $0\le\min(u_j,u)\le u\in L^1(m)$ implies
$\int_{\T}\min(u_j,u)\,dm\to\int_{\T}u\,dm$ by dominated convergence
on almost-everywhere convergent further subsequences. Integrating
$|u_j-u|=u_j+u-2\min(u_j,u)$ proves the asserted $L^1$ limit.
\end{proof}

\begingroup
\hyphenpenalty=10000\relax

\newpage \section{Constructing the counterexample}\label{sec:13}
\endgroup

The two approximation theorems now give a way to find outer canonical
functions with arbitrarily small radial-zero capacity near any given
proper outer canonical function. The intermediate canonical function belongs
to a space defined by a finite atomic measure. It connects the two
approximations in Dirichlet norm. The estimate in Section~\ref{subsec:04-continuity} then bounds the
\(L^1(m)\) distance between the boundary logarithms of the given outer
canonical function and its outer approximant.

We then choose these nearby functions as centers of nested closed
balls in the metric \(\dlog\). A strict capacity bound at a center
lets us choose a radius for which the bound holds throughout the ball.
Completeness supplies a point in all the balls. The shrinking capacity
bounds give this point a polar radial zero set, and the initial ball
keeps its value at the origin below \(1\). The canonical projection
formula then gives a positive distance from \(1\) to the subspace
generated by this same function.

\subsection{Density at an outer target}\label{subsec:13-density}

Openness and density here are relative to the normalized outer
canonical class \(\X\), equipped with the logarithmic metric
\(\dlog\). The two approximation theorems supply nearby functions
for each proper target. The constant function \(1\) already
belongs to every strict capacity sublevel with a positive
threshold, since its radial zero set is empty. For proper targets,
the strict bounds below on the approximants' values at zero keep them in the
proper class and leave room for the later closed neighborhoods.

\begin{proposition}[Density of small radial-zero capacity]
\label{prop:small-capacity-density}
For every \(\varepsilon>0\), the strict sublevel
\begin{equation}\label{eq:strict-capacity-sublevels}
 \mathcal U_\varepsilon=\{f\in\X:\Phi(f)<\varepsilon\},
 \qquad \Phi(f)=C_R(\Zrad(f)),
\end{equation}
is open and dense in \((\X,\dlog)\). For a proper target \(f\in\Xpr\)
with \(a=f(0)\) and any \(R>0\), one can choose
\(g=G_{\D_K}\in\Xpr\) for a nonpolar Carleson set \(K\).
The function \(g\) extends continuously to \(\overline{\DD}\) and
satisfies \(\{z\in\overline{\DD}:g(z)=0\}=\Zrad(g)=K\). It also
satisfies
\[
 [g]=\D_K,\qquad \dlog(g,f)<R,\qquad
 \Phi(g)=C_R(K)<\varepsilon,\qquad
 \frac a2<g(0)<\frac{1+a}{2}.
\]
\end{proposition}

\begin{proof}
Fix \(f\in\Xpr\), \(R>0\), and \(\varepsilon>0\), and put
\(a=f(0)\in(0,1)\). We first decide how close the eventual outer
function \(g\) must be to \(f\) in Dirichlet norm. In
Proposition~\ref{prop:outer-continuity}, the quantity \(T_f(t)\)
measures the integral of \(-\log|f^*|\) where \(|f^*|<t\).
Choose \(0<t<1\) so that \(T_f(t)<R/8\). This makes the tail
contribution to the metric bound less than \(R/4\). With \(t\)
fixed, set
\[
 \eta=\min\!\left\{\frac a2,\ \frac{1-a}{2},\
          \frac{R}{2(1+2/a+2/t)}\right\}>0.
\]
The first two entries in the minimum allow for the two strict
margins at the origin. The last makes the part of the metric estimate
proportional to the norm error less than \(R/2\). The choice of
\(\eta\) uses only \(f\) and \(R\).

We allow half of this bound for each approximation. Apply Theorem~\ref{thm:A} with a bound of \(\eta/2\) on the norm error and any positive bound on the total mass. It gives a nonzero finite positive measure \(\mu\) of finite support. Its canonical function \(G_\mu\) in \(M_\mu=S_\mu H^2\cap\D\) satisfies \(\|G_\mu-f\|_{\D}<\eta/2\).
Keep this measure fixed and write \(N\) for the number of distinct atoms.

Theorem~\ref{thm:B} gives a sequence with this same \(G_\mu\) as
its norm limit. By \eqref{eq:finite-recovery-simultaneous-tolerances} with
\(\delta=\eta/2\), choose an integer \(\ell>5N/\varepsilon\)
large enough that \(\|H_{\mu,\ell}-G_\mu\|_{\D}<\eta/2\). Set \(K=K_{\mu,\ell}\) and
\(g=H_{\mu,\ell}\) at this index. By Theorem~\ref{thm:B},
\([g]=\D_K\), and \(g\) extends continuously to
\(\overline{\DD}\) with
\(\{z\in\overline{\DD}:g(z)=0\}=\Zrad(g)=K\). Hence
\[
 \Phi(g)=C_R(K)\le\frac{5N}{\ell}<\varepsilon.
\]

The triangle inequality through the fixed atomic target gives
\[
 \|g-f\|_{\D}\le\|g-G_\mu\|_{\D}+\|G_\mu-f\|_{\D}<\eta.
\]
The intermediate \(G_\mu\) has nontrivial inner factor \(S_\mu\).
Since \(f,g\in\X\) and \(\|g-f\|_{\D}<\eta\le a/2\),
Proposition~\ref{prop:outer-continuity} applies with \(f\) fixed
and the previously chosen \(t\). Inequality~\eqref{eq:outer-continuity} gives
\[
 \dlog(f,g)\le(1+2/a+2/t)\|g-f\|_{\D}+2T_f(t)<3R/4<R.
\]
Evaluation at zero has norm one in \(\D\). Therefore
\(|g(0)-a|<\eta\), and the first two entries defining \(\eta\)
yield the two strict bounds on \(g(0)\). The upper bound and 
the identity \eqref{eq:certificate} give
\[
 \operatorname{dist}_{\D}(1,[g])^2=1-g(0)^2
 >1-\left(\frac{1+a}{2}\right)^2>0.
\]
The only point of \(\X\) outside \(\Xpr\) is \(1\), which
already belongs to \(\mathcal U_\varepsilon\) because its radial
zero set is empty. Thus every relative metric ball about any
point of \(\X\) contains a point of \(\mathcal U_\varepsilon\).

For openness, let \(h\in\mathcal U_\varepsilon\). The strict
inequality \(\Phi(h)<\varepsilon\) leaves room for an increase in
capacity under perturbation. Corollary~\ref{cor:capacity-ball}
uses that room to give a closed Dirichlet ball of positive radius
centered at \(h\) on which \(\Phi<\varepsilon\). The closed relative
\(\dlog\) ball of the same radius lies in that Dirichlet ball
because \(\|g-h\|_{\D}\le\dlog(g,h)\) for \(g\in\X\). In
particular \(h\) has an open relative neighborhood in
\(\mathcal U_\varepsilon\). The closed relative ball itself is
what we will use when taking a limit of centers.
\end{proof}

By the Baire category theorem, the intersection of the open dense
sets \(\mathcal U_{2^{-n}}\) is dense in \((\X,\dlog)\). Its
members have radial-zero capacity zero. Every canonical function in a ball on which the values at zero
are uniformly bounded away from one is noncyclic.
The nested-ball argument below selects trace-zero generators
\(f_n\) and gives the additional estimate
\(\dlog(f_n,f)\le r_0\,4^{-n}\) for the distance to their limit.
B\'en\'eteau, Ivrii, Manolaki, and Seco \cite[Theorem~4.1]{BIMS20} use category
for universality of optimal approximants on a prescribed compact polar set. Here the dense sets come from the two approximations,
and the sets defining the centers vary.

\subsection{Nested closed balls}\label{subsec:13-balls}

We use the nested closed balls from the usual proof of the Baire
category theorem. The constant-function example in
Section~\ref{subsec:03-neighborhood} illustrates that radial-zero
capacity can increase at a norm limit in the ambient Dirichlet
space. In the present construction, a bound on every point of
each closed ball applies directly to its eventual common point,
including when that point lies on the boundary. Density supplies
each new center, and Corollary~\ref{cor:capacity-ball} supplies a
radius with the required capacity bound throughout the ball.
Both the radii and the capacity bounds tend to zero.

Write \(B_{\dlog}(g,r)\) and \(\overline B_{\dlog}(g,r)\) for the
open and closed balls in \(\X\). The first ball will keep the
distance from \(1\) to the generated spaces bounded away from zero.
Fix \(f_0\in\Xpr\), put \(a_0=f_0(0)\), and choose
\begin{equation}\label{eq:local-initial-radius}
 0<r_0<\tfrac14\min\{a_0,1-a_0\}.
\end{equation}
For \(B_0=\overline B_{\dlog}(f_0,r_0)\), evaluation at zero and
\eqref{eq:metric-origin} give
\[
 0<a_0-r_0\le g(0)\le a_0+r_0<1\qquad(g\in B_0).
\]
The upper endpoint is strictly less than \(1\), so every \(g\)
in \(B_0\) belongs to \(\Xpr\). The projection identity
\eqref{eq:certificate} then gives
\(1-(a_0+r_0)^2>0\) as a lower bound for
\(\dist_{\D}(1,[g])^2\) throughout \(B_0\).

For \(n\ge1\), suppose \(f_{n-1}\) and \(r_{n-1}>0\) have been chosen with
\(B_{n-1}\subset B_0\), which holds at the initial step. We choose
the next center inside the open ball of half the previous radius.
This leaves room for a closed ball around it. Use
Proposition~\ref{prop:small-capacity-density} at \(f_{n-1}\) with
\(R=r_{n-1}/2\) and \(\varepsilon=2^{-n}\) to choose
\[
 f_n\in\mathcal U_{2^{-n}}\cap
       B_{\dlog}(f_{n-1},r_{n-1}/2),
 \qquad b_n:=2^{-n}-\Phi(f_n)>0 .
\]
Here \(\mathcal U_\varepsilon\) is the strict sublevel in
\eqref{eq:strict-capacity-sublevels}. Since \(f_{n-1}\in B_0\)
is proper, the proposition lets us choose \(f_n\) as the canonical
generator of a trace-zero space.

The positive margin \(b_n\) leaves room for two errors. The first
comes from replacing the zero set of the trace of \(f_n\) by the
set \(E_{t_n}(f_n)\) where its modulus is at most \(t_n\).
The second comes from changing \(f_n\) to a nearby function.
With \(f_n\) and \(b_n\) fixed, use the near-zero limit in
\eqref{eq:capacity-perturbation} to choose \(t_n>0\) so that
\[
 C_R(E_{t_n}(f_n))\le \Phi(f_n)+b_n/4 .
\]
With \(t_n\) chosen, we next select the radius. The term \(r_{n-1}/4\) leaves room
inside the previous ball, \(2^{-n}\) imposes an absolute size bound,
and the last term makes the squared perturbation term at most
\(b_n/4\). For each \(n\ge1\), define
\begin{equation}\label{eq:local-ball-radii}
 r_n=\min\{r_{n-1}/4,\ 2^{-n},\ t_n\sqrt{b_n}/2\},
 \qquad B_n=\overline B_{\dlog}(f_n,r_n).
\end{equation}
The radius is positive. For every \(g\in B_n\), the bound
\(\|g-f_n\|_{\D}\le\dlog(g,f_n)\le r_n\) holds, including
when \(\dlog(g,f_n)=r_n\).
Applying the squared perturbation bound in
\eqref{eq:capacity-perturbation} gives
\begin{equation}\label{eq:local-ball-capacity}
 \Phi(g)
 \le C_R(E_{t_n}(f_n))+\frac{r_n^2}{t_n^2}
 \le \Phi(f_n)+\frac{b_n}{4}+\frac{b_n}{4}
 =2^{-n}-\frac{b_n}{2}<2^{-n}.
\end{equation}
Thus half the original margin remains on the whole closed ball.
To show that this ball also lies inside the previous one, add the
distance from \(g\) to \(f_n\) to the distance between the centers.
The triangle inequality gives, for every \(g\in B_n\),
\[
 \dlog(g,f_{n-1})
 \le \dlog(g,f_n)+\dlog(f_n,f_{n-1})
 < r_{n-1}/4+r_{n-1}/2<r_{n-1}.
\]
The capacity estimate and this inclusion give the inclusions below.
Iterating \(r_n\le r_{n-1}/4\) gives the radius bound.
\begin{equation}\label{eq:local-nested-balls}
 B_n\subset\mathcal U_{2^{-n}}\cap
      B_{\dlog}(f_{n-1},r_{n-1})\subset B_{n-1},
 \qquad 0<r_n\le r_0\,4^{-n}.
\end{equation}
In particular \(B_n\subset B_0\), so the argument can be repeated
at the next step.

\subsection{Counterexamples near noncyclic outer canonical functions}\label{subsec:13-limit}

Completeness gives a function in all the closed balls. The initial ball controls its value at zero, and the later balls force the capacity of its radial zero set to be zero.

\begin{theorem}[Local strengthening]\label{thm:local}
Let \(f_0\in\Xpr\), put \(a_0=f_0(0)\), and let \(r_0\) satisfy
\eqref{eq:local-initial-radius}. There is \(f\in\Xpr\) such that
\[
 \dlog(f,f_0)\le r_0,\qquad C_R(\Zrad(f))=0.
\]
In particular \(\Zrad(f)\) has logarithmic capacity zero. The
value \(a=f(0)\) satisfies
\begin{equation}\label{eq:local-origin-bounds}
 0<a_0-r_0\le a\le a_0+r_0<1,\qquad
 a_0e^{-r_0}\le a\le a_0e^{r_0}.
\end{equation}
Moreover,
\begin{equation}\label{eq:local-distance-certificate}
 P_{[f]}1=af,\qquad
 \dist_{\D}(1,[f])^2=1-a^2\ge1-(a_0+r_0)^2>0 .
\end{equation}
For every polynomial \(p\in\mathbb C[z]\),
\[
 \|1-pf\|_{\D}^2=1-a^2+\|(p-a)f\|_{\D}^2 ,
\]
and the constant polynomial \(p=a\) attains the minimum.

For every neighborhood \(\mathcal O\) of \(f_0\) relative to
\((\X,\dlog)\) and every \(\eta>0\), one can choose such an
\(f\in\mathcal O\cap\Xpr\) with \(|f(0)-a_0|<\eta\).
\end{theorem}

\begin{proof}
Take the centers and balls just constructed. For \(k\ge n\),
nesting gives \(f_k\in B_n\), so \(\dlog(f_k,f_n)\le r_n\).
By \eqref{eq:local-nested-balls}, \(r_n\to0\), which makes
\((f_n)\) Cauchy. Theorem~\ref{thm:outer-completeness} gives a
limit \(f\in\X\) in the metric \(\dlog\). Here both terms of the
metric matter. In the completeness proof, the limit of the boundary
logarithms reconstructs the Hardy outer function, and strong
Dirichlet convergence retains the unit norm and the canonical
moments \eqref{eq:canonical-moments}.

Fix \(n\ge0\) and let \(k\to\infty\) with \(k\ge n\). The centers
\(f_k\) all lie in \(B_n\), and closedness of this ball puts \(f\)
in \(B_n\) as well. This gives the error bound
\begin{equation}\label{eq:local-center-error}
 \dlog(f,f_n)\le r_n\le r_0\,4^{-n}\qquad(n\ge0).
\end{equation}
The common point of the balls is unique because their diameters are
bounded by \(2r_n\). For every \(n\ge1\), apply the capacity estimate
\eqref{eq:local-ball-capacity} to \(g=f\) in \(B_n\).
It follows that \(0\le\Phi(f)<2^{-n}\) for every \(n\ge1\), and
therefore \(\Phi(f)=0\). This step uses the bound on every point
of each \(B_n\). By definition, \(\Phi(f)\) is the \(C_R\)-capacity
of the radial zero set. The Borel description preceding
\eqref{eq:radial-zero-capacity} applies to this set, and
Section~\ref{subsec:03-pairing} gives the equivalence of its
\(C_R\)-polarity and zero logarithmic capacity.

Membership in \(B_0\) gives \(\dlog(f,f_0)\le r_0\). With
\(a=f(0)\), the first line of \eqref{eq:metric-origin}
bounds \(|a-a_0|\), and the second bounds \(|\log a-\log a_0|\).
Exponentiating the latter
bound gives the multiplicative interval in
\eqref{eq:local-origin-bounds}. The additive interval gives
\(a\le a_0+r_0<1\), so \(f\in\Xpr\).

The canonical moments of \(f\) identify \(af\) as the projection
of \(1\) onto \([f]\) by \eqref{eq:certificate}. The residual
\(1-af\) has squared norm \(1-a^2\), which is at least
\(1-(a_0+r_0)^2\). To approximate \(1\) by \(pf\), decompose
\(1-pf\) into this residual and \((a-p)f\in[f]\).
Orthogonality gives \eqref{eq:polynomial-error}, and choosing the
constant polynomial \(p=a\) makes the second term vanish.

For the neighborhood assertion, choose an open metric ball about
\(f_0\) contained in \(\mathcal O\). Repeat the construction with an
initial radius \(s<\min\{r_0,\eta\}\) smaller than the radius of that
open ball. Then \(\overline B_{\dlog}(f_0,s)\subset\mathcal O\)
and the resulting \(f\) has \(|f(0)-a_0|\le s<\eta\). Choosing
\(s<\varepsilon\) also gives \(\dlog(f,f_0)<\varepsilon\)
for any prescribed \(\varepsilon>0\).
\end{proof}

The residual vectors themselves converge in Dirichlet norm.
Put \(a_n=f_n(0)\) and \(a=f(0)\). In the difference of the
residuals, the change in the coefficient contributes at most
\(|a_n-a|\) because \(f_n\) has unit norm. The change in the
function contributes at most \(a\|f_n-f\|_{\D}\), with \(a\le1\).
Evaluation at zero and \eqref{eq:local-center-error} bound each
contribution by \(r_n\). Thus
\begin{equation}\label{eq:local-residual-vectors}
 \|(1-a_nf_n)-(1-af)\|_{\D}
 \le |a_n-a|+a\|f_n-f\|_{\D}\le 2r_n .
\end{equation}
The geometric error in \eqref{eq:local-center-error} measures the
distance from each selected center to the limit. To find a center,
the proof chooses a logarithmic threshold, an atomic measure, and
a recovery index for the current target. The subsequent capacity
threshold puts an upper bound on its ball radius. The closed balls are nested in $\X$, but the compact sets
defining their centers may vary from step to step.

\subsection{The counterexample and the distance from 1 to the space it generates}\label{subsec:13-numerical}

The local theorem can start from any proper outer canonical
function. We first use a concrete function from the boundary
construction. To obtain a numerical lower bound on the approximation
error, we then use a function whose value at zero is prescribed.

\begin{proof}[Proof of Theorem~\ref{thm:main}]
Take the proper outer canonical function \(e_0\) constructed in
Section~\ref{subsec:03-seed} from a nonpolar Carleson set. Its
value at zero, \(a_0=e_0(0)\), lies in \((0,1)\), so we can choose
\(0<r_0<\frac14\min\{a_0,1-a_0\}\). Apply
Theorem~\ref{thm:local} with \(f_0=e_0\). It gives \(f\in\Xpr\),
a Hardy outer canonical function with \(C_R(\Zrad(f))=0\).
The equivalence of polar sets in Section~\ref{subsec:03-pairing}
gives zero logarithmic capacity for this radial zero set.

To check noncyclicity, use the projection identity at the resulting
function itself. With \(a=f(0)\in(0,1)\), \eqref{eq:certificate}
gives
\[
 P_{[f]}1=af,\qquad
 \operatorname{dist}_{\D}(1,[f])^2=1-a^2>0.
\]
The positive distance excludes \(1\) from \([f]\), which makes
\([f]\) a proper subspace of \(\D\). Equation~\eqref{eq:polynomial-error}
gives the exact error for each polynomial multiple of \(f\).
This gives a counterexample to Question~12 of Brown and Shields
\cite[Question~12, p.~296]{BS84}, formulated using the radial zero
set defined in \cite[Notation~3(a), p.~292]{BS84}.
\end{proof}

Brown and Shields also asked whether every outer \(f\in\D\)
generates \(\D_E\) for some \(E\subset\T\)
\cite[Question~11, p.~296]{BS84}. For the function just constructed,
suppose that \([f]=\D_E\). Membership of \(f\) forces
\(E\setminus\Zrad(f)\) to be polar by the radial trace agreement.
Since \(\Zrad(f)\) is polar, \(E\) is polar. Its trace condition
then holds for every Dirichlet function, giving \([f]=\D\), a
contradiction. Thus this \(f\) also gives a negative answer to
Question~11. Furthermore, \(\Zapp(f)\) must be uncountable by
the countable-set theorem of Hedenmalm--Shields and
Richter--Sundberg \cite[Theorem~1.2]{EFKR06}, since \(f\) is
outer and noncyclic. The capacity conclusion of
Theorem~\ref{thm:main} concerns the radial set.

The following bounds on \(f(0)\) give a numerical instance of the
distance formula in \cite[Theorem~2.3]{BFKSS18} at the counterexample.

\begin{corollary}[Numerical specialization]\label{cor:numerical}
The function in Theorem~\ref{thm:main} may be chosen with
\begin{equation}\label{eq:numerical-origin}
 \tfrac12 e^{-1/16}\le a=f(0)\le\tfrac12 e^{1/16}<1.
\end{equation}
For this choice,
\begin{equation}\label{eq:numerical-gap}
 \operatorname{dist}_{\D}(1,[f])^2
=1-a^2\ge 1-\tfrac14 e^{1/8}>\tfrac{175}{256}.
\end{equation}
For every \(p\in\mathbb C[z]\),
\[
 \|1-pf\|_{\D}^2=1-a^2+\|(p-a)f\|_{\D}^2  \ge 1-\tfrac14 e^{1/8}.
\]
\end{corollary}

\begin{proof}
To bound \(1-a^2\) from below, we start at a known value at zero
and use the multiplicative bounds in Theorem~\ref{thm:local}.
Take the first member of the family in Theorem~\ref{thm:G} at
prescribed value at zero equal to \(1/2\) and direction \(1\). This gives a function
\(f_0\in\Xpr\) with \(f_0(0)=1/2\). The radius \(r_0=1/16\)
is admissible because the upper bound for the initial radius is
\(\frac14\min\{f_0(0),1-f_0(0)\}=1/8\).
The local theorem gives a counterexample whose value at zero lies
between \(\tfrac12e^{-1/16}\) and \(\tfrac12e^{1/16}\).

For \(0<x<1\), the exponential series is strictly smaller than
the geometric series \(1/(1-x)\). At \(x=1/16\) this gives
\(e^{1/16}<16/15\), so the upper bound for \(a\) in
\eqref{eq:numerical-origin} is strictly below one. Squaring that
endpoint gives \(a^2\le\tfrac14 e^{1/8}\). At \(x=1/8\) the same
comparison gives \(e^{1/8}<8/7<81/64\), and therefore
\[
 1-\tfrac14 e^{1/8}>1-\tfrac{81}{256}=\tfrac{175}{256}.
\]
The distance identity \eqref{eq:certificate} now gives
\eqref{eq:numerical-gap}. By \eqref{eq:polynomial-error},
every polynomial multiple has this same squared distance
\(1-a^2\) as a lower bound for its squared error. The constant
polynomial \(p=a\) attains the squared distance exactly.
\end{proof}

\subsection{Dense radial zeros and small evaluation norms on the generated subspace}\label{subsec:13-small-values}

Richter and Sundberg constructed a nonzero invariant subspace
\(\D_E\) in which every function has limit inferior zero at
every boundary point as \(z\) approaches that point from within
\(\DD\) \cite[Theorem~4.3, pp.~876--877]{RS94}.
In their construction, \(E\) has positive logarithmic capacity
in every open arc. The following corollary gives a dense
\(G_\delta\) family in \((\Xpr,\dlog)\) whose members have polar
radial zero sets. For each member \(f\) of this family and each \(\zeta\in\T\),
the corollary gives a sequence \(z_n\in\DD\) tending to \(\zeta\)
such that \(h(z_n)\to0\) uniformly for \(h\) in the unit ball of \([f]\).

\begin{corollary}[Dense polar radial zeros]\label{cor:dense-radial-zero-set}
There is a dense \(G_\delta\) subset of \((\Xpr,\dlog)\) whose
members \(f\) are noncyclic outer functions with a dense polar
\(G_\delta\) radial zero set in \(\T\). At every point of this set,
\(f\) tends to zero along every approach from within \(\DD\).
For each such \(f\) and every \(\zeta\in\T\), there are
points \(z_n\in\DD\) tending to \(\zeta\) for which
\[
 \sup_{\substack{h\in[f]\\\|h\|_{\D}\le1}} |h(z_n)|\longrightarrow0.
\]
In particular, \(\Zapp(f)=\T\).
There is also a dense \(G_\delta\) subset of \(\T\), contained in
\(\Zrad(f)\), at whose points the sequence can be chosen on the
radius to that point.
\end{corollary}

\begin{proof}
We first show that an atom of small mass can be added at any
prescribed \(\xi\in\T\). Fix a nonzero finitely supported positive
measure \(\mu\), and put \(M_t=M_{\mu+t\delta_\xi}\subset M_\mu\).
Write \(q_\xi(z)=1-\bar\xi z\). The functions
\(F_{t,\xi}=q_\xi S_{t,\xi}\) of Lemma~\ref{lem:atomic-witness}
converge locally uniformly to \(q_\xi\) and have squared norms
\(3+2t\to3=\|q_\xi\|_{\D}^2\). Bounded point evaluations
identify the weak limit, so the convergence is strong in \(\D\).
For each polynomial \(p\), \(pG_\mu F_{t,\xi}\in M_t\) by
Hardy divisibility and converges to \(pq_\xi G_\mu\) in \(\D\),
since \(p\) and \(G_\mu\) are multipliers. Since \(q_\xi\) is
cyclic, these limits approximate \(G_\mu\) in \(\D\), so
\(\operatorname{dist}_{\D}(G_\mu,M_t)\to0\).
Now \(P_{M_t}1=P_{M_t}P_{M_\mu}1\) and
\(P_{M_\mu}1=G_\mu(0)G_\mu\). The latter vector is nonzero, so
normalizing the convergent projections gives
\(G_{\mu+t\delta_\xi}\to G_\mu\) in \(\D\).

For an arc \(I\) in a countable base of \(\T\) and \(n\ge2\), let
\(\mathcal V_{I,n}\) be the set of \(f\in\Xpr\) for which some
\(\xi\in I\) and \(1-1/n<r<1\) satisfy
\[
 B_f(r\xi)<1/n,\qquad B_f(z)=\frac{|f(z)|}{\sqrt{1-|z|^2}}.
\]
These sets are relatively open by continuity of evaluation.
Given \(e\in\Xpr\) and a relative \(\dlog\)-neighborhood of \(e\),
Proposition~\ref{prop:outer-continuity} gives \(\varepsilon>0\)
such that all \(g\in\Xpr\) within \(3\varepsilon\) of \(e\) in
\(\D\) belong to the neighborhood. Apply Theorem~\ref{thm:A} to
obtain \(\|G_\mu-e\|_{\D}<\varepsilon\). Choose \(\xi\in I\)
and \(t>0\) so that \(\nu=\mu+t\delta_\xi\) satisfies
\(\|G_\nu-G_\mu\|_{\D}<\varepsilon\). By Hardy point evaluation
on \(G_\nu/S_\nu\), whose \(H^2\) norm is at most one,
\[
 \frac{|G_\nu(r\xi)|}{\sqrt{1-r^2}}
 \le \frac{\exp(-t(1+r)/(1-r))}{1-r^2}\longrightarrow0
 \qquad (r\uparrow1).
\]
Fix a radius \(r>1-1/n\) at which this bound is less than
\(1/(2n)\). Theorem~\ref{thm:B} and continuity of evaluation at
this fixed point give \(g\in\Xpr\) with
\(\|g-G_\nu\|_{\D}<\varepsilon\) and \(B_g(r\xi)<1/n\).
Thus \(\mathcal V_{I,n}\) is dense.

The proper class is open in the complete metric space of
Theorem~\ref{thm:outer-completeness}. Intersect all
\(\mathcal V_{I,n}\) with the open dense sets
\(\mathcal U_{2^{-n}}\cap\Xpr\) of
Proposition~\ref{prop:small-capacity-density}. Their intersection
is a dense \(G_\delta\) subset of \(\Xpr\), and each member has
polar radial zero set. For any \(\zeta\in\T\), take arcs from
the base containing \(\zeta\) with diameters tending to zero.
The corresponding points \(z_n\) tend to \(\zeta\) and satisfy
\(B_f(z_n)<1/n\). Lemma~\ref{lem:completed-canonical-quotient}
and the Hardy evaluation estimate give
\[
 \sup_{\substack{h\in[f]\\\|h\|_{\D}\le1}} |h(z)|\le B_f(z).
\]
Taking \(z=z_n\) and using \(B_f(z_n)<1/n\) shows that the norm of
point evaluation at \(z_n\) on \([f]\) tends to zero.

For fixed \(f\) in this intersection, the sets
\[
 W_n=\bigcup_{1-1/n<r<1}\{\xi\in\T:B_f(r\xi)<1/n\}
\]
are open in \(\T\) and meet every base arc. Their intersection
is dense and supplies radial sequences as claimed.
On these sequences \(f\) tends to zero. By Richter and Sundberg
\cite[Theorem~5.1(b)--(d), p.~877]{RS94}, the radial limit of \(|f|\) at each
\(\xi\in\T\) exists and equals
\(\limsup_{\DD\ni z\to\xi}|f(z)|\). Consequently
\(\Zrad(f)\) contains \(\bigcap_{n\ge2}W_n\) and is the
zero set of the latter upper limit as a function of \(\xi\).
This function is upper semicontinuous, so its zero set is
\(G_\delta\) in \(\T\).
\end{proof}

\paragraph{Strong projection convergence and continuous functions.}
Choose \(f\) from Corollary~\ref{cor:dense-radial-zero-set} and
write \(M=[f]\). Let \(A(\DD)\) denote the disk algebra. Then
\(M\cap A(\DD)=\{0\}\). Indeed, fix \(h\in M\cap A(\DD)\). For each \(\zeta\in\T\), choose
a sequence \(z_n\to\zeta\) as in the corollary. Then \(h(z_n)\to0\),
so continuity gives \(h(\zeta)=0\). The maximum modulus principle
gives \(h=0\) on \(\DD\).

Proposition~\ref{prop:small-capacity-density} gives nonpolar
Carleson sets \(K_j\) whose canonical functions
\(g_j=G_{\D_{K_j}}\in A(\DD)\) tend to \(f\) in \(\D\).
Their polynomial multiples show that \(\D_{K_j}\cap A(\DD)\)
is dense in \(\D_{K_j}\). The projection convergence theorem
of Luo and Richter \cite[Corollary~4.5]{LR15} gives
\[
 P_{\D_{K_j}}y\longrightarrow P_My
 \quad\text{in }\D\quad (y\in\D).
\]

\section*{Acknowledgments}

\noindent 
\textit{Funding.}
X. F. is supported by the National Science and Technology
Council (NSTC) of Taiwan (Grant No.~114-2115-M-A49-003-MY3).
F. G. is supported by the National Natural Science Foundation of China
(Grant No. 12601229).

\medskip
\textit{AI Statement.}
The proof was developed during a long chain of discussions with artificial-intelligence tools.
The authors also used AI tools for language editing and
\LaTeX\ formatting. All mathematical content was independently written and checked by the
authors, who take full responsibility.

\clearpage

\section*{Notation}
{\addcontentsline{toc}{section}{Notation}}
 
\begin{longtable}{@{}L{.19\textwidth}L{.56\textwidth}L{.18\textwidth}@{}}
\toprule Symbol & Meaning & Definition\\\midrule\endfirsthead
\toprule Symbol & Meaning & Definition\\\midrule\endhead
$\DD,\T,m,A_0$ & Unit disk, circle, normalized arclength and area & Section~\ref{subsec:01-problem}\\
$\D,D(f)$ & Dirichlet space and derivative energy & \eqref{eq:full-norm}\\
$[f]$ & Closed polynomial module in the Dirichlet norm & \eqref{eq:module}\\
$f^*$ & Hardy boundary values as an arclength a.e. class & Section~\ref{subsec:01-problem}\\
$\Zrad(f)$ & Radial zero set defined by complex radial limits & \eqref{eq:radial-zero}\\
$P_M,G_M$ & Projection onto \(M\) and its canonical function when \(P_M1\ne0\) & Section~\ref{subsec:02-projection}\\
$k_a,K(r)$ & Evaluation kernel and its radial diagonal & Section~\ref{subsec:02-norm}\\
$D_\xi(F),P[v]$ & Local Dirichlet energy and Poisson integral of \(v\) & Section~\ref{subsec:02-energy}\\
$\ell_e,\sigma_e$ & Local Dirichlet integral and associated measure & \eqref{eq:source}\\
$w,D_w$ & Canonical boundary weight and Poisson-weighted energy & Section~\ref{subsec:02-multiplier}\\
$\mu,S_\mu,M_\mu,G_\mu$ & Singular measure, associated singular inner function, invariant subspace, and canonical function & \eqref{eq:atomic-target}, Section~\ref{subsec:12-assembly}\\
$m_\mu,E_\mu,N$ & Total mass, support, and number of atoms for finite support & Theorem~\ref{thm:B}\\
$V,\widetilde u$ & Harmonic trace space and quasi-continuous representative & \eqref{eq:V}, Section~\ref{subsec:03-trace}\\
$C_R$ & Obstacle capacity & Section~\ref{subsec:03-norm}\\
$\D_K$ & Trace-zero space & Section~\ref{subsec:03-ideal}\\
$\Phi(f),E_t(f)$ & Radial-zero capacity and near-zero trace level set & Section~\ref{subsec:03-neighborhood}\\
$\X,\Xpr$ & Outer canonical class and its proper part & Section~\ref{subsec:04-class}\\
$\dlog$ & Dirichlet distance plus \(L^1\) boundary-log distance & Section~\ref{subsec:01-logarithm-history}, \eqref{eq:dlog}\\
$T_f(t)$ & Boundary-log tail of the fixed target \(f\) & Proposition~\ref{prop:outer-continuity}\\
$K_{\mu,\ell},H_{\mu,\ell}$ & Recovery compactum and canonical function for fixed \(\mu\) & Theorem~\ref{thm:B}\\
$\mathcal H_w$ & Hardy functions with finite weighted Dirichlet energy & Section~\ref{subsec:05-model}\\
$\rho_e,Q_e$ & Boundary density and extended-valued quadratic form & \eqref{eq:cost-density}\\
$c_e,T_e$ & Positive lower bound for \(\ell_e\) and weighted trace on \([e]\) & \eqref{eq:cost-floor}, \eqref{eq:cost-trace}\\
$M_A,q_A,e_A,b(A)$ & Finite-zero space, complementary projection of \(1\), canonical function, and logarithmic sum & \eqref{eq:finite-zero-projection}\\
$V_\nu$ & Total mass of the density and its quadrature & Section~\ref{subsec:06-quadrature}\\
$W_I$ & Local energy weight of the inner function \(I\) & \eqref{eq:finite-product}\\
$V^*$ & Dual of the harmonic trace space & \eqref{eq:boundary-dual-norm}\\
$A^2(\Omega)$ & Bergman space with area \(dx\,dy/\pi\) & Section~\ref{subsec:08-companion}\\
$k,\beta$ & Derivative $(ze)'$ and Bergman function with distributional trace $e^*\overline{k^*}$ & \eqref{eq:companion-identity}\\
$b_\rho,B_\rho$ & Modulus lower bound and continuation area bound & \eqref{eq:common-nonzero-collar}, \eqref{eq:common-continuation-bounds}\\
$\lambda_j,T$ & Boundary logarithm measure and limiting mass at \(\xi\) & Section~\ref{subsec:09-logarithm}\\
$\Zapp(h),E^{(t)}$ & Unrestricted approach-zero set and arclength neighborhood of \(E\) & Section~\ref{subsec:09-cyclicity}\\
$I_j,N_j$ & Angular interval and its number of cells & \eqref{eq:fixed-packed-scales}\\ 
$\varphi_b,c(b),B(b)$ & Disk involution, squared value at zero, and squared energy evaluation norm & Section~\ref{sec:11}\\
$P_\ell,P$ & Trial products for strong recovery & Section~\ref{subsec:12-products}\\
$K_\ell,H_\ell$ & Union compactum and its canonical function & Section~\ref{sec:12}\\
\bottomrule\end{longtable}

\end{document}